\documentclass[11pt]{article}
\usepackage[T1]{fontenc}
\usepackage{lmodern}
\usepackage[margin=1in]{geometry}
\usepackage{amsmath,amssymb,amsthm,mathtools,mathrsfs}
\usepackage{array,booktabs}
\usepackage{float}
\usepackage{needspace}
\usepackage{tikz}
\usetikzlibrary{arrows.meta}
\usepackage{microtype}
\usepackage[nocompress]{cite}
\usepackage[colorlinks=true,linkcolor=blue,citecolor=blue,urlcolor=blue]{hyperref}
\newcommand{\R}{\mathbb R}
\newcommand{\Z}{\mathbb Z}
\newcommand{\N}{\mathbb N}
\newcommand{\Sphere}{\mathbb S}
\newcommand{\norm}[1]{\left\lVert#1\right\rVert}

\DeclareMathOperator{\dist}{dist}
\DeclareMathOperator{\supp}{supp}
\newtheorem{theorem}{Theorem}[section]
\newtheorem{proposition}[theorem]{Proposition}
\newtheorem{lemma}[theorem]{Lemma}
\newtheorem{corollary}[theorem]{Corollary}
\theoremstyle{definition}
\newtheorem{assumption}[theorem]{Assumption}
\theoremstyle{remark}
\newtheorem{remark}[theorem]{Remark}
\numberwithin{equation}{section}

\def\ns{\noalign{\smallskip}}

\def\ds{\displaystyle}

\allowdisplaybreaks[2]
\makeatletter
\newcounter{proofstepanchor}
\newcommand{\proofsteplabel}[2]{%
\begingroup
\refstepcounter{proofstepanchor}%
\def\@currentlabel{#2}%
\label{#1}%
\endgroup}
\makeatother
\title{\bf Observability of wave equations with time- and space-dependent potentials: a counterexample under GCC and a positive criterion on the disk}
\author{
Zhonghua Liao\thanks{School of Mathematics, Sichuan University, Chengdu, 610064, China. The research of this author is supported by the NSF of China under grants 12401586.  {\small\it E-mail:} {\small\tt zhonghualiao@yeah.net}. }
\quad  Qi L\"{u}\thanks{School of Mathematics, Sichuan University, Chengdu, 610064, China. The research of this author is supported by the NSF of China under grants 12625102. {\small\it E-mail:} {\small\tt lu@scu.edu.cn}.}  \quad Xu Zhang\thanks{School of Mathematics, Sichuan University, Chengdu, 610064, China. The research of this author is partially supported by the National Key R\&D Program of China under grant 2024YFA1013400 and the New Cornerstone Investigator Program. {\small\it E-mail:} {\small\tt zhang\_xu@scu.edu.cn}.}}

\date{}

\begin{document}

\maketitle

\begin{abstract}
We show that the geometric control condition (GCC) does not guarantee internal observability for a wave
equation with a
smooth time-dependent potential. On a bounded connected domain with
analytic boundary, we construct a nonempty open observation set such
that the pair satisfies GCC and the free Dirichlet wave is observable.
For this same domain and observation set, we construct a smooth
$2\pi$-periodic potential for which there exists a nonzero smooth
Dirichlet wave that vanishes on the observed cylinder for all time.
The wave is the unique mild solution in the space
$X_\Omega=L^2(\Omega)\times H^{-1}(\Omega)$ on every finite time
interval, so no finite observation horizon yields an observability
inequality for this fixed potential.

The construction starts from a real-valued local counterexample to
unique continuation across a noncharacteristic hyperplane. This local
statement is of independent interest, since earlier constructions for
the wave operator use complex-valued waves and potentials. The global
counterexample is then obtained by joining flat oscillatory solutions
along the edges of a polygonal observation set, repairing their sum at
the vertices without changing its zero set, and using a spatial
transition to a common time-independent pair near the boundary. The
planar construction also yields counterexamples in every space
dimension $n\ge3$.

For comparison, and to identify a geometry that supplies both
propagation and continuation, we show on the unit disk that an
observation region containing a central disk and a neighborhood of an
embedded cross-cut yields observability for every weak wave at
sufficiently large observation time, with constants uniform over
bounded measurable potentials. The proof uses a cut, logarithmic polar
coordinates, and a Carleman weight whose favorable flux is confined to
the observed artificial boundary.
\end{abstract}


\clearpage
\begingroup
\small
\tableofcontents
\endgroup
\clearpage

\section{Introduction}

Let $\Omega\subset\mathbb R^2$ be a bounded domain with smooth boundary, let
$\omega\subset\Omega$ be open, and let $T>0$. We consider the Dirichlet wave
equation
\begin{equation}\label{eq:intro-wave}
\begin{cases}
p_{tt}-\Delta p+V(t,x)p=0 &\text{in }(0,T)\times\Omega,\\
p=0 &\text{on }(0,T)\times\partial\Omega,\\
p(0)=p_0,\quad p_t(0)=p_1 &\text{in }\Omega.
\end{cases}
\end{equation}
For $V\in L^\infty((0,T)\times\Omega)$, the natural state space for
an $L^2$ observation of the displacement is
$X_\Omega\triangleq L^2(\Omega)\times H^{-1}(\Omega)$. The corresponding
weak solution lies in
$C([0,T];L^2(\Omega))\cap C^1([0,T];H^{-1}(\Omega))$; see, for example,
\cite[Theorem~2.2(ii)]{DZZ}.

The internal observability problem asks
whether there exists $C>0$ such that
\begin{equation}\label{eq:intro-observation}
\|(p_0,p_1)\|_{X_\Omega}^2
\le C\int_0^T\!\int_\omega |p(t,x)|^2\,\mathrm{d}x\,\mathrm{d}t
\end{equation}
for every initial state. By the Hilbert Uniqueness Method, \eqref{eq:intro-observation} is
equivalent to exact internal controllability in
$H_0^1(\Omega)\times L^2(\Omega)$ with an $L^2$ control
\cite{Lions,LionsBook,TucsnakWeiss}. Therefore, its validity depends on the
observation geometry, the observation time, and the lower-order coefficients.

For the free wave equation, the geometric control condition (GCC) gives a
fundamental sufficient condition for observability: every unit-speed generalized
ray must meet the observation region within a common finite time. The internal
result appears in Bardos--Lebeau--Rauch \cite{BLRinterior}, while their boundary
theory is developed in \cite{BLR}; sharp geometric questions are also addressed
by Burq--G\'erard \cite{BurqGerard}. In a domain with boundary, the relevant
rays include reflection, glancing, and gliding. Their treatment uses the
propagation theory of Melrose--Sj\"ostrand \cite{MS1,MS2}; see also
\cite[Chapter~24]{HormanderIII}. Le Rousseau--Lebeau--Terpolilli--Tr\'elat
\cite{LRLT} extend this observation framework to moving observation sets. In
particular, their Proposition~1.2 and Theorems~1.8 and~1.8$'$ provide the free
displacement estimate in the space $X_\Omega$. Note that their coefficients are
independent of time, even though the observation region may move. For the
one-dimensional free wave on the torus, Niu--Wang--Xiang
\cite{NiuWangXiang} characterize observability from measurable spacetime
sets by GCC, formulated in a measure-theoretic sense, together with an
observable symmetry condition.

The GCC theory, however, concerns only the principal symbol.
A zero-order potential leaves the principal symbol and its rays unchanged,
so this observation alone does not settle \eqref{eq:intro-observation}.
Indeed, the usual proof combines a propagation estimate, which controls high
frequencies up to a weaker remainder, with a uniqueness argument, which excludes
solutions that vanish in the observed cylinder. The sensitivity of high-frequency
waves to autonomous potentials is studied quantitatively in \cite{DehmanErvedoza}.
For a time-dependent potential, a compactness argument still requires a
unique-continuation statement for the limiting nonautonomous equation.
Consequently, ray geometry and unique continuation play distinct roles in the
observability problem.

For time-independent potentials, spectral arguments and elliptic unique continuation can remove the remaining
obstruction. Laurent--L\'eautaud \cite{LaurentLeautaud} obtain uniform estimates
under GCC, including explicit dependence on the norm of a bounded potential
$V=V(x)$. In the partially analytic setting, the local
results of Tataru \cite{Tataru1995} and Robbiano--Zuily \cite{RobbianoZuily},
and the quantitative theory of Laurent--L\'eautaud \cite{LaurentLeautaudUCP},
establish continuation results using analyticity in the
appropriate variables. These results cannot be applied to an arbitrary merely
measurable dependence on $t$.

A different approach is provided by global Carleman estimates, which can
handle time-dependent lower-order terms. Zhang \cite{Zhang2000,ZhangLowerOrder}
proves explicit observability estimates with nonsmooth time-dependent
coefficients under specified geometric assumptions. Duyckaerts--Zhang--Zuazua
\cite{DZZ} develop estimates for systems with potentials in
$L^\infty(0,T;L^r(\Omega))$, including bounded measurable potentials. Their
internal wave estimate uses a noncritical pseudoconvex spatial weight and
an observation collar of its favorable boundary region. They also study the
growth of the observability constant as the potential increases. The sharpness
constructions for systems in that paper concern this quantitative growth, and
should be distinguished from the existence of a nonzero scalar wave whose
observation vanishes identically.

More recent Carleman constructions have substantially refined the observation
regions. Shao \cite{Shao2019} obtains boundary observability on domains with
timelike moving boundaries for wave equations with time-dependent lower-order
terms. Jena \cite{Jena2021} derives interior observability by restricting
an outflow collar to the exterior of a null cone, and also treats observation
centers inside the domain. His subsequent work \cite{Jena2022} treats interior
control on time-dependent domains. The principal statements in these two
papers are formulated for smooth lower-order coefficients, without assuming
analyticity in time. Fu--Liao \cite{FuLiao2022} obtain interior estimates for
variable spatial principal coefficients, with a bounded time-dependent
potential and appropriately regular first-order terms. Fu--Liao--L\"u
\cite{FuLiaoLu} further establish internal and boundary estimates for
hyperbolic systems with rough potentials, with explicit constants and refined
observation regions.

Thus, these results rely on specific weights and their associated sign
conditions, rather than on GCC alone. This raises the question whether GCC
itself, without additional weight conditions, can guarantee observability
for arbitrary bounded time-dependent potentials. The present paper answers
this question negatively.

\medskip

We now state the contributions of this paper. The main contribution is a
counterexample showing that GCC for the free equation does not imply
observability once an arbitrary smooth time-dependent potential is present.
We construct a bounded connected domain $\Omega\subset\mathbb R^2$ with
analytic boundary, a nonempty open set $\omega\Subset\Omega$ satisfying
GCC, and real-valued functions $p,V\in C^\infty(\R\times\overline\Omega)$
such that
\[
\Box p+Vp=0,\qquad
p|_{\R\times\partial\Omega}=0,\qquad
p|_{\R\times\omega}=0,\qquad \|p(0)\|_2>0.
\]
The pair $(p,V)$ is $2\pi$-periodic in time, every derivative of each is
bounded on the whole time axis, and the wave is the unique mild solution in
$X_\Omega$ on every finite interval. Consequently, for this one fixed
potential no finite observation horizon yields an observability inequality,
even though the free wave is observable on the same domain and observation
set.

The counterexample rests on a local statement of independent interest:
a counterexample to local unique continuation across a noncharacteristic
hyperplane for the wave equation with a smooth \emph{real-valued} scalar
potential. This local statement is new in the real-valued setting.
Alinhac--Baouendi \cite{AlinhacBaouendi1979,AlinhacBaouendi} constructed
local nonuniqueness for operators of principal type under smooth
zeroth-order perturbations; for the wave operator their construction uses
a complex-valued wave and a complex-valued potential. Guisset--Shao
\cite{GuissetShao2024} extended this approach to inverse-square singular
potentials and to larger interface regions, again in the complex-valued
setting. In contrast, the local counterexample in this paper is real-valued,
periodic on the entire time axis, and smooth up to the interface. The key
point is to handle the zeros of the oscillatory solution: the construction
confines $\Box f$ to a region where $f$ is strictly positive, so that
$V=-\Box f/f$ is smooth and real there and extends by zero elsewhere.
Taking the real part of a complex counterexample does not suffice, since a
complex potential couples the real and imaginary parts of the wave.

The passage from the local counterexample to the global one requires
additional work. The global construction begins with flat oscillatory
solutions along the edges of a polygonal observation set, repairs their
residuals where the solutions meet without changing the zero set, and joins
the resulting pair to a common time-independent pair whose solution is
positive just inside the outer boundary. The repair uses a finite-chain
correction whose denominator is a sum of squares. This denominator remains
nonzero in the repair region, even at zeros of the oscillatory solution.
A spatial transition then brings the edge and vertex pairs to $(p,V)=(1,0)$,
where they can be glued to a static Dirichlet pair. The observation polygons
satisfy GCC by a triangular-lattice interception argument, and the analytic
boundary contains a concave arc that prevents rays from remaining in a thin
boundary collar.

The planar counterexample also provides the starting point for
counterexamples in every higher space dimension. For each fixed integer
$n\ge3$, a dimension-dependent axisymmetric conjugation and a positive
static completion turn the planar local solution into a smooth
$n$-dimensional wave that vanishes on a ball and equals one outside a
compact set. Combining these local waves with a finite-horizon ball
packing and an analytic boundary from whose thin collar every generalized
ray escapes gives a bounded star-shaped domain in $\R^n$ and a finite
union of mutually separated congruent interior observation balls
satisfying GCC. A single smooth real-valued $2\pi$-periodic potential
then admits a nonzero Dirichlet wave invisible for all time, so both
observability and exact internal controllability fail on every finite
interval.

For comparison, and to identify a geometry that supplies both propagation
and continuation, we also prove a positive observability result on the unit
disk, Theorem~\ref{thm:disk}. If the observation region contains a central
disk and a tubular neighborhood of a proper embedded cross-cut joining an
inner circle to the outer boundary, then \eqref{eq:intro-observation} holds
for every sufficiently large observation time and every bounded measurable
time-dependent potential. The constants are uniform when the potential ranges
over a bounded subset of $L^\infty$. The cross-cut may wind around the center
and backtrack in the radial direction, and in particular the condition does
not require its tube to be a collar of the outflow region of one point
multiplier on the original disk; this is the classical geometry associated
with the multiplier method, and the disk theorem shows that it can be
relaxed.

The proof of the disk theorem cuts inside the observed tube and chooses a
continuous logarithmic polar branch on the complement. This cut produces an
artificial boundary contained in the observation region. A unitary conjugation
then retains a scalar divergence-form spatial principal part, for which we
construct a Carleman weight whose favorable boundary flux is confined to the
artificial boundary. The resulting estimate allows us to apply the auxiliary
optimal control argument and the modified energy estimates of
\cite{FuLiaoLu} to obtain observability for weak solutions. The bounded
potential is handled without differentiating its coefficients.

Several questions remain open: a general characterization of observation
sets that ensure observability for arbitrary bounded time-dependent
potentials; the optimal observation time in the disk theorem; and the
optimal dependence of the control cost on the potential bound.

The remainder of this paper is organized as follows.
Section~\ref{sec:main-results} states the main results.
The proof of the counterexample occupies
Sections~\ref{ce:front-section}--\ref{ce:completion-section}:
Section~\ref{ce:front-section} constructs the oscillatory solutions.
Section~\ref{ce:preliminaries} constructs the GCC observation set and
establishes the finite-chain repair and the two general criteria for
merging solutions and reaching a positive static solution.
Section~\ref{ce:completion-section} verifies both criteria,
assembles the local solutions, and completes the proof of
Theorem~\ref{thm:counter}. Section~\ref{disk:section}, which proves the
disk observability theorem, appears before the counterexample sections
for technical convenience; it is independent of them and can be read
separately. Finally, Section~\ref{sec:general-dimensions} extends the
planar counterexample to every dimension $n\ge3$.

\section{Main results}\label{sec:main-results}

For a smooth bounded domain $\Omega$, let $A_\Omega\triangleq-\Delta_D$.
We equip $X_\Omega=L^2(\Omega)\times H^{-1}(\Omega)$ with
\begin{equation}\label{eq:state-norm}
\|(a,b)\|_{X_\Omega}^2\triangleq
\|a\|_{L^2(\Omega)}^2+
\|A_\Omega^{-1/2}b\|_{L^2(\Omega)}^2,
\qquad \Box\triangleq\partial_t^2-\Delta_x.
\end{equation}
Thus $H^{-1}(\Omega)$ is understood as the dual of $H_0^1(\Omega)$
equipped with the Dirichlet-energy norm. We use
$\|\cdot\|_2$ and $\|\cdot\|_\infty$ for the Lebesgue norms on the
indicated domain. We also use Euclidean tensor norms and set
$|\partial^r f|\triangleq\max_{|\alpha|=r}|\partial^\alpha f|$
for Cartesian spacetime derivatives, where $r\in\mathbb N_0$ and
$\alpha\in\mathbb N_0^{1+n}$ in $n$ space dimensions.

Unless otherwise stated, we write $a\asymp b$ if there is a constant $C\ge1$, depending only on
the fixed geometry and dimension, such that $C^{-1}b\le a\le Cb$.
Dependence on additional fixed parameters is written in parentheses using
the constant's own symbol, as in $C=C(r)$ or
$\lambda_0=\lambda_0(\mu)$, or stated locally.

The positive constants $c,C$ may change from line to line.
For $G>0$, the notation $F=O(G)$ means $|F|\le CG$, while
$F=o(G)$ means $|F|/G\to0$ in the limit under consideration.
The symbols $O(\cdot)$ and $o(\cdot)$ denote bounds uniform on the specified regions, with the parameter dependence stated at each use.
We write $E\Subset U$ when $\overline E$ is compact and contained in $U$,
and $\mathbf1_E$ for the indicator function of $E$.

Throughout, generalized rays are understood in the compressed
bicharacteristic sense. When an interior observation reaches the physical
boundary, GCC is formulated using its relatively open extension to
$\overline\Omega$, while observation integrals are always taken over the
interior. This convention includes gliding rays; see
\cite[Definition~1.6 and Remark~1.7]{LRLT}.

Let
\[
\mathbb D\triangleq\{x\in\R^2:|x|<1\},\qquad 0<\rho<\tfrac14.
\]
A proper cross-cut is a $C^3$ embedding
$\Gamma:[0,1]\to\overline{\mathbb D\setminus B_\rho(0)}$ satisfying
\begin{equation}\label{eq:proper-crosscut}
\Gamma(0)\in\partial B_\rho(0),\quad
\Gamma(1)\in\partial\mathbb D,\quad
\Gamma((0,1))\subset\mathbb D\setminus\overline{B_\rho(0)}.
\end{equation}
For a set $E$ and $\delta>0$ we write
$\mathcal N_\delta(E)\triangleq\{x:\operatorname{dist}(x,E)<\delta\}$.


\begin{theorem}[A permanently invisible wave under GCC]\label{thm:counter}
There exist a bounded connected domain $\Omega\subset\R^2$ with analytic
boundary and a nonempty open set $\omega\Subset\Omega$ such that the
pair $(\Omega,\omega)$ satisfies GCC and the free Dirichlet wave is
observable in $X_\Omega$. For this domain and observation set, there
exist real-valued functions $p,V\in C^\infty(\R\times\overline\Omega)$
satisfying
\begin{equation}\label{eq:counterexample}
\begin{cases}\ds
\Box p+Vp=0 \qquad \mbox{ in }\R\times\overline\Omega,\\
p|_{\R\times\partial\Omega}=0, \qquad p|_{\R\times\omega}=0,\qquad
\|p(0)\|_2>0.
\end{cases}
\end{equation}
The pair $(p,V)$ is $2\pi$-periodic in time and satisfies, for every
fixed spacetime multi-index $\alpha$,
\begin{equation}\label{eq:counter-bounds}
\sup_{\R\times\overline\Omega}
\bigl(|\partial^\alpha p|+|\partial^\alpha V|\bigr)<\infty.
\end{equation}
The constructed wave is the unique mild solution in $X_\Omega$ on each
finite time interval. In particular, for this fixed potential, no
finite observation horizon yields an observability inequality.
\end{theorem}
\begin{remark}[Scope of the counterexample]\label{rem:scope}
The domain and observation set in Theorem~\ref{thm:counter} are
specifically constructed. For a prescribed domain and a prescribed time-dependent
potential, a general characterization of which observation sets
satisfying GCC yield observability remains open.
\end{remark}

The construction of Theorem~\ref{thm:counter} uses the following local
statement, which is of independent interest: a counterexample to local
unique continuation across a noncharacteristic hyperplane for the wave
equation with a smooth \emph{real-valued} scalar potential. In local
coordinates $(t,z,y)$, write $\Box=\partial_t^2-\partial_z^2-\partial_y^2$.
Here flatness at $y=0$ means that every derivative vanishes on that plane;
it allows a smooth extension by zero across the interface.

\begin{theorem}[Failure of local unique continuation]\label{ce:front}
For every $Z>0$ there exist $Y_*>0$ and a pair of real-valued functions
$f,V\in C^\infty(\R\times(-Z,Z)\times(-Y_*,Y_*))$ such that
\begin{equation}\label{ce:front-equation}
\Box f+Vf=0,\qquad
(f,V)=0\quad(y\le0),\qquad
f(0,0,y)>0\quad(y>0).
\end{equation}
The pair $(f,V)$ is $2\pi$-periodic in $t$ and flat at $y=0$.
\end{theorem}
\begin{remark}[The local counterexample with a real-valued potential]\label{rem:local-ucp}
Theorem~\ref{ce:front} is also of independent interest. Alinhac--Baouendi
\cite{AlinhacBaouendi1979,AlinhacBaouendi} constructed local counterexamples
to unique continuation for operators of principal type under smooth
zeroth-order perturbations. For the wave operator, their construction
uses a complex-valued wave and a complex-valued potential; it is local near
one point of the interface and does not produce a real-valued pair.
Guisset--Shao \cite{GuissetShao2024} extended this approach to
inverse-square singular potentials and to larger interface regions along
which the underlying null geodesic family remains regular, again in the
complex-valued setting.

The real-valued conclusion of Theorem~\ref{ce:front} requires more
than taking real parts in the Alinhac--Baouendi type complex geometric-optics constructions. Writing $u=v+iw$ and
$a=A+iB$, with $v,w,A,B$ real, gives
\[
(\Box+a)u=0\qquad\Longrightarrow\qquad
(\Box+A)v=Bw.
\]
The right-hand side need not vanish at a zero of $v$, so division by
$v$ need not produce a smooth real scalar potential. At the level of
the construction, the factors $e^{i\lambda\varphi}$ have no zeros, and
amplitude comparison restricts cancellation between adjacent beams to
selected hypersurfaces, where their residuals are corrected before
dividing by the wave; see \cite[Section~1.4]{GuissetShao2024}.
Real oscillations introduce additional nodal sets within individual
bands, and the original correction does not impose the required
divisibility at these new zeros. Arbitrarily high-order smallness of
the residual in the frequency does not imply smooth divisibility by
the real wave. The real-valued problem is explicitly left open in the
setting of \cite[Remark~1.5]{GuissetShao2024}. Our local construction
instead confines $\Box f$ to a region where $f$ is strictly positive,
so $-\Box f/f$ is real and smooth there and extends by zero elsewhere.

Passing from local counterexamples to an invisible Dirichlet wave on
a GCC geometry requires additional matching. Each local wave comes
with its own potential, whereas the assembled wave must solve one
equation with a common smooth potential. For example, if
$\Box p_j+V_jp_j=0$ for $j=1,2$, then $q=p_1+p_2$ would require
\[
W=\frac{V_1p_1+V_2p_2}{p_1+p_2}\qquad(q\ne0).
\]
At a cancellation point $p_1=-p_2\ne0$, the numerator is
$(V_1-V_2)p_1$ and need not vanish. Even its vanishing at such points
would not, by itself, ensure smooth divisibility by $q$. Cutoffs also
introduce residuals whose quotients need not extend smoothly; this
issue is discussed in
\cite[Section~3.1, before Corollary~3.5]{GuissetShao2024}.

The corrections in the earlier complex-valued local counterexamples control
cancellation between frequency bands within one local construction;
see \cite[Sections~2.2--2.3]{GuissetShao2024}. They do not provide the
estimates for the new cancellation set formed by joining two rotated
copies at a vertex. For those local waves, we do not know how to verify
the junction nondegeneracy, all-order interface matching, and static
entrance conditions of Propositions~\ref{ce:merging-criterion} and~\ref{ce:activation-criterion}.

Our real-valued local construction is designed to provide this
additional structure. On a neighborhood of the support of
$\Box f_i$, the corresponding wave satisfies
$f_i\ge5E_i/8$, where $E_i$ is its positive amplitude envelope.
Near the cancellation set of $f_1+f_2$, the wave producing the
residual is therefore positive and the other wave must be negative;
the portions of the two residual sources requiring repair have
disjoint supports.
A modulus bound for complex waves does not give this separation,
since opposite phases can cancel even when both waves are nonzero
on their residual supports. Together with the phase and derivative estimates of Subsection~\ref{ce:vertex-verification}, the sign relations in \eqref{ce:owner-signs} yield a quantitative lower bound for
the finite-chain repair denominator. The residual decay then
absorbs the losses from division and differentiation, so the repair
coefficients and all their derivatives are flat across both
interfaces. The finite-chain lower bound and flatness of the repair follow from \eqref{ce:incoming-unweighted} and \eqref{ce:merging-correction-flatness}, respectively. A separate choice of the transition
curve verifies the entrance condition in
Proposition~\ref{ce:activation}, allowing completion to $(1,0)$ and
matching to a common static Dirichlet pair. GCC is established
independently by the observation geometry.
\end{remark}

Starting from Theorem~\ref{ce:front}, we construct the global
counterexample in Theorem~\ref{thm:counter} using a polygonal observation
set. This set is placed inside an explicit domain with analytic boundary
and a single concave indentation; it consists of the interiors of finitely
many mutually separated convex polygons, and it satisfies GCC. Along each
straight edge, a rotated and translated local solution has its zero side
facing the polygon interior. Section~\ref{ce:front-section} constructs two
versions of the local solution with different amplitude weights, both
satisfying the qualitative conclusion of Theorem~\ref{ce:front}. An even
number of edges permits these two versions to alternate around each
polygon, while a common small exterior turning angle allows uniform vertex
estimates.

The main difficulty is to merge the solutions at the finitely many vertices
while maintaining a smooth potential and a wave that vanishes in each
observed polygon. A local repair corrects the equation while preserving the
zero set of the sum of the two solutions meeting at a vertex. A common
spatial transition then brings the solution--potential pairs along the
edges and near the vertices to $(p,V)=(1,0)$. They can then be glued to a
prescribed time-independent pair whose solution is positive in a collar
just inside the physical boundary and vanishes on the physical boundary.
The pairs are equal on an open overlap, so the glued pair is smooth and
satisfies the Dirichlet condition.

The planar local assembly also yields the following higher-dimensional
extension of Theorem~\ref{thm:counter}. Its proof is deferred to
Section~\ref{sec:general-dimensions}. Here $\Delta_n$ denotes the Euclidean
Laplacian in $\R^n$, and $X_{\Omega_n}$ carries the norm
\eqref{eq:state-norm}.
\begin{theorem}[Counterexamples for general $n\ge3$]\label{thm:counter-high-dim}
For $n\ge3$, there exist a bounded star-shaped domain
$\Omega_n\subset\R^n$ with real-analytic boundary and a nonempty open
set $\omega_n\Subset\Omega_n$, a finite union of mutually separated
congruent Euclidean balls, such that $(\Omega_n,\omega_n)$ satisfies
GCC. There also exist real-valued $p_n,V_n\in C^\infty(\R\times\overline\Omega_n)$
such that
\begin{enumerate}
\item[(i)] the pair $(p_n,V_n)$ is $2\pi$-periodic and satisfies
\begin{equation}\label{eq:main-wave}
\partial_t^2p_n-\Delta_n p_n+V_np_n=0
\quad\text{in }\R\times\overline\Omega_n,
\qquad
p_n|_{\R\times\partial\Omega_n}=0,
\qquad
p_n|_{\R\times\omega_n}=0;
\end{equation}
\item[(ii)] $\|p_n(t)\|_{L^2(\Omega_n)}>0$ for every $t\in\R$;
\item[(iii)] for every spacetime multi-index $\alpha$,
\begin{equation}\label{eq:main-bounds}
\sup_{\R\times\overline\Omega_n}
\bigl(|\partial^\alpha p_n|+|\partial^\alpha V_n|\bigr)<\infty.
\end{equation}
\end{enumerate}
On each finite interval, $p_n$ is the unique mild solution in
$X_{\Omega_n}$ with its own initial data. For this fixed
$(\Omega_n,\omega_n,V_n)$, no $T>0$ admits \eqref{eq:intro-observation}
with $\Omega=\Omega_n$ and $\omega=\omega_n$.
\end{theorem}
The free Dirichlet wave on this same geometric pair is observable for sufficiently large times by the GCC theorem \cite[Theorems~1.8 and~1.8$'$]{LRLT}. For the fixed smooth nonautonomous potential in Theorem~\ref{thm:counter-high-dim}, however, there is a nonzero wave that vanishes throughout the observed cylinder.


\begin{theorem}[Observation through an embedded cross-cut]\label{thm:disk}
Suppose that an open set $\omega\subset\mathbb D$ satisfies\vspace{-2mm}
\begin{equation}\label{eq:crosscut-observation}
B_{2\rho}(0)\cup
\bigl(\mathcal N_\delta(\Gamma([0,1]))\cap\mathbb D\bigr)
\subset\omega\vspace{-2mm}
\end{equation}
for some proper cross-cut $\Gamma$ and some $\delta>0$. Then there
exists $T_*=T_*(\rho,\Gamma,\delta)>0$ such that the following holds.
For every $T>T_*$ and every $L>0$, there is a constant
$C=C(T,\rho,\Gamma,\delta,L)>0$ such that, for every
$V\in L^\infty((0,T)\times\mathbb D)$ with
$\|V\|_\infty\le L$, every weak Dirichlet solution of $\Box p+Vp=0$
satisfies\vspace{-2mm}
\begin{equation}\label{eq:disk-observability}
\|(p(0),p_t(0))\|_{X_{\mathbb D}}^2
\le C\int_0^{T}\int_\omega|p(t,x)|^2\,\mathrm{d}x\,\mathrm{d}t.\vspace{-2mm}
\end{equation}
\end{theorem}
\begin{remark}[Comparison with the point-multiplier geometry]\label{rem:multiplier}
Theorem~\ref{thm:disk} does not require the observation set to contain the
outflow collar of a point multiplier. For example, take a radial segment
as $\Gamma$ and choose $\delta$ small. The closure of the observed tube
then meets the unit circle in an arc of length less than $\pi$, and the
central disk does not meet that circle. For any $x_0\in\R^2$, the outflow
set is\vspace{-3mm}
\[
\Gamma_+(x_0)\triangleq
\{x\in\partial\mathbb D:(x-x_0)\cdot x>0\}
=\{x\in\partial\mathbb D:1-x_0\cdot x>0\}.\vspace{-2mm}
\]
It has length strictly larger than $\pi$ for finite $x_0$, and hence its
entire interior collar cannot be contained in this observation set.
The improvement therefore concerns this classical geometry on the original
disk; the proof still constructs a pseudoconvex weight after cutting.
\end{remark}
%

By the duality between observability and exact controllability
\cite{Lions,LionsBook,TucsnakWeiss}, Theorems~\ref{thm:disk},
\ref{thm:counter}, and~\ref{thm:counter-high-dim} yield the following
corollary. Its proof is given after the proof of
Theorem~\ref{thm:counter-high-dim}.

\begin{corollary}[Consequences for exact controllability]\label{cor:control}
Under \eqref{eq:crosscut-observation}, for every $T>T_*$ and every
$V\in L^\infty((0,T)\times\mathbb D)$, the system\vspace{-2mm}
\begin{equation}\label{eq:controlled}
y_{tt}-\Delta y+Vy=\mathbf1_\omega f,\qquad
y|_{(0,T)\times\partial\mathbb D}=0\vspace{-2mm}
\end{equation}
is exactly controllable in $H_0^1(\mathbb D)\times L^2(\mathbb D)$ by
controls $f\in L^2((0,T)\times\omega)$.

Conversely, for every $n\ge2$, let $(\Omega_n,\omega_n,V_n)$ denote
the data of Theorem~\ref{thm:counter} when $n=2$ and of
Theorem~\ref{thm:counter-high-dim} when $n\ge3$. Then for every $T>0$,
the system
\begin{equation}\label{eq:control}
y_{tt}-\Delta_n y+V_ny=\mathbf1_{\omega_n}f,
\qquad y|_{(0,T)\times\partial\Omega_n}=0
\end{equation}
is not exactly controllable in
$H_0^1(\Omega_n)\times L^2(\Omega_n)$ by controls
$f\in L^2((0,T)\times\omega_n)$.
\end{corollary}

\section{Observability on the disk}
\label{disk:section}

We prove Theorem~\ref{thm:disk} in three steps. Cut out the smooth cell in \eqref{disk:cell-map} and transform the equation by \eqref{disk:operator-conjugation}. Then use \eqref{disk:weight-conditions} to prove the Carleman estimate. Finally combine the displacement estimate \eqref{disk:middle-displacement} with the weak energy bounds \eqref{disk:weak-energy-state-bounds} and absorb the remainder.

On the disk, the equation is\vspace{-2mm}
\begin{equation}\label{disk:equation}
\begin{cases}
p_{tt}-\Delta p+V(t,x)p=0 & \mbox{ in }(0,T)\times \mathbb D,\\
p =0 & \mbox{ on }(0,T)\times\partial\mathbb D.
\end{cases}
\end{equation}

\subsection{The physical cut and exact conjugation}\label{disk:physical-cut}

The cut removes only observed points and opens the remaining annulus
into a smooth cell. Logarithmic polar coordinates then give a single
angular coordinate on that cell, and a unitary conjugation preserves
the $L^2$ quantity in the observation estimate. We proceed in three
steps: first cut the annulus along the observed cross-cut, then round
the corners to obtain a smooth cell, and finally pass to logarithmic
polar coordinates.

\medskip\noindent\textbf{Step 1. An observed cut through the annulus.}
Set
\[
\widehat\omega\triangleq B_{2\rho}(0)\cup
\bigl(\mathcal N_\delta(\Gamma([0,1]))\cap\overline{\mathbb D}\bigr),
\qquad \widehat\omega\cap\mathbb D\subset\omega.
\]
Closures of subsets of $\overline{\mathbb D}$ are taken in
$\overline{\mathbb D}$. By Sard's theorem, choose $r_0\in(\rho,2\rho)$ such that
$r_0^2$ is a regular value of $a\mapsto|\Gamma(a)|^2$. Thus
\[
|\Gamma(a)|=r_0\quad\Longrightarrow\quad
\Gamma(a)\cdot\Gamma'(a)\ne0\qquad(a\in[0,1]).
\]
By \eqref{eq:proper-crosscut} and the intermediate value theorem, the level set
$\{a\in[0,1]:|\Gamma(a)|=r_0\}$ is nonempty and compact, so we may define
\[
a_0\triangleq\max\{a\in[0,1]:|\Gamma(a)|=r_0\},\qquad
A_0\triangleq\{x\in\mathbb R^2:r_0\le|x|\le1\}.
\]
Then
\[
\Gamma(a_0)\cdot\Gamma'(a_0)>0,\qquad
r_0<|\Gamma(a)|<1\quad(a_0<a<1).
\]
Taking the last crossing of the inner circle leaves the remaining arc
inside the annulus, even when the original curve backtracks radially.

We replace this arc by a smooth cross-cut and thicken it inside the
observed region. Smooth approximation and the tubular neighborhood theorem give
$\varepsilon>0$ and a map
\begin{equation}\label{disk:cut-map-properties}
\begin{aligned}
&F\in C^\infty([0,1]\times[-\varepsilon,\varepsilon];\mathbb R^2),
\qquad \det DF(a,b)\ne0,\\
&F(a,b)=F(\widetilde a,\widetilde b)\quad\Longrightarrow\quad
(a,b)=(\widetilde a,\widetilde b),\\
&F(0,0)=\Gamma(a_0),\qquad F(1,0)=\Gamma(1),\\
&|F(0,b)|=r_0,\qquad |F(1,b)|=1
\quad(|b|\le\varepsilon),\\
&r_0<|F(a,b)|<1
\quad(0<a<1,\ |b|\le\varepsilon),\\
&F([0,1]\times[-\varepsilon,\varepsilon])
\subset A_0\cap\mathcal N_\delta(\Gamma([a_0,1])).
\end{aligned}
\end{equation}
Define
\[
\mathcal V\triangleq F([0,1]\times(-\varepsilon,\varepsilon)),\qquad
K_{\rm obs}\triangleq\overline B_{r_0}(0)\cup\overline{\mathcal V}.
\]
The central curve $a\mapsto F(a,0)$, $a\in[0,1]$, is injective and satisfies
\[
|\partial_aF(a,0)|>0,\qquad
F(j,0)\cdot\partial_aF(j,0)>0\quad(j\in\{0,1\}).
\]
Choose $d_0>0$ such that
$\mathcal N_{2d_0}(K_{\rm obs})\cap\overline{\mathbb D}\subset\widehat\omega$,
and set
\[
U_*\triangleq\mathcal N_{d_0}(K_{\rm obs})\cap\overline{\mathbb D}.
\]
Thus
\begin{equation}\label{disk:observed-neighborhood}
\overline B_{r_0}(0)\cup\overline{\mathcal V}
=K_{\rm obs}\subset U_*\subset\overline{U_*}\subset\widehat\omega.
\end{equation}

\medskip\noindent\textbf{Step 2. A smooth cell with observed artificial boundary.}
Set
\[
G_0\triangleq\mathbb D\setminus K_{\rm obs}.
\]
After local $C^\infty$ regularization at
$\{F(j,\sigma\varepsilon):j=0,1,\ \sigma=\pm1\}$, the smooth
Jordan--Schoenflies theorem extends the smooth boundary curve to a
diffeomorphism from the disk onto its interior. This gives a map $J$,
with the prescribed matching along the physical boundary arc, satisfying
\begin{equation}\label{disk:cell-map}
\begin{aligned}
&J\in C^\infty(\overline{\mathbb D};\mathbb R^2),\qquad
\det DJ(x)\ne0\quad(x\in\overline{\mathbb D}),\\
&J(x)=J(y)\quad\Longrightarrow\quad x=y
\qquad(x,y\in\overline{\mathbb D}),\\
&G_{\rm ph}\triangleq J(\mathbb D),\qquad
G_{\rm ph}\subset G_0,\qquad
\overline{G_0\setminus G_{\rm ph}}\subset U_*.
\end{aligned}
\end{equation}
Define
\[
\Sigma_{\rm ph}\triangleq\partial G_{\rm ph}\cap\partial\mathbb D,
\qquad
\Sigma_{\rm art}\triangleq
\overline{\partial G_{\rm ph}\setminus\Sigma_{\rm ph}}.
\]
\begin{center}
\begin{tabular}{@{}ll@{}}\toprule
Object & Meaning\\\midrule
$K_{\rm obs}$ & Observed central disk and thickened cut\\
$G_0=\mathbb D\setminus K_{\rm obs}$ & Region left after removing the observed core\\
$G_{\rm ph}$ & Smooth cell left after rounding the corners\\
$\Sigma_{\rm art}$ & Artificial boundary inside the observed region\\\bottomrule
\end{tabular}
\end{center}
The conditions in \eqref{disk:cell-map} imply
\[
\begin{aligned}
&\mathbb D\setminus G_{\rm ph}
\subset K_{\rm obs}\cup\overline{G_0\setminus G_{\rm ph}}\subset U_*,\qquad
\Sigma_{\rm art}
\subset K_{\rm obs}\cup\overline{G_0\setminus G_{\rm ph}}\subset U_*,\\
&T_x(\partial G_{\rm ph})
=\{\xi\in\mathbb R^2:\xi\cdot x=0\}
\quad(x\in\Sigma_{\rm ph}).
\end{aligned}
\]
Consequently, for some $\varepsilon_*>0$,
\begin{equation}\label{disk:cell-collar}
\begin{aligned}
&\overline{G_{\rm ph}}\subset A_0,\qquad
\mathbb D\setminus G_{\rm ph}\subset U_*\cap\mathbb D\subset\omega,\qquad
\mathcal N_{\varepsilon_*}(\Sigma_{\rm art})\cap\mathbb D
\subset U_*\cap\mathbb D\subset\omega.
\end{aligned}
\end{equation}

\medskip\noindent\textbf{Step 3. Logarithmic polar coordinates.}
The map $J$ gives a simply connected open neighborhood of
$\overline{G_{\rm ph}}$ that does not contain the origin. On this
neighborhood choose a smooth angular branch $\Theta_G$ satisfying
\[
(\cos\Theta_G(x),\sin\Theta_G(x))=\frac{x}{|x|},\qquad
\nabla\Theta_G(x)=\frac{(-x_2,x_1)}{|x|^2}.
\]
Define
\begin{equation}\label{disk:covering}
\begin{aligned}
&\pi(s,\vartheta)\triangleq e^s(\cos\vartheta,\sin\vartheta),\qquad
\kappa(s)\triangleq e^{-2s},\\
&\Phi_G(x)\triangleq(\log|x|,\Theta_G(x)),\qquad
G\triangleq\Phi_G(G_{\rm ph}).
\end{aligned}
\end{equation}
Throughout this section, write $z=(s,\vartheta)$ and
$\mathrm{d}z=\mathrm{d}s\,\mathrm{d}\vartheta$.
The maps $\Phi_G$ and $\pi$ are mutually inverse on neighborhoods of
$\overline{G_{\rm ph}}$ and $\overline G$. Moreover,
\[
\begin{aligned}
&\det D\Phi_G(x)=|x|^{-2},\qquad
G\subset(\log r_0,0)\times\mathbb R,\\
&\sup_{(s,\vartheta)\in G}|\vartheta|<\infty,\qquad
\partial G\in C^\infty,\\
&C^{-1}|x-y|\le|\Phi_G(x)-\Phi_G(y)|\le C|x-y|
\quad(x,y\in\overline{G_{\rm ph}})
\end{aligned}
\]
for some $C>0$. The following lemma shows that the passage to
logarithmic polar coordinates preserves the relevant function spaces
and transforms the Laplacian into divergence form.
\begin{lemma}[Unitary logarithmic conjugation]\label{disk:conjugation}
The map $\mathcal U f\triangleq e^s f\circ\pi$ is unitary from
$L^2(G_{\rm ph})$ onto $L^2(G)$. It induces a bounded isomorphism
$\mathcal U_1\triangleq\mathcal U|_{H_0^1(G_{\rm ph})}$ from
$H_0^1(G_{\rm ph})$ onto $H_0^1(G)$, and, by duality,
a bounded isomorphism $\mathcal U_{-1}\triangleq(\mathcal U_1^{-1})'$
between the corresponding negative spaces, where the prime denotes the
dual operator; on the $L^2$ pivot,
$\mathcal U_{-1}$ equals $\mathcal U$.
In the sense of sesquilinear forms,
\begin{equation}\label{disk:operator-conjugation}
\mathcal U_{-1}(-\Delta)\mathcal U_1^{-1}
=-\operatorname{div}(\kappa\nabla)-\kappa.
\end{equation}
Consequently, if $p$ solves $\Box p+Vp=0$ on
$(0,T)\times G_{\rm ph}$ with homogeneous Dirichlet boundary condition, then
$\widetilde p\triangleq\mathcal U p$ satisfies
\begin{equation}\label{disk:transformed}
\begin{aligned}
&\widetilde p_{tt}-\operatorname{div}(\kappa\nabla\widetilde p)=B\widetilde p,\qquad \mbox{ for }
B(t,s,\vartheta)\triangleq\kappa(s)-V(t,\pi(s,\vartheta)),
\end{aligned}
\end{equation}
in the Dirichlet transposition sense on $(0,T)\times G$.
\end{lemma}

\begin{proof}[Proof of Lemma~\ref{disk:conjugation}]
We prove the three assertions in order: unitarity and the $H_0^1$
isomorphism, the transformation formula for the Laplacian, and the
transformed equation.

\medskip
\noindent\textbf{Step 1: Unitarity and the $H_0^1$ isomorphism.}
The pullback metric and Jacobian of $\pi$ are
$e^{2s}(\mathrm{d}s^2+\mathrm{d}\vartheta^2)$ and $e^{2s}$, respectively.
In particular,
\[
\|\mathcal Uf\|_{L^2(G)}^2
=\int_G e^{2s}|f\circ\pi|^2\,\mathrm{d}z
=\int_{G_{\rm ph}}|f(x)|^2\,\mathrm{d}x
=\|f\|_{L^2(G_{\rm ph})}^2.
\]
The inverse is
\[
\mathcal U^{-1}\widetilde p(x)
=|x|^{-1}\widetilde p(\log|x|,\Theta_G(x)),
\qquad x\in G_{\rm ph}.
\]
To compare the Dirichlet energies, write
$W\triangleq p\circ\pi=e^{-s}\widetilde p$. Since
$\mathrm{d}x=e^{2s}\mathrm{d}s\,\mathrm{d}\vartheta$ and
$|\nabla_x p|^2\circ\pi=e^{-2s}(|W_s|^2+|W_\vartheta|^2)$, we have
\[
\int_{G_{\rm ph}}|\nabla p|^2\mathrm{d}x
=\int_G(|W_s|^2+|W_\vartheta|^2)\,\mathrm{d}z.
\]
Moreover,
\[
W_s=e^{-s}(\widetilde p_s-\widetilde p),\qquad
W_\vartheta=e^{-s}\widetilde p_\vartheta,\qquad
\kappa=e^{-2s},
\]
and hence
\[
\int_{G_{\rm ph}}|\nabla p|^2\mathrm{d}x
=\int_G\kappa\bigl(|\widetilde p_s-\widetilde p|^2
+|\widetilde p_\vartheta|^2\bigr)\,\mathrm{d}z.
\]
Since $r_0\le e^s\le1$, the Poincar\'e inequality shows that the two
quadratic forms
\[
\int_G\kappa\bigl(|\widetilde p_s-\widetilde p|^2
+|\widetilde p_\vartheta|^2\bigr)\,\mathrm{d}z
\quad\text{and}\quad
\int_G\kappa|\nabla\widetilde p|^2\,\mathrm{d}z
\]
are comparable. Density of smooth functions with compact support
therefore gives the two $H_0^1$ bounds, and $\mathcal U_1$ is an
isomorphism by the open mapping theorem.

\medskip
\noindent\textbf{Step 2: The transformation formula for the Laplacian.}
By construction, $\mathcal U_{-1}=(\mathcal U_1^{-1})'$. For
$f\in L^2(G_{\rm ph})$ and $\Phi\in H_0^1(G)$, the identity
$\mathcal U_1^{-1}\Phi=\mathcal U^{-1}\Phi$ gives
\[
\langle\mathcal U_{-1}f,\Phi\rangle
=(f,\mathcal U_1^{-1}\Phi)
=(\mathcal U f,\Phi),
\]
where the parentheses denote $L^2$ inner products. Thus
$\mathcal U_{-1}$ equals $\mathcal U$ on the $L^2$ pivot.

For smooth functions, $\kappa_s=-2\kappa$ and
\[
\mathcal U(-\Delta)\mathcal U^{-1}\widetilde p
=-\kappa(\widetilde p_{ss}+\widetilde p_{\vartheta\vartheta}
-2\widetilde p_s+\widetilde p)
=-\operatorname{div}(\kappa\nabla\widetilde p)-\kappa\widetilde p.
\]
For $\widetilde p,\Phi\in C_c^\infty(G)$, the corresponding weak identity is
\[
\begin{aligned}
\langle\mathcal U_{-1}(-\Delta)\mathcal U_1^{-1}\widetilde p,\Phi\rangle
& =\int_G\kappa\bigl[(\widetilde p_s-\widetilde p)
(\overline\Phi_s-\overline\Phi)
+\widetilde p_\vartheta\overline\Phi_\vartheta\bigr]\,\mathrm{d}z\\
& =\int_G\kappa\nabla\widetilde p\cdot\nabla\overline\Phi
-\kappa\partial_s(\widetilde p\overline\Phi)
+\kappa\widetilde p\overline\Phi\,\mathrm{d}z\\
& =\int_G\kappa\nabla\widetilde p\cdot\nabla\overline\Phi
+(\kappa_s+\kappa)\widetilde p\overline\Phi\,\mathrm{d}z\\
& =\int_G\kappa\nabla\widetilde p\cdot\nabla\overline\Phi
-\kappa\widetilde p\overline\Phi\,\mathrm{d}z.
\end{aligned}
\]
The last equality uses $\kappa_s=-2\kappa$. Since both sides are
continuous in $\widetilde p\in H_0^1(G)$ and $\Phi\in H_0^1(G)$,
density of smooth functions proves
\eqref{disk:operator-conjugation} on $H_0^1(G)$.

\medskip
\noindent\textbf{Step 3: The transformed equation.}
Let $p$ solve $\Box p+Vp=0$ on $(0,T)\times G_{\rm ph}$ with homogeneous
Dirichlet boundary condition, and set $\widetilde p\triangleq\mathcal U p$.
The operators $\mathcal U$ and $\mathcal U_1$ are independent of time.
Thus, in the Dirichlet transposition sense,
\[
\begin{aligned}
0&=\mathcal U_{-1}\bigl(p_{tt}-\Delta p+Vp\bigr) =\widetilde p_{tt}
-\operatorname{div}(\kappa\nabla\widetilde p)
-\kappa\widetilde p+(V\circ\pi)\widetilde p,
\end{aligned}
\]
and therefore
\[
\widetilde p_{tt}-\operatorname{div}(\kappa\nabla\widetilde p)
=(\kappa-V\circ\pi)\widetilde p=B\widetilde p.
\]
The boundedness of $\mathcal U$ and $\mathcal U_{-1}$ also ensures
$
\widetilde p\in C([0,T];L^2(G))
\cap C^1([0,T];H^{-1}(G))$. 
\end{proof}

\subsection{A spatial weight with observed positive flux}\label{disk:weight-section}

We construct a Carleman weight on $G$ whose favorable boundary flux is
supported on the artificial boundary
$\pi^{-1}(\Sigma_{\rm art})\subset\pi^{-1}(U_*)$. This is the key
geometric property of the disk argument: the weight is chosen so that,
on the physical boundary, its normal derivative has the wrong sign,
while on the artificial boundary it has the right sign and the
artificial boundary lies inside the observation region.

For a smooth $d$ define
\begin{equation}\label{disk:quadratic-form}
\mathfrak Q_d(z,\xi)\triangleq2\kappa^2\bigl[
(d_{ss}-d_s)\xi_1^2+2(d_{s\vartheta}-d_\vartheta)\xi_1\xi_2
+(d_{\vartheta\vartheta}+d_s)\xi_2^2\bigr].
\end{equation}
Choose $C_0$ so that $\vartheta+C_0\ge1$ on $\overline G$.
For $\varpi>0$, define
\begin{equation}\label{disk:weight}
d(s,\vartheta)\triangleq
\exp\!\bigl(\varpi[r_0^2e^{-2s}+(\vartheta+C_0)^2]\bigr).
\end{equation}
Let $\nu$ be the outer normal of $G$. The part of the boundary with
strictly positive flux is
\[
\Gamma_+\triangleq\{z\in\partial G:\kappa\nabla d\cdot\nu>0\}.
\]
We call this the positive-flux region. The next lemma proves that the
weight satisfies the pointwise positivity $\mathfrak Q_d\ge4\kappa|\xi|^2$
and that the positive-flux region is observed.

\begin{lemma}[A weight with observed positive flux]\label{disk:weight-lemma}
There exists a constant $\varpi_0=\varpi_0(G,r_0,C_0)>0$, depending
only on these fixed geometric data, such that, for every
$\varpi\ge\varpi_0$, the weight $d$ satisfies
\begin{equation}\label{disk:weight-conditions}
\mathfrak Q_d(z,\xi)\ge4\kappa|\xi|^2,\qquad
E\triangleq\kappa|\nabla d|^2\ge4d>0.
\end{equation}

There exists $\varepsilon_1>0$ such that the set
\begin{equation}\label{disk:positive-collar}
\omega_1\triangleq
\{z\in G:\operatorname{dist}(z,\Gamma_+)<\varepsilon_1\}
\end{equation}
is an open collar of $\Gamma_+$ in $G$ and satisfies
\[
\overline{\omega_1}^{\,G}\subset\pi^{-1}(U_*).
\]
Here the superscript $G$ denotes closure relative to $G$; for example,
$\overline{\omega_1}^{\,G}=\overline{\omega_1}\cap G$, where the closure
on the right is taken in $\R^2$.
\end{lemma}

\begin{proof}[Proof of Lemma~\ref{disk:weight-lemma}]
We prove the two assertions separately: first the pointwise bounds
\eqref{disk:weight-conditions}, then the collar construction
\eqref{disk:positive-collar}.

\medskip
\noindent\textbf{Step 1: The pointwise bounds.}
Throughout this step set
$R\triangleq r_0^2e^{-2s}$ and $Y\triangleq\vartheta+C_0$, so that
$r_0^2\le R\le1$ and $Y\ge1$ on $\overline G$. The derivatives of
\eqref{disk:weight} are
\[
\begin{aligned}
&d_s=-2\varpi Rd,\;\; d_\vartheta=2\varpi Yd,\;\; d_{ss}=(4\varpi R+4\varpi^2R^2)d,\;\;
d_{s\vartheta}=-4\varpi^2RYd,\;\;
d_{\vartheta\vartheta}=(2\varpi+4\varpi^2Y^2)d.
\end{aligned}
\]
Consequently,
$\mathfrak Q_d=2\kappa^2d\,\xi^{\mathsf T}M_\varpi\xi$, where
\[
M_\varpi\triangleq
\begin{pmatrix}
6\varpi R+4\varpi^2R^2&-2\varpi Y(1+2\varpi R)\\
-2\varpi Y(1+2\varpi R)&2\varpi(1-R)+4\varpi^2Y^2
\end{pmatrix},
\]
and
\[
\det M_\varpi
=4\varpi^2\bigl[R(3+2\varpi R)(1-R)+(2\varpi R-1)Y^2\bigr].
\]
For $\varpi\ge r_0^{-2}$,
\[
\begin{aligned}
2\varpi R-1&\ge\varpi R\ge r_0^2\varpi,\\
\det M_\varpi
&\ge4\varpi^2(2\varpi R-1)Y^2
\ge4r_0^2\varpi^3,\qquad
\operatorname{tr}M_\varpi\le C\varpi^2.
\end{aligned}
\]
The first diagonal entry and determinant are positive, so $M_\varpi$ is
positive definite. Its eigenvalues therefore satisfy
\[
\lambda_{\min}(M_\varpi)
=\frac{\det M_\varpi}{\lambda_{\max}(M_\varpi)}
\ge\frac{\det M_\varpi}{\operatorname{tr}M_\varpi}
\ge c\varpi.
\]
Since $\kappa\ge1$ and $d\ge1$ on $\overline G$,
\[
\mathfrak Q_d(z,\xi)
=2\kappa^2d\,\xi^{\mathsf T}M_\varpi\xi
\ge2c\varpi\kappa^2d|\xi|^2
\ge4\kappa|\xi|^2
\qquad(c\varpi\ge2).
\]

For the second inequality, note that
\[
E=\kappa|\nabla d|^2
=4\kappa\varpi^2d^2(R^2+Y^2)
\ge4\kappa\varpi^2d^2Y^2
\ge4d
\qquad(\varpi\ge1),
\]
which proves the second inequality in \eqref{disk:weight-conditions}.

\medskip
\noindent\textbf{Step 2: The favorable collar.}
On the physical boundary arc $\pi^{-1}(\Sigma_{\rm ph})$ of $G$ we have
$s=0$ and $\nu=(1,0)$,
including the common endpoints by smooth tangency. Hence
\[
\kappa\nabla d\cdot\nu=\kappa d_s
=-2\varpi\kappa r_0^2d<0
\]
there, because $d>0$ and $\kappa>0$. Thus $\Gamma_+$ does not meet the
physical arc, and its closure is a compact subset of the artificial arc
$\pi^{-1}(\Sigma_{\rm art})$, separated from
$\pi^{-1}(\Sigma_{\rm ph})$ by a positive distance.
By \eqref{disk:cell-collar}, $\Sigma_{\rm art}$ is compactly contained
in $U_*$. Since $\pi$ is bi-Lipschitz, choose $\varepsilon_1>0$ such that
\begin{equation}\label{disk:artificial-collar-choice}
\{z\in\overline G:
\operatorname{dist}(z,\pi^{-1}(\Sigma_{\rm art}))\le4\varepsilon_1\}
\subset\pi^{-1}(U_*).
\end{equation}
Set
$
\omega_1=\mathcal N_{\varepsilon_1}(\Gamma_+)\cap G.
$
Since $\Gamma_+\subset\pi^{-1}(\Sigma_{\rm art})$, this collar
satisfies \eqref{disk:positive-collar}.
\end{proof}

\subsection{Carleman estimate for the equation \eqref{disk:equation}}

\subsection{Carleman estimate}\label{disk:carleman-section}

We first prove the strong estimate \eqref{disk:carleman-estimate}, then
obtain the weak estimate \eqref{disk:weak-carleman-free}. Fix
$\varpi\ge\varpi_0$, the weight $d$, and the collar $\omega_1$ as in
Lemma~\ref{disk:weight-lemma}. The spatial cutoff removes the
contribution near $\Gamma_+$; the time weight makes the integrated
time-face difference nonpositive, so the strong estimate requires no
zero endpoint data.

In this subsection $\tau$ is a translated physical time coordinate. Set
\[
P_G\triangleq\partial_\tau^2-\operatorname{div}(\kappa\nabla),\qquad
R_G\triangleq\frac12\max_{\overline{G\setminus\omega_1}}\sqrt E.
\]
To separate the bulk quadratic terms from the spacetime boundary fluxes,
we use the scalar multiplier identity in
\cite[Lemma~2.3, equations~(2.9)--(2.11)]{FuLiaoLu}. Here $\Psi$ may
depend on time, so we retain the additional product-rule term
$2\Psi_\tau vv_\tau$.

\medskip
\noindent\textbf{Step 1: A pointwise weighted identity.}
The following lemma is the algebraic core of the estimate. It expresses
$\theta^2|P_Gu|^2$ as a sum of a positive bulk term and spacetime
divergences, up to lower-order errors.

\begin{lemma}[Pointwise weighted identity]\label{disk:identity-lemma}
For real-valued $u,\ell,\Psi\in C^\infty(\mathbb R\times G)$,
set
$\theta\triangleq e^\ell$ and $v\triangleq\theta u$.
The following pointwise identity holds:
\begin{align}
&\theta^2|P_Gu|^2-I_1^2-I_2^2
+2\operatorname{div}\mathbf J+2\mathcal M_\tau
\label{disk:identity}\\
& =2[\ell_{\tau\tau}+\operatorname{div}(\kappa\nabla\ell)
+\Psi]v_\tau^2
-8\kappa\nabla\ell_\tau\cdot\nabla v\,v_\tau
+2\sum_{i,j=1}^2c_\ell^{ij}v_iv_j
-2\kappa\nabla\Psi\cdot\nabla v\,v
+\mathcal B_\ell v^2+2\Psi_\tau vv_\tau,\notag
\end{align}
Here $\partial_1\triangleq\partial_s$,
$\partial_2\triangleq\partial_\vartheta$,
$\ell_i\triangleq\partial_i\ell$ and $v_i\triangleq\partial_i v$,
with $i,j\in\{1,2\}$. The coefficient and flux terms are defined by
\begin{equation}\label{disk:ABC}
\begin{cases}
A_\ell\triangleq\kappa|\nabla\ell|^2
-\operatorname{div}(\kappa\nabla\ell)-\ell_\tau^2
+\ell_{\tau\tau}-\Psi,\\[4pt]
c_\ell^{ij}\triangleq2\kappa\partial_i(\kappa\ell_j)
-\delta_{ij}\operatorname{div}(\kappa^2\nabla\ell)
+\kappa\delta_{ij}(\ell_{\tau\tau}-\Psi),\\[4pt]
\mathcal B_\ell\triangleq2\bigl[A_\ell\Psi
-(A_\ell\ell_\tau)_\tau
+\operatorname{div}(A_\ell\kappa\nabla\ell)\bigr],\\[4pt]
\mathbf J\triangleq 2\kappa^2(\nabla\ell\cdot\nabla v)\nabla v
+\kappa A_\ell v^2\nabla\ell-\kappa\Psi v\nabla v -\kappa^2|\nabla v|^2\nabla\ell
-2\kappa\ell_\tau v_\tau\nabla v
+\kappa v_\tau^2\nabla\ell, \\[4pt]
\mathcal M\triangleq\ell_\tau(v_\tau^2+\kappa|\nabla v|^2)
-2\kappa(\nabla\ell\cdot\nabla v)v_\tau
+\Psi vv_\tau-A_\ell\ell_\tau v^2.
\end{cases}
\end{equation}
The square terms are
\[
I_1\triangleq P_Gv-A_\ell v,\qquad
I_2\triangleq-2\ell_\tau v_\tau
+2\kappa\nabla\ell\cdot\nabla v-\Psi v,
\]
and satisfy $\theta P_Gu=I_1+I_2$.
\end{lemma}

\medskip
\noindent\textbf{Step 2: The strong Carleman estimate.}
We now choose the weight parameters and prove the strong estimate. The
proof has four parts: the pointwise inequality, the bulk estimate, the
boundary fluxes, and the conversion from $v$ back to $u$.

\begin{proposition}[Carleman estimate on the cut domain]\label{disk:carleman}
Fix $S>R_G$ and $\alpha\in(R_G/S,1)$, and define
\[
\varphi\triangleq e^{\mu(d-\alpha\tau^2)},\quad
\ell\triangleq\lambda\varphi,\quad \theta\triangleq e^\ell.
\]
For every sufficiently large $\mu$ there exist constants $C=C(\mu)>0$ and
$\lambda_0=\lambda_0(\mu)>0$ such that
\begin{equation}\label{disk:carleman-estimate}
\lambda\mu\int_{-S}^S\!\int_G\theta^2
(|u_\tau|^2+|\nabla u|^2+\lambda^2\mu^2|u|^2)
\le C(\mu)\int_{-S}^S\!\int_G\theta^2|P_Gu|^2
\end{equation}
for $\lambda\ge\lambda_0(\mu)$ whenever
$u\in H^1((-S,S)\times G)$ has zero lateral trace,
$P_Gu\in L^2$, and $u=0$ on $(-S,S)\times\omega_1$.
\end{proposition}

\begin{proof}[Proof of Proposition~\ref{disk:carleman}]
By smooth approximation, we may assume that $u$ is smooth up to the
boundary, vanishes near $\Gamma_+$, and is spatially supported where
$\sqrt E/2\le R_\sharp<\alpha S$, with
$R_\sharp\triangleq(R_G+\alpha S)/2$. It suffices to treat real-valued
$u$, since
\begin{equation}\label{disk:complex-splitting}
|u|^2=|\operatorname{Re}u|^2+|\operatorname{Im}u|^2,\qquad
|P_Gu|^2=|P_G\operatorname{Re}u|^2+|P_G\operatorname{Im}u|^2,
\end{equation}
and \eqref{disk:complex-splitting} also holds for every first derivative of $u$.

\medskip
\noindent\textit{Part 1: The pointwise inequality.}
Put $v\triangleq\theta u$ and define
\[
\Psi\triangleq\ell_{\tau\tau}-\operatorname{div}(\kappa\nabla\ell)
+2(1+\alpha)\lambda\mu\varphi.
\]
Use the coefficients $A_\ell,c_\ell^{ij},\mathcal B_\ell$ and the terms
$\mathbf J,\mathcal M,I_1,I_2$ defined in
Lemma~\ref{disk:identity-lemma}. Equation~\eqref{disk:identity}, with
$I_1^2+I_2^2\ge0$, implies
\[
\begin{aligned}
&\theta^2|P_Gu|^2+2\operatorname{div}\mathbf J+2\mathcal M_\tau\\
& \ge
2[\ell_{\tau\tau}+\operatorname{div}(\kappa\nabla\ell)+\Psi]v_\tau^2
-8\kappa\nabla\ell_\tau\cdot\nabla v\,v_\tau
+2\sum_{i,j=1}^2c_\ell^{ij}v_iv_j
-2\kappa\nabla\Psi\cdot\nabla v\,v
+\mathcal B_\ell v^2+2\Psi_\tau vv_\tau.
\end{aligned}
\]

\medskip
\noindent\textit{Part 2: The bulk estimate.}
The weight derivatives are
\begin{equation}\label{disk:carleman-weight-derivatives}
\ell_i=\lambda\mu\varphi d_i,\quad
\ell_{ij}=\lambda\mu^2\varphi d_id_j+\lambda\mu\varphi d_{ij},\quad
\ell_\tau=-2\lambda\mu\alpha\tau\varphi,\quad \ell_{\tau\tau}=4\lambda\mu^2\alpha^2\tau^2\varphi
-2\lambda\mu\alpha\varphi.
\end{equation}
The quadratic terms in \eqref{disk:identity} satisfy
\begin{equation}\label{disk:bulk-quadratic}
\begin{aligned}
&2[\ell_{\tau\tau}+\operatorname{div}(\kappa\nabla\ell)+\Psi]v_\tau^2 -8\kappa\nabla\ell_\tau\cdot\nabla v\,v_\tau
+2\sum_{i,j=1}^2c_\ell^{ij}v_iv_j\\
& =4(1-\alpha)\lambda\mu\varphi(v_\tau^2+\kappa|\nabla v|^2) +4\lambda\mu^2\varphi
(2\alpha\tau v_\tau+\kappa\nabla d\cdot\nabla v)^2 +2\lambda\mu\varphi
\bigl[\mathfrak Q_d(z,\nabla v)-4\kappa|\nabla v|^2\bigr]\\
&\quad\ge4(1-\alpha)\lambda\mu\varphi(v_\tau^2+\kappa|\nabla v|^2).
\end{aligned}
\end{equation}
Recall $E=\kappa|\nabla d|^2$ from \eqref{disk:weight-conditions}.
For the zeroth-order coefficient set $F\triangleq E-4\alpha^2\tau^2$.
The definitions in \eqref{disk:ABC}, combined with
\eqref{disk:carleman-weight-derivatives}, give the exact formulas
\[
A_\ell=\lambda^2\mu^2\varphi^2F-2(1+\alpha)\lambda\mu\varphi,
\]
and
\[
\mathcal B_\ell=4(1+\alpha)\lambda\mu\varphi A_\ell
+4\lambda\mu\alpha\tau\varphi\partial_\tau A_\ell
+2\lambda\mu\varphi\kappa\nabla d\cdot\nabla A_\ell.
\]
Since $\kappa_s=-2\kappa$,
\[
\kappa\nabla d\cdot\nabla E
=\mathfrak Q_d(z,\nabla d)\ge4E.
\]
For fixed $\mu$, the remainders denoted by $O(\lambda^2)$ below are
bounded in absolute value by $C\lambda^2$, with $C=C(\mu)$ independent
of $\lambda$, uniformly on $[-S,S]\times\overline G$.
The coefficient $\mathcal B_\ell$ therefore satisfies
\begin{equation}\label{disk:bulk-zero-coefficient}\begin{aligned}
\mathcal B_\ell
&=4\lambda^3\mu^4\varphi^3F^2 +2\lambda^3\mu^3\varphi^3
\bigl[2(1+\alpha)F-16\alpha^3\tau^2
+\kappa\nabla d\cdot\nabla E\bigr]\\
&\quad-4(1+\alpha)\lambda^2\mu^2\varphi^2
\bigl[2(1+\alpha)+\mu F\bigr]\\
&\ge4\lambda^3\mu^4\varphi^3F^2
+2\lambda^3\mu^3\varphi^3
\bigl[(6+2\alpha)E-(16\alpha^3+8(1+\alpha)\alpha^2)\tau^2\bigr]
+O(\lambda^2)\\
&=4\lambda^3\mu^4\varphi^3F^2
+2\lambda^3\mu^3\varphi^3
\bigl[4(1-\alpha)E+2(1+3\alpha)F\bigr]+O(\lambda^2)\\
&=4\lambda^3\mu^4\varphi^3
\left(F+\frac{1+3\alpha}{2\mu}\right)^2 +8(1-\alpha)\lambda^3\mu^3\varphi^3E
-(1+3\alpha)^2\lambda^3\mu^2\varphi^3+O(\lambda^2)\\
&\ge8(1-\alpha)\lambda^3\mu^3\varphi^3E
-(1+3\alpha)^2\lambda^3\mu^2\varphi^3+O(\lambda^2).
\end{aligned}\end{equation}
Since $E\ge4\min d>0$, choosing $\mu$ and then $\lambda$ large gives
$\mathcal B_\ell\ge c\lambda^3\mu^3\varphi^3$.
Moreover $|\nabla\Psi|+|\Psi_\tau|\le C(\mu)\lambda\varphi$.
Young's inequality gives\vspace{-2mm}
\begin{equation}\label{disk:bulk-young}
2|\kappa\nabla\Psi\cdot\nabla v\,v|+2|\Psi_\tau vv_\tau|
\le(1-\alpha)\lambda\mu\varphi
(v_\tau^2+\kappa|\nabla v|^2)+C(\mu)\lambda\varphi v^2.\vspace{-2mm}
\end{equation}
For sufficiently large $\lambda$,
$C(\mu)\lambda\varphi\le(c/2)\lambda^3\mu^3\varphi^3$,
so \eqref{disk:bulk-quadratic}, \eqref{disk:bulk-zero-coefficient},
and \eqref{disk:bulk-young} yield\vspace{-2mm}
\begin{equation}\label{disk:bulk}
\theta^2|P_Gu|^2+2\operatorname{div}\mathbf J+2\mathcal M_\tau
\ge c\lambda\mu\varphi(v_\tau^2+\kappa|\nabla v|^2)
+c\lambda^3\mu^3\varphi^3v^2.\vspace{-2mm}
\end{equation}

\noindent\textit{Part 3: Boundary fluxes.}
On the lateral boundary, $v=v_\tau=0$ and
$\nabla v=\nu\partial_\nu v$, so\vspace{-2mm}
\[
\mathbf J\cdot\nu
=\lambda\mu\varphi\kappa(\kappa\nabla d\cdot\nu)
|\partial_\nu v|^2.\vspace{-2mm}
\]
This vanishes near $\Gamma_+$ and is nonpositive elsewhere.
At the left time face,\vspace{-2mm}
\begin{align*}
\mathcal M(-S)
&={}2\lambda\mu\alpha S\varphi
(v_\tau^2+\kappa|\nabla v|^2)
-2\lambda\mu\varphi(\kappa\nabla d\cdot\nabla v)v_\tau
+\Psi vv_\tau\\
&\quad+2\lambda^3\mu^3\alpha S\varphi^3(4\alpha^2S^2-E)v^2
+4(1+\alpha)\lambda^2\mu^2\alpha S\varphi^2v^2.
\end{align*}
Here $2|(\kappa\nabla d\cdot\nabla v)v_\tau|
\le\sqrt E(v_\tau^2+\kappa|\nabla v|^2)$ and
$|\Psi vv_\tau|\le v_\tau^2+C(\mu)\lambda^2v^2$. The strict margins\vspace{-2mm}
\[
2\alpha S-\sqrt E\ge2(\alpha S-R_\sharp)>0,\qquad
4\alpha^2S^2-E\ge4(\alpha^2S^2-R_\sharp^2)>0\vspace{-2mm}
\]
hold on the spatial support. At either time face,\vspace{-2mm}
\[
\begin{aligned}
-\sigma\mathcal M(\sigma S)
&\ge\bigl[\lambda\mu\varphi(2\alpha S-\sqrt E)-1\bigr]
(v_\tau^2+\kappa|\nabla v|^2)\\
&\quad+\bigl[2\lambda^3\mu^3\alpha S\varphi^3
(4\alpha^2S^2-E)-C(\mu)\lambda^2\bigr]v^2 \ge0
\qquad(\sigma\in\{-1,1\})
\end{aligned}
\]
for large $\lambda$. Consequently, \eqref{disk:bulk} implies
\begin{equation}\label{disk:integrated-bulk}
\begin{aligned}
&c\lambda\mu\int_{-S}^S\!\int_G
\varphi(v_\tau^2+\kappa|\nabla v|^2)
+c\lambda^3\mu^3\int_{-S}^S\!\int_G\varphi^3v^2\\
&\le\int_{-S}^S\!\int_G\theta^2|P_Gu|^2
+2\int_{-S}^S\!\int_{\partial G}\mathbf J\cdot\nu
+2\int_G[\mathcal M(S)-\mathcal M(-S)] \le\int_{-S}^S\!\int_G\theta^2|P_Gu|^2.
\end{aligned}
\end{equation}

\medskip
\noindent\textit{Part 4: From $v$ back to $u$.}
For each time or spatial derivative,
\begin{equation}\label{disk:weighted-derivative-conversion}
\theta\partial u=\partial v-(\partial\ell)v,\qquad
\theta^2|\partial u|^2
\le2|\partial v|^2+C\lambda^2\mu^2\varphi^2v^2.
\end{equation}
Since $0<\min\varphi\le\max\varphi<\infty$ for fixed $\mu$,
\eqref{disk:integrated-bulk} and \eqref{disk:weighted-derivative-conversion}
prove \eqref{disk:carleman-estimate}.
\end{proof}

\begin{proof}[Proof of Theorem~\ref{thm:disk}]
We first localize the solution in space and time. A weak Carleman
estimate controls its displacement on a central time interval up to an
exponentially small global remainder; a weak energy estimate controls
the remainder. The proof has four steps: choose the cutoffs, apply the
weak Carleman estimate, control the displacement near the middle time,
and absorb the remainder.

\medskip
\noindent\textbf{Step 1. Choose the time interval and cutoffs.}
Set $T_*\triangleq2R_G$ and fix $T>T_*$. Choose $R_G<S<T/2$ and
$\alpha\in(R_G/S,1)$, and set $\tau\triangleq t-T/2$.
Fix a sufficiently large $\mu$ in Proposition~\ref{disk:carleman}.
Choose $0<s_0<s_1<S$ such that\vspace{-2mm}
\[
\alpha s_0^2<\frac12\min_{\overline G}d,\qquad
\alpha s_1^2>R_G^2.\vspace{-2mm}
\]
The first inequality makes the weight strictly positive on the central
time slab $|\tau|\le s_0$; the second makes it strictly negative on the
outer slabs $s_1\le|\tau|\le S$. This gap is what produces the
exponentially small remainder below.

Take $\eta\in C_c^\infty(-S,S)$ with $0\le\eta\le1$ and
$\eta=1$ on $[-s_1,s_1]$. Choose
$\chi\in C^\infty(\overline G;[0,1])$ satisfying
\[
\chi(z)=
\begin{cases}
0,&\operatorname{dist}(z,\pi^{-1}(\Sigma_{\rm art}))\le2\varepsilon_1,\\
1,&\operatorname{dist}(z,\pi^{-1}(\Sigma_{\rm art}))\ge3\varepsilon_1.
\end{cases}
\]
The two closed sets are separated by at least $\varepsilon_1$, so such
a smooth cutoff exists. The choice \eqref{disk:artificial-collar-choice}
of $\varepsilon_1$ gives\vspace{-2mm}
\[
\begin{aligned}
&\chi=0\quad\hbox{on }\omega_1,\qquad
\chi=1\quad\hbox{on }G\setminus\pi^{-1}(U_*),\\
&\operatorname{supp}_{\overline G}\nabla\chi
\cup\operatorname{supp}_{\overline G}\operatorname{div}(\kappa\nabla\chi)
\subset\pi^{-1}(U_*).
\end{aligned}
\]
Thus the spatial commutators are supported in the observed collar,
while the solution is unmodified outside it.

\medskip
\noindent\textbf{Step 2. A weighted estimate for weak solutions.}
We apply the weak Carleman estimate of \cite{FuLiaoLu} to the
transformed operator $P_G-B$. Take the scalar case of
\cite[Theorem~4.1, equation~(4.1)]{FuLiaoLu}, with zero potential
$\widetilde q=0$, coefficients $h^{jk}=\kappa\delta_{jk}$, and spatial
weight $d$. The second, potential-dependent term on its right-hand side
vanishes, giving constants $C>0$ and $\lambda_0>0$ such that
\begin{equation}\label{disk:weak-carleman-free}
\lambda\|\theta w\|_{L^2((-S,S)\times G)}^2
\le C\|\theta P_Gw\|_{H^{-1}((-S,S)\times G)}^2,
\qquad \lambda\ge\lambda_0.
\end{equation}
Here $w\in C([-S,S];L^2(G))$ has homogeneous lateral Dirichlet data in
the transposition sense, $w(-S)=w(S)=0$, $w=0$ on
$(-S,S)\times\omega_1$, and $P_Gw\in H^{-1}((-S,S)\times G)$; these are
the conditions in \cite[equation~(3.2)]{FuLiaoLu}. The negative space is
the dual of $H_0^1((-S,S)\times G)$, and the right-hand side is the
norm of the product $\theta P_Gw$.

In that paper's notation, take $\Omega=G$, $\phi=d$, and
$\widetilde T=S$. The matrix $h^{jk}=\kappa\delta_{jk}$ is smooth,
symmetric, and uniformly elliptic. The first inequality in
\eqref{disk:weight-conditions} gives $s_0=4$ in Condition~1.1(i); the
second gives both the absence of critical points and equation~(1.10)
there. Moreover, $\Gamma_0=\Gamma_+$ and $\widetilde R_1=R_G$ in its
equation~(2.2), while $S>R_G$ allows the choice $\alpha S>R_G$.
Proposition~\ref{disk:carleman} verifies the corresponding strong
estimate, so the auxiliary optimal control argument in its Section~3
and proof of Theorem~4.1 applies here.

Since $L^2$ embeds continuously into this negative space,
\eqref{disk:weak-carleman-free} and $\|B\|_\infty\le r_0^{-2}+L$ imply
\begin{equation}\label{disk:weak-carleman-potential}
\lambda\|\theta w\|_2^2
\le C\|\theta(P_G-B)w\|_{H^{-1}((-S,S)\times G)}^2
\end{equation}
whenever $\lambda\ge C\cdot(1+\|B\|_\infty^2)$ and $\lambda\ge\lambda_0$.
Indeed, these choices of $\lambda$ give
\[
\begin{aligned}
\lambda\|\theta w\|_2^2
&\le C\|\theta(P_G-B)w\|_{H^{-1}((-S,S)\times G)}^2
+C\|B\|_\infty^2\|\theta w\|_2^2\\
&\le C\|\theta(P_G-B)w\|_{H^{-1}((-S,S)\times G)}^2
+\frac\lambda2\|\theta w\|_2^2.
\end{aligned}
\]
All constants here are independent of $V$ within the stated bound.

\medskip
\noindent\textbf{Step 3. Control the displacement near the middle time.}
No zero trace of $p$ is assumed on $\Sigma_{\rm art}$, so we localize
before applying \eqref{disk:weak-carleman-potential}. Set
$\widetilde p\triangleq\mathcal U(p|_{G_{\rm ph}})$ and
$Z\triangleq\eta\chi\widetilde p$. The localized transposition equation is
\begin{equation}\label{disk:localized-equation}
(P_G-B)Z
=2\partial_\tau(\eta'\chi\widetilde p)-\eta''\chi\widetilde p
+\eta\operatorname{div}(\kappa\nabla\chi)\widetilde p
-2\operatorname{div}(\eta\kappa\widetilde p\nabla\chi).
\end{equation}
The spatial cutoff gives $Z$ homogeneous lateral data, including on
$\pi^{-1}(\Sigma_{\rm art})$. The temporal cutoff gives $Z(\pm S)=0$,
and the remaining hypotheses of \eqref{disk:weak-carleman-potential}
hold because $\widetilde p\in C([-S,S];L^2(G))$ and the source in
\eqref{disk:localized-equation} is in spacetime $H^{-1}$.
For any $L^2$ function $F$, the weighted distributional derivative
satisfies
\[
\theta\partial_jF=\partial_j(\theta F)-\partial_j\theta\,F,
\qquad j\in\{\tau,s,\vartheta\}.
\]
Since $|\partial\theta|\le C\lambda\theta$ on the fixed cylinder,
\[
\begin{aligned}
\|\theta\partial_jF\|_{H^{-1}}
&\le\|\partial_j(\theta F)\|_{H^{-1}}
+\|(\partial_j\theta)F\|_{H^{-1}}
\le C\|\theta F\|_2+C\lambda\|\theta F\|_2
\le C\lambda\|\theta F\|_2
\qquad(\lambda\ge1).
\end{aligned}
\]
The spatial source terms are supported in $\pi^{-1}(U_*)$, while the
temporal source terms have $|\tau|\ge s_1$. Applying
\eqref{disk:weak-carleman-potential} to $Z$ therefore gives
\begin{equation}\label{disk:localized-carleman}
\begin{aligned}
&\lambda\|\theta Z\|_2^2\le C\lambda^2\Bigg[
\int_{-S}^S\!\int_{G\cap\pi^{-1}(U_*)}
\theta^2|\widetilde p|^2\,\mathrm{d}z\,\mathrm{d}\tau +
\int_{s_1\le|\tau|\le S}\!\int_{G\setminus\omega_1}
\theta^2|\widetilde p|^2\,\mathrm{d}z\,\mathrm{d}\tau\Bigg].
\end{aligned}
\end{equation}
By \eqref{disk:weight-conditions} and the definition of $R_G$, we have
$d(z)\le E(z)/4\le R_G^2$ on $G\setminus\omega_1$.
Thus, on the temporal source region and the central cylinder, respectively,
\[
\begin{alignedat}{2}
&d(z)-\alpha\tau^2\le R_G^2-\alpha s_1^2<0
&&\qquad(s_1\le|\tau|\le S,\ z\in G\setminus\omega_1),\\
&d(z)-\alpha\tau^2\ge\tfrac12\min_{\overline G}d>0
&&\qquad(|\tau|\le s_0,\ z\in G).
\end{alignedat}
\]
The strict weight gap is
\[
\begin{aligned}
\min_{\substack{|\tau|\le s_0\\z\in\overline G}}\varphi
&\ge e^{\frac{\mu}{2}\min_{\overline G}d}>1,\quad
\max_{\substack{s_1\le|\tau|\le S\\z\in\overline{G\setminus\omega_1}}}
\varphi
 \le e^{\mu(R_G^2-\alpha s_1^2)}<1.
\end{aligned}
\]
On $|\tau|\le s_0$, $Z=\chi\widetilde p$, with $\chi=1$ on
$G\setminus\pi^{-1}(U_*)$. Thus \eqref{disk:localized-carleman} gives\vspace{-2mm}
\[
\begin{aligned}
&\int_{-s_0}^{s_0}\!\int_{G\setminus\pi^{-1}(U_*)}
|\widetilde p|^2 \le
e^{-2\lambda e^{\frac{\mu}{2}\min_{\overline G}d}}
\|\theta Z\|_2^2\\
&\le C\lambda e^{C\lambda}
\int_{-S}^S\!\int_{G\cap\pi^{-1}(U_*)}|\widetilde p|^2 +C\lambda
\exp\!\left\{-2\lambda\left[
e^{\frac{\mu}{2}\min_{\overline G}d}
-e^{\mu(R_G^2-\alpha s_1^2)}\right]\right\}
\int_{-S}^S\!\int_G|\widetilde p|^2.
\end{aligned}\vspace{-2mm}
\]
Let $I\triangleq(T/2-s_0,T/2+s_0)$. Since
$\mathbb D\setminus G_{\rm ph}\subset U_*\cap\mathbb D\subset\omega$
and $\mathcal U$ is unitary,\vspace{-2mm}
\[
\begin{aligned}
&\int_I\|p(t)\|_2^2\,\mathrm{d}t
\le\int_I\!\int_\omega|p|^2\,\mathrm{d}x\,\mathrm{d}t +\int_{-s_0}^{s_0}\!\int_{G\setminus\pi^{-1}(U_*)}
|\widetilde p|^2\,\mathrm{d}z\,\mathrm{d}\tau,\\
&\int_{-S}^S\!\int_{G\cap\pi^{-1}(U_*)}|\widetilde p|^2
\le\int_0^T\!\int_\omega|p|^2\,\mathrm{d}x\,\mathrm{d}t,\\
&\int_{-S}^S\!\int_G|\widetilde p|^2
\le\int_0^T\|p(t)\|_2^2\,\mathrm{d}t.
\end{aligned}\vspace{-2mm}
\]
Consequently, for constants $c=c(T)>0$ and $C=C(T)>0$,\vspace{-2mm}
\begin{equation}\label{disk:middle-displacement}
\int_I\|p(t)\|_2^2\,\mathrm{d}t
\le Ce^{C\lambda}\int_0^T\!\int_\omega|p|^2\,\mathrm{d}x\,\mathrm{d}t
+C\lambda e^{-c\lambda}\int_0^T\|p(t)\|_2^2\,\mathrm{d}t.\vspace{-2mm}
\end{equation}
This estimates the displacement on a central time interval: the
observation contributes a factor $e^{C\lambda}$, while the global
remainder is exponentially small.

\smallskip
\noindent\textbf{Step 4. Control the state norm and absorb the remainder.}
It remains to bound the initial state by the displacement on $I$.
Recall $A_{\mathbb D}=-\Delta_D$ and define
\[
\mathcal E(t)\triangleq\|p(t)\|_2^2+
\|A_{\mathbb D}^{-1/2}p_t(t)\|_2^2.
\]
For smooth solutions,
\begin{equation}\label{disk:weak-energy-derivative}
\begin{aligned}
\mathcal E'(t)
&=2\operatorname{Re}(p_t,p)
+2\operatorname{Re}(A_{\mathbb D}^{-1}p_t,p_{tt})
=2\operatorname{Re}(p_t,p)
-2\operatorname{Re}(A_{\mathbb D}^{-1}p_t,A_{\mathbb D}p+Vp)
=-2\operatorname{Re}(Vp,A_{\mathbb D}^{-1}p_t),\\
|\mathcal E'(t)|
&\le2L\|A_{\mathbb D}^{-1/2}\|_{L^2\to L^2}
\|p\|_2\|A_{\mathbb D}^{-1/2}p_t\|_2
\le CL\mathcal E(t).
\end{aligned}
\end{equation}
Moreover, in distributions in time,
\begin{equation}\label{disk:weak-energy-second-derivative}
\frac{\mathrm{d}^2}{\mathrm{d}t^2}
\|A_{\mathbb D}^{-1/2}p\|_2^2
=2\|A_{\mathbb D}^{-1/2}p_t\|_2^2
-2\|p\|_2^2-2\operatorname{Re}(Vp,A_{\mathbb D}^{-1}p).
\end{equation}
Choose $\xi\in C_c^\infty(I)$ independently of $V$, with
$0\le\xi\le1$ and $\xi=1$ on $[T/2-s_0/2,T/2+s_0/2]$.
The compact support of $\xi$ in $I$ implies\vspace{-2mm}
\[
\begin{aligned}
\int_I\xi\,\frac{\mathrm{d}^2}{\mathrm{d}t^2}
\|A_{\mathbb D}^{-1/2}p\|_2^2\,\mathrm{d}t
&=-\int_I\xi'\,\frac{\mathrm{d}}{\mathrm{d}t}
\|A_{\mathbb D}^{-1/2}p\|_2^2\,\mathrm{d}t
=\int_I\xi''\|A_{\mathbb D}^{-1/2}p\|_2^2\,\mathrm{d}t.
\end{aligned}\vspace{-2mm}
\]
Therefore,\vspace{-2mm}
\[
\begin{aligned}
\int_{T/2-s_0/2}^{T/2+s_0/2}\mathcal E(t)\,\mathrm{d}t
&\le\int_I\xi\mathcal E\,\mathrm{d}t =2\int_I\xi\|p\|_2^2\,\mathrm{d}t
+\operatorname{Re}\int_I\xi(Vp,A_{\mathbb D}^{-1}p)\,\mathrm{d}t +\frac12\int_I\xi''
\|A_{\mathbb D}^{-1/2}p\|_2^2\,\mathrm{d}t\\
& \le\left[2+\|A_{\mathbb D}^{-1}\|_{L^2\to L^2}
\left(L+\tfrac12\|\xi''\|_\infty\right)\right]
\int_I\|p\|_2^2\,\mathrm{d}t\\
& \le C (1+L)\int_I\|p\|_2^2\,\mathrm{d}t.
\end{aligned}\vspace{-2mm}
\]
Here $C=C(\xi)>0$ depends only on the fixed cutoff $\xi$ and the disk,
and is independent of $V$ and $L$.
Gronwall's inequality ensures
\[
e^{-CL|t-s|}\mathcal E(s)\le\mathcal E(t)
\le e^{CL|t-s|}\mathcal E(s)
\qquad(s,t\in[0,T]).
\]
In particular,\vspace{-2mm}
\[
s_0e^{-CLT}\mathcal E(0)
\le\int_{T/2-s_0/2}^{T/2+s_0/2}\mathcal E(t)\,\mathrm{d}t
\le C(1+L)\int_I\|p(t)\|_2^2\,\mathrm{d}t.\vspace{-2mm}
\]
Consequently,\vspace{-3mm}
\begin{equation}\label{disk:weak-energy-state-bounds}
\mathcal E(0)\le C(T,L)\int_I\|p(t)\|_2^2\,\mathrm{d}t,
\qquad
\int_0^T\|p(t)\|_2^2\,\mathrm{d}t\le C(T,L)\mathcal E(0).\vspace{-2mm}
\end{equation}
By approximation, \eqref{disk:weak-energy-derivative},
\eqref{disk:weak-energy-second-derivative}, and their consequence
\eqref{disk:weak-energy-state-bounds} also hold for weak solutions.

\smallskip\noindent\textit{Absorbing the remainder.}
Combining \eqref{disk:middle-displacement} and
\eqref{disk:weak-energy-state-bounds} gives\vspace{-2mm}
\[
\mathcal E(0)
\le C(T,L)e^{C\lambda}\int_0^T\!\int_\omega|p|^2\,\mathrm{d}x\,\mathrm{d}t
+C(T,L)\lambda e^{-c\lambda}\mathcal E(0).\vspace{-2mm}
\]
Choose $\lambda$ large enough to satisfy the weak Carleman threshold,
the potential absorption condition, and
$C(T,L)\lambda e^{-c\lambda}\le1/2$. Then\vspace{-2mm}
\[
\begin{aligned}
\tfrac12\mathcal E(0)
&\le\bigl[1-C(T,L)\lambda e^{-c\lambda}\bigr]\mathcal E(0)
\le C(T,L)e^{C\lambda}
\int_0^T\!\int_\omega|p|^2\,\mathrm{d}x\,\mathrm{d}t,
\end{aligned}\vspace{-2mm}
\]
which proves \eqref{eq:disk-observability}.
\end{proof}

\section{A real-valued local counterexample to unique continuation}\label{ce:front-section}

To prove Theorem~\ref{ce:front}, fix $Z>0$. We construct the solution
first for $y>0$ and show at the end that it extends smoothly by zero
across $y=0$. We use the geometric-optics framework of
Alinhac--Baouendi \cite[pp.~XXII.4--XXII.6]{AlinhacBaouendi1979}; see
also Guisset--Shao \cite[Sections~1.4 and~2]{GuissetShao2024}. The
periodic correction below gives the solution and potential a common
time period.

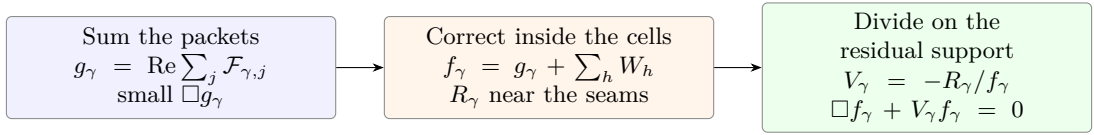
\begin{figure}[!htbp]
\centering
\begin{tikzpicture}[font=\footnotesize,>={Stealth[length=1.5mm]}]
\node[draw=black!45,rounded corners=2pt,fill=blue!6,text width=4.1cm,
align=center,minimum height=1.35cm] (packets) at (0,0)
{Sum the packets\\$g_\gamma=\operatorname{Re}\sum_j\mathcal F_{\gamma,j}$\\small $\Box g_\gamma$};
\node[draw=black!45,rounded corners=2pt,fill=orange!8,text width=4.1cm,
align=center,minimum height=1.35cm] (correction) at (5,0)
{Correct inside the cells\\$f_\gamma=g_\gamma+\sum_hW_h$\\$R_\gamma$ near the seams};
\node[draw=black!45,rounded corners=2pt,fill=green!7,text width=4.1cm,
align=center,minimum height=1.35cm] (potential) at (10,0)
{Divide on the residual support\\$V_\gamma=-R_\gamma/f_\gamma$\\$\Box f_\gamma+V_\gamma f_\gamma=0$};
\draw[->] (packets.east)--(correction.west);
\draw[->] (correction.east)--(potential.west);
\end{tikzpicture}
\caption{The three stages of the local construction. The positivity needed
for the last division is proved in Lemma~\ref{ce:division-region}.}
\label{fig:front-construction-roadmap}
\end{figure}

The construction has three stages, summarized in
Figure~\ref{fig:front-construction-roadmap}. For each amplitude weight
$\gamma\in\{0,24\log2\}$, we first sum oscillatory packets to obtain an
approximate solution $g_\gamma$ of the free wave equation. We then add a
periodic causal correction $\sum_{h\ge5}W_h$, supported in narrow phase
bands, to cancel the cell residuals and produce a corrected field
$f_\gamma$. Finally, on the support of the residual, where $f_\gamma$
is bounded below by a fixed fraction of its envelope, we divide by
$f_\gamma$ to define the potential $V_\gamma$.\vspace{-2mm}
\[
\begin{aligned}
g_\gamma&=\operatorname{Re}\sum_j\mathcal F_{\gamma,j}
&&\text{(sum the packets)},\\
f_\gamma&=g_\gamma+\sum_{h\ge5}W_h
&&\text{(cancel the cell residuals)},\\
V_\gamma&=-R_\gamma/f_\gamma,\qquad R_\gamma=\Box f_\gamma
&&\text{(divide on the residual support)}.
\end{aligned}\vspace{-2mm}
\]
After proving $f_\gamma>0$ in the last step, we extend $V_\gamma$ by
zero off the residual support and verify smoothness and flatness at
$y=0$.
The normal width $Y_*$ will be decreased finitely many times, always
by a choice common to the two values of $\gamma$. Define\vspace{-2mm}
\begin{equation}\label{ce:front-parameters}
\begin{aligned}
&\beta\triangleq2^{1/16},\quad
\lambda_j\triangleq2^j,\quad K_j\triangleq\lfloor j^2\beta^j\rfloor,
\quad j\in\mathbb Z,\ j\ge4,\quad u\triangleq t+z,\quad v\triangleq t-z.
\end{aligned}\vspace{-2mm}
\end{equation}
Here $j$ labels an oscillatory packet, $\lambda_j$ is its frequency,
and $K_j$ is the number of recursive correction orders retained in it.
The factor $\beta$ relates the packet index to its normal height:
a packet with index $j$ is supported where $y$ is comparable to
$\beta^{-j}$. The integer frequencies will give the common period $2\pi$.

The choice of truncation order has the sizes\vspace{-3mm}
\[
K_j\asymp j^2\beta^j,\qquad
jK_j\asymp j^3\beta^j,\qquad
\lambda_j^{-K_j}=\exp(-jK_j\log2).\vspace{-3mm}
\]
On $y\asymp\beta^{-j}$, this factor is smaller than every fixed power
$y^N$ and still tends to zero after multiplication by a fixed $y^{-N}$.
The choice of $K_j$ gives this decay; including the growth of the
coefficient derivatives yields the residual bound in
Lemma~\ref{ce:sum-residual}.

The coordinates $u$ and $v$ satisfy\vspace{-3mm}
\[
t=(u+v)/2,\qquad z=(u-v)/2,\qquad
\partial_u=\tfrac12(\partial_t+\partial_z),\qquad
\partial_v=\tfrac12(\partial_t-\partial_z),\vspace{-3mm}
\]
so $\Box=4\partial_u\partial_v-\partial_y^2$ and
$\partial_uz=1/2$, $\partial_vz=-1/2$.
Initially set $Y_*\triangleq\beta^{-7}$. The packet sums are fixed
from the outset. The periodic corrections are fixed after the first
decrease of $Y_*$. Every later decrease restricts the same functions
to a smaller cylinder. The final width in
\eqref{ce:final-front-width} is independent of derivative order.

\vspace{-2mm}

\subsection{Construction of the solution}\label{ce:front-construction}

The construction proceeds in three steps: choose amplitudes that decay
toward $y=0$, cancel the residual recursively by a finite transport
expansion, and then add the periodic cell correction.

\medskip\noindent\textbf{Packet amplitudes and normal bands.}
The profile functions are
\begin{equation}\label{ce:cutoff-functions}
\begin{aligned}
&a(s)\triangleq\begin{cases}e^{-1/s},&s>0,\\0,&s\le0,\end{cases}
\quad b(s)\triangleq\frac{a(s)}{a(s)+a(1-s)},\\
&\rho(s)\triangleq
\begin{cases}
b((s-\beta^{-1})/(1-\beta^{-1})),&\beta^{-1}<s\le1,\\
b((\beta-s)/(\beta-1)),&1\le s<\beta,\\
0,&\text{otherwise}.
\end{cases}
\end{aligned}
\end{equation}
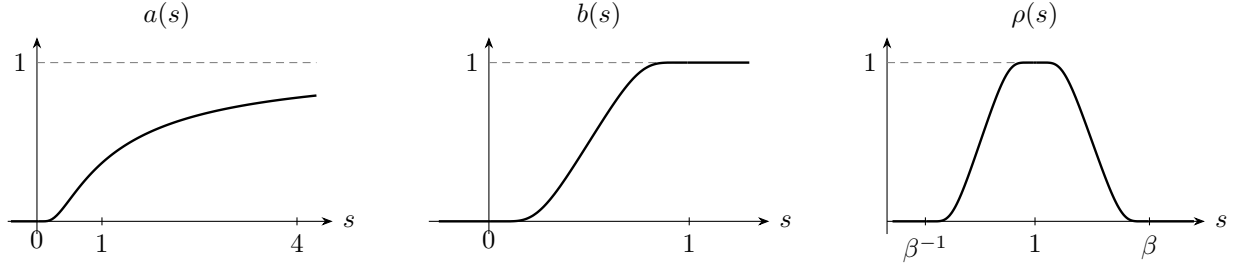
\begin{figure}[H]
\centering
\begin{minipage}[t]{.31\textwidth}
\centering
\begin{tikzpicture}[baseline=0pt,x=.86cm,y=2.1cm,font=\small,
>={Stealth[length=1.5mm]}]
\node at (2,1.30) {$a(s)$};
\draw[->,thin] (-.45,0)--(4.55,0) node[right] {$s$};
\draw[->,thin] (0,-.08)--(0,1.16);
\draw[gray,densely dashed] (0,1)--(4.3,1);
\node[left] at (0,1) {$1$};
\node[below] at (0,0) {$0$};
\draw (1,.025)--(1,-.025) node[below] {$1$};
\draw (4,.025)--(4,-.025) node[below] {$4$};
\draw[line width=.9pt] (-.4,0)--(.02,0);
\draw[line width=.9pt,domain=.02:4.3,samples=121]
plot (\x,{exp(-1/\x)});
\end{tikzpicture}
\end{minipage}\hfill
\begin{minipage}[t]{.31\textwidth}
\centering
\begin{tikzpicture}[baseline=0pt,x=2.65cm,y=2.1cm,font=\small,
>={Stealth[length=1.5mm]}]
\node at (.55,1.30) {$b(s)$};
\draw[->,thin] (-.30,0)--(1.40,0) node[right] {$s$};
\draw[->,thin] (0,-.08)--(0,1.16);
\draw[gray,densely dashed] (0,1)--(1,1);
\node[left] at (0,1) {$1$};
\node[below] at (0,0) {$0$};
\draw (1,.025)--(1,-.025) node[below] {$1$};
\draw[line width=.9pt] (-.25,0)--(.01,0);
\draw[line width=.9pt,domain=.01:.99,samples=161]
plot (\x,{exp(-1/\x)/(exp(-1/\x)+exp(-1/(1-\x)))});
\draw[line width=.9pt] (.99,1)--(1.30,1);
\end{tikzpicture}
\end{minipage}\hfill
\begin{minipage}[t]{.31\textwidth}
\centering
\begin{tikzpicture}[baseline=0pt,x=1.45cm,y=2.1cm,font=\small,
>={Stealth[length=1.5mm]}]
\pgfmathsetmacro{\cutoffbeta}{2^(1/16)}
\node at (1,1.30) {$\rho(s)$};
\draw[->,thin] (-.35,0)--(2.55,0) node[right] {$s$};
\draw[->,thin] (-.35,-.08)--(-.35,1.16);
\draw[gray,densely dashed] (-.35,1)--(1,1);
\node[left] at (-.35,1) {$1$};
\draw (0,.025)--(0,-.025) node[below] {$\beta^{-1}$};
\draw (1,.025)--(1,-.025) node[below] {$1$};
\draw ({1+\cutoffbeta},.025)--({1+\cutoffbeta},-.025)
node[below] {$\beta$};
\draw[line width=.9pt] (-.30,0)--(.01,0);
\draw[line width=.9pt,domain=.01:.99,samples=161]
plot (\x,{exp(-1/\x)/(exp(-1/\x)+exp(-1/(1-\x)))});
\draw[line width=.9pt] (.99,1)--({1+.01*\cutoffbeta},1);
\draw[line width=.9pt,domain=.01:.99,samples=161]
plot ({1+\cutoffbeta*\x},
{exp(-1/(1-\x))/(exp(-1/\x)+exp(-1/(1-\x)))});
\draw[line width=.9pt] ({1+.99*\cutoffbeta},0)--(2.45,0);
\end{tikzpicture}
\end{minipage}
\caption{Graphs of the functions in \eqref{ce:cutoff-functions}.
The horizontal axis in the third panel is magnified near $s=1$;
here $\beta=2^{1/16}$. Dashed lines indicate the level $1$.}
\label{fig:cutoff-functions}\vspace{-2mm}
\end{figure}
Since $K_j$ grows with the packet index, the estimates below depend on
the derivative order.
\begin{lemma}[Gevrey bounds of order two]\label{ce:cutoff-bounds}
The functions $a,b,\rho\in C^\infty(\R)$. There exists a
constant $C>0$, independent of $r$, such that\vspace{-2mm}
\[
\|h^{(r)}\|_{L^\infty(\R)}\le C^{r+1}(r!)^2,
\qquad r\in\mathbb N_0,\quad h\in\{a,b,\rho\}.\vspace{-2mm}
\]
Moreover, $b(s)+b(1-s)=1$ for every $s\in\R$.
\end{lemma}
For integer packet indices $j\ge4$ and $y>0$, define
\begin{equation}\label{ce:front-amplitudes}
A_{\gamma,j}(y)\triangleq
e^{-1/y-j^2-\gamma j}\rho(\beta^jy),\quad
E_\gamma(y)\triangleq\sum_{j\ge4}A_{\gamma,j}(y).
\end{equation}
The coefficients $A_{\gamma,j}$ are the nonnegative amplitude weights of
the packets. Their sum $E_\gamma$ is the envelope used to measure the
size of the oscillatory function.

\begin{figure}[H]
\centering
\begin{tikzpicture}[x=1.14cm,y=1.05cm,font=\small,
>={Stealth[length=1.6mm]},line join=round]
\definecolor{amplitudeblue}{RGB}{40,94,138}
\definecolor{amplituderust}{RGB}{166,86,43}
\definecolor{amplitudegreen}{RGB}{43,117,94}
\pgfmathsetmacro{\bandbeta}{2^(1/16)}
\pgfmathsetmacro{\bandleft}{22*(1/\bandbeta-.74)}
\pgfmathsetmacro{\bandright}{22*(\bandbeta-.74)}
\fill[amplitudeblue!8] (\bandleft,0) rectangle (\bandright,3.65);
\draw[amplitudeblue!45,densely dashed]
(\bandleft,0)--(\bandleft,3.65)
(\bandright,0)--(\bandright,3.65);
\draw[->,thin] (-.55,0)--(11.7,0) node[right] {$y$};
\draw[->,thin] (-.55,-.08)--(-.55,4.25)
node[above] {$A_{\gamma,j}$};
\node[below left] at (-.55,0) {$0$};
\draw[white,line width=3pt] (.12,0)--(.36,0);
\draw[thin] (.10,-.07)--(.20,.07) (.28,-.07)--(.38,.07);
\node at (.80,.30) {$\cdots$};
\foreach \offset in {-4,-3,-2,-1,0,1,2,3,4,5}{
\pgfmathsetmacro{\amplitudecentre}{\bandbeta^(-\offset)}
\pgfmathsetmacro{\amplitudeheight}{1.5*exp(-.19*\offset)}
\def\amplitudecolor{gray!65}
\def\amplitudewidth{.7pt}
\ifnum\offset=-1
\def\amplitudecolor{amplitudegreen}
\def\amplitudewidth{1.1pt}
\fi
\ifnum\offset=0
\def\amplitudecolor{amplituderust}
\def\amplitudewidth{1.1pt}
\fi
\ifnum\offset=1
\def\amplitudecolor{amplitudeblue}
\def\amplitudewidth{1.1pt}
\fi
\draw[\amplitudecolor,line width=\amplitudewidth,
domain=.01:.99,samples=81]
plot ({22*(\amplitudecentre*(1/\bandbeta+(1-1/\bandbeta)*\x)-.74)},
{\amplitudeheight*exp(-1/\x)/(exp(-1/\x)+exp(-1/(1-\x)))});
\draw[\amplitudecolor,line width=\amplitudewidth,
domain=.01:.99,samples=81]
plot ({22*(\amplitudecentre*(1+(\bandbeta-1)*\x)-.74)},
{\amplitudeheight*exp(-1/(1-\x))/(exp(-1/\x)+exp(-1/(1-\x)))});
}
\node[amplitudeblue] at (3.55,2.30) {$A_{\gamma,m+1}$};
\draw[amplitudeblue,thin] (4.05,2.10)--(\bandleft,1.28);
\node[amplituderust] at (5.40,2.95) {$A_{\gamma,m}$};
\draw[amplituderust,thin] (5.48,2.73)--(5.72,1.55);
\node[amplitudegreen] at (7.55,3.45) {$A_{\gamma,m-1}$};
\draw[amplitudegreen,thin] (7.18,3.23)--(\bandright,1.86);
\node at (5.55,4.05)
{$E_\gamma=A_{\gamma,m+1}+A_{\gamma,m}+A_{\gamma,m-1}\quad\text{on }J_m$};
\draw[thin] (\bandleft,.06)--(\bandleft,-.06)
node[below] {$\beta^{-(m+1)}$};
\draw[thin] (\bandright,.06)--(\bandright,-.06)
node[below] {$\beta^{-(m-1)}$};
\draw[amplitudeblue!70,thin]
(\bandleft,-.65)--(\bandleft,-.77)--(\bandright,-.77)--(\bandright,-.65);
\node[below,amplitudeblue] at ({(\bandleft+\bandright)/2},-.77) {$J_m$};
\draw[->,gray!80] (2.65,1.20)--(1.10,.68);
\node[above,align=center,font=\footnotesize] at (1.75,1.35)
{Decreasing peak heights\\toward $y=0$};
\end{tikzpicture}
\caption{Schematic graphs of several functions $A_{\gamma,j}$ for a large index $m$.
Peak heights are shown schematically, and the horizontal axis is broken near the origin.
On the shaded interval $J_m$, only $A_{\gamma,m-1}$, $A_{\gamma,m}$,
and $A_{\gamma,m+1}$ can be nonzero; at any fixed $y\in J_m$,
at most two of them contribute to $E_\gamma(y)$.
The omitted functions continue toward $y=0$ with decreasing peak heights.}
\label{fig:packet-amplitudes}
\end{figure}
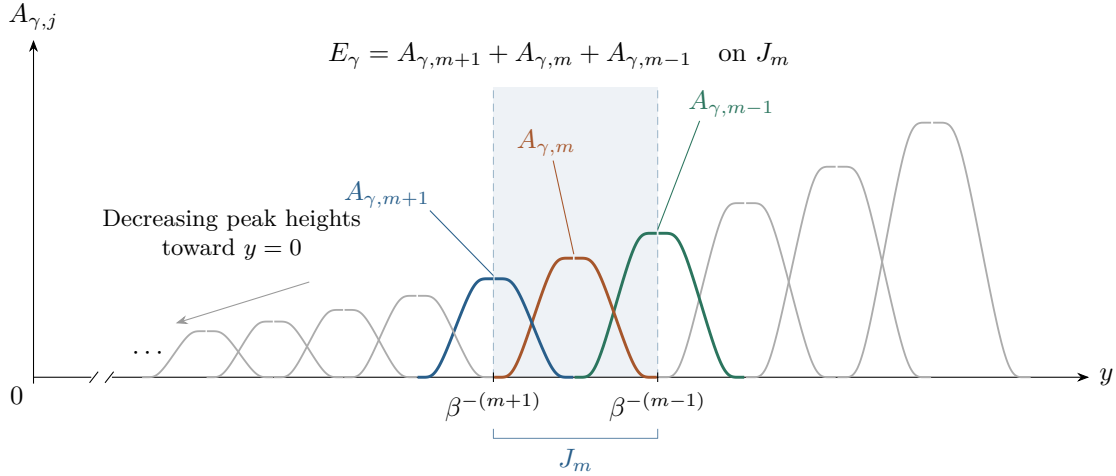

Set $J_m\triangleq[\beta^{-1},\beta]\beta^{-m}$ for $m\in\mathbb Z$.
Figure~\ref{fig:packet-amplitudes} illustrates the supports of the functions
$A_{\gamma,j}$ and the decrease of their peak heights toward $y=0$.

For every integer $m\ge7$, on $\beta^{-(m+1)}<y\le\beta^{-m}$, we have
$A_{\gamma,j}(y)=0$ for every $j\ge4$ with $j\notin\{m,m+1\}$. Thus
\begin{equation}\label{ce:two-packet-envelope}
E_\gamma(y)=A_{\gamma,m}(y)+A_{\gamma,m+1}(y).
\end{equation}
For $\beta^{-(m+1)}<y\le\beta^{-m}$, the two cutoff arguments satisfy
\[
\begin{aligned}
\rho(\beta^my)&=b\!\left(\frac{\beta^my-\beta^{-1}}{1-\beta^{-1}}\right),\qquad
\rho(\beta^{m+1}y) =b\!\left(1-\frac{\beta^my-\beta^{-1}}{1-\beta^{-1}}\right).
\end{aligned}
\]
Since $\rho(\beta^my)+\rho(\beta^{m+1}y)=1$, it follows that
\begin{equation}\label{ce:packet-envelope}
e^{-\beta^{m+1}-(m+1)^2-\gamma(m+1)}
\le E_\gamma(y)\le e^{-\beta^m-m^2-\gamma m}.
\end{equation}
For $m\ge5$ and $y\in J_m$, we similarly have
$A_{\gamma,j}(y)=0$ for every $j\ge4$ with $j\notin\{m-1,m,m+1\}$, so
\begin{equation}\label{ce:band-envelope-lower}
E_\gamma(y)=A_{\gamma,m-1}(y)+A_{\gamma,m}(y)+A_{\gamma,m+1}(y)
\ge e^{-1/y-(m+1)^2-\gamma(m+1)}.
\end{equation}
The lower bound follows from nonnegativity and the identities
\[
\sum_{j=m-1}^{m+1}\rho(\beta^jy)=1,\qquad
e^{-j^2-\gamma j}\ge e^{-(m+1)^2-\gamma(m+1)}
\quad(m-1\le j\le m+1).
\]
The half-open intervals $\beta^{-(m+1)}<y\le\beta^{-m}$ with
$m\ge7$ partition $(0,\beta^{-7}]$, which contains every normal
interval used below. On each such interval, $E_\gamma$
is the two-term sum in \eqref{ce:two-packet-envelope}. The intervals $J_m$
overlap and contain the supports used in the correction.

The indices used below have the following roles:
\begin{center}
\small
\begin{tabular}{@{}p{.19\textwidth}p{.35\textwidth}p{.40\textwidth}@{}}
\toprule
Notation&Range or location&Role\\
\midrule
$j$&$j\ge4$, height $y\asymp\beta^{-j}$&packet index; frequency $\lambda_j$\\
$m$&band $J_m$, $m\ge5$&normal-band index; neighboring packets have indices $m-1,m,m+1$\\
$h$&$h\ge5$, band $J_h$&running index when the corrections $W_h$ are summed\\
$A_{\gamma,j},E_\gamma$&nonnegative weights and their sum&measure amplitudes before oscillatory cancellation\\
\bottomrule
\end{tabular}
\end{center}
The condition $m\ge5$ ensures $m-1\ge4$. The half-open working bands
start at $m=7$ because the initial width is $Y_*\le\beta^{-7}$.

\medskip\noindent\textbf{Recursive cancellation of the wave residual.}
For each integer $j\ge4$, we construct an approximate solution of
the free wave equation as a finite sum in inverse powers of $\lambda_j$.
The summation index is $\ell\in\{0,\ldots,K_j\}$, where $K_j$ is fixed by
\eqref{ce:front-parameters}. This upper limit is chosen so that the
residual is sufficiently small.
Define the coefficients recursively by
\begin{align}
&H_{\gamma,j,0}\triangleq A_{\gamma,j},\quad
H_{\gamma,j,\ell+1}(z,y)\triangleq\frac{\mathrm i}{2}
\int_0^z(\partial_\sigma^2+\partial_y^2)
H_{\gamma,j,\ell}(\sigma,y)\,\mathrm{d}\sigma,\quad
\mathcal F_{\gamma,j}\triangleq e^{\mathrm i\lambda_ju}
\sum_{\ell=0}^{K_j}\lambda_j^{-\ell}H_{\gamma,j,\ell}.
\label{ce:approximate-solution}
\end{align}
The recursive integral is chosen so that the leading transport term
cancels. To see this, compute the residual of $\mathcal F_{\gamma,j}$.
For a coefficient $H=H(z,y)$ independent of $t$, direct differentiation gives
\[
\begin{aligned}
\partial_t^2(e^{\mathrm i\lambda_j u}H)
&=-\lambda_j^2e^{\mathrm i\lambda_j u}H,\qquad
\partial_z^2(e^{\mathrm i\lambda_j u}H)
=e^{\mathrm i\lambda_j u}(-\lambda_j^2H+2\mathrm i\lambda_j\partial_zH+\partial_z^2H).
\end{aligned}
\]
The two terms of order $\lambda_j^2$ cancel in their difference. Thus
\[
\Box\bigl(e^{\mathrm i\lambda_ju}H_{\gamma,j,\ell}\bigr)
=e^{\mathrm i\lambda_ju}
\bigl(-2\mathrm i\lambda_j\partial_zH_{\gamma,j,\ell}
-\Delta_{z,y}H_{\gamma,j,\ell}\bigr).
\]
The recursion in \eqref{ce:approximate-solution} gives
\[
-2\mathrm i\partial_zH_{\gamma,j,\ell+1}
=\Delta_{z,y}H_{\gamma,j,\ell},\qquad
\partial_zH_{\gamma,j,0}=0.
\]
Hence
\[
\begin{aligned}
e^{-\mathrm i\lambda_ju}\Box\mathcal F_{\gamma,j}
&=\sum_{\ell=0}^{K_j-1}\lambda_j^{-\ell}
\bigl(-2\mathrm i\partial_zH_{\gamma,j,\ell+1}
-\Delta_{z,y}H_{\gamma,j,\ell}\bigr) -\lambda_j^{-K_j}\Delta_{z,y}H_{\gamma,j,K_j}.
\end{aligned}
\]
Each term in parentheses is zero, so
\begin{equation}\label{ce:seed-residual}
\Box\mathcal F_{\gamma,j}
=-\lambda_j^{-K_j}e^{\mathrm i\lambda_ju}
\Delta_{z,y}H_{\gamma,j,K_j}.
\end{equation}

The recurrence in \eqref{ce:approximate-solution} cancels successive
errors and gives the residual \eqref{ce:seed-residual}.

\medskip\noindent\textbf{Control of the recursive coefficients.}
Throughout this section, $C=C(Z)$ and $c=c(Z)>0$ unless additional
fixed parameters are indicated. We write $C(r)=C(Z,r)$ and
$C(r,B)=C(Z,r,B)$, with the same convention for other derivative
indices. These constants, including those in $O(\cdot)$, are independent
of $j,m$ and the spacetime point, and uniform for both values of $\gamma$.
Constants obtained on a cylinder remain valid when its normal width
is later restricted.
\begin{lemma}[Packet coefficients and their estimates]\label{ce:coefficient-bounds}
The coefficients defined in \eqref{ce:approximate-solution} have the
following properties. The estimates in (i) and (iii) hold on the band
$J_m$ for every integer $m\ge5$, uniformly in $t\in\R$ and
$|z|\le Z+1$.
\begin{enumerate}
\item[(i)] There exists a constant $C\ge1$, independent of $Z$, $m$, $j$,
and $r$, such that for every $j\in\{m-1,m,m+1\}$ and every integer
$r\ge0$,
\begin{equation}\label{ce:growing-amplitude}
|A_{\gamma,j}^{(r)}|
\le e^{Cm}C^{r+1}(r!)^2y^{-2r}E_\gamma.
\end{equation}

\item[(ii)] For every integer $j\ge4$ and every integer $\ell\ge1$, the
coefficients defined in \eqref{ce:approximate-solution} satisfy
\begin{equation}\label{ce:coefficient-formula}
H_{\gamma,j,\ell}(z,y)=\left(\frac{\mathrm i}{2}\right)^\ell
\sum_{\substack{1\le a\le\ell\\a\equiv\ell\pmod2}}
\frac a\ell\binom{\ell}{(\ell-a)/2}
\frac{z^a}{a!}A_{\gamma,j}^{(\ell+a)}(y).
\end{equation}
In particular, $H_{\gamma,j,\ell}(0,y)=0$ for every $\ell\ge1$.

\item[(iii)] For every integer $r\ge0$ there exists a constant
$C(r)=C(Z,r)\ge1$ such that, for every $j\in\{m-1,m,m+1\}$ and every
$1\le\ell\le K_j$,
\begin{align}
\log\frac{|\partial^r(e^{\mathrm i\lambda_ju}
\lambda_j^{-\ell}H_{\gamma,j,\ell})|}{E_\gamma}
& \le-\tfrac7{16}j\ell\log2
+10\ell\log(j+2)+C(r)\ell+C(r)m.
\label{ce:coefficient-estimate}
\end{align}
\end{enumerate}
\end{lemma}
\begin{proof}[Proof outline of Lemma~\ref{ce:coefficient-bounds}]
Differentiate the common exponential and cutoff separately to obtain
\eqref{ce:growing-amplitude}. Induction in the recurrence gives
\eqref{ce:coefficient-formula}; combining factorial growth with
$\lambda_j^{-\ell}$ gives \eqref{ce:coefficient-estimate}.
The complete calculation is given in \hyperref[app:Coefficients]{Appendix A}.
\end{proof}

\medskip\noindent\textbf{Sum the packets toward the interface.}
Although \eqref{ce:coefficient-estimate} with $\ell=K_j$ and
\eqref{ce:seed-residual} make the residual of each
$\mathcal F_{\gamma,j}$ very small, any fixed such function vanishes
for $0<y<\beta^{-j-1}$.
We need a real approximate solution with nonzero values arbitrarily
close to the interface $y=0$, so we define
\begin{equation}\label{ce:approximate-sum}
g_\gamma\triangleq\operatorname{Re}\sum_{j\ge4}\mathcal F_{\gamma,j}.
\end{equation}
By \eqref{ce:coefficient-formula} and the flatness of $\rho$ at its
support endpoints, a nonzero summand must have
$\beta^{-j-1}<y<\beta^{1-j}$, so at each fixed $y>0$ at most two
summands are nonzero.
On every compact subset of $\{y>0\}$, the variable $y$ has a positive
lower bound, whereas $\sup J_j=\beta^{1-j}\to0$ as $j\to\infty$.
Hence only finitely many indices $j$ occur there.
Each $\mathcal F_{\gamma,j}$ is itself a finite sum over
$0\le\ell\le K_j$, so $g_\gamma$ is well defined and smooth for $y>0$.
Moreover, \eqref{ce:coefficient-formula} and
\eqref{ce:packet-envelope} give
\[
g_\gamma(0,0,y)
=\sum_{j\ge4}A_{\gamma,j}(y)
=E_\gamma(y)>0,\qquad 0<y<Y_*.
\]

\begin{lemma}[Residual of the approximate solution]\label{ce:sum-residual}
There exists $c>0$ such that, for every integer $m\ge5$
and every fixed Cartesian derivative order $r\ge0$, the approximate
solution \eqref{ce:approximate-sum} satisfies, on the full band $J_m$,
uniformly for $t\in\R$ and $|z|\le Z+1$,
\begin{equation}\label{ce:seed-reserve}
\frac{|\partial^r\Box g_\gamma|}{E_\gamma}
\le C(r)\exp\{-cm^3\beta^m+C(r)m^2\beta^m\log(m+2)\}.
\end{equation}
\end{lemma}

\begin{proof}[Proof of Lemma~\ref{ce:sum-residual}]
For each fixed derivative order, the choice of $K_j$ gives a negative
exponent of size $jK_j$ and positive terms of size at most
$K_j\log(j+2)$. We express these in the band index $m$.
Fix $r\ge0$, $m\ge5$, $y\in J_m$, and $|z|\le Z+1$.
We apply \eqref{ce:coefficient-estimate} with $\ell=K_j$ for
$j\in\{m-1,m,m+1\}$.
Since $j\ge4$, we have
\[
\tfrac12j^2\beta^j\le j^2\beta^j-1
\le K_j\le j^2\beta^j.
\]
Since $m\ge2$ and $|j-m|\le1$,
\[
\frac m2\le j\le2m,\qquad
\beta^{-1}\beta^m\le\beta^j\le\beta\beta^m,
\]
so
\[
\frac{m^2\beta^m}{8\beta}\le K_j\le4\beta m^2\beta^m,
\qquad jK_j\ge\frac{m^3\beta^m}{16\beta}.
\]
Consequently,
\[
\begin{aligned}
-\frac7{16}jK_j\log2
&\le-\frac{7\log2}{256\beta}m^3\beta^m,\\
10K_j\log(j+2)+C(r)K_j+C(r)m
&\le4\beta m^2\beta^m
\bigl[20\log(m+2)+C(r)\bigr]+C(r)m\\
&\le C(r)m^2\beta^m\log(m+2).
\end{aligned}
\]
Here $\log(j+2)\le\log(m+3)\le2\log(m+2)$.
In particular, \eqref{ce:coefficient-estimate} gives, uniformly in
$j\in\{m-1,m,m+1\}$,
\[
\frac{|\partial^r(e^{\mathrm i\lambda_ju}
\lambda_j^{-K_j}H_{\gamma,j,K_j})|}{E_\gamma}
\le C(r)\exp\{-cm^3\beta^m
+C(r)m^2\beta^m\log(m+2)\},
\qquad c=\frac{7\log2}{256\beta}>0.
\]

On $J_m$, the locally finite sum \eqref{ce:approximate-sum} and
\eqref{ce:seed-residual} satisfy
\[
\begin{aligned}
\Box g_\gamma
&=\operatorname{Re}\sum_{j\ge4}\Box\mathcal F_{\gamma,j}
=-\operatorname{Re}\sum_{j=m-1}^{m+1}
\lambda_j^{-K_j}e^{\mathrm i\lambda_ju}
\Delta_{z,y}H_{\gamma,j,K_j}.
\end{aligned}
\]
Also,
\[
\begin{aligned}
&(\partial_z-\partial_t)
\bigl(e^{\mathrm i\lambda_ju}H_{\gamma,j,K_j}\bigr)
=e^{\mathrm i\lambda_ju}\partial_zH_{\gamma,j,K_j},\\
&e^{\mathrm i\lambda_ju}\Delta_{z,y}H_{\gamma,j,K_j}
=\bigl[\partial_z^2-2\partial_z\partial_t
+\partial_t^2+\partial_y^2\bigr]
\bigl(e^{\mathrm i\lambda_ju}H_{\gamma,j,K_j}\bigr).
\end{aligned}
\]
For $|\alpha|=r$, equation \eqref{ce:seed-residual} therefore implies
\[
\begin{aligned}
\frac{|\partial^\alpha\Box g_\gamma|}{E_\gamma}
&\le5\sum_{j=m-1}^{m+1}\max_{|\nu|=r+2}
\frac{|\partial^\nu(e^{\mathrm i\lambda_ju}
\lambda_j^{-K_j}H_{\gamma,j,K_j})|}{E_\gamma}\\
&\le15C(r+2)\exp\{-cm^3\beta^m
+C(r+2)m^2\beta^m\log(m+2)\}\\
&\le C(r)\exp\{-cm^3\beta^m
+C(r)m^2\beta^m\log(m+2)\},
\end{aligned}
\]
which is \eqref{ce:seed-reserve}.
\end{proof}

For every fixed $r,N\in\mathbb N_0$, \eqref{ce:seed-reserve} also gives,
for all sufficiently large $m$,
\[
y^{-N}\frac{|\partial^r\Box g_\gamma|}{E_\gamma}
\le C(r)\exp\!\left[-\frac c2m^3\beta^m+N(m+1)\log\beta\right]
\longrightarrow0.
\]
Here $C(r)\log(m+2)/m\to0$ absorbs the positive term in the exponent,
and $y^{-1}\le\beta^{m+1}$ on $J_m$. Thus the residual bound can pay
any fixed power of the normal scale introduced in the correction.

Although Lemma~\ref{ce:sum-residual} makes $\Box g_\gamma$ small,
$g_\gamma$ can have zeros and \eqref{ce:seed-reserve} does not ensure that
$-\Box g_\gamma/g_\gamma$ is smooth there.
The following proposition adds periodic corrections to obtain $f_\gamma$
and confines its residual to regions where a positive lower bound for
$f_\gamma$ will allow us to define a smooth real potential $V_\gamma$
and obtain our final solution of $\Box f_\gamma+V_\gamma f_\gamma=0$.

The periodic correction is constructed separately for each fixed
$\gamma\in\{0,24\log2\}$. We suppress this dependence in the notation
$W_h,w_m,h_m^{(k_0)},h_m,T_m$ used below. These functions may differ
for the two weights, whereas the widths, characteristic cells, and
cutoff functions are common, and all estimates are uniform in $\gamma$.
\begin{proposition}[Periodic correction and residual bounds]\label{ce:reserve}
After decreasing $Y_*$ by a choice depending only on $Z$, there exist
real smooth $2\pi$-periodic functions $W_h$, indexed by integers
$h\ge5$, on $\R\times(-Z,Z)\times(0,Y_*)$, whose support in $y$ lies
in the band $J_h$. Set
\[
f_\gamma\triangleq g_\gamma+\sum_{h\ge5}W_h,
\qquad
R_\gamma\triangleq\Box f_\gamma.
\]
Then the following properties hold.
\begin{enumerate}
\item[(i)] The residual is confined to narrow phase bands: relative
to $\R\times(-Z,Z)\times(0,Y_*)$,
\[
\operatorname{supp}R_\gamma\subset
\bigcup_{h\ge5}
\Bigl\{(t,z,y):y\in J_h,\quad
\operatorname{dist}(2^{h-1}(t+z),2\pi\mathbb Z)\le\pi/72\Bigr\}.
\]
The index $h$ may depend on the point $(t,z,y)$.

\item[(ii)] There exists a constant $c>0$ such that, for every integer
$m\ge5$ and every fixed $\alpha\in\mathbb N_0^3$, uniformly in
$t\in\R$, $|z|<Z$, and $y\in J_m\cap(0,Y_*)$,
\begin{equation}\label{ce:reserve-estimate}
\frac{|\partial^\alpha\Box g_\gamma|+|\partial^\alpha R_\gamma|
+\sum_h|\partial^\alpha W_h|}{E_\gamma}
\le C(\alpha)\exp\{-cm^3\beta^m+C(\alpha) m^2\beta^m\log(m+2)\}.
\end{equation}

\item[(iii)] For every fixed $A,N\in\R$, $P\ge0$, and every
multi-index $\alpha$, uniformly in $t\in\R$, $|z|<Z$, and
$y\in J_m\cap(0,Y_*)$ with $m\ge5$,
\[
e^{A\beta^m+Am^2+Am}y^{-N}
\frac{|\partial^\alpha\Box g_\gamma|+|\partial^\alpha R_\gamma|
+\sum_h|\partial^\alpha W_h|}{E_\gamma}
=o(y^P)\qquad(y\downarrow0).
\]
\end{enumerate}
\end{proposition}

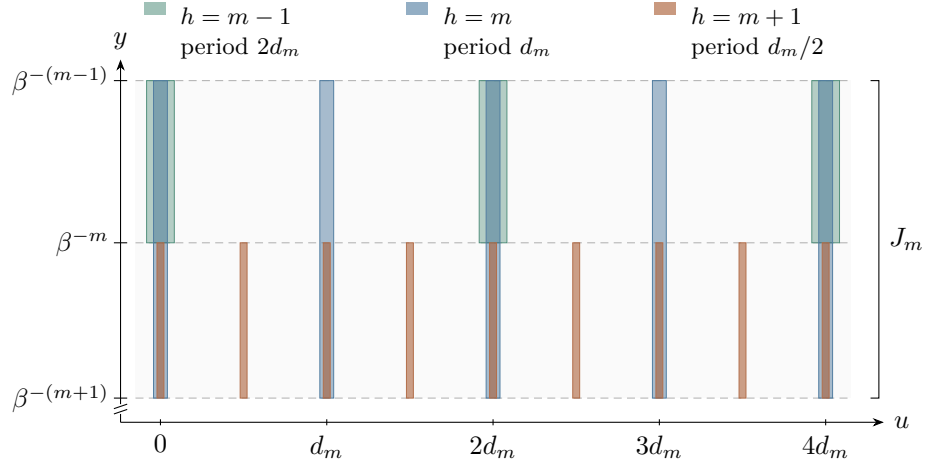
\begin{figure}[H]
\centering
\begin{tikzpicture}[x=2.2cm,y=1cm,font=\small,
>={Stealth[length=1.5mm]}]
\definecolor{residualgreen}{RGB}{43,117,94}
\definecolor{residualblue}{RGB}{40,94,138}
\definecolor{residualrust}{RGB}{166,86,43}
\pgfmathsetmacro{\residualbeta}{2^(1/16)}
\pgfmathsetmacro{\residualmid}{4.2/(1+\residualbeta)}
\fill[black!2] (-.15,0) rectangle (4.15,4.2);
\draw[gray!60,densely dashed] (-.15,0)--(4.15,0);
\draw[gray!60,densely dashed] (-.15,\residualmid)--(4.15,\residualmid);
\draw[gray!60,densely dashed] (-.15,4.2)--(4.15,4.2);
\begin{scope}
\clip (-.15,0) rectangle (4.15,4.2);
\foreach \residualcentre in {0,2,4}{
\filldraw[fill=residualgreen,fill opacity=.34,
draw=residualgreen!80,line width=.2pt]
({\residualcentre-1/12},\residualmid)
rectangle ({\residualcentre+1/12},4.2);
}
\foreach \residualcentre in {0,1,2,3,4}{
\filldraw[fill=residualblue,fill opacity=.40,
draw=residualblue!80,line width=.2pt]
({\residualcentre-1/24},0)
rectangle ({\residualcentre+1/24},4.2);
}
\foreach \residualcentre in {0,.5,1,1.5,2,2.5,3,3.5,4}{
\filldraw[fill=residualrust,fill opacity=.56,
draw=residualrust!85,line width=.2pt]
({\residualcentre-1/48},0)
rectangle ({\residualcentre+1/48},\residualmid);
}
\end{scope}
\draw[->,thin] (-.24,-.32)--(4.35,-.32) node[right] {$u$};
\draw[->,thin] (-.24,-.32)--(-.24,4.48) node[above] {$y$};
\draw[white,line width=2.4pt] (-.24,-.21)--(-.24,-.08);
\draw[thin] (-.28,-.21)--(-.20,-.15)
(-.28,-.14)--(-.20,-.08);
\foreach \residualtick/\residuallabel in
{0/0,1/{d_m},2/{2d_m},3/{3d_m},4/{4d_m}}{
\draw[thin] (\residualtick,-.28)--(\residualtick,-.36)
node[below] {$\residuallabel$};
}
\draw[thin] (-.28,0)--(-.20,0)
node[left=3pt] {$\beta^{-(m+1)}$};
\draw[thin] (-.28,\residualmid)--(-.20,\residualmid)
node[left=3pt] {$\beta^{-m}$};
\draw[thin] (-.28,4.2)--(-.20,4.2)
node[left=3pt] {$\beta^{-(m-1)}$};
\draw[thin] (4.27,0)--(4.32,0)--(4.32,4.2)--(4.27,4.2);
\node[right] at (4.32,2.1) {$J_m$};
\fill[residualgreen,fill opacity=.45] (-.10,5.07) rectangle (.03,5.24);
\node[anchor=north west,align=left,font=\footnotesize] at (.06,5.30)
{$h=m-1$\\period $2d_m$};
\fill[residualblue,fill opacity=.50] (1.48,5.07) rectangle (1.61,5.24);
\node[anchor=north west,align=left,font=\footnotesize] at (1.64,5.30)
{$h=m$\\period $d_m$};
\fill[residualrust,fill opacity=.60] (2.97,5.07) rectangle (3.10,5.24);
\node[anchor=north west,align=left,font=\footnotesize] at (3.13,5.30)
{$h=m+1$\\period $d_m/2$};
\end{tikzpicture}
\caption{A schematic enclosure of $\operatorname{supp}R_\gamma$ from
Proposition~\ref{ce:reserve}, shown in $\operatorname{int}J_m$ for large $m$,
with $d_m=2\pi/2^{m-1}$.
The colors correspond to $h=m-1,m,m+1$; the middle family spans
all of $J_m$ and overlaps the other two.
The true phase half-width is $\pi/(72\,2^{h-1})$;
the widths are enlarged sixfold for visibility.
Shading indicates where residual support may occur.}
\label{fig:residual-support}
\end{figure}

Part (i) confines the residual to narrow phase bands. Lemma~\ref{ce:division-region} gives $f_\gamma\ge5E_\gamma/8$ near these bands, allowing division there; (ii)--(iii) ensure that $V_\gamma=-R_\gamma/f_\gamma$ extends smoothly by zero across $y=0$.

\smallskip\noindent\textit{The equation to solve inside one cell.}
The characteristic coordinates reduce the required correction to
\[
(4\partial_r\partial_s-\partial_y^2)w=h,
\qquad w(0,s,y)=w(r,0,y)=0.
\]
Indeed, with $T=(r+s)/2$ and $X=(r-s)/2$,
\[
4\partial_r\partial_s-\partial_y^2
=\partial_T^2-\partial_X^2-\partial_y^2.
\]
Write the two characteristic displacements as $a,b\ge0$. Then
\[
\Delta T=\frac{a+b}{2},\qquad \Delta X=\frac{a-b}{2},\qquad
(\Delta T)^2-(\Delta X)^2=ab.
\]
For $a,b>0$, the integral averages the source over the normal
interval $|\Delta y|\le\sqrt{ab}$ inside the light cone, with weight
$(ab-(\Delta y)^2)^{-1/2}$. The support bound \eqref{ce:causal-support} will keep the correction
inside the prescribed normal band, and the derivative bound \eqref{ce:causal-derivative-bound} will be uniform as the cell width tends to zero.
\begin{lemma}[Characteristic solution formula and estimates]\label{ce:characteristic-solver}
Let $d,L>0$, and let $h\in C^\infty([0,d]\times[0,L]\times\R;\R)$.
Assume that there are two fixed real numbers $y_-<y_+$, independent of
$(r,s)$, such that
\[
h(r,s,y)=0\qquad\text{for all }(r,s)\in[0,d]\times[0,L]
\text{ and }y\notin[y_-,y_+].
\]
Write
\[
\operatorname{supp}_y h\triangleq
\overline{\{y\in\R:\ h(r,s,y)\ne0
\text{ for some }(r,s)\in[0,d]\times[0,L]\}}\subset[y_-,y_+].
\]
For $g\in C^\infty(\R)$ and $q\ge0$, define
\begin{equation}\label{ce:angular-mean}
\mathcal A_qg(y)\triangleq\frac1\pi\int_{-1}^1
\frac{g(y+\sqrt q\,\theta)}{\sqrt{1-\theta^2}}\,\mathrm d\theta,
\end{equation}
and set
\begin{equation}\label{ce:causal-inverse}
w(r,s,y)\triangleq\frac14\int_0^r\int_0^s
\mathcal A_{ab}[h(r-a,s-b,\cdot)](y)\,\mathrm db\,\mathrm da.
\end{equation}
Then the following properties hold.
\begin{enumerate}
\item[(i)] The function $w$ is real and smooth up to the boundary of
$[0,d]\times[0,L]\times\R$, and solves the characteristic problem
\[
\begin{cases}
(4\partial_r\partial_s-\partial_y^2)w=h,
&0<r<d,\quad0<s<L,\quad y\in\R,\\
w(0,s,y)=w(r,0,y)=0.
\end{cases}
\]

\item[(ii)] The support of $w$ in the $y$ variable satisfies
\begin{equation}\label{ce:causal-support}
\operatorname{supp}_y w
\subset[y_- -\sqrt{dL},\,y_+ +\sqrt{dL}].
\end{equation}
The same support inclusion holds for every derivative of $w$.

\item[(iii)] For every $g\in C^\infty(\R)$ and all
$\nu,k\in\mathbb N_0$, the angular mean is smooth for $q\ge0$ and
\begin{equation}\label{ce:angular-bound}
|\partial_q^\nu\partial_y^k\mathcal A_qg(y)|
\le\frac1{4^\nu\nu!}
\sup_{|y'-y|\le\sqrt q}|g^{(2\nu+k)}(y')|.
\end{equation}
At $q=0$, the derivative on the left equals
$g^{(2\nu+k)}(y)/(4^\nu\nu!)$ before taking its absolute value.

\item[(iv)] Fix $D>0$. For each multi-index
$\alpha\in\mathbb N_0^3$, if $d\le D$, then
\begin{equation}\label{ce:causal-derivative-bound}
|\partial_{r,s,y}^\alpha w(r,s,y)|
\le C_\alpha(D,L)\max_{|\mu|\le3|\alpha|}
\sup_{\substack{0\le r'\le r,\ 0\le s'\le s\\
|y'-y|\le\sqrt{rs}}}
|\partial_{r,s,y}^\mu h(r',s',y')|.
\end{equation}
The constant $C_\alpha(D,L)$ is independent of $d$, $h$, $y_-$, $y_+$,
and the evaluation point.
\end{enumerate}
\end{lemma}
Proposition~\ref{ce:reserve} uses Lemma~\ref{ce:characteristic-solver}(i) to cancel the source, (ii) to bound propagation by $\sqrt{dL}$, and (iv) to bound derivatives without negative powers of $d$; (iii) proves (iv). The complete proof is in \hyperref[app:Characteristic]{Appendix A}.

\begin{proof}[Proof of Proposition~\ref{ce:reserve}]
The construction proceeds in three stages: cancel
$\chi_m\Box g_\gamma$ inside each characteristic cell, cut off the
cell solution and repeat it periodically, and verify that the resulting
residual and its derivatives satisfy the required bounds:
\[
-\chi_m\Box g_\gamma\ \longmapsto\ w_m
\ \longmapsto\ W_m=\eta_mw_m
\ \longmapsto\ f_\gamma=g_\gamma+\sum_mW_m.
\]
The first stage gives the exact residual formula
\eqref{ce:complete-front} and proves its phase localization.
Step~\ref{step:reserve-residual} bounds every derivative of the
correction and residual; Step~\ref{step:reserve-correction} shows
that \eqref{ce:reserve-estimate} absorbs the later division losses.

\medskip\noindent\textbf{Step 1: Construct the periodic correction and locate the residual.}

\smallskip\noindent\emph{Localize the source in the normal variable.}
We use the partition functions $\chi_m$ to split the source between
neighboring normal bands. The smaller interval $I_m$ will contain the
source, while the larger band $J_m$ leaves room for its propagation
under the wave equation. Put
\[
\delta_\beta\triangleq
\min\{\beta^{-3/4}-\beta^{-1},\,\beta-\beta^{3/4}\}>0,
\]
and make the first decrease of the normal width by setting
\begin{equation}\label{ce:geometric-front-width}
Y_*\triangleq Y_{\rm geom}\triangleq
\min\left\{\beta^{-7},\,
\beta^{-1}\left(\frac{\delta_\beta}
{2\sqrt{4\pi(4Z+2)}}\right)^{1/7}\right\}.
\end{equation}
All correction indices are integers $m\ge5$. Call an index active if
$J_m\cap(0,Y_{\rm geom})\ne\varnothing$. The cell construction below
is performed only for active indices; for every other index set
$w_m=W_m=0$. For every $y>0$, define
\[
\psi(s)\triangleq a(s-\beta^{-3/4})a(\beta^{3/4}-s),\qquad
\chi_m(y)\triangleq
\frac{\psi(\beta^my)}{\sum_{h\in\Z}\psi(\beta^hy)}.
\]
Their restrictions form a smooth partition
$\sum_m\chi_m(y)=1$ on $(0,Y_*)$. Set
$I_m\triangleq[\beta^{-3/4},\beta^{3/4}]\beta^{-m}$. Then
$\operatorname{supp}\chi_m\subset I_m$ and
$|\partial_y^\nu\chi_m|\le C(\nu)\beta^{m\nu}$ for $\nu\in\mathbb N_0$.
Moreover,
\[
0\le\chi_m\le1,\qquad
\chi_m(y)=1\quad\text{for }\quad
\beta^{-m-1/4}\le y\le\beta^{-m+1/4}.
\]
At each $y\in(0,Y_*)$, at most two adjacent functions $\chi_h$ are
nonzero, and their sum is $1$; see Figure~\ref{fig:partition-of-unity}.

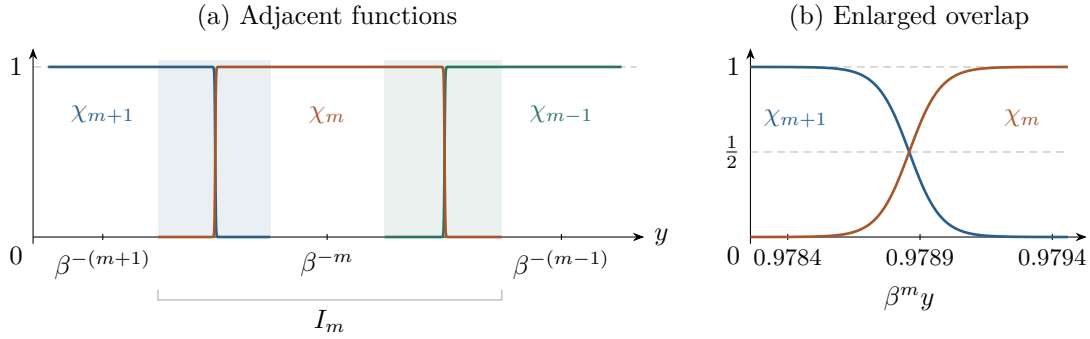
\begin{figure}[H]
\centering
\begin{tikzpicture}[x=1cm,y=2.25cm,font=\small,
>={Stealth[length=1.5mm]}]
\definecolor{partitionblue}{RGB}{40,94,138}
\definecolor{partitionrust}{RGB}{166,86,43}
\definecolor{partitiongreen}{RGB}{43,117,94}
\fill[partitionblue!10] (-2.237837,0) rectangle (-0.754039,1.04);
\fill[partitiongreen!10] (0.762250,0) rectangle (2.311742,1.04);
\draw[gray!55,densely dashed] (-3.889841,1)--(4.095162,1);
\draw[->,thin] (-3.889841,0)--(4.195162,0) node[right] {$y$};
\draw[->,thin] (-3.889841,-.04)--(-3.889841,1.11);
\node[left] at (-3.889841,1) {$1$};
\node[below left] at (-3.889841,0) {$0$};
\draw[thin] (-2.967770,.018)--(-2.967770,-.018) node[below] {$\beta^{-(m+1)}$};
\draw[thin] (0.000000,.018)--(0.000000,-.018) node[below] {$\beta^{-m}$};
\draw[thin] (3.099165,.018)--(3.099165,-.018) node[below] {$\beta^{-(m-1)}$};
\node at (0,1.29) {(a) Adjacent functions};
\draw[partitionblue,line width=1pt,line join=round] plot coordinates {
(-3.689841,1.0000000) (-1.534074,0.9999977) (-1.514792,0.9997641) (-1.507374,0.9986073)
(-1.502923,0.9959710) (-1.499956,0.9918424) (-1.495504,0.9767188) (-1.492536,0.9537924)
(-1.489568,0.9103620) (-1.485116,0.7780048) (-1.476211,0.2941649) (-1.471758,0.1255323)
(-1.468789,0.0658546) (-1.465820,0.0334439) (-1.462851,0.0166889) (-1.458398,0.0057890)
(-1.455428,0.0028426) (-1.450974,0.0009740) (-1.437609,0.0000384) (-0.754039,0.0000000)
};
\draw[partitiongreen,line width=1pt,line join=round] plot coordinates {
(0.762250,0.0000000) (1.508013,0.0000278) (1.517307,0.0002359) (1.523504,0.0009768)
(1.528152,0.0028287) (1.531251,0.0057352) (1.535899,0.0164441) (1.538999,0.0328653)
(1.542098,0.0646045) (1.545197,0.1230966) (1.549847,0.2891852) (1.559146,0.7738895)
(1.563796,0.9085911) (1.566897,0.9529321) (1.569997,0.9763381) (1.574648,0.9917495)
(1.577749,0.9959419) (1.582400,0.9986075) (1.590153,0.9997677) (1.608764,0.9999969)
(3.895162,1.0000000)
};
\draw[partitionrust,line width=1pt,line join=round] plot coordinates {
(-2.237837,0.0000000) (-1.525175,0.0000195) (-1.510341,0.0006849) (-1.505890,0.0019852)
(-1.502923,0.0040290) (-1.498472,0.0115908) (-1.494020,0.0328653) (-1.491052,0.0646045)
(-1.488084,0.1230966) (-1.483632,0.2891852) (-1.473243,0.8300064) (-1.468789,0.9341454)
(-1.464336,0.9763381) (-1.461367,0.9882560) (-1.458398,0.9942110) (-1.453943,0.9980100)
(-1.449489,0.9993190) (-1.436124,0.9999732) (-1.333572,1.0000000) (1.497172,0.9999977)
(1.517307,0.9997641) (1.525054,0.9986073) (1.529702,0.9959710) (1.532800,0.9918424)
(1.537449,0.9767188) (1.540548,0.9537924) (1.543647,0.9103620) (1.548297,0.7780048)
(1.557596,0.2941649) (1.562246,0.1255323) (1.565346,0.0658546) (1.568447,0.0334439)
(1.571547,0.0166889) (1.576198,0.0057890) (1.579299,0.0028426) (1.583951,0.0009740)
(1.597907,0.0000384) (2.311742,0.0000000)
};
\node[partitionblue] at (-2.967770,.72) {$\chi_{m+1}$};
\node[partitionrust] at (0,.72) {$\chi_m$};
\node[partitiongreen] at (3.099165,.72) {$\chi_{m-1}$};
\draw[gray!65,thin] (-2.237837,-.31)--(-2.237837,-.37)--(2.311742,-.37)--(2.311742,-.31);
\node[below] at (0.036952,-.37) {$I_m$};
\draw[->,thin] (5.6,0)--(10.05,0);
\draw[->,thin] (5.6,-.04)--(5.6,1.11);
\draw[gray!55,densely dashed] (5.6,1)--(9.8,1);
\draw[gray!55,densely dashed] (5.6,.5)--(9.8,.5);
\node[left] at (5.6,1) {$1$};
\node[left] at (5.6,.5) {$\frac12$};
\node[below left] at (5.6,0) {$0$};
\node at (7.7,1.29) {(b) Enlarged overlap};
\draw[thin] (6.093525,.018)--(6.093525,-.018) node[below,font=\footnotesize] {0.9784};
\draw[thin] (7.843525,.018)--(7.843525,-.018) node[below,font=\footnotesize] {0.9789};
\draw[thin] (9.593525,.018)--(9.593525,-.018) node[below,font=\footnotesize] {0.9794};
\node at (7.7,-.36) {$\beta^m y$};
\draw[partitionblue,line width=1pt,line join=round] plot coordinates {
(5.600000,0.9999569) (6.109250,0.9995035) (6.392750,0.9980732) (6.503000,0.9967390)
(6.597500,0.9948849) (6.676250,0.9925630) (6.749750,0.9894641) (6.823250,0.9850942)
(6.886250,0.9799608) (6.944000,0.9737575) (7.001750,0.9657015) (7.049000,0.9573818)
(7.096250,0.9471547) (7.143500,0.9346411) (7.185500,0.9212561) (7.222250,0.9075344)
(7.259000,0.8917027) (7.295750,0.8735376) (7.332500,0.8528281) (7.369250,0.8293887)
(7.400750,0.8070158) (7.432250,0.7824776) (7.469000,0.7511003) (7.495250,0.7269102)
(7.526750,0.6960188) (7.584500,0.6346705) (7.642250,0.5686071) (7.778750,0.4069327)
(7.847000,0.3311107) (7.883750,0.2933840) (7.915250,0.2631343) (7.952000,0.2304768)
(7.983500,0.2048272) (8.020250,0.1776550) (8.057000,0.1533889) (8.093750,0.1319026)
(8.135750,0.1105251) (8.172500,0.0943686) (8.214500,0.0785189) (8.256500,0.0651372)
(8.303750,0.0526316) (8.345750,0.0434500) (8.387750,0.0358080) (8.435000,0.0287555)
(8.482250,0.0230575) (8.534750,0.0180146) (8.592500,0.0137125) (8.655500,0.0101689)
(8.723750,0.0073467) (8.876000,0.0035465) (9.070250,0.0013948) (9.338000,0.0003836)
(9.800000,0.0000410)
};
\draw[partitionrust,line width=1pt,line join=round] plot coordinates {
(5.600000,0.0000431) (6.109250,0.0004965) (6.392750,0.0019268) (6.503000,0.0032610)
(6.597500,0.0051151) (6.676250,0.0074370) (6.749750,0.0105359) (6.823250,0.0149058)
(6.886250,0.0200392) (6.944000,0.0262425) (7.001750,0.0342985) (7.049000,0.0426182)
(7.096250,0.0528453) (7.143500,0.0653589) (7.185500,0.0787439) (7.222250,0.0924656)
(7.259000,0.1082973) (7.295750,0.1264624) (7.332500,0.1471719) (7.369250,0.1706113)
(7.400750,0.1929842) (7.432250,0.2175224) (7.469000,0.2488997) (7.495250,0.2730898)
(7.526750,0.3039812) (7.584500,0.3653295) (7.642250,0.4313929) (7.778750,0.5930673)
(7.847000,0.6688893) (7.883750,0.7066160) (7.915250,0.7368657) (7.952000,0.7695232)
(7.983500,0.7951728) (8.020250,0.8223450) (8.057000,0.8466111) (8.093750,0.8680974)
(8.135750,0.8894749) (8.172500,0.9056314) (8.214500,0.9214811) (8.256500,0.9348628)
(8.303750,0.9473684) (8.345750,0.9565500) (8.387750,0.9641920) (8.435000,0.9712445)
(8.482250,0.9769425) (8.534750,0.9819854) (8.592500,0.9862875) (8.655500,0.9898311)
(8.723750,0.9926533) (8.876000,0.9964535) (9.070250,0.9986052) (9.338000,0.9996164)
(9.800000,0.9999590)
};
\node[partitionblue] at (6.2,.70) {$\chi_{m+1}$};
\node[partitionrust] at (9.2,.70) {$\chi_m$};
\end{tikzpicture}
\caption{Three adjacent partition functions (left) and an enlarged overlap
(right); at most two functions are nonzero at each point, and their sum is $1$.}
\label{fig:partition-of-unity}
\end{figure}

\smallskip\noindent\emph{Solve the same source problem in every cell.}
Choose the cell period so that all three packet frequencies present
on $I_m$ are integer multiples of the cell frequency. Let
\begin{equation}\label{ce:correction-scales}
\Lambda_m\triangleq2^{m-1},\qquad d_m\triangleq2\pi/\Lambda_m.
\end{equation}
For every correction index $m\ge5$, we have
$d_m=2\pi\,2^{1-m}\le\pi/8<1$.
For each cell index $k_0\in\mathbb Z$, in the cell
$k_0d_m<u<(k_0+1)d_m$, set
\begin{equation}\label{ce:cell-coordinates}
\begin{gathered}
r\triangleq u-k_0d_m,\qquad s\triangleq v-k_0d_m+2Z+1,\\
s=r-2z+2Z+1,\qquad
0<r<d_m,\quad 1<s<4Z+2\qquad(|z|<Z).
\end{gathered}
\end{equation}
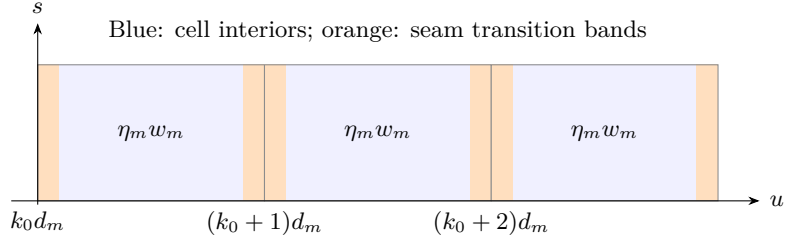
\begin{figure}[!htbp]
\centering
\begin{tikzpicture}[x=1cm,y=1cm,font=\footnotesize,>={Stealth[length=1.5mm]}]
\foreach \k in {0,1,2}{
\fill[blue!6] ({3*\k},0) rectangle ({3*\k+3},1.8);
\fill[orange!25] ({3*\k},0) rectangle ({3*\k+.28},1.8);
\fill[orange!25] ({3*\k+2.72},0) rectangle ({3*\k+3},1.8);
\draw[black!45] ({3*\k},0) rectangle ({3*\k+3},1.8);
\node at ({3*\k+1.5},.9) {$\eta_m w_m$};}
\draw[->] (-.35,0)--(9.55,0) node[right] {$u$};
\node[below] at (0,0) {$k_0d_m$};
\node[below] at (3,0) {$(k_0+1)d_m$};
\node[below] at (6,0) {$(k_0+2)d_m$};
\draw[->] (0,-.05)--(0,2.35) node[above] {$s$};
\node[align=center] at (4.5,2.25) {Blue: cell interiors; orange: seam transition bands};
\end{tikzpicture}
\caption{Periodic characteristic cells in the $u$ direction; seam widths are schematic. The function $w_m$ cancels the source inside each cell, $\eta_m$ smooths the seams, and only the two terms in \eqref{ce:complete-front} remain.}
\label{fig:characteristic-cells}
\end{figure}

For $(r,s)\in[0,d_m]\times[0,4Z+2]$ and $y\in\R$, define
\[
h_m^{(k_0)}(r,s,y)\triangleq-\chi_m(y)\Box g_\gamma
\left(k_0d_m+\tfrac12(r+s-2Z-1),\tfrac12(r-s+2Z+1),y\right).
\]
Here we use the packet sum $g_\gamma$ and the functions $\chi_m$
defined for all $y>0$, including $y\ge Y_{\rm geom}$, and extend
$\chi_m\Box g_\gamma$ by zero only to $y\le0$.
For each fixed $m$, this extension is smooth because
$\operatorname{supp}\chi_m\subset I_m\Subset(0,\infty)$.

We now show that the source is independent of the cell index $k_0$.
Set
\[
t_0\triangleq\tfrac12(r+s-2Z-1),\qquad
z_0\triangleq\tfrac12(r-s+2Z+1).
\]
Then $t_0+z_0=r$, and the spacetime point in the definition of
$h_m^{(k_0)}$ has
\[
t=k_0d_m+t_0,\qquad z=z_0,\qquad u=t+z=k_0d_m+r.
\]
The coefficient formula \eqref{ce:coefficient-formula} gives
$\operatorname{supp}_y\Delta_{z,y}H_{\gamma,j,K_j}\subset J_j$.
Since $I_m\cap J_j=\varnothing$ whenever $|j-m|\ge2$, only
$j=m-1,m,m+1$ can contribute on $I_m$.
Thus \eqref{ce:seed-residual} and
$\Box g_\gamma=\operatorname{Re}\sum_j\Box\mathcal F_{\gamma,j}$ give,
for $y\in I_m$,
\[
h_m^{(k_0)}(r,s,y)
=\chi_m(y)\operatorname{Re}\sum_{j=m-1}^{m+1}
\lambda_j^{-K_j}e^{\mathrm i\lambda_jk_0d_m}
e^{\mathrm i\lambda_jr}
\bigl(\Delta_{z,y}H_{\gamma,j,K_j}\bigr)(z_0,y).
\]
For each of these three indices,
\[
\frac{\lambda_j}{\Lambda_m}=2^{j-m+1}\in\{1,2,4\},\qquad
e^{\mathrm i\lambda_jk_0d_m}
=e^{2\pi\mathrm i k_0(\lambda_j/\Lambda_m)}=1.
\]
Moreover, $z_0$ is independent of $k_0$, and the coefficients
$H_{\gamma,j,K_j}$ have no $t$ dependence. Consequently,
\begin{equation}\label{ce:cell-source-common}
\begin{aligned}
h_m^{(k_0)}(r,s,y)
&=\chi_m(y)\operatorname{Re}\sum_{j=m-1}^{m+1}
\lambda_j^{-K_j}e^{\mathrm i\lambda_jr}
\bigl(\Delta_{z,y}H_{\gamma,j,K_j}\bigr)(z_0,y)
=h_m^{(0)}(r,s,y).
\end{aligned}
\end{equation}
Outside $I_m$, both sources vanish because $\chi_m=0$, including
under the extension to $y\le0$. Hence \eqref{ce:cell-source-common} holds on the
whole rectangle $[0,d_m]\times[0,4Z+2]\times\R$.
Denote this common source by $h_m$. Since
$\operatorname{supp}_y h_m\subset I_m$, apply
Lemma~\ref{ce:characteristic-solver} with
\[
d=d_m,\quad L=4Z+2,\quad h=h_m,\quad
y_-=\beta^{-m-3/4},\quad y_+=\beta^{-m+3/4}.
\]
Let $w_m$ be the smooth solution given by
\eqref{ce:causal-inverse}. By part (i) of the lemma, it satisfies
\[
\begin{cases}
(4\partial_r\partial_s-\partial_y^2)w_m=h_m,\\[3pt]
w_m(0,s,y)=w_m(r,0,y)=0.
\end{cases}
\]

\smallskip\noindent\emph{Keep the solution inside $J_m$.}
We check the source region and compare the propagation distance with
the gap from $I_m$ to $\partial J_m$.
Every source point in the formula has
$0\le r'\le d_m\le1$, $0\le s'\le4Z+2$, and therefore
$z'=(r'-s'+2Z+1)/2\in[-Z-1/2,Z+1]$.
In particular, $|z'|\le Z+1$, so the residual estimate of
Lemma~\ref{ce:sum-residual} holds on the full source region.
The integral only uses points with $|y'-y|\le\sqrt{(4Z+2)d_m}$, and
\[
\sqrt{(4Z+2)d_m}\beta^m=\sqrt{4\pi(4Z+2)}\,2^{-7m/16}\longrightarrow0.
\]
For an active index, $\beta^{-m}<\beta Y_{\rm geom}$, and hence
\eqref{ce:geometric-front-width} gives
\[
\begin{aligned}
\beta^m\sqrt{(4Z+2)d_m}
&=\sqrt{4\pi(4Z+2)}\,(\beta^{-m})^7
\le\sqrt{4\pi(4Z+2)}\,(\beta Y_{\rm geom})^7
\le\frac{\delta_\beta}{2}.
\end{aligned}
\]
Thus $\sqrt{(4Z+2)d_m}\le\delta_\beta\beta^{-m}/2$ for every active
correction index.
Lemma~\ref{ce:characteristic-solver}(ii), applied to the source
interval $I_m$, gives
\begin{equation}\label{ce:cell-correction-support}
\operatorname{supp}_y w_m
\subset[\beta^{-m-3/4}-\sqrt{(4Z+2)d_m},\,\beta^{-m+3/4}+\sqrt{(4Z+2)d_m}].
\end{equation}
Since $\operatorname{dist}(I_m,\partial J_m)=\delta_\beta\beta^{-m}$,
the enlarged interval remains inside $\operatorname{int}J_m$, with
\[
\operatorname{dist}(\operatorname{supp}_y w_m,\partial J_m)
\ge\delta_\beta\beta^{-m}-\sqrt{(4Z+2)d_m}
\ge\frac{\delta_\beta}{2}\beta^{-m}.
\]
Lemma~\ref{ce:characteristic-solver}(ii) gives the support inclusion
\eqref{ce:causal-support} for every derivative of $w_m$.

\smallskip\noindent\emph{Cut off the cell solution and repeat it periodically.}
The values of $w_m$ on opposite cell edges need not match, so we make
the correction zero near both edges before repeating it. Choose the
smooth periodic cutoff
\[
\eta(\theta)\triangleq
b\left(\frac{144}{\pi}\operatorname{dist}(\theta,2\pi\Z)-1\right),
\qquad \eta_m(u)\triangleq\eta(\Lambda_mu).
\]
Since $b=0$ on $(-\infty,0]$ and $b=1$ on $[1,\infty)$,
\[
\eta(\theta)=
\begin{cases}
0,&\operatorname{dist}(\theta,2\pi\mathbb Z)\le\pi/144,\\
1,&\operatorname{dist}(\theta,2\pi\mathbb Z)\ge\pi/72.
\end{cases}
\]
The periodic profile and a magnified neighborhood of $2k\pi$ are
shown in Figure~\ref{fig:periodic-eta}.
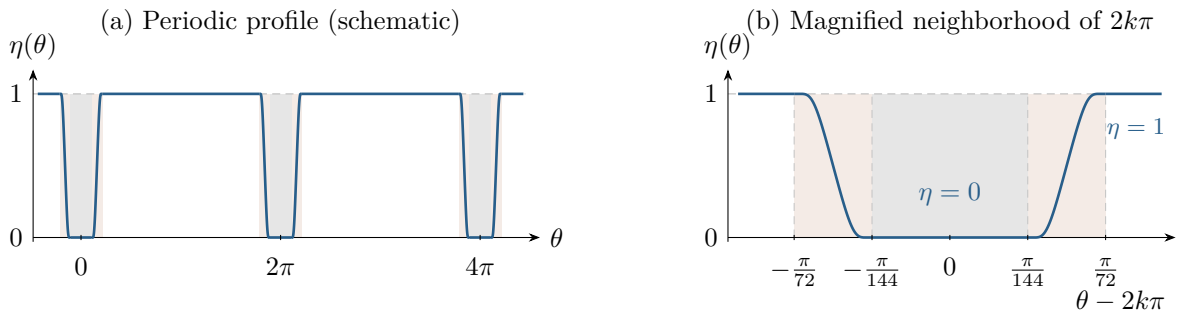
\begin{figure}[H]
\centering
\definecolor{etablue}{RGB}{40,94,138}
\definecolor{etarust}{RGB}{166,86,43}
\pgfmathdeclarefunction{etastep}{1}{%
\pgfmathparse{exp(-1/(#1))/(exp(-1/(#1))+exp(-1/(1-(#1))))}%
}
\begin{minipage}[t]{.47\textwidth}
\centering
\begin{tikzpicture}[baseline=(current bounding box.north),
x=.88cm,y=1.9cm,font=\small,
>={Stealth[length=1.5mm]}]
\node at (3,1.50) {(a) Periodic profile (schematic)};
\foreach \etacentre in {0,3,6}{
\fill[etarust!12] ({\etacentre-.32},0)
rectangle ({\etacentre+.32},1);
\fill[gray!20] ({\etacentre-.16},0)
rectangle ({\etacentre+.16},1);
}
\draw[gray!55,densely dashed] (-.72,1)--(6.68,1);
\draw[->,thin] (-.72,0)--(6.90,0) node[right] {$\theta$};
\draw[->,thin] (-.72,-.03)--(-.72,1.17)
node[above] {$\eta(\theta)$};
\node[left] at (-.72,0) {$0$};
\node[left] at (-.72,1) {$1$};
\draw[etablue,line width=1pt] (-.65,1)--(-.32,1)
(.32,1)--(2.68,1) (3.32,1)--(5.68,1)
(6.32,1)--(6.65,1);
\foreach \etacentre in {0,3,6}{
\draw[etablue,line width=1pt]
({\etacentre-.16},0)--({\etacentre+.16},0);
\draw[etablue,line width=1pt] ({\etacentre-.32},1)
-- plot[domain=.01:.99,samples=81]
({\etacentre-.32+.16*\x},{1-etastep(\x)})
-- ({\etacentre-.16},0);
\draw[etablue,line width=1pt] ({\etacentre+.16},0)
-- plot[domain=.01:.99,samples=81]
({\etacentre+.16+.16*\x},{etastep(\x)})
-- ({\etacentre+.32},1);
}
\foreach \etacentre/\etalabel in {0/0,3/{2\pi},6/{4\pi}}{
\draw[thin] (\etacentre,.018)--(\etacentre,-.018)
node[below=3pt] {$\etalabel$};
}
\end{tikzpicture}
\end{minipage}\hfill
\begin{minipage}[t]{.49\textwidth}
\centering
\begin{tikzpicture}[baseline=(current bounding box.north),
x=1.03cm,y=1.9cm,font=\small,
>={Stealth[length=1.5mm]}]
\node at (0,1.50) {(b) Magnified neighborhood of $2k\pi$};
\fill[etarust!12] (-2,0) rectangle (-1,1);
\fill[etarust!12] (1,0) rectangle (2,1);
\fill[gray!20] (-1,0) rectangle (1,1);
\draw[gray!55,densely dashed] (-2.85,1)--(2.75,1);
\foreach \etabound in {-2,-1,1,2}{
\draw[gray!45,densely dashed] (\etabound,0)--(\etabound,1);
}
\draw[->,thin] (-2.85,0)--(2.90,0);
\draw[->,thin] (-2.85,-.03)--(-2.85,1.17)
node[above] {$\eta(\theta)$};
\node[left] at (-2.85,0) {$0$};
\node[left] at (-2.85,1) {$1$};
\draw[etablue,line width=1pt] (-2.72,1)--(-2,1)
(-1,0)--(1,0) (2,1)--(2.72,1);
\draw[etablue,line width=1pt] (-2,1)
-- plot[domain=.01:.99,samples=81] ({-2+\x},{1-etastep(\x)})
-- (-1,0);
\draw[etablue,line width=1pt] (1,0)
-- plot[domain=.01:.99,samples=81] ({1+\x},{etastep(\x)})
-- (2,1);
\node[etablue] at (0,.30) {$\eta=0$};
\node[etablue,font=\footnotesize] at (2.39,.76) {$\eta=1$};
\foreach \etabound/\etalabel in
{-2/{-\frac\pi{72}},-1/{-\frac\pi{144}},0/0,
1/{\frac\pi{144}},2/{\frac\pi{72}}}{
\draw[thin] (\etabound,.018)--(\etabound,-.018)
node[below=3pt] {$\etalabel$};
}
\node[anchor=north east] at (2.90,-.30) {$\theta-2k\pi$};
\end{tikzpicture}
\end{minipage}
\caption{The $2\pi$-periodic cutoff $\eta$.
The narrow neighborhoods of $2\pi\Z$ are widened horizontally in (a).
In (b), the gray region has $\eta=0$; the shaded side strips are the
smooth transitions to $\eta=1$.}
\label{fig:periodic-eta}
\end{figure}
These constant values make $\eta$ smooth also at the multiples of
$\pi$, where the distance function is not smooth. Using the cell
coordinates $r,s$ in \eqref{ce:cell-coordinates}, define
\[
W_m(t,z,y)\triangleq
\begin{cases}
\eta(\Lambda_mr)w_m(r,s,y),
&k_0d_m<t+z<(k_0+1)d_m,\\
0,&t+z\in d_m\mathbb Z.
\end{cases}
\]
The function is zero on a neighborhood of every seam and hence is
smooth across the seams. Since $\Lambda_mk_0d_m=2\pi k_0$,
$\eta(\Lambda_mr)=\eta_m(u)$ in each cell.

\smallskip\noindent\emph{Compute the residual left by the cutoff.}
Within each cell, $\partial_u=\partial_r$ and $\partial_v=\partial_s$.
Since $\eta_m$ depends only on $u$,
\[
\partial_r\partial_s(\eta_mw_m)
=\eta_m\partial_r\partial_sw_m+\eta_m'\partial_sw_m,
\qquad
\partial_y^2(\eta_mw_m)=\eta_m\partial_y^2w_m.
\]
Consequently,
\[
\Box W_m
=\eta_m(4\partial_r\partial_s-\partial_y^2)w_m
+4\eta_m'\partial_sw_m
=-\eta_m\chi_m\Box g_\gamma+4\eta_m'\partial_sw_m.
\]
On $0<y<Y_{\rm geom}$, every inactive $\chi_m$ vanishes and the active
partition functions sum to $1$. Thus the function $f_\gamma$ in the
statement satisfies
\begin{equation}\label{ce:complete-front}
\begin{aligned}
R_\gamma
&=\sum_m\chi_m\Box g_\gamma
-\sum_m\eta_m\chi_m\Box g_\gamma
+4\sum_m\eta_m'\partial_sw_m =\sum_m(1-\eta_m)\chi_m\Box g_\gamma
+4\sum_m\eta_m'\partial_sw_m.
\end{aligned}
\end{equation}
Here $\eta_m'=\partial_u\eta_m$.
The first term on the last line is the part of the original source
not canceled after the cutoff; the second comes from differentiating
the cutoff. Both vanish where $\eta_m=1$.
The factor $\chi_m$ was included in the source equation rather than
multiplied onto $w_m$, so no derivatives of $\chi_m$ appear in this
residual identity.

We next check periodicity. For the cell containing $t+z$,
\[
\begin{aligned}
r(t+d_m)&=t+d_m+z-(k_0+1)d_m=r(t),\\
s(t+d_m)&=t+d_m-z-(k_0+1)d_m+2Z+1=s(t).
\end{aligned}
\]
Thus
\[
W_m(t+d_m,z,y)=W_m(t,z,y),\qquad
g_\gamma(t+2\pi,z,y)=g_\gamma(t,z,y),
\]
where the second equality uses $\lambda_j\in\mathbb N$.
Since $2\pi/d_m=\Lambda_m\in\mathbb N$, every $W_m$ is also
$2\pi$-periodic. The locally finite sum defining $f_\gamma$ therefore
has period $2\pi$.

The correction functions are now fixed. A later decrease of $Y_*$
will only restrict their sum to the smaller cylinder.

Finally, we locate the support of the residual. The definition of
$\eta$ gives
\[
\operatorname{supp}(1-\eta_m)\cup\operatorname{supp}\eta_m'
\subset\{u:\operatorname{dist}(\Lambda_mu,2\pi\mathbb Z)
\le\pi/72\}.
\]
Since the bands $J_m$ are locally finite for $y>0$, formula
\eqref{ce:complete-front} and the support inclusion
\eqref{ce:cell-correction-support} for $w_m$ imply
\[
\operatorname{supp}R_\gamma\subset
\bigcup_{m\ge5}
\{(t,z,y):y\in J_m,\quad
\operatorname{dist}(\Lambda_m(t+z),2\pi\mathbb Z)\le\pi/72\}.
\]
Since $\Lambda_m=2^{m-1}$, this proves (i).

\medskip\noindent\textbf{Step 2: Bound the derivatives uniformly in the band index.}\proofsteplabel{step:reserve-residual}{2}
We transfer the residual bound through the source $h_m$, the cell
solution $w_m$, and the two joined terms $W_m,T_m$. The only new
losses are the envelope comparison at nearby normal points and
derivatives of the cutoffs; for each fixed derivative order, these
contribute at most $e^{C(\alpha)m}$.

\smallskip\noindent\emph{Compare the envelopes at the source and evaluation points.}
It suffices to consider active indices and $y\in J_m$, since all
other correction terms vanish. For an active index, $w_m$ and its
derivatives vanish outside $J_m$, and every nonzero source term has
$y'\in I_m$.
For $x\in J_m$, we have $\sum_{j\ge4}\rho(\beta^jx)=1$, so
\[
e^{-1/x-(m+1)^2-\gamma(m+1)}
\le E_\gamma(x)\le e^{-1/x-(m-1)^2-\gamma(m-1)}.
\]
Since $|y'-y|\le\sqrt{(4Z+2)d_m}$, it follows that
\[
\begin{aligned}
\left|\frac1{y'}-\frac1y\right|
&\le C\beta^{2m}\sqrt{(4Z+2)d_m}\le C2^{-3m/8},\\
\frac{E_\gamma(y')}{E_\gamma(y)}
&\le\exp\left(\left|\frac1{y'}-\frac1y\right|+4m+2\gamma\right)
\le e^{Cm}.
\end{aligned}
\]

\smallskip\noindent\emph{Transfer the source estimate to the cell solution.}
Lemma~\ref{ce:characteristic-solver}(iv) is uniform for
$0<d_m\le1$: it uses at most $3|\alpha|$ source derivatives and no
factor $d_m^{-1}$. With $D=1$ and $L=4Z+2$, it gives, for every
$\alpha\in\mathbb N_0^3$,
\[
|\partial_{r,s,y}^\alpha w_m(r,s,y)|
\le C(\alpha)\max_{|\nu|\le3|\alpha|}
\sup_{\substack{0\le r'\le r,\ 0\le s'\le s\\
|y'-y|\le\sqrt{(4Z+2)d_m}}}
|\partial_{r,s,y}^\nu h_m(r',s',y')|.
\]
Here $\sqrt{rs}\le\sqrt{(4Z+2)d_m}$, and $C(\alpha)$ depends only on
$Z$ and $\alpha$, uniformly in the correction index $m$.
On the nonzero source region, $y'\in I_m\subset J_m$ and
$-Z-1/2\le z'\le Z+1$. At the corresponding source point,
\[
\begin{aligned}
\partial_r^a\partial_s^b\partial_y^d h_m
&=-\sum_{p=0}^d\binom dp\chi_m^{(p)}
\left(\frac{\partial_t+\partial_z}{2}\right)^a
\left(\frac{\partial_t-\partial_z}{2}\right)^b
\partial_y^{d-p}\Box g_\gamma.
\end{aligned}
\]
Thus \eqref{ce:seed-reserve} and
$|\partial_y^j\chi_m|\le C(j)\beta^{mj}$ imply
\[
\frac{|\partial_{r,s,y}^\nu h_m(r',s',y')|}{E_\gamma(y')}
\le C(\nu)\exp\{-cm^3\beta^m
+C(\nu)m^2\beta^m\log(m+2)+C(\nu)m\}.
\]
The envelope comparison therefore yields
\begin{equation}\label{ce:cell-solution-bound}
\frac{|\partial_{r,s,y}^\alpha w_m(r,s,y)|}{E_\gamma(y)}
\le C(\alpha)\exp\{-cm^3\beta^m
+C(\alpha)m^2\beta^m\log(m+2)+C(\alpha)m\}.
\end{equation}

\smallskip\noindent\emph{Differentiate the joined correction and its residual.}
Within each cell,
$\partial_t=\partial_r+\partial_s$ and
$\partial_z=\partial_r-\partial_s$. Set
\[
T_m\triangleq(1-\eta_m)\chi_m\Box g_\gamma
+4\eta_m'\partial_sw_m,
\qquad R_\gamma=\sum_mT_m.
\]
For each fixed integer $j\ge0$,
$|\eta_m^{(j)}|\le C(j)\Lambda_m^j\le C(j)e^{C(j)m}$.
For Cartesian multiindices $\alpha$,
\[
\begin{aligned}
\partial^\alpha W_m
&=\sum_{\mu\le\alpha}\binom\alpha\mu
(\partial^\mu\eta_m)\partial^{\alpha-\mu}w_m,\\
\partial^\alpha T_m
&=\sum_{\mu+\nu+\sigma=\alpha}
\frac{\alpha!}{\mu!\nu!\sigma!}
\partial^\mu(1-\eta_m)\partial^\nu\chi_m
\partial^\sigma\Box g_\gamma +4\sum_{\mu\le\alpha}\binom\alpha\mu
(\partial^\mu\eta_m')\partial^{\alpha-\mu}\partial_sw_m.
\end{aligned}
\]
Estimate \eqref{ce:cell-solution-bound}, up to order $|\alpha|+1$,
therefore implies
\begin{equation}\label{ce:joined-correction-bound}
\frac{|\partial^\alpha W_m|+|\partial^\alpha T_m|}{E_\gamma(y)}
\le C(\alpha)\exp\{-cm^3\beta^m
+C(\alpha)m^2\beta^m\log(m+2)+C(\alpha)m\}.
\end{equation}
Estimate \eqref{ce:joined-correction-bound} extends across the seams
because $W_m$ and $\eta_m'\partial_sw_m$ vanish on a neighborhood of
each seam.
If $y\in J_m$, a nonzero term with band index $h$ must satisfy
$|h-m|\le2$. Since $m,h\ge5$,
\[
h^3\beta^h\ge\frac{\beta^{-2}}8m^3\beta^m,\qquad
h^2\beta^h\log(h+2)+h\le Cm^2\beta^m\log(m+2).
\]
Consequently,
\[
\begin{aligned}
&\frac{|\partial^\alpha R_\gamma|
+\sum_h|\partial^\alpha W_h|}{E_\gamma(y)}
\le\sum_{|h-m|\le2}
\frac{|\partial^\alpha T_h|+|\partial^\alpha W_h|}{E_\gamma(y)}
\le C(\alpha)\exp\{-cm^3\beta^m
+C(\alpha)m^2\beta^m\log(m+2)\}.
\end{aligned}
\]
Together with \eqref{ce:seed-reserve}, this is
\eqref{ce:reserve-estimate}, with $c$ decreased if necessary.
This proves (ii).

\medskip\noindent\textbf{Step 3: Absorb the weights needed for division.}\proofsteplabel{step:reserve-correction}{3}
The negative exponent $-cm^3\beta^m$ dominates both
$m^2\beta^m\log(m+2)$ and every fixed weight in (iii).
To verify the stated $o(y^P)$ bound, fix $A,N\in\R$, $P\ge0$, and
$\alpha\in\mathbb N_0^3$. On $J_m$,\vspace{-3mm}
\[
y^{-N-P}\le y^{-|N|-P}\le\beta^{(|N|+P)(m+1)},\vspace{-2mm}
\]
and (ii) gives
\[
\begin{aligned}
&e^{A\beta^m+Am^2+Am}y^{-N-P}
\frac{|\partial^\alpha\Box g_\gamma|+|\partial^\alpha R_\gamma|
+\!\sum_h|\partial^\alpha W_h|}{E_\gamma}\\
& \le C(\alpha)\exp\bigl\{-cm^3\beta^m\!
+C(\alpha)m^2\beta^m\log(m\!+\!2) +|A|(\beta^m\!+m^2\!+m)\! +(|N|+P)(m+1)\log\beta\bigr\}\\
& \le C(\alpha)e^{-cm^3\beta^m/2}\longrightarrow0.
\end{aligned}
\]
Thus the weighted quantity in (iii) is $o(y^P)$.
\end{proof}

\begin{proof}[Proof outline of Lemma~\ref{ce:characteristic-solver}]
The operator \eqref{ce:angular-mean} averages over a light-cone section, and \eqref{ce:causal-inverse} is the characteristic solution formula. Differentiating the average proves (iii), and its Fourier multiplier $J_0(\zeta\sqrt q)$ verifies (i). The integration range gives (ii), while differentiation on a fixed square gives (iv).
The complete calculation is given in \hyperref[app:Characteristic]{Appendix A}.
\end{proof}

\subsection{Derivative estimates and the potential}\label{ce:front-estimates}

In this subsection we complete the construction of the pair $(f_\gamma,V_\gamma)$.
We first bound the derivatives of $f_\gamma$ and prove that $f_\gamma$
is bounded below by a fixed fraction of its envelope on the residual
support. This lower bound allows us to define
$V_\gamma=-R_\gamma/f_\gamma$ as a smooth real function. The rapid
decay of $R_\gamma/E_\gamma$ then gives a smooth zero extension of the
potential and all its derivatives across $y=0$.

We retain the amplitudes $A_{\gamma,j}$ and the envelope $E_\gamma$
from \eqref{ce:front-amplitudes}, and the bands
$J_m=[\beta^{-1},\beta]\beta^{-m}$.

\medskip\noindent\textbf{Amplitude derivatives and high-order terms.}
The first step is to control the derivatives of the amplitude
coefficients. The next lemma gives both the fixed-order growth of
$A_{\gamma,j}$ and the decay of the high-order tail of the transport
expansion.
\begin{lemma}[Amplitude asymptotics and packet tails]\label{ce:packet-bounds}
For every integer $m\ge5$, the following estimates hold on the full
band $J_m$. All estimates involving spacetime variables are uniform in
$t\in\R$ and $|z|\le Z+1$.
\begin{enumerate}
\item[(i)] For each fixed $r\ge1$, the amplitude derivatives satisfy
the following estimate after summing over the active packets:
\begin{equation}\label{ce:normal-amplitude}
\sum_{j\ge4}|A_{\gamma,j}^{(r)}-y^{-2r}A_{\gamma,j}|
\le C(r)E_\gamma y^{-2r}m^{2r}\beta^{-m}.
\end{equation}

\item[(ii)] Given an integer derivative order $r\ge0$ and a constant
$B>0$, there exists an integer $L=L(r,B)\ge1$ such that
\[
\left|\partial^r\operatorname{Re}
\sum_{j\ge4}e^{\mathrm i\lambda_ju}
\sum_{\ell=L}^{K_j}\lambda_j^{-\ell}H_{\gamma,j,\ell}\right|
\le C(r,B)E_\gamma e^{-Bm}.
\]
The inner sum is zero when $K_j<L$. The integer $L$ only separates
the terms with $\ell<L$ from those with $\ell\ge L$ in the estimate;
it does not alter the sums defining $\mathcal F_{\gamma,j}$.
\end{enumerate}
\end{lemma}
\begin{proof}[Proof outline of Lemma~\ref{ce:packet-bounds}]
Separate the highest derivative of the common exponential and bound the
remainder by \eqref{ce:normal-amplitude-half-band}. Sum
\eqref{ce:coefficient-estimate} over $\ell\ge L$ to obtain the
prescribed decay $e^{-Bm}$. The index $L$ is used only in the estimates
and does not alter the construction.
The complete calculation is given in \hyperref[app:PacketTail]{Appendix A}.
\end{proof}
\noindent\textbf{Three scales for the corrected solution.}
For $0<y<Y_*$, choose the unique integer $m\ge7$ satisfying
$\beta^{-(m+1)}<y\le\beta^{-m}$, and put
\[
\lambda\triangleq2^m,\qquad k\triangleq y^{-2},\qquad
\tau\triangleq\frac{k^2}{\lambda}=\frac1{\lambda y^4}.
\]
The two active packet frequencies are $\lambda$ and $2\lambda$.
Since $\beta=2^{1/16}$,
\[
\lambda\asymp y^{-16},\qquad k=y^{-2},\qquad
y^{12}\le\tau<2y^{12},\qquad \tau\le\tfrac12.
\]
The table lists the derivative scales, with all derivative orders fixed:
\begin{center}
\begin{tabular}{c@{\qquad}l}
derivative or term & factor relative to $E_\gamma$\\[2pt]
$\partial_u^a\partial_v^b\partial_y^h f_\gamma$
& $C(a,b,h)\lambda^a\tau^bk^h$\\
$\displaystyle\sum_{j=m}^{m+1}|A_{\gamma,j}^{(r)}|$
& $C(r)k^r$\\
$\displaystyle\partial_u^a\partial_y^h
(f_\gamma-\sum_jA_{\gamma,j}\cos(\lambda_ju))$
& $C(a,h)\lambda^a k^h\tau$
\end{tabular}
\end{center}
Thus the $u$ derivatives have the large factor $\lambda$, whereas the
$v$ derivatives carry the small factor $\tau$. The leading cosine
terms do not depend on $v$. The nonzero $v$ derivatives of $f_\gamma$
arise from the higher packet coefficients and the cell corrections.
Also, $C(r)$ may grow with $r$: the fixed-order amplitude bound does
not replace the Gevrey bound \eqref{ce:growing-amplitude} when the
derivative order grows with $K_j$.

The next lemma proves the first and third rows, equivalently
\eqref{ce:front-jets}; the second bound supplies the small factor
$\tau$ in the remainder.

\begin{lemma}[Derivative bounds for the corrected solution]\label{ce:front-derivative-bounds}
For every fixed triple of nonnegative integers $a,b,h$, uniformly in
$t\in\R$, $|z|<Z$, and $0<y<Y_*$, the corrected solution satisfies
\begin{equation}\label{ce:front-jets}
\begin{aligned}
&|\partial_u^a\partial_v^b\partial_y^h f_\gamma|
\le C(a,b,h)E_\gamma\lambda^{a-b}y^{-2h-4b},\\
&\left|\partial_u^a\partial_y^h\!\left(f_\gamma-
\sum_jA_{\gamma,j}\cos(\lambda_ju)\right)\right|
\le C(a,h)E_\gamma\lambda^{a-1}y^{-2h-4}.
\end{aligned}
\end{equation}
\end{lemma}

\begin{proof}[Proof of Lemma~\ref{ce:front-derivative-bounds}]
Fix $a,b,h$ and the interval
$\beta^{-(m+1)}<y\le\beta^{-m}$. We first estimate the terms with
$\ell<L$ for a fixed $L$. We then choose $L$ large enough that
Lemma~\ref{ce:packet-bounds}  controls the
remaining terms; the cell corrections are estimated separately by
Proposition~\ref{ce:reserve}.
Recall from \eqref{ce:approximate-solution},
\eqref{ce:approximate-sum}, and Proposition~\ref{ce:reserve} that
\[
f_\gamma=g_\gamma+\sum_{q\ge5}W_q.
\]
On this interval, only the indices $j=m,m+1$ contribute to $g_\gamma$,
and only $W_m,W_{m+1}$ can be nonzero. Since $\lambda_m=\lambda$ and
$\lambda_{m+1}=2\lambda$, formula~\eqref{ce:coefficient-formula} gives
the real trigonometric expansion
\[
\begin{aligned}
f_\gamma(t,z,y)
&=A_{\gamma,m}(y)\cos(\lambda u)
+A_{\gamma,m+1}(y)\cos(2\lambda u)\\
&\quad+\sum_{\ell=1}^{K_m}
\frac{\cos(\lambda u+\ell\pi/2)}{(2\lambda)^\ell}
\sum_{\substack{1\le d\le\ell\\d\equiv\ell\pmod2}}
\frac d\ell\binom{\ell}{(\ell-d)/2}
\frac{z^d}{d!}A_{\gamma,m}^{(\ell+d)}(y)\\
&\quad+\sum_{\ell=1}^{K_{m+1}}
\frac{\cos(2\lambda u+\ell\pi/2)}{(4\lambda)^\ell}
\sum_{\substack{1\le d\le\ell\\d\equiv\ell\pmod2}}
\frac d\ell\binom{\ell}{(\ell-d)/2}
\frac{z^d}{d!}A_{\gamma,m+1}^{(\ell+d)}(y)\\
&\quad+W_m(t,z,y)+W_{m+1}(t,z,y).
\end{aligned}
\]
Since $Y_*\le\beta^{-7}$ and $m\ge7$, the definitions give
\[
\beta^{2m}\le y^{-2}<\beta^{2m+2},\qquad
2^{-3m/4}\le\frac1{\lambda y^4}<\beta^4 2^{-3m/4}\le\tfrac12,
\qquad \lambda\ge1,\quad 0<y<1.
\]

For every fixed integer $s\ge1$, \eqref{ce:normal-amplitude} gives
\begin{equation}\label{ce:amplitude-fixed-order}
\begin{aligned}
&|A_{\gamma,m}^{(s)}|+|A_{\gamma,m+1}^{(s)}|\\
& \le |A_{\gamma,m}^{(s)}-y^{-2s}A_{\gamma,m}|
+|A_{\gamma,m+1}^{(s)}-y^{-2s}A_{\gamma,m+1}|
+y^{-2s}(A_{\gamma,m}+A_{\gamma,m+1})\\
& \le E_\gamma y^{-2s}\bigl[1+C(s)m^{2s}\beta^{-m}\bigr]
\le C(s)E_\gamma y^{-2s},
\end{aligned}
\end{equation}
where $\sup_{m\ge7}m^{2s}\beta^{-m}<\infty$.
For $s=0$, the corresponding identity is
$A_{\gamma,m}+A_{\gamma,m+1}=E_\gamma$.

\medskip
\noindent\textit{Step 1: Bound the finitely many terms with $\ell<L$.}
\proofsteplabel{step:front-low-orders}{1}
Let $L\ge2$ be fixed for the moment.
The constants in this step may depend on $a,b,h,L$.
For $j=m,m+1$, $1\le\ell<L$, $1\le d\le\ell$, and $b\le d$, we have
\[
\begin{aligned}
&\partial_u^a\partial_v^b\partial_y^h
\left(e^{\mathrm i\lambda_ju}\lambda_j^{-\ell}
\frac{z^d}{d!}A_{\gamma,j}^{(\ell+d)}(y)\!\right) \!=\left(\!-\frac12\right)^b
e^{\mathrm i\lambda_ju}\lambda_j^{-\ell}
\sum_{p=0}^{\min\{a,d-b\}}\!\binom ap
(\mathrm i\lambda_j)^{a-p}2^{-p}
\frac{z^{d-b-p}}{(d\!-\!b\!-\!p)!}
A_{\gamma,j}^{(\ell+d+h)}(y).
\end{aligned}
\]
Here $\partial_uz=1/2$, $\partial_vz=-1/2$, and $\partial_vu=0$.
The derivative is zero if $b>d$.
Since $\lambda\le\lambda_j\le2\lambda$ and $|z|<Z$,
the amplitude bound \eqref{ce:amplitude-fixed-order} implies
\[
\begin{aligned}
\left|
\partial_u^a\partial_v^b\partial_y^h
\left(e^{\mathrm i\lambda_ju}\lambda_j^{-\ell}
\frac{z^d}{d!}A_{\gamma,j}^{(\ell+d)}\right)\right|
& \le C\lambda^{a-\ell}|A_{\gamma,j}^{(\ell+d+h)}|
\le CE_\gamma\lambda^{a-\ell}y^{-2(\ell+d+h)}\\
& \le CE_\gamma\lambda^ay^{-2h}
(\lambda y^4)^{-\ell}
=CE_\gamma\lambda^{a-\ell}y^{-2h-4\ell}.
\end{aligned}
\]
The coefficients in \eqref{ce:coefficient-formula} satisfy
\[
2^{-\ell}\frac d\ell
\binom{\ell}{(\ell-d)/2}\le1
\qquad(1\le d\le\ell,\ d\equiv\ell\pmod2).
\]
Thus, for $1\le\ell<L$,
\[
\left|\partial_u^a\partial_v^b\partial_y^h
\left(e^{\mathrm i\lambda_ju}\lambda_j^{-\ell}
H_{\gamma,j,\ell}\right)\right|
\le CE_\gamma\lambda^ay^{-2h}(\lambda y^4)^{-\ell}
\qquad(j=m,m+1,\ b\le\ell),
\]
and this derivative is zero if $b>\ell$.
For $\ell=0$, $H_{\gamma,j,0}=A_{\gamma,j}$, so
\[
\left|
\partial_u^a\partial_v^b\partial_y^h
\left(e^{\mathrm i\lambda_ju}A_{\gamma,j}\right)\right|
\le C(a,h)E_\gamma\lambda^ay^{-2h}\mathbf1_{\{b=0\}},
\qquad j=m,m+1.
\]
Consequently,
\begin{equation}\label{ce:front-low-order-sum}
\begin{aligned}
&\left|\partial_u^a\partial_v^b\partial_y^h
\operatorname{Re}\left(e^{\mathrm i\lambda_ju}
\sum_{\ell=0}^{\min\{K_j,L-1\}}
\lambda_j^{-\ell}H_{\gamma,j,\ell}\right)\right|\\
& \le CE_\gamma\lambda^ay^{-2h}
\left[\mathbf1_{\{b=0\}}+
\sum_{\ell=\max\{1,b\}}^{L-1}(\lambda y^4)^{-\ell}\right]
\le CE_\gamma\lambda^{a-b}y^{-2h-4b},\qquad j=m,m+1,
\end{aligned}
\end{equation}
because, for every integer $q\ge1$,
\[
\sum_{\ell=q}^{\infty}(\lambda y^4)^{-\ell}
=\frac{(\lambda y^4)^{-q}}{1-(\lambda y^4)^{-1}}
\le2(\lambda y^4)^{-q}.
\]
Repeating the summation in \eqref{ce:front-low-order-sum} with lower
index $\ell=1$ gives
\begin{equation}\label{ce:front-low-order-remainder}
\left|\partial_u^a\partial_y^h
\operatorname{Re}\left(e^{\mathrm i\lambda_ju}
\sum_{\ell=1}^{\min\{K_j,L-1\}}
\lambda_j^{-\ell}H_{\gamma,j,\ell}\right)\right|
\le C(a,h,L)E_\gamma\lambda^{a-1}y^{-2h-4},\qquad j=m,m+1.
\end{equation}

\medskip
\noindent\textit{Step 2: Absorb the tail and the cell corrections.}
\proofsteplabel{step:front-tail}{2}
Set $r=a+b+h$ and $B=(b+1)\log2$.
The absolute-value estimate \eqref{ce:tail-absolute-bound},
applied to the finitely many derivative orders $0,\ldots,r$, gives an
integer $L>r+1$ such that
\begin{equation}\label{ce:front-tail-bound}
\sum_{j\ge4}\sum_{\ell=L}^{K_j}
\max_{|\alpha|\le r}
\left|\partial^\alpha
\left(e^{\mathrm i\lambda_ju}\lambda_j^{-\ell}
H_{\gamma,j,\ell}\right)\right|
\le C(r,B)E_\gamma e^{-Bm}.
\end{equation}
Estimate \eqref{ce:tail-absolute-bound} bounds a geometric series of
absolute values, so the same $L$ can be used for all these derivative
orders. Also, Proposition~\ref{ce:reserve}(ii) gives
\begin{equation}\label{ce:front-correction-bound}
\sum_{q\ge5}\max_{|\alpha|\le r}|\partial^\alpha W_q|
\le C(r)E_\gamma
\exp\{-cm^3\beta^m+C(r)m^2\beta^m\log(m+2)\}
\le C(r,B)E_\gamma e^{-Bm}.
\end{equation}
The last inequality in \eqref{ce:front-correction-bound} follows, for
sufficiently large $m$, from
\[
-cm^3\beta^m+C(r)m^2\beta^m\log(m+2)
\le-\tfrac c2m^3\beta^m\le-Bm;
\]
the finitely many remaining values of $m$ are covered by the constant.
Since
$
\partial_u=\tfrac12(\partial_t+\partial_z)$ and
$\partial_v=\tfrac12(\partial_t-\partial_z)$, 
\eqref{ce:front-tail-bound} and \eqref{ce:front-correction-bound}
also hold, up to constants depending on $r$, for the mixed
derivatives used in \eqref{ce:front-jets}.
Moreover,
\[
e^{-Bm}=2^{-(b+1)m}
\le2^{-3(b+1)m/4}\le(\lambda y^4)^{-(b+1)}
\le\min\{(\lambda y^4)^{-b},(\lambda y^4)^{-1}\}.
\]
The definitions \eqref{ce:approximate-sum} and
\eqref{ce:approximate-solution}, together with
\eqref{ce:front-low-order-sum}, \eqref{ce:front-tail-bound}, and
\eqref{ce:front-correction-bound}, therefore give
\[
\begin{aligned}
|\partial_u^a\partial_v^b\partial_y^hf_\gamma|
&\le C(a,b,h,L)E_\gamma\lambda^{a-b}y^{-2h-4b}
+C(r,B)E_\gamma e^{-Bm} \le C(a,b,h)E_\gamma\lambda^{a-b}y^{-2h-4b}.
\end{aligned}
\]
For the second inequality, choose $b=0$ in Step~\ref{step:front-tail},
so that $r=a+h$ and $B=\log2$.
The identity
\[
f_\gamma-\sum_{j\ge4}A_{\gamma,j}\cos(\lambda_ju)
=\operatorname{Re}\sum_{j\ge4}e^{\mathrm i\lambda_ju}
\sum_{\ell=1}^{K_j}\lambda_j^{-\ell}H_{\gamma,j,\ell}
+\sum_{q\ge5}W_q
\]
together with \eqref{ce:front-low-order-remainder},
\eqref{ce:front-tail-bound}, and \eqref{ce:front-correction-bound},
yields
\[
\begin{aligned}
\left|\partial_u^a\partial_y^h
\left(f_\gamma-\sum_{j\ge4}A_{\gamma,j}\cos(\lambda_ju)\right)\right| &\le C(a,h,L)E_\gamma\lambda^{a-1}y^{-2h-4}
+C(a,h)E_\gamma e^{-m\log2}\\ &
\le C(a,h)E_\gamma\lambda^{a-1}y^{-2h-4}.
\end{aligned}
\]
Both bounds in \eqref{ce:front-jets} follow.
\end{proof}
\noindent\textbf{Leading derivative terms and normal extraction.}
Still working with $t\in\R$, $|z|<Z$, and $0<y<Y_*$, let $m\ge7$ be the
unique integer satisfying $\beta^{-(m+1)}<y\le\beta^{-m}$, and retain
$\lambda=2^m$, $u=t+z$, and $v=t-z$.
Write $a_0\triangleq A_{\gamma,m}$ and $a_1\triangleq A_{\gamma,m+1}$;
primes denote derivatives in $y$.

The derivative bounds alone do not show that $f_\gamma$ is positive on
the residual support. We now keep the leading terms explicitly.
The $u$ derivatives come from the leading cosines. A $v$ derivative
annihilates those cosines, so its first contribution comes from
$\ell=1$; two $v$ derivatives first see $\ell=2$. For the normal
derivative, $\partial_y-y^{-2}$ cancels the derivative of the common
factor $e^{-1/y}$. The packet cutoff derivatives remain, together with
a smaller contribution from the corrections.
\begin{lemma}[Leading derivative terms]\label{ce:leading-derivatives}
Let $t\in\R$, $|z|<Z$, and $0<y<Y_*$, and let $m\ge7$ be the unique
integer such that $\beta^{-(m+1)}<y\le\beta^{-m}$. Put
$\lambda=2^m$, $u=t+z$, and $v=t-z$, and write
$a_0\triangleq A_{\gamma,m}$ and $a_1\triangleq A_{\gamma,m+1}$,
with primes denoting derivatives in $y$. Then, on the spacetime region
of Lemma~\ref{ce:front-derivative-bounds}, the following expansions
hold.

\smallskip
\noindent\emph{Characteristic derivatives.}
The derivatives along $u$ and $v$ satisfy
\begin{equation}\label{ce:front-extraction}
\begin{cases}
\ds f_{\gamma,u}=-\lambda(a_0\sin(\lambda u)+2a_1\sin(2\lambda u))
+O(E_\gamma y^{-4}), \\[2pt]
\ns\ds f_{\gamma,v}=\frac1{4\lambda}
(a_0''\sin(\lambda u)+\tfrac12a_1''\sin(2\lambda u))
+O(E_\gamma\lambda^{-2}y^{-8}),\\[2pt]
\ns\ds f_{\gamma,vv}=-\frac1{16\lambda^2}
(a_0^{(4)}\cos(\lambda u)+\tfrac14a_1^{(4)}\cos(2\lambda u))
+O(E_\gamma\lambda^{-3}y^{-12}).
\end{cases}
\end{equation}

\smallskip
\noindent\emph{Normal derivative.}
The normal derivative, after subtracting the contribution of the
common exponential factor, is
\begin{equation}\label{ce:normal-extraction}
\begin{aligned}
f_{\gamma,y}-y^{-2}f_\gamma
&=e^{-1/y-m^2-\gamma m}\beta^m\rho'(\beta^my)\cos(\lambda u)\\
&\quad+e^{-1/y-(m+1)^2-\gamma(m+1)}
\beta^{m+1}\rho'(\beta^{m+1}y)\cos(2\lambda u)
+O(E_\gamma\lambda^{-1}y^{-6}).
\end{aligned}
\end{equation}
\end{lemma}
\begin{proof}[Proof of Lemma~\ref{ce:leading-derivatives}]
The lemma has two parts: the characteristic derivatives along $u$ and
$v$, and the normal derivative after subtracting the contribution of
the common exponential factor. We prove them in this order.

\medskip
\noindent\textit{Part 1: The $u$ and $v$ derivatives.}

\smallskip
\noindent\emph{The $u$ derivative.}
By \eqref{ce:front-jets},
\[
\begin{aligned}
\left|f_{\gamma,u}
+\lambda\bigl(a_0\sin(\lambda u)+2a_1\sin(2\lambda u)\bigr)\right|
&=\left|\partial_u\Big(f_\gamma-
\sum_{j=m}^{m+1}A_{\gamma,j}\cos(\lambda_ju)\Big)\right| \le CE_\gamma\lambda^{1-1}y^{-4}=CE_\gamma y^{-4}.
\end{aligned}
\]
This gives the first line of \eqref{ce:front-extraction}: the
derivative in $u$ is dominated by the leading cosines, and the
remainder is $O(E_\gamma y^{-4})$.

\smallskip
\noindent\emph{The $v$ derivatives.}
The leading cosines do not depend on $v$, so their $v$ derivatives
vanish. The first nonzero contribution comes from the transport
coefficient with $\ell=1$, and the second from $\ell=2$. For
$j=m,m+1$, \eqref{ce:coefficient-formula} and $\partial_vz=-1/2$,
$\partial_vu=0$ give
\[
\begin{aligned}
&\operatorname{Re}\bigl(e^{\mathrm i\lambda_ju}
\lambda_j^{-1}H_{\gamma,j,1}\bigr)
=-\frac{zA_{\gamma,j}''}{2\lambda_j}\sin(\lambda_ju),\qquad
\partial_v\operatorname{Re}\bigl(e^{\mathrm i\lambda_ju}
\lambda_j^{-1}H_{\gamma,j,1}\bigr)
=\frac{A_{\gamma,j}''}{4\lambda_j}\sin(\lambda_ju),\\
&\operatorname{Re}\bigl(e^{\mathrm i\lambda_ju}
\lambda_j^{-2}H_{\gamma,j,2}\bigr)
=-\frac{z^2A_{\gamma,j}^{(4)}}{8\lambda_j^2}\cos(\lambda_ju),\qquad
\partial_v^2\operatorname{Re}\bigl(e^{\mathrm i\lambda_ju}
\lambda_j^{-2}H_{\gamma,j,2}\bigr)
=-\frac{A_{\gamma,j}^{(4)}}{16\lambda_j^2}\cos(\lambda_ju).
\end{aligned}
\]
For general $b=1,2$ and $1\le\ell\le K_j$,
\[
\begin{aligned}
&\partial_v^b\operatorname{Re}\bigl(e^{\mathrm i\lambda_ju}
\lambda_j^{-\ell}H_{\gamma,j,\ell}\bigr) =\frac{(-1)^b\cos(\lambda_ju+\ell\pi/2)}
{2^{\ell+b}\lambda_j^\ell}
\sum_{\substack{\max\{1,b\}\le d\le\ell\\d\equiv\ell\pmod2}}
\frac d\ell\binom\ell{(\ell-d)/2}
\frac{z^{d-b}}{(d-b)!}A_{\gamma,j}^{(\ell+d)}.
\end{aligned}
\]
In particular,
\[
\partial_v^b\operatorname{Re}\bigl(e^{\mathrm i\lambda_ju}
\lambda_j^{-\ell}H_{\gamma,j,\ell}\bigr)=0
\qquad(0\le\ell<b).
\]
Thus only the levels $\ell\ge b$ contribute to the $b$-th $v$
derivative.

We now separate the levels below a fixed threshold from the remaining
tail. Let $L\in\mathbb N$ satisfy
\[
L\ge\max\left\{4,\max_{r=0,1,2}L(r,3\log2)\right\},
\]
where $L(r,3\log2)$ is from
Lemma~\ref{ce:packet-bounds}\ref{cond:packet-tail}.
Step~\ref{step:front-low-orders} of the proof of
Lemma~\ref{ce:front-derivative-bounds} gives, for
$b\le\ell\le\min\{K_j,L-1\}$,
\[
\begin{aligned}
\left|\partial_v^b\operatorname{Re}\bigl(e^{\mathrm i\lambda_ju}
\lambda_j^{-\ell}H_{\gamma,j,\ell}\bigr)\right|
& \le C(L)\lambda^{-\ell}\sum_{d=b}^{\ell}
|A_{\gamma,j}^{(\ell+d)}|
\le C(L)E_\gamma\lambda^{-\ell}
\sum_{d=b}^{\ell}y^{-2(\ell+d)}\\
& \le C(L)E_\gamma\lambda^{-\ell}y^{-4\ell}
=C(L)E_\gamma(\lambda y^4)^{-\ell}.
\end{aligned}
\]
The absolute-value estimates in the proof of
Lemma~\ref{ce:packet-bounds} give
\[
\partial_v^b
=2^{-b}(\partial_t-\partial_z)^b
=2^{-b}\sum_{a=0}^b(-1)^a\binom ba
\partial_t^{b-a}\partial_z^a,
\]
and hence
\[
\begin{aligned}
\sum_{j=m}^{m+1}\sum_{\ell=L}^{K_j}
\left|\partial_v^b\bigl(e^{\mathrm i\lambda_ju}
\lambda_j^{-\ell}H_{\gamma,j,\ell}\bigr)\right|
&\le2^{-b}\sum_{a=0}^b\binom ba
\sum_{j=m}^{m+1}\sum_{\ell=L}^{K_j}
\left|\partial_t^{b-a}\partial_z^a
\bigl(e^{\mathrm i\lambda_ju}\lambda_j^{-\ell}
H_{\gamma,j,\ell}\bigr)\right|\\
& \le CE_\gamma e^{-3m\log2}=CE_\gamma2^{-3m}.
\end{aligned}
\]
By Proposition~\ref{ce:reserve}(ii),
\[
\sum_{q\ge5}|\partial_v^bW_q|
\le CE_\gamma
\exp\{-cm^3\beta^m+Cm^2\beta^m\log(m+2)\}
\le CE_\gamma2^{-3m}.
\]
Combining the last two estimates, the tail of the transport expansion
and the cell corrections together satisfy
\[
\begin{aligned}
&\left|\partial_v^b\left(
\operatorname{Re}\sum_{j=m}^{m+1}e^{\mathrm i\lambda_ju}
\sum_{\ell=L}^{K_j}\lambda_j^{-\ell}H_{\gamma,j,\ell}
+\sum_{q\ge5}W_q\right)\right|\\
& \le\sum_{j=m}^{m+1}\sum_{\ell=L}^{K_j}
\left|\partial_v^b\bigl(e^{\mathrm i\lambda_ju}
\lambda_j^{-\ell}H_{\gamma,j,\ell}\bigr)\right|
+\sum_{q\ge5}|\partial_v^bW_q|
\le CE_\gamma2^{-3m}
\le CE_\gamma(\lambda y^4)^{-3}.
\end{aligned}
\]

We can now isolate the leading level. For $b=1,2$,
\[
\begin{aligned}
&\partial_v^bf_\gamma
-\operatorname{Re}\sum_{j=m}^{m+1}\partial_v^b
\bigl(e^{\mathrm i\lambda_ju}\lambda_j^{-b}H_{\gamma,j,b}\bigr)\\
& =\operatorname{Re}\sum_{j=m}^{m+1}
\sum_{\ell=b+1}^{\min\{K_j,L-1\}}
\partial_v^b\bigl(e^{\mathrm i\lambda_ju}
\lambda_j^{-\ell}H_{\gamma,j,\ell}\bigr) +\partial_v^b\left(
\operatorname{Re}\sum_{j=m}^{m+1}e^{\mathrm i\lambda_ju}
\sum_{\ell=L}^{K_j}\lambda_j^{-\ell}H_{\gamma,j,\ell}
+\sum_{q\ge5}W_q\right),
\end{aligned}
\]
and therefore
\[
\begin{aligned}
&\left|\partial_v^bf_\gamma
-\operatorname{Re}\sum_{j=m}^{m+1}\partial_v^b
\bigl(e^{\mathrm i\lambda_ju}\lambda_j^{-b}H_{\gamma,j,b}\bigr)\right|\\
& \le C(L)E_\gamma\sum_{\ell=b+1}^{L-1}(\lambda y^4)^{-\ell}
+CE_\gamma2^{-3m}\\
& \le CE_\gamma\left(
\frac{(\lambda y^4)^{-b-1}}{1-(\lambda y^4)^{-1}}
+(\lambda y^4)^{-3}\right)
\le CE_\gamma(\lambda y^4)^{-b-1}
=CE_\gamma\lambda^{-b-1}y^{-4b-4},
\end{aligned}
\]
where $(\lambda y^4)^{-1}\le1/2$. Thus the remainder in the $b$-th
$v$ derivative is smaller than the leading level by a factor
$(\lambda y^4)^{-1}$.

It remains to identify the leading contributions explicitly. For
$j=m,m+1$,
\[
\begin{aligned}
&\operatorname{Re}\sum_{j=m}^{m+1}\partial_v
\bigl(e^{\mathrm i\lambda_ju}\lambda_j^{-1}H_{\gamma,j,1}\bigr)
=\frac1{4\lambda}
\bigl(a_0''\sin(\lambda u)+\tfrac12a_1''\sin(2\lambda u)\bigr),\\
&\operatorname{Re}\sum_{j=m}^{m+1}\partial_v^2
\bigl(e^{\mathrm i\lambda_ju}\lambda_j^{-2}H_{\gamma,j,2}\bigr)
=-\frac1{16\lambda^2}
\bigl(a_0^{(4)}\cos(\lambda u)+\tfrac14a_1^{(4)}\cos(2\lambda u)\bigr).
\end{aligned}
\]
Combining these identities with the remainder estimate above gives
the second and third lines of \eqref{ce:front-extraction}. This
completes Part 1.

\medskip
\noindent\textit{Part 2: The normal derivative.}
For the leading cosine sum, the cancellation of $\partial_y e^{-1/y}$
is exact. For the complete function, we retain the correction
remainder. Specifically,
\[
\begin{aligned}
f_{\gamma,y}-y^{-2}f_\gamma
&=\sum_{j=m}^{m+1}
\bigl(A_{\gamma,j}'-y^{-2}A_{\gamma,j}\bigr)\cos(\lambda_ju) +(\partial_y-y^{-2})\left(f_\gamma-
\sum_{j=m}^{m+1}A_{\gamma,j}\cos(\lambda_ju)\right),
\end{aligned}
\]
where
$
A_{\gamma,j}'-y^{-2}A_{\gamma,j}
=e^{-1/y-j^2-\gamma j}\beta^j\rho'(\beta^jy).
$
The first sum is the leading normal contribution displayed in
\eqref{ce:normal-extraction}. For the remainder, \eqref{ce:front-jets}
gives
\[
\begin{aligned}
&\bigg|(\partial_y-y^{-2})\Big(f_\gamma-
\sum_{j=m}^{m+1}A_{\gamma,j}\cos(\lambda_ju)\Big)\bigg|\\
& \le\bigg|\partial_y\Big(f_\gamma-
\sum_{j=m}^{m+1}A_{\gamma,j}\cos(\lambda_ju)\Big)\bigg|
+y^{-2}\bigg|f_\gamma-
\sum_{j=m}^{m+1}A_{\gamma,j}\cos(\lambda_ju)\bigg|
\le CE_\gamma\lambda^{-1}y^{-6}.
\end{aligned}
\]
This proves \eqref{ce:normal-extraction} and completes Part 2.
\end{proof}
To put the leading terms on a common scale, apply
\eqref{ce:normal-amplitude} with $r=2$ and
\eqref{ce:weighted-cutoff-derivatives} with $q=1$:
\begin{equation}\label{ce:common-null-amplitude-bounds}
\begin{aligned}
\sum_{j=m}^{m+1}|A_{\gamma,j}''-y^{-4}A_{\gamma,j}|
&\le CE_\gamma y^{-4}m^4\beta^{-m},\\
e^{-1/y}\sum_{j=m}^{m+1}e^{-j^2-\gamma j}
\beta^j|\rho'(\beta^jy)|
&\le CE_\gamma m^2\beta^m
\le CE_\gamma y^{-2}m^4\beta^{-m}.
\end{aligned}
\end{equation}
The last inequality in \eqref{ce:common-null-amplitude-bounds} uses
$y^{-2}\ge\beta^{2m}$ and $m^2\le m^4$.
Substitution into \eqref{ce:front-extraction} and
\eqref{ce:normal-extraction}, respectively, gives
\begin{equation}\label{ce:common-null}
\begin{aligned}
f_{\gamma,v}&=\frac1{4\lambda y^4}
(a_0\sin(\lambda u)+\tfrac12a_1\sin(2\lambda u))
+O(E_\gamma\lambda^{-1}y^{-4}(m^4\beta^{-m}+\lambda^{-1}y^{-4})),\\
f_{\gamma,y}&=y^{-2}(a_0\cos(\lambda u)+a_1\cos(2\lambda u))
+O(E_\gamma y^{-2}(m^4\beta^{-m}+\lambda^{-1}y^{-4})).
\end{aligned}
\end{equation}

\medskip\noindent\textbf{Choose one width for the division region.}
All estimates so far hold on $0<y<Y_{\rm geom}$, where the corrected
solution is already fixed. Choose $C_*=C_*(Z)\ge1$ to dominate the
error constants in the second inequality of \eqref{ce:front-jets}
for $a=0,2,4$ and $h=0$, and in the estimate for $f_{\gamma,y}$ in
\eqref{ce:common-null}, for both values of $\gamma$.
This constant is fixed on the present cylinder before any further
restriction. Define
\[
\Phi(s)\triangleq\beta s\bigl(\log_\beta(1/s)\bigr)^4+2s^{12},
\qquad 0<s\le\beta^{-7}.
\]
On $\beta^{-(m+1)}<y\le\beta^{-m}$, we have
\[
m\le\log_\beta(1/y),\qquad \beta^{-m}<\beta y,\qquad
\frac1{\lambda y^4}=\frac{\beta^{-16m}}{y^4}
<\beta^{16}y^{12}=2y^{12}.
\]
Consequently,
\[
m^4\beta^{-m}+\frac1{\lambda y^4}
\le\beta y\bigl(\log_\beta(1/y)\bigr)^4+2y^{12}
=\Phi(y).
\]
Since $\Phi(s)\to0$ as $s\downarrow0$, we can now make the final
choice of the normal width:
\begin{equation}\label{ce:final-front-width}
0<Y_*\le Y_{\rm geom},\qquad
C_*\sup_{0<s\le Y_*}\Phi(s)\le\frac18.
\end{equation}
The margin $1/8$ leaves $\cos(\pi/12)-1/8>5/8$, which will give the
lower bound needed for division. Fix one such $Y_*=Y_*(Z)$ for the
rest of the proof, and restrict the same corrected solution to
$0<y<Y_*$. This collects the geometric support condition and the
finite-order smallness condition.
All previously established constants remain valid. For a higher
derivative order, Lemma~\ref{ce:packet-bounds}\ref{cond:packet-tail}
may require a larger $L$ and the estimates may require larger
constants, but neither changes the function or decreases $Y_*$.

\medskip\noindent\textbf{The division region.}
To locate the division region, recall from \eqref{ce:correction-scales}
that $\Lambda_h=2^{h-1}$ is the frequency of the periodic cutoff
$\eta_h(u)=\eta(\Lambda_hu)$, whose period in $u$ is
$d_h=2\pi/\Lambda_h$. Formula~\eqref{ce:complete-front} shows that
the residual is supported where
$\operatorname{dist}(\Lambda_hu,2\pi\Z)\le\pi/72$ and $y\in J_h$
for some correction index $h\ge5$. These support bands are locally
finite in $y>0$, so the union of the indicated closed phase sets is
relatively closed. The following lemma gives uniform bounds on a
slightly larger phase region.

\begin{lemma}[Bounds near the residual support]\label{ce:division-region}
For the final width chosen in \eqref{ce:final-front-width}, there is
a constant $c>0$ with the following property.
Let $t\in\R$, $|z|<Z$, and $0<y<Y_*$, and choose
the unique integer $m\ge7$ satisfying
$\beta^{-(m+1)}<y\le\beta^{-m}$. Put $\lambda=2^m$ and $u=t+z$.
If, for some integer $h\ge5$,
\[
y\in J_h,\qquad
\operatorname{dist}(\Lambda_hu,2\pi\Z)\le\pi/48,
\qquad \Lambda_h=2^{h-1},
\]
then
\begin{equation}\label{ce:source-margin}
f_\gamma\ge\tfrac58E_\gamma,\qquad
f_{\gamma,y}\ge cE_\gamma y^{-2},\qquad
f_{\gamma,uu}\le-cE_\gamma\lambda^2,\qquad
f_{\gamma,uuuu}\ge cE_\gamma\lambda^4.
\end{equation}
In particular, \eqref{ce:source-margin} holds on
$\operatorname{supp}R_\gamma$.
\end{lemma}

\begin{proof}[Proof of Lemma~\ref{ce:division-region}]
If $A_{\gamma,j}(y)\ne0$ and $y\in J_h$, then
$j\in\{h-1,h,h+1\}$. Consequently,
\[
\frac{\lambda_j}{\Lambda_h}=2^{j-h+1}\in\{1,2,4\},\qquad
\operatorname{dist}(\lambda_ju,2\pi\Z)
\le2^{j-h+1}\operatorname{dist}(\Lambda_hu,2\pi\Z)
\le\pi/12.
\]
Thus every nonzero leading amplitude has a cosine factor at least
$\cos(\pi/12)$. Since
$A_{\gamma,m},A_{\gamma,m+1}\ge0$ and
$A_{\gamma,m}+A_{\gamma,m+1}=E_\gamma$,
\[
A_{\gamma,m}\cos(\lambda u)
+2^{2q}A_{\gamma,m+1}\cos(2\lambda u)
\ge E_\gamma\cos(\pi/12),\qquad q=0,1,2.
\]
The second estimate in \eqref{ce:front-jets}, for $a=0,2,4$ and no
$y$ derivatives, and the estimate for $f_{\gamma,y}$ in
\eqref{ce:common-null}, with the common constant $C_*$ fixed before
\eqref{ce:final-front-width}, give
\[
\begin{aligned}
f_\gamma
&\ge E_\gamma\left(\cos(\pi/12)-\frac{C_*}{\lambda y^4}\right),\\
f_{\gamma,y}
&\ge E_\gamma y^{-2}\left[\cos(\pi/12)
-C_*\left(m^4\beta^{-m}+\frac1{\lambda y^4}\right)\right],\\
-f_{\gamma,uu}
&\ge E_\gamma\lambda^2
\left(\cos(\pi/12)-\frac{C_*}{\lambda y^4}\right),\\
f_{\gamma,uuuu}
&\ge E_\gamma\lambda^4
\left(\cos(\pi/12)-\frac{C_*}{\lambda y^4}\right).
\end{aligned}
\]
The final choice \eqref{ce:final-front-width} gives
\[
C_*\left(m^4\beta^{-m}+\frac1{\lambda y^4}\right)
\le C_*\Phi(y)
\le\tfrac18.
\]
Since $\cos(\pi/12)>3/4$, the bounds in \eqref{ce:source-margin}
follow with $c=1/2$.
Proposition~\ref{ce:reserve}(i) places
$\operatorname{supp}R_\gamma$ in the smaller phase region
$\operatorname{dist}(\Lambda_hu,2\pi\Z)\le\pi/72$, so
\eqref{ce:source-margin} holds on that support.
\end{proof}

\begin{proof}[Proof of Theorem~\ref{ce:front}]
We construct $V_\gamma$ and verify that both $f_\gamma$ and $V_\gamma$
extend flatly by zero at $y=0$, that
$f_\gamma(0,0,y)=E_\gamma(y)>0$ for $y>0$, and that the pair is
$2\pi$-periodic. The proof has four steps: define the potential on
$y>0$, extend it flatly across $y=0$, extend the solution and retain a
nonzero value, and verify periodicity.

\medskip
\noindent\textbf{Step 1: Define a smooth potential for $y>0$.}
\proofsteplabel{step:front-potential}{1}
Set $Q_+\triangleq\R\times(-Z,Z)\times(0,Y_*)$, and take supports
relative to this cylinder. The two open sets
$
U_{\rm div}\triangleq\{f_\gamma>E_\gamma/2\}$ and
$U_{\rm zero}\triangleq Q_+\setminus\operatorname{supp}R_\gamma
$
cover $Q_+$, because Lemma~\ref{ce:division-region} gives
$\operatorname{supp}R_\gamma\subset U_{\rm div}$. Define
\[
V_\gamma=
\begin{cases}
-R_\gamma/f_\gamma,&\text{on }U_{\rm div},\\
0,&\text{on }U_{\rm zero}.
\end{cases}
\]
On the overlap, $R_\gamma=0$ and $f_\gamma>E_\gamma/2>0$, so
$-R_\gamma/f_\gamma=0$. The definitions therefore give one smooth
real function. Moreover, $\Box f_\gamma=R_\gamma$, and hence
\[
\begin{cases}
\ds\Box f_\gamma+V_\gamma f_\gamma
=R_\gamma-(R_\gamma/f_\gamma)f_\gamma=0
&\text{on }U_{\rm div},\\
\ns\ds \Box f_\gamma+V_\gamma f_\gamma =R_\gamma=0
& \text{on }U_{\rm zero}.
\end{cases}
\]
Thus
\begin{equation}\label{ce:front-potential-equation}
\Box f_\gamma+V_\gamma f_\gamma=0.
\end{equation}
Since $f_\gamma>0$ on $\operatorname{supp}R_\gamma$, the definitions give
\[
\{V_\gamma\ne0\}=\{R_\gamma\ne0\},\qquad
\operatorname{supp}V_\gamma
=\overline{\{V_\gamma\ne0\}}
=\overline{\{R_\gamma\ne0\}}
=\operatorname{supp}R_\gamma.
\]

\medskip
\noindent\textbf{Step 2: Extend the potential flatly across $y=0$.}
\proofsteplabel{step:front-potential-flatness}{2}
We prove $|\partial^\alpha V_\gamma|\le C(\alpha,N)y^N$ for every
fixed $\alpha$ and $N$. The factor $E_\gamma$ in the residual will
cancel the factor $E_\gamma^{-1}$ from differentiating
$f_\gamma^{-1}$. By Proposition~\ref{ce:reserve}, on the support band
$J_m$, for every fixed integer $r\ge0$,
\[
|\partial^rR_\gamma|\le C(r)E_\gamma
\exp\{-cm^3\beta^m+C(r)m^2\beta^m\log(m+2)\}.
\]
Since $2^m\asymp y^{-16}$ on $J_m$, for each fixed multi-index $\mu$,
\eqref{ce:front-jets} gives an exponent $M(\mu)\ge0$ with
\[
|\partial^\mu f_\gamma|\le C(\mu)E_\gamma y^{-M(\mu)}.
\]
On $U_{\rm div}=\{f_\gamma>E_\gamma/2\}$, differentiating
$f_\gamma f_\gamma^{-1}=1$ gives, for $|\nu|\ge1$,
\begin{equation}\label{ce:inverse-derivative-recursion}
\partial^\nu(f_\gamma^{-1})
=-f_\gamma^{-1}\sum_{0<\mu\le\nu}\binom\nu\mu
(\partial^\mu f_\gamma)\partial^{\nu-\mu}(f_\gamma^{-1}).
\end{equation}
Set $N(0)=0$ and, inductively,\vspace{-2mm}
\[
N(\nu)\triangleq\max_{0<\mu\le\nu}
\{M(\mu)+N(\nu-\mu)\}.\vspace{-2mm}
\]
The bound $|f_\gamma^{-1}|\le2/E_\gamma$ starts the induction, and
identity \eqref{ce:inverse-derivative-recursion} yields
\[
\begin{aligned}
|\partial^\nu(f_\gamma^{-1})|
&\le\frac2{E_\gamma}\sum_{0<\mu\le\nu}\binom\nu\mu
C(\mu)E_\gamma y^{-M(\mu)}
C(\nu-\mu)E_\gamma^{-1}y^{-N(\nu-\mu)} \le C(\nu)E_\gamma^{-1}y^{-N(\nu)}.
\end{aligned}
\]
For each fixed multi-index $\alpha$, on this open set,
\[
\begin{aligned}
|\partial^\alpha V_\gamma|
&\le\sum_{\nu\le\alpha}\binom\alpha\nu
|\partial^{\alpha-\nu}R_\gamma|
|\partial^\nu(f_\gamma^{-1})| \le C(\alpha)y^{-N(\alpha)}
\exp\{-cm^3\beta^m+C(\alpha)m^2\beta^m\log(m+2)\},
\end{aligned}
\]
where $N(\alpha)\ge0$ depends only on $\alpha$ and the fixed
parameters.
For every fixed integer $N\ge0$, the exponent satisfies
\[
\frac{C(\alpha)m^2\beta^m\log(m+2)
+(N(\alpha)+N)(m+1)\log\beta}{m^3\beta^m}
\longrightarrow0.
\]
Since $y^{-1}\le\beta^{m+1}$ on $J_m$, it follows, for all
sufficiently large $m$, that
\[
\begin{aligned}
y^{-N}|\partial^\alpha V_\gamma|
&\le C(\alpha)\exp\bigl\{-cm^3\beta^m
+C(\alpha)m^2\beta^m\log(m+2) +(N(\alpha)+N)(m+1)\log\beta\bigr\}\\
&\le C(\alpha)e^{-cm^3\beta^m/2}\le C(\alpha).
\end{aligned}
\]
Increasing the constant covers the finitely many remaining bands.
Thus
\begin{equation}\label{ce:potential-flatness}
|\partial^\alpha V_\gamma|\le C(\alpha,N)y^N
\qquad(0<y<Y_*,\ |z|<Z).
\end{equation}
Estimate \eqref{ce:potential-flatness} also holds where $V_\gamma$ is
identically zero. Hence $V_\gamma$ extends smoothly and flatly by zero
to $y\le0$.

\medskip
\noindent\textbf{Step 3: Extend the solution and retain a nonzero value.}
\proofsteplabel{step:front-solution-flatness}{3}
By \eqref{ce:packet-envelope} and \eqref{ce:front-jets}, for every
multiindex $\alpha$ there is a fixed $M(\alpha)>0$ such that, on
$\beta^{-(m+1)}<y\le\beta^{-m}$,
\[
\begin{aligned}
y^{-N}|\partial^\alpha f_\gamma|
&\le C(\alpha)y^{-N-M(\alpha)}E_\gamma(y)
\le C(\alpha)\beta^{(N+M(\alpha))(m+1)}
e^{-\beta^m-m^2-\gamma m}\rightarrow0
\qquad(N\in\mathbb N_0).
\end{aligned}
\]
Thus $f_\gamma$ also extends flatly by zero to $y\le0$.
Equation \eqref{ce:front-potential-equation} holds across $y=0$ by
smoothness. After zero extension, $V_\gamma=R_\gamma=0$ for $y\le0$,
so $\{V_\gamma\ne0\}=\{R_\gamma\ne0\}$ still holds, and their supports
also agree in $\R\times(-Z,Z)\times(-Y_*,Y_*)$.
At $(t,z)=(0,0)$, all coefficients $H_{\gamma,j,\ell}$ with $\ell\ge1$
vanish because they contain a positive power of $z$, while all cell
corrections vanish because the periodic cutoff is zero at $u=0$.
Therefore
\[
f_\gamma(0,0,y)=\sum_{j\ge4}A_{\gamma,j}(y)
=E_\gamma(y)>0\qquad(0<y<Y_*).
\]

\medskip
\noindent\textbf{Step 4: Verify periodicity.}
\proofsteplabel{step:front-periodicity}{4}
The oscillation frequencies $\lambda_j=2^j$ are integers. By
Proposition~\ref{ce:reserve}, each cell correction has period
$d_m=2\pi/\Lambda_m$, where $\Lambda_m=2^{m-1}\in\mathbb N$.
Hence $f_\gamma$ is $2\pi$-periodic in $t$, and so is $V_\gamma$ by
its definition. The construction applies to both values of $\gamma$.
Taking $(f,V)=(f_0,V_0)$ proves Theorem~\ref{ce:front}.
\end{proof} 

\section{Observation geometry and gluing criteria}\label{ce:preliminaries}

\begin{figure}[H]
\centering
\begin{tikzpicture}[
block/.style={draw=black!70,rounded corners=1.5pt,fill=black!3,
align=center,text width=3.8cm,minimum height=1.2cm,inner sep=4pt,font=\small},
flow/.style={-{Stealth[length=1.6mm]},line width=.5pt}]
\node[block] (geometry) at (0,0)
{Scaled GCC polygons\\Subsection~\ref{ce:geometry-section}};
\node[block] (edges) at (5.1,0)
{Local solutions along edges\\rotation and translation\\Section~\ref{ce:front-section}};
\node[block] (vertices) at (10.2,0)
{Merge at the vertices\\Proposition~\ref{ce:merging-criterion}};
\node[block] (static) at (10.2,-2)
{Join to a static pair\\Proposition~\ref{ce:activation-criterion}};
\node[block] (exterior) at (5.1,-2)
{Glue to the Dirichlet\\exterior pair};
\node[block] (counterexample) at (0,-2)
{Global counterexample\\Theorem~\ref{thm:counter}};
\draw[flow] (geometry.east) -- (edges.west);
\draw[flow] (edges.east) -- (vertices.west);
\draw[flow] (vertices.south) -- (static.north);
\draw[flow] (static.west) -- (exterior.east);
\draw[flow] (exterior.west) -- (counterexample.east);
\end{tikzpicture}
\caption{Road map for the proof of Theorem~\ref{thm:counter}.}
\label{fig:counterexample-roadmap}
\end{figure}
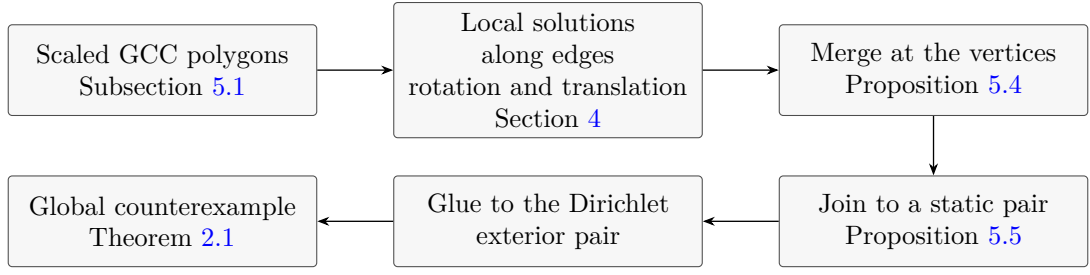

This section provides the geometric and analytic ingredients for the
proof of Theorem~\ref{thm:counter}. The construction has four stages,
summarized in Figure~\ref{fig:counterexample-roadmap}. First, we select
a polygonal observation set satisfying GCC, and we prove that the free
wave is observable on it. Second, we establish a finite-chain repair
lemma that corrects the equation without changing the zero set of the
wave. Third, we merge two edge solutions at a vertex, using the repair
lemma to handle the region where both solutions are present. Fourth, we
join the resulting local pair to a positive static solution, so that it
can be matched to a Dirichlet exterior pair. The verification of these
criteria and the final assembly are carried out in
Section~\ref{ce:completion-section}.
Throughout this section, $\Box=\partial_t^2-\Delta_x$. The repair in
Lemma~\ref{ce:repair} uses a positive multiplier
$p=e^{\Phi_{\rm rep}}q$, which preserves the zero set because
$e^{\Phi_{\rm rep}}>0$. This positivity is the reason why the repair
can be applied at the zeros of the oscillatory wave.

\subsection{Construction of observation geometry}\label{ce:geometry-section}

With the local solutions constructed in Section~\ref{ce:front-section},
we now construct the observation geometry for the global counterexample
in Theorem~\ref{thm:counter}. The construction has two tasks: rays must
leave a thin boundary collar, and the observation polygons must
intercept every sufficiently long interior segment.

We begin with a bounded connected domain whose boundary is analytic and
contains a single concave indentation:
\begin{equation}\label{ce:flower}
\Omega_0\triangleq
\{(r\cos\theta,r\sin\theta):0\le r<1+\tfrac23\cos\theta\}.
\end{equation}
Let $d(x)\triangleq\operatorname{dist}(x,\partial\Omega_0)$.
For $0<\delta<1/12$, define the triangular lattice with spacing $\delta/10$ by
\[
\Lambda_\delta\triangleq
\left\{\frac{\delta}{10}\left(m+\tfrac n2,\tfrac{\sqrt3\,n}{2}\right):
m,n\in\Z\right\}.
\]
Joining nearest neighbours gives equilateral triangles of side $\delta/10$.
Next, define the regular $64$-gon centred at the origin, with
inradius $12/25$, by
\begin{equation}\label{ce:lattice-polygon}
Q_{64}\triangleq\bigcap_{k=0}^{63}
\{x\in\R^2:x\cdot n_k\le12/25\},\qquad
n_k\triangleq\left(\cos\tfrac{k\pi}{32},\sin\tfrac{k\pi}{32}\right),
\end{equation}
where $n_k$ are the outward unit normals to its sides. Each inequality
specifies the half-plane containing the origin.
Retain the lattice centres at distance at least $2\delta$ from the
boundary, and define their polygons and the observation set by
\begin{equation}\label{ce:polygons}
\begin{aligned}
&\mathcal J_\delta\triangleq
\{q\in\Lambda_\delta\cap\Omega_0:d(q)\ge2\delta\},\qquad P_q^0\triangleq q+\frac{\delta}{10}Q_{64}\quad(q\in\mathcal J_\delta),\qquad
\omega_0\triangleq
\bigcup_{q\in\mathcal J_\delta}\operatorname{int}P_q^0.
\end{aligned}
\end{equation}
Each observation polygon has inradius $r_\delta\triangleq6\delta/125$ and
exterior angle $\alpha_0\triangleq\pi/32$. We enumerate the family as
$P_1^0,\ldots,P_J^0$, where $J\triangleq\#\mathcal J_\delta$.

\begin{figure}[H]
\centering
\begin{tikzpicture}[x=1.55cm,y=1.55cm,font=\small,line join=round,
>={Stealth[length=1.5mm]}]
\definecolor{observationblue}{RGB}{205,226,239}
\definecolor{observationedge}{RGB}{56,103,134}
\definecolor{curvaturearc}{RGB}{162,77,39}
\def\Runit{.48/cos(180/64)}
\def\unitpolygon{(2.8125:{\Runit})--(8.4375:{\Runit})--(14.0625:{\Runit})--(19.6875:{\Runit})--(25.3125:{\Runit})--(30.9375:{\Runit})--(36.5625:{\Runit})--(42.1875:{\Runit})
--(47.8125:{\Runit})--(53.4375:{\Runit})--(59.0625:{\Runit})--(64.6875:{\Runit})--(70.3125:{\Runit})--(75.9375:{\Runit})--(81.5625:{\Runit})--(87.1875:{\Runit})
--(92.8125:{\Runit})--(98.4375:{\Runit})--(104.0625:{\Runit})--(109.6875:{\Runit})--(115.3125:{\Runit})--(120.9375:{\Runit})--(126.5625:{\Runit})--(132.1875:{\Runit})
--(137.8125:{\Runit})--(143.4375:{\Runit})--(149.0625:{\Runit})--(154.6875:{\Runit})--(160.3125:{\Runit})--(165.9375:{\Runit})--(171.5625:{\Runit})--(177.1875:{\Runit})
--(182.8125:{\Runit})--(188.4375:{\Runit})--(194.0625:{\Runit})--(199.6875:{\Runit})--(205.3125:{\Runit})--(210.9375:{\Runit})--(216.5625:{\Runit})--(222.1875:{\Runit})
--(227.8125:{\Runit})--(233.4375:{\Runit})--(239.0625:{\Runit})--(244.6875:{\Runit})--(250.3125:{\Runit})--(255.9375:{\Runit})--(261.5625:{\Runit})--(267.1875:{\Runit})
--(272.8125:{\Runit})--(278.4375:{\Runit})--(284.0625:{\Runit})--(289.6875:{\Runit})--(295.3125:{\Runit})--(300.9375:{\Runit})--(306.5625:{\Runit})--(312.1875:{\Runit})
--(317.8125:{\Runit})--(323.4375:{\Runit})--(329.0625:{\Runit})--(334.6875:{\Runit})--(340.3125:{\Runit})--(345.9375:{\Runit})--(351.5625:{\Runit})--(357.1875:{\Runit})
}
\node[font=\footnotesize] at (0,1.44) {(a) Domain $\Omega_0$};
\begin{scope}[shift={(-.666667,0)}]
\def\fr{(1+2*cos(\t)/3)}
\path[fill=black!2,draw=black,line width=.75pt,
samples=300,domain=0:360,variable=\t]
plot ({\fr*cos(\t)},{\fr*sin(\t)}) -- cycle;
\draw[curvaturearc,line width=2pt,
samples=30,domain=169.26:190.74,variable=\t]
plot ({\fr*cos(\t)},{\fr*sin(\t)});
\fill (0,0) circle[radius=1pt];
\node[anchor=south west,font=\footnotesize] at (.03,.03) {$0$};
\end{scope}
\node[curvaturearc,anchor=east,font=\footnotesize] at (-1.16,.50)
{$\kappa<0$};
\draw[curvaturearc,line width=.4pt] (-1.14,.46)--(-1,0);
\draw[->,curvaturearc,line width=.8pt] (-.98,-.04)--(-.40,-.23);
\node[font=\footnotesize] at (0,-1.40)
{$r=1+\tfrac23\cos\theta$};

\begin{scope}[shift={(3.30,0)}]
\node[font=\footnotesize] at (0,1.44)
{(b) Lattice cell, enlarged};
\foreach \px/\py/\lab in {-.5/-.3/q_0,.5/-.3/q_1,0/.5660254/q_2} {
\begin{scope}[shift={(\px,\py)}]
\path[fill=observationblue,draw=observationedge,line width=.6pt]
\unitpolygon -- cycle;
\end{scope}
}
\draw[black!55,densely dashed,line width=.5pt]
(-.5,-.3)--(.5,-.3)--(0,.5660254)--cycle;
\foreach \px/\py/\lab in {-.5/-.3/q_0,.5/-.3/q_1} {
\begin{scope}[shift={(\px,\py)}]
\fill (0,0) circle[radius=.8pt];
\node[font=\scriptsize,anchor=north] at (0,-.04) {$\lab$};
\end{scope}
}
\fill (0,.5660254) circle[radius=.8pt];
\node[font=\scriptsize,anchor=south] at (0,.6060254) {$q_2$};
\fill[curvaturearc] (0,0) circle[radius=1pt];
\draw[->,curvaturearc,line width=.8pt] (0,0)--(.99,-.12);
\draw[black!45,line width=.4pt]
(-.5,-.81)--(-.5,-1.01) (.5,-.81)--(.5,-1.01);
\draw[<->,line width=.45pt] (-.5,-.96)--(.5,-.96)
node[midway,fill=white,inner sep=1pt,font=\footnotesize] {$\delta/10$};
\node[font=\footnotesize] at (0,-1.40) {$q_i\in\Lambda_\delta$};
\end{scope}

\begin{scope}[shift={(6.25,0)}]
\node[font=\footnotesize] at (0,1.44)
{(c) Regular $64$-gon};
\begin{scope}[scale=1.65]
\path[fill=observationblue,draw=observationedge,line width=.65pt]
\unitpolygon -- cycle;
\foreach \k in {0,...,63} {
\fill[observationedge]
({(\k+.5)*5.625}:{\Runit}) circle[radius=.55pt];
}
\fill (0,0) circle[radius=.8pt];
\draw[->,line width=.45pt] (0,0)--(.48,0)
node[midway,below,font=\scriptsize] {$r_\delta$};
\node[font=\footnotesize] at (0,.20) {$64$ sides};
\end{scope}
\node[font=\footnotesize] at (0,-1.12) {$r_\delta=6\delta/125$};
\node[font=\footnotesize] at (0,-1.40) {$\alpha_0=\pi/32$};
\end{scope}
\end{tikzpicture}
\caption{Construction of the observation geometry.
Panel (a) shows the exact boundary of $\Omega_0$, with its concave arc
highlighted and an arrow pointing into the interior. Panel (b) enlarges
one equilateral lattice cell: a straight segment starting in the central
gap meets a blue observation polygon. The dashed segments join the
polygon centres. Panel (c) further enlarges one polygon and marks all
$64$ vertices. Panels (b) and (c) are local magnifications of the family
in \eqref{ce:polygons}.}
\label{fig:polygonal-geometry}
\end{figure}
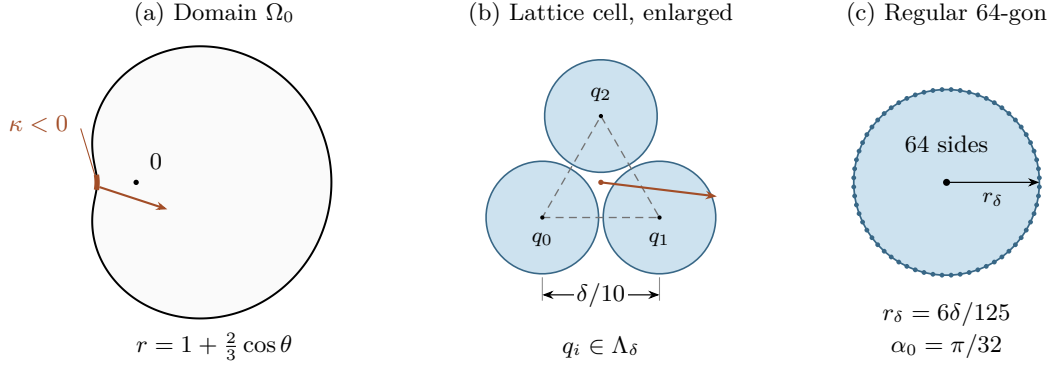

In Figure~\ref{fig:polygonal-geometry}, the concave arc prevents
boundary confinement and the triangular arrangement intercepts interior
straight rays. The gaps between polygons and their distance from the
boundary leave room for the local construction.

For $a>0$, write $\Omega_a=a\Omega_0$, $\omega_a=a\omega_0$, and let
$\mathcal S_a(t)$ be the free Dirichlet wave group on $\Omega_a$, with
$\pi_1(f,g)=f$. The number $\delta_*>0$ below is a uniform geometric
threshold for boundary escape; Step~\ref{step:geometry-1} proves its
existence. The next proposition states the four properties of the
construction: escape from the boundary collar, the geometry of the
polygons, interception of interior rays, and free observability after
rescaling.
\begin{proposition}[Observation geometry and free observability]\label{ce:geometry}
There exists $\delta_*=\delta_*(\Omega_0)\in(0,1/12)$ such that the
following assertions hold.
\begin{enumerate}
\item[(i)]\proofsteplabel{cond:geometry-escape}{(i)} For every $0<\delta<\delta_*$, every generalized unit-speed
ray in $\Omega_0$ enters $\{x\in\Omega_0:d(x)>4\delta\}$ before time
$34$.

\item[(ii)]\proofsteplabel{cond:geometry-polygons}{(ii)} For every $0<\delta<1/12$, the family
$\{P_q^0:q\in\mathcal J_\delta\}$ is finite and nonempty.
Its members are mutually separated, congruent regular $64$-gons
compactly contained in $\Omega_0$, with inradius $6\delta/125$ and
exterior angle $\pi/32$. Their common edge length is
$2r_\delta\tan(\pi/64)>0$. In particular,
$\omega_0\Subset\Omega_0$.

\item[(iii)]\proofsteplabel{cond:geometry-interception}{(iii)} Let $0<\delta<\delta_*$. For every $x\in\Omega_0$ with
$d(x)\ge3\delta$ and every unit vector $v$, the segment
$x+[0,\delta/2]v$ meets $\omega_0$.
Every generalized unit-speed ray enters $\omega_0$ before time
$34+\delta/2<35$. Consequently, $(\Omega_0,\omega_0)$ satisfies GCC.

\item[(iv)]\proofsteplabel{cond:geometry-free}{(iv)} Let $0<\delta<\delta_*$. For $a>0$, set
$\Omega_a\triangleq a\Omega_0$ and $\omega_a\triangleq a\omega_0$,
and let $\mathcal S_a$ be the free Dirichlet wave group on $\Omega_a$.
The pair $(\Omega_a,\omega_a)$ satisfies GCC before time $35a$.
There exists $C(\delta)>0$, independent of $a$, such that
\begin{equation}\label{ce:free-scaled} 
\|z_0\|_{L^2(\Omega_a)\times H^{-1}(\Omega_a)}^2
\le\frac{C(\delta)}a\int_0^{35a}
\|\mathbf1_{\omega_a}\pi_1\mathcal S_a(t)z_0\|_2^2\,\mathrm dt,\qquad
z_0\in L^2(\Omega_a)\times H^{-1}(\Omega_a).
\end{equation}
where $\pi_1(f,g)\triangleq f$ is the displacement projection.
\end{enumerate}
\end{proposition}
Following the order in Figure~\ref{fig:polygonal-geometry}, we prove
assertions~\ref{cond:geometry-escape}--\ref{cond:geometry-free}.
\begin{proof}[Proof of Proposition~\ref{ce:geometry}]
The proposition has four assertions. We prove them in the order in
which they appear: escape from the boundary collar, the geometry of
the observation polygons, interception of interior rays together with
GCC, and free observability after rescaling. The first three assertions
are purely geometric; the last one combines GCC with a scaling
argument.

\medskip
\noindent\textbf{Step 1. Uniform escape from a boundary collar:
proof of \ref{cond:geometry-escape}.}
\proofsteplabel{step:geometry-1}{1}

We use compactness. Rays confined to arbitrarily thin collars would
converge to a ray lying entirely on the boundary, which the concave
part of the boundary makes impossible. The proof has three parts:
compute the geometry of the boundary, establish compactness of the
compressed bicharacteristics, and derive the contradiction at the
concave contact.

\smallskip\noindent\textit{Geometry of the boundary.}
The counterclockwise parametrization
\[
q(\theta)=\left(1+\frac23\cos\theta\right)(\cos\theta,\sin\theta)
\]
satisfies
\[
|q_\theta|^2=\frac{13}{9}+\frac43\cos\theta\ge\frac19,
\qquad
L_b\triangleq\int_0^{2\pi}|q_\theta|\,\mathrm d\theta
\le\frac{2\pi\sqrt{13}}3<8,
\]
where the length bound follows from Cauchy--Schwarz.
Its signed curvature, with the inward normal convention, is
\[
\kappa(\theta)=
\frac{\frac{17}{9}+2\cos\theta}
{\left(\frac{13}{9}+\frac43\cos\theta\right)^{3/2}},
\qquad \kappa(\pi)=-3.
\]
The only zeros of $\kappa$ satisfy $\cos\theta=-17/18$, and there
\[
\kappa_\theta
=\frac{-2\sin\theta}
{\left(\frac{13}{9}+\frac43\cos\theta\right)^{3/2}}\ne0.
\]
Thus tangencies have order two where $\kappa\ne0$ and order three at
its zeros; there are no contacts of infinite order.
Figure~\ref{fig:collar-escape-mechanism} illustrates the obstruction
to a ray remaining on the boundary.

\begin{figure}[H]
\centering
\begin{minipage}[t]{.48\textwidth}
\centering
\begin{tikzpicture}[x=1.65cm,y=1.65cm,font=\footnotesize,
>={Stealth[length=1.5mm]},line join=round]
\def\fr{(1+2*cos(\t)/3)}
\path[fill=blue!5,draw=black!80,samples=240,domain=0:360,variable=\t]
plot ({\fr*cos(\t)},{\fr*sin(\t)}) -- cycle;
\draw[->,orange!80!black,line width=1pt,samples=50,
domain=115:179,variable=\t]
plot ({\fr*cos(\t)},{\fr*sin(\t)});
\fill[orange!80!black] (-1/3,0) circle[radius=1.3pt];
\node[anchor=east,orange!80!black] at (-.48,-.02) {$q(\pi)$};
\draw[orange!80!black,line width=.4pt] (-.46,-.02)--(-.35,-.02);
\node at (.65,.15) {$\Omega_0$};
\node[align=center] at (.62,-1.48)
{(a) A hypothetical boundary ray};
\node[anchor=east,align=right] at (-.42,.84) {boundary\\ray};
\draw[black!55,line width=.4pt] (-.40,.73)--(-.23,.48);
\end{tikzpicture}
\end{minipage}\hfill
\begin{minipage}[t]{.48\textwidth}
\centering
\begin{tikzpicture}[x=1.05cm,y=1.65cm,font=\footnotesize,
>={Stealth[length=1.5mm]},line join=round]
\path[fill=blue!5] (-1.55,.85)--(1.55,.85)
-- plot[samples=80,domain=1.55:-1.55,variable=\t]
(\t,{-.42*\t*\t}) -- cycle;
\draw[black!80,line width=.8pt,samples=80,domain=-1.55:1.55,variable=\t]
plot (\t,{-.42*\t*\t});
\draw[->,orange!80!black,line width=.9pt] (-1.62,0)--(1.65,0);
\fill (0,0) circle[radius=1pt];
\node[anchor=north] at (0,-.07) {$q(\pi)$};
\node[orange!80!black,anchor=south] at (1.1,.04) {tangent ray};
\node at (-.85,.58) {interior};
\node[anchor=west] at (1.45,-1.02) {$\partial\Omega_0$};
\draw[black!65,line width=.4pt] (1.40,-.96)--(1.20,-.61);
\node[align=center] at (0,-1.48)
{(b) Straight continuation at $q(\pi)$};
\end{tikzpicture}
\end{minipage}
\caption{The obstruction to a ray lying entirely on the boundary.
In (a), such a hypothetical ray would follow the boundary and reach
$q(\pi)$. In (b), the tangent line lies locally inside $\Omega_0$ on
both sides of the contact, so a generalized ray must enter the interior.
Panel (a) uses the exact boundary; panel (b) is a rotated local schematic.}
\label{fig:collar-escape-mechanism}
\end{figure}
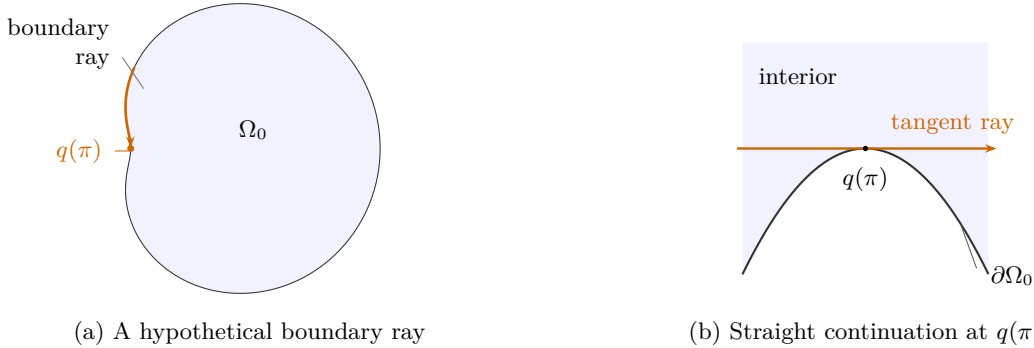

\noindent\textit{Compactness of the compressed
bicharacteristics.}
We apply compactness to the compressed bicharacteristics before
projecting to their positions. Let $\pi_b$ denote the compression map,
which identifies incoming and reflected covectors at the boundary.
For the symbol $h(t,x;\tau,\xi)=\tau^2-|\xi|^2$, the normalization
$\tau=1/2$, $|\xi|=1/2$ gives
\[
H_h t=2\tau=1,\qquad |H_h x|=2|\xi|=1.
\]
Thus the compressed lift of every unit-speed generalized ray on
$[0,T]$ takes values in the compact set
\[
\mathcal K_T\triangleq
\pi_b\bigl\{(t,x;1/2,\xi):t\in[0,T],\,
x\in\overline{\Omega_0},\ |\xi|=1/2\bigr\}.
\]
The compact smooth manifold with boundary $\overline{\Omega_0}$, with
its Euclidean metric, is a special case of Vasy's setting; at its
boundary, \cite[Definition~1.1(iii)]{VasyCorners} imposes the
Melrose--Sj\"ostrand rule. By
\cite[Propositions~5.4--5.5]{VasyCorners} (see also
\cite[Proposition~24.3.12 and the preceding
discussion]{HormanderIII}), every sequence of such lifts has a
uniformly convergent subsequence whose limit is again a compressed
generalized bicharacteristic. The continuous position projection
therefore gives uniform convergence to another unit-speed generalized
ray. Moreover, the absence of infinite-order contacts gives uniqueness
through each compressed state \cite[Corollary~24.3.10]{HormanderIII}.
If the initial states converge, every subsequential limit has that
initial state and hence is the same curve. Compactness then gives
uniform convergence of the full sequence, proving continuous
dependence as well.

\smallskip\noindent\textit{Derivation of the contradiction.}
We show that, for all sufficiently small $\delta>0$, every ray reaches
$d>4\delta$ at some time in $[0,16]$. Otherwise, there are
$\delta_n\downarrow0$ and unit-speed generalized rays $x_n$ such that
\[
0\le d(x_n(t))\le4\delta_n\qquad(0\le t\le16).
\]
Applying \cite[Proposition~5.5]{VasyCorners} with $T=16$ and then
projecting to position gives, after passing to a subsequence,
\[
x_n\longrightarrow x\quad\text{uniformly on }[0,16],
\qquad d(x(t))=0\quad(0\le t\le16),
\]
where $x$ is a unit-speed generalized ray. Such a ray cannot have a
transverse reflection, whose adjacent pieces lie in the interior. Its
velocity is therefore tangent to the boundary, with magnitude one and
continuous direction. It follows the boundary in one fixed
orientation. Since $16>2L_b$, it reaches $q(\pi)=(-1/3,0)$ at some
$t_*\in(0,16)$.

The tangent line at $q(\pi)$ is $(-1/3,y)$. For
$0<|y|<1/\sqrt3$, its polar radius $r=\sqrt{1/9+y^2}$ satisfies
$1/3<r<2/3$, and
\begin{equation}\label{ce:concave-tangent-interior}
r^2-r+\frac29
=\left(r-\frac13\right)\left(r-\frac23\right)<0
\quad\Longrightarrow\quad
r<1-\frac2{9r}.
\end{equation}
Since $\cos\theta=-1/(3r)$ on this line, the last inequality in
\eqref{ce:concave-tangent-interior} is exactly the defining inequality
\eqref{ce:flower}, so $(-1/3,y)\in\Omega_0$. At this strictly concave
contact, a generalized ray continues along the straight tangent line
locally \cite[Definition~24.3.7(i)]{HormanderIII}. Hence
$x(t)\in\Omega_0$ for $0<|t-t_*|$ sufficiently small, contradicting
$d(x(t))=0$.

This contradiction gives a uniform $\delta_*>0$ for escape at some
time at most $16$. Decreasing $\delta_*$ if necessary so that
$\delta_*<1/12$ proves assertion~\ref{cond:geometry-escape}, since
$16<34$.

\medskip
\noindent\textbf{Step 2. Separation and geometry of the polygons:
proof of \ref{cond:geometry-polygons}.}
Fix $0<\delta<1/12$. Since $\Omega_0$ is bounded, it contains only
finitely many lattice points, so $\mathcal J_\delta$ is finite.
Moreover, $d(0)=1/3>2\delta$, hence $0\in\mathcal J_\delta$.

By \eqref{ce:lattice-polygon} and \eqref{ce:polygons}, the polygons
are congruent regular $64$-gons with inradius $r_\delta=6\delta/125$
and exterior angle $\pi/32$. Direct calculation gives their
circumradius and common edge length as
\[
R_\delta=\frac{6\delta}{125\cos(\pi/64)}<\frac{\delta}{20},
\qquad 2r_\delta\tan(\pi/64)>0,
\]
respectively. The lattice definition gives $|q-q'|\ge\delta/10$ for
distinct centres, so
\[
\operatorname{dist}(P_q^0,P_{q'}^0)
\ge\frac{\delta}{10}-2R_\delta>0.
\]
Since $d(q)\ge2\delta>R_\delta$ and
$P_q^0\subset\overline B(q,R_\delta)$, every polygon is compactly
contained in $\Omega_0$. The finite union
$\overline{\omega_0}=\bigcup_{q\in\mathcal J_\delta}P_q^0$ is
therefore also compact and contained in $\Omega_0$. This proves
$\omega_0\Subset\Omega_0$ and all the asserted geometric properties.

\medskip
\noindent\textbf{Step 3. Interception of interior rays and GCC:
proof of \ref{cond:geometry-interception}.}
\proofsteplabel{step:geometry-3}{3}
Fix $0<\delta<\delta_*$. The proof has two parts: a lattice
interception lemma, and the derivation of GCC from the escape estimate
and the lattice interception.

\smallskip\noindent\textit{Lattice interception.}
Every segment avoiding all disks $B(q,r_\delta)$,
$q\in\Lambda_\delta$, has length at most $\delta/5$. Indeed, on a
lattice-edge line its length is at most
$\delta/10-2r_\delta<\delta/5$. Otherwise, if it crosses three or more
lattice triangles, a middle triangle contains a complete chord $ST$
with $S\in AB$ and $T\in AC$. Avoidance of the vertex disks gives
$|AS|,|AT|\le\delta/10-r_\delta$. Since $\angle BAC=\pi/3$,
\[
\left|\frac{S+T}{2}-A\right|^2
=\frac{|AS|^2+|AT|^2+|AS|\,|AT|}{4}
\le\frac34\left(\frac{\delta}{10}-r_\delta\right)^2
<r_\delta^2,
\]
contradicting disk avoidance. Thus the segment lies in at most two
triangles, each of diameter $\delta/10$, proving the assertion.

\smallskip\noindent\textit{Interior segments.}
If $d(x)\ge3\delta$ and $|v|=1$, then
$x+[0,\delta/2]v\subset\Omega_0$ and every point $p$ of this segment
satisfies $d(p)\ge5\delta/2$. The segment has length
$\delta/2>\delta/5$, so it meets some $B(q,r_\delta)$ at a point $p$.
Since $|p-q|<r_\delta<d(p)$, $q\in\Omega_0$, and\vspace{-2mm}
\begin{equation}\label{ce:geometry-interior-hit}
d(q)\ge d(p)-|p-q|>\frac{5\delta}{2}-r_\delta>2\delta.\vspace{-2mm}
\end{equation}
Hence $q\in\mathcal J_\delta$ and
$B(q,r_\delta)\subset\operatorname{int}P_q^0\subset\omega_0$, which
proves the first assertion.

\smallskip\noindent\textit{GCC.}
By Step~\ref{step:geometry-1}, every generalized unit-speed ray has a
time $t_0<34$ with $d(x(t_0))>4\delta$. During the next $\delta/2$
units of time it cannot reach the boundary, so it follows a straight
segment. The interior interception established by
\eqref{ce:geometry-interior-hit} gives\vspace{-2mm}
\[
t_{\rm hit}\le t_0+\frac{\delta}{2}<34+\frac{\delta}{2}<35.\vspace{-2mm}
\]
If the intersection occurs at time zero, openness gives a positive
observation time as well. Thus every generalized ray meets $\omega_0$
in $(0,34+\delta/2)$, proving GCC.
Figure~\ref{fig:gcc-interception} illustrates this argument.

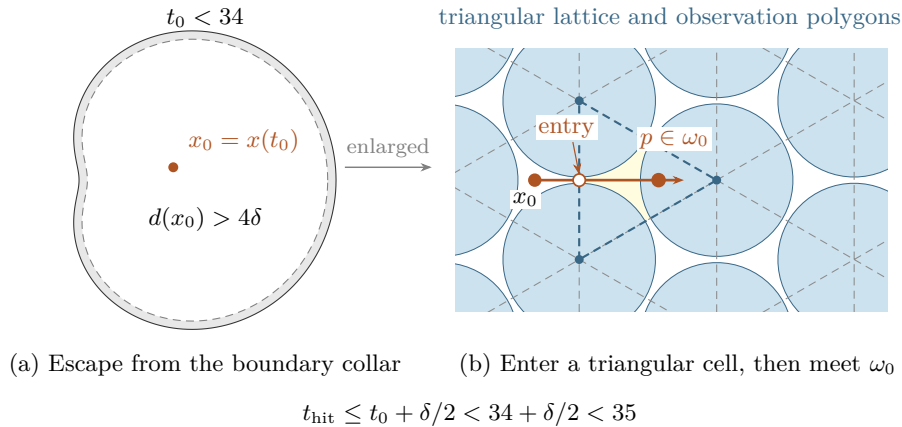
\begin{figure}[H]
\centering
\begin{tikzpicture}[x=1cm,y=1cm,font=\footnotesize,
>={Stealth[length=1.5mm]},line join=round]
\definecolor{gccblue}{RGB}{205,226,239}
\definecolor{gccedge}{RGB}{56,103,134}
\definecolor{gccray}{RGB}{174,83,34}
\begin{scope}[shift={(-.5,0)},scale=1.7]
\def\fr{(1+2*cos(\t)/3)}
\def\speed{sqrt(13/9+4*cos(\t)/3)}
\def\nx{(-(cos(\t)+2*cos(2*\t)/3)/\speed)}
\def\ny{(-sin(\t)*(1+4*cos(\t)/3)/\speed)}
\path[fill=black!9,even odd rule,samples=240,domain=0:360,variable=\t]
plot ({\fr*cos(\t)},{\fr*sin(\t)}) -- cycle
plot ({\fr*cos(\t)+.065*\nx},{\fr*sin(\t)+.065*\ny}) -- cycle;
\draw[black!80,samples=240,domain=0:360,variable=\t]
plot ({\fr*cos(\t)},{\fr*sin(\t)});
\draw[black!50,densely dashed,samples=240,domain=0:360,variable=\t]
plot ({\fr*cos(\t)+.065*\nx},{\fr*sin(\t)+.065*\ny});
\fill[gccray] (.4,.10) circle[radius=1.1pt];
\node[anchor=south west,gccray] at (.44,.13) {$x_0=x(t_0)$};
\node at (.65,-.30) {$d(x_0)>4\delta$};
\end{scope}
\node at (.6,2.15) {$t_0<34$};
\node[align=center] at (.6,-2.45) {(a) Escape from the boundary collar};
\draw[->,black!50,line width=.55pt] (2.45,.17)--(3.60,.17)
node[midway,above,font=\scriptsize] {enlarged};

\begin{scope}[shift={(5.55,0)},scale=2.1]
\begin{scope}
\clip (-.78,-.83) rectangle (1.94,.83);
\path[fill=yellow!15] (0,-.5)--(0,.5)--({sqrt(3)/2},0)--cycle;
\foreach \j in {-2,-1,0,1,2,3} {
\foreach \i in {-2,-1,0,1,2} {
\path[fill=gccblue,draw=gccedge,line width=.35pt,
samples=65,domain=0:360,variable=\a]
plot ({\j*sqrt(3)/2+(.48/cos(180/64))*cos(\a+180/64)},
{\i+\j/2-.5+(.48/cos(180/64))*sin(\a+180/64)}) -- cycle;
}
}
\foreach \j in {-2,-1,0,1,2,3} {
\foreach \i in {-2,-1,0,1,2} {
\draw[black!45,densely dashed,line width=.4pt]
({\j*sqrt(3)/2},{\i+\j/2-.5})--
({\j*sqrt(3)/2},{\i+\j/2+.5})
({\j*sqrt(3)/2},{\i+\j/2-.5})--
({(\j+1)*sqrt(3)/2},{\i+\j/2})
({\j*sqrt(3)/2},{\i+\j/2-.5})--
({(\j+1)*sqrt(3)/2},{\i+\j/2-1});
}
}
\draw[gccedge,densely dashed,line width=.8pt]
(0,-.5)--(0,.5)--({sqrt(3)/2},0)--cycle;
\foreach \p in {(0,-.5),(0,.5),({sqrt(3)/2},0)}
\fill[gccedge] \p circle[radius=.8pt];
\end{scope}
\draw[->,gccray,line width=1pt] (-.28,0)--(.66,0);
\fill[gccray] (-.28,0) circle[radius=1.2pt];
\node[anchor=north,fill=white,inner sep=1pt] at (-.34,-.06) {$x_0$};
\draw[gccray,fill=white,line width=.65pt] (0,0) circle[radius=1.1pt];
\node[gccray,anchor=south,fill=white,inner sep=1pt] at (-.08,.25)
{entry};
\draw[->,gccray,line width=.45pt] (-.04,.24)--(0,.035);
\fill[gccray] (.50,0) circle[radius=1.3pt];
\node[gccray,anchor=south,fill=white,inner sep=1pt] at (.60,.18)
{$p\in\omega_0$};
\end{scope}
\node[gccedge] at (6.75,2.15) {triangular lattice and observation polygons};
\node[align=center] at (6.85,-2.45)
{(b) Enter a triangular cell, then meet $\omega_0$};
\node at (4.1,-3.10)
{$t_{\rm hit}\le t_0+\delta/2<34+\delta/2<35$};
\end{tikzpicture}
\caption{Why the observation set satisfies GCC.
(a) Every generalized ray escapes the boundary collar before time $34$.
(b) Its next interior segment, of length $\delta/2$, must hit an
observation polygon. The enlarged local view shows the ray entering
the highlighted triangular cell before meeting the blue polygon at
$p$. Dashed lines form the triangular lattice, of spacing $\delta/10$.
The boundary in (a) is exact; the collar width is schematic.}
\label{fig:gcc-interception}
\end{figure}

\noindent\textbf{Step 4. Free observability and spatial dilation:
proof of \ref{cond:geometry-free}.}
Step~\ref{step:geometry-3} gives GCC for $(\Omega_0,\omega_0)$ before
time $34+\delta/2<35$. The boundary is analytic and has no contacts of
infinite order, as verified in Step~\ref{step:geometry-1}. Hence
\cite[Theorem~1.8$'$]{LRLT} gives
\begin{equation}\label{ce:free-reference}
\|u(0)\|_{L^2(\Omega_0)}^2+\|u_t(0)\|_{H^{-1}(\Omega_0)}^2
\le C(\delta)\int_0^{35}
\|\mathbf1_{\omega_0}u(t)\|_2^2\,\mathrm dt
\end{equation}
for every free Dirichlet wave $u$ on $\Omega_0$ with initial data in
$L^2(\Omega_0)\times H^{-1}(\Omega_0)$.

For $a>0$, the dilation $x(t)\mapsto a^{-1}x(a\tau)$ preserves
unit-speed generalized rays, so $(\Omega_a,\omega_a)$ satisfies GCC
before time $35a$. Let $z(t)=\pi_1\mathcal S_a(t)z_0$ and set
$
\widetilde z(\tau,\xi)=z(a\tau,a\xi).
$
Then $\widetilde z$ is a free Dirichlet wave on $\Omega_0$. In
dimension two, changing variables and using the dual definition of the
$H^{-1}$ norm give
\[
\begin{aligned}
\|\widetilde z(0)\|_{L^2(\Omega_0)}
&=a^{-1}\|z(0)\|_{L^2(\Omega_a)},\\
\|\partial_\tau\widetilde z(0)\|_{H^{-1}(\Omega_0)}
&=a^{-1}\|z_t(0)\|_{H^{-1}(\Omega_a)},\\
\int_0^{35}\|\mathbf1_{\omega_0}\widetilde z(\tau)\|_2^2\,\mathrm d\tau
&=a^{-3}\int_0^{35a}
\|\mathbf1_{\omega_a}z(t)\|_2^2\,\mathrm dt.
\end{aligned}
\]
Applying \eqref{ce:free-reference} to $\widetilde z$ and multiplying
by $a^2$ proves \eqref{ce:free-scaled}, with $C(\delta)$ independent
of $a$.
\end{proof}

\subsection{Finite-chain repair preserving zeros}\label{ce:repair-section}

We now establish the repair lemma. Its role in the global construction
is to merge two solutions at a vertex without changing the zero set of
the sum, even at points where the sum vanishes. The key idea is to
replace division by the oscillatory seed $q$ with division by a finite
sum of squares, which remains positive on the region where the repair
is needed.

For spacetime covectors $\zeta,\xi$ and a function $q$, set
\begin{equation}\label{ce:lorentz}
\mathcal B(\zeta,\xi)\triangleq\zeta_t\xi_t-\zeta_x\cdot\xi_x,\quad
\zeta^\sharp\triangleq(\zeta_t,-\zeta_x),
\qquad D_qh\triangleq\mathcal B(\mathrm{d}q,\mathrm{d}h).
\end{equation}
Given $\Box q=N+Gq$, we seek a real exponent $\Phi_{\rm rep}$ for
which $p=e^{\Phi_{\rm rep}}q$ satisfies $\Box p+Vp=0$. The positive
multiplier preserves the zero set. The sequence
$Q,D_qQ,\ldots,D_q^rQ$ consists of repeated derivatives along the
Lorentz gradient of $q$; we call it a finite chain.
\begin{lemma}[Repair preserving the zero set]\label{ce:repair}
Let $n\in\mathbb N$, let $O\subset\R^{1+n}$ be open, and let
$q,N,G\in C^\infty(O;\R)$ satisfy $\Box q=N+Gq$.
For $r\in\mathbb N_0$, set
\begin{equation}\label{ce:repair-chain}
Q\triangleq\mathcal B(\mathrm{d}q,\mathrm{d}q),\quad
\mathcal C_{-1}\triangleq q,\quad
\mathcal C_l\triangleq D_q^lQ\ (0\le l\le r),\quad
F\triangleq\sum_{l=0}^r \mathcal C_l^2.
\end{equation}
Assume that $F>0$ on $\operatorname{supp}_O N$, where
$\operatorname{supp}_O N$ is the closure of $\{N\ne0\}$ in $O$.
For $0\le l\le r$, define
\begin{equation}\label{ce:repair-b}
\beta_l\triangleq
\begin{cases}
-\dfrac{N\mathcal C_l}{2F},&F>0,\\[3pt]
0,&F=0,
\end{cases}
\end{equation}
and set
\begin{equation}\label{ce:repair-exponent}
\Phi_{\rm rep}\triangleq\sum_{l=0}^r\sum_{a=0}^l
(-1)^a\mathcal C_{l-1-a}D_q^a\beta_l,
\qquad
{H_{\rm rep}}\triangleq2\sum_{l=0}^r(-1)^lD_q^{l+1}\beta_l.
\end{equation}
Then the functions $\Phi_{\rm rep}$ and
\begin{equation}\label{ce:repair-output}
V\triangleq-G-{H_{\rm rep}}-\Box\Phi_{\rm rep}-\mathcal B(\mathrm{d}\Phi_{\rm rep},\mathrm{d}\Phi_{\rm rep})
\end{equation}
are real and smooth on $O$, and $p=e^{\Phi_{\rm rep}}q$ satisfies
$\Box p+Vp=0$ on $O$.
\end{lemma}

\begin{proof}
The coefficient definition gives
$2\sum_{l=0}^r\beta_l\mathcal C_l=-N$. Leibniz expansion then yields
$2D_q\Phi_{\rm rep}=-N+qH_{\rm rep}$. Substituting into
$\Box(e^{\Phi_{\rm rep}}q)$ gives the stated $V$. No step divides by
$q$, so the identity holds at zeros as well.
Appendix~\ref{app:repair-details} gives the coefficient sum and
product expansions.
\end{proof}

The next corollary collects the properties of the construction that
will be used in the vertex merging and in the completion to a static
pair. In particular, it records the identity
$2D_q\Phi_{\rm rep}=-N+qH_{\rm rep}$, the positivity of the multiplier,
the locality of the repair, and the preservation of time periodicity.

\begin{corollary}[Properties of the explicit repair]\label{ce:repair-properties}
Under the hypotheses of Lemma~\ref{ce:repair}, the construction
\eqref{ce:repair-chain}--\eqref{ce:repair-output} has the following
properties.
\begin{enumerate}
\item[(i)]\proofsteplabel{cond:repair-smooth}{(i)} The functions
$\beta_0,\ldots,\beta_r,\Phi_{\rm rep},{H_{\rm rep}}$ are real and
smooth on $O$, and
\begin{equation}\label{ce:repair-identity}
2D_q\Phi_{\rm rep}=-N+q{H_{\rm rep}}\qquad\hbox{throughout }O.
\end{equation}
\item[(ii)]\proofsteplabel{cond:repair-zeros}{(ii)} The multiplier
$e^{\Phi_{\rm rep}}$ is strictly positive, and
$p^{-1}(0)=q^{-1}(0)$.
\item[(iii)]\proofsteplabel{cond:repair-locality}{(iii)} On every open
subset of $O$ where $N\equiv0$, all the functions
$\beta_0,\ldots,\beta_r,\Phi_{\rm rep},{H_{\rm rep}}$ vanish and
$(p,V)=(q,-G)$.
\item[(iv)]\proofsteplabel{cond:repair-periodic}{(iv)} If $T>0$, the
translation $\tau_T(t,x)=(t+T,x)$ satisfies $\tau_T(O)=O$, and
$q,N,G$ are $T$-periodic, then all the constructed functions are
$T$-periodic.
\end{enumerate}
\end{corollary}

\begin{proof}[Proof of Corollary~\ref{ce:repair-properties}]
The smoothness in \ref{cond:repair-smooth} was established in the
proof of Lemma~\ref{ce:repair}, and its identity is
\eqref{ce:repair-multiplier-identity}. Since $\Phi_{\rm rep}$ is real,
$e^{\Phi_{\rm rep}}>0$ at every point of $O$, which proves
\ref{cond:repair-zeros}. If $U\subset O$ is open and $N|_U=0$, the
definition \eqref{ce:repair-b} gives $\beta_l|_U=0$, whether or not
$F$ vanishes. Every derivative of $\beta_l$ is therefore zero on $U$,
as are $\Phi_{\rm rep},{H_{\rm rep}}$ and every derivative of
$\Phi_{\rm rep}$. The definition $p=e^{\Phi_{\rm rep}} q$ and
\eqref{ce:repair-output} then give $(p,V)=(q,-G)$ on $U$, proving
\ref{cond:repair-locality}.

For \ref{cond:repair-periodic}, the periodicity of $q$ implies
periodicity of the coefficients of $D_q$ and the commutation rule
\begin{equation}\label{ce:repair-translation-rule}
D_q(f\circ\tau_T)=(D_qf)\circ\tau_T
\qquad(f\in C^\infty(O)).
\end{equation}
It follows inductively that $Q,\mathcal C_l$ and $F$ are $T$-periodic.
The two branches in \eqref{ce:repair-b} are preserved by $\tau_T$, so
each $\beta_l$ is $T$-periodic as well. The commutation rule
\eqref{ce:repair-translation-rule} applies to its iterates; hence
$\Phi_{\rm rep},{H_{\rm rep}}$ are periodic. Finally, coordinate
derivatives commute with time translation, and the definitions of
$p=e^{\Phi_{\rm rep}} q$ and $V$ in \eqref{ce:repair-output} show that
$p,V$ have the same period.
\end{proof}

For the iterates $\mathcal C_l=D_q^lQ$, $l\in\{1,2,3\}$, we record the
first three identities using the actual derivatives of the vector
field $D_q$. Here $\partial^kq$ denotes the $k$th total differential
in the spacetime variables, viewed as a symmetric $k$-linear form.
The operator $D_q$ acts on the coordinate components of covectors, so
that $D_q\mathrm{d}q=\partial^2q((\mathrm{d}q)^\sharp,\cdot)$.
The chain derivative identities \eqref{ce:chain-expansion} are
\begin{equation}\label{ce:chain-expansion}
\!\!\!\begin{cases}
\ds\mathcal C_1 =2\partial^2q((\mathrm{d}q)^\sharp,(\mathrm{d}q)^\sharp),\\
\ns\ds	\mathcal C_2 =2\partial^3q((\mathrm{d}q)^\sharp,(\mathrm{d}q)^\sharp,(\mathrm{d}q)^\sharp)
+4\mathcal B(D_q\mathrm{d}q,D_q\mathrm{d}q),\\
\ns\ds	\mathcal C_3 =2\partial^4q((\mathrm{d}q)^\sharp,(\mathrm{d}q)^\sharp,
(\mathrm{d}q)^\sharp,(\mathrm{d}q)^\sharp) \!+\!14\partial^3q((D_q\mathrm{d}q)^\sharp,(\mathrm{d}q)^\sharp,(\mathrm{d}q)^\sharp)
\!+\!8\partial^2q((D_q\mathrm{d}q)^\sharp,(D_q\mathrm{d}q)^\sharp).
\end{cases}
\end{equation}
Indeed, $D_q((\mathrm{d}q)^\sharp)=(D_q\mathrm{d}q)^\sharp$, and
\begin{equation}\label{ce:chain-differentiation}
\begin{cases}
\ds D_q\bigl[2\partial^3q((\mathrm{d}q)^\sharp,(\mathrm{d}q)^\sharp,(\mathrm{d}q)^\sharp)\bigr] =2\partial^4q((\mathrm{d}q)^\sharp,(\mathrm{d}q)^\sharp,
(\mathrm{d}q)^\sharp,(\mathrm{d}q)^\sharp)
+6\partial^3q((D_q\mathrm{d}q)^\sharp,(\mathrm{d}q)^\sharp,(\mathrm{d}q)^\sharp),\\
\ns\ds D_q\bigl[4\mathcal B(D_q\mathrm{d}q,D_q\mathrm{d}q)\bigr] =8\partial^3q((D_q\mathrm{d}q)^\sharp,(\mathrm{d}q)^\sharp,(\mathrm{d}q)^\sharp)
+8\partial^2q((D_q\mathrm{d}q)^\sharp,(D_q\mathrm{d}q)^\sharp).
\end{cases}
\end{equation}
The identities \eqref{ce:chain-differentiation} establish
\eqref{ce:chain-expansion}. In particular, $D_q$ cannot be replaced by
differentiation in a fixed direction when estimating the higher chain
terms.

\subsection{Merging two solutions}\label{ce:incoming-section}

We now apply the repair to the sum of two solutions near intersecting
straight interfaces. Each given pair vanishes on one side of its
interface; only the region where both pairs are present needs
correction. The main result of this subsection is
Proposition~\ref{ce:merging-criterion}, which states sufficient
conditions under which two edge pairs can be merged into a single
smooth pair.

Let $U\subset\R^2$ be an open neighborhood of $v$, and choose linearly
independent unit normals $n_1,n_2$. Set
\[
y_i(x)\triangleq n_i\cdot(x-v),\quad i\in\{1,2\},\qquad
D\triangleq U\cap\{y_1>0,y_2>0\}.
\]
Let $(f_i,W_i)\in C^\infty(\R\times U;\R)^2$ satisfy
\[
\begin{cases}
\ds\Box f_i+W_if_i=0\quad\hbox{on }\R\times U,\\
\ns\ds (f_i,W_i)=(0,0)\quad\hbox{on }\R\times(U\cap\{y_i\le0\}).
\end{cases}
\]
In the correction formulas, use
\begin{equation}\label{ce:merging-data}
q\triangleq f_1+f_2,\qquad
R\triangleq W_2f_1+W_1f_2
=\Box q+(W_1+W_2)q.
\end{equation}
Thus $R$ is precisely the residual of the sum.

\begin{proposition}[Merging two solutions]\label{ce:merging-criterion}
Let $r\in\mathbb N_0$. Suppose that there are envelopes
$E_i\in C^\infty(U;[0,\infty))$, $i\in\{1,2\}$, with
$E_i=0$ where $y_i\le0$ and $E_i>0$ where $y_i>0$.
All supports below are taken relative to $\R\times D$.
The estimates in \ref{cond:merging-potential}--\ref{cond:merging-chain}
and in the conclusion are understood for every compact
$I\subset\R$ and compact $K\subset U$, with positive constants
allowed to depend on $I$ and $K$. Assume the following conditions
for $i\in\{1,2\}$.

\begin{enumerate}
\item[(i)]\proofsteplabel{cond:merging-amplitude}{(i)}
\emph{Nonvanishing on the residual support.}
There are constants $c,C>0$ such that
\[
|f_i|\le CE_i\quad\hbox{on }\R\times U,
\qquad
|f_i|\ge cE_i\quad\hbox{on }\operatorname{supp}W_i.
\]

\item[(ii)]\proofsteplabel{cond:merging-potential}{(ii)}
\emph{Decay of the potentials relative to the envelopes.}
For every $\alpha\in\mathbb N_0^3$ and every $P\in\mathbb N$,
\begin{equation}\label{ce:criterion-potential-decay}
|\partial^\alpha W_i|\le C_{I,K,\alpha,P}E_i^P
\quad\hbox{on }I\times(K\cap D).
\end{equation}

\item[(iii)]\proofsteplabel{cond:merging-chain}{(iii)}
\emph{Quantitative lower bound for the finite chain.}
Choose $0<a<b<c$ and a cutoff
$\vartheta\in C^\infty(\R;[0,1])$ equal to one on $[-a,a]$
and zero outside $(-b,b)$. On $\R\times D$, set
\[
\delta\triangleq\frac{E_1E_2}{E_1+E_2},
\qquad
\chi_*\triangleq\vartheta(q/\delta),
\qquad
F_0\triangleq\sum_{l=0}^r \mathcal C_l^2,
\]
where $\mathcal C_l$ is defined by \eqref{ce:repair-chain}.
For some fixed $\kappa>0$,
\begin{equation}\label{ce:criterion-nondegeneracy}
F_0\ge c_{I,K}E_i^\kappa
\quad\hbox{on }(I\times K)\cap\operatorname{supp}(\chi_*W_i).
\end{equation}
\end{enumerate}

Then there are real smooth functions $p,V,\Phi_{\rm rep}$ on
$\R\times U$ satisfying
\[
\Box p+Vp=0,
\qquad
p=e^{\Phi_{\rm rep}}(f_1+f_2),
\qquad
\Phi_{\rm rep}=0\quad\hbox{on }\R\times(U\setminus D),
\]
and
\[
\begin{aligned}
&(p,V)=(f_i,W_i)
\quad\hbox{where }y_{3-i}\le0<y_i,\qquad i\in\{1,2\},\\
&(p,V)=(0,0)\quad\hbox{where }y_1\le0,\ y_2\le0.
\end{aligned}
\]
Moreover, for every $\alpha\in\mathbb N_0^3$ and every
$N\in\mathbb N_0$,
\begin{equation}\label{ce:merging-correction-flatness}
|\partial^\alpha\Phi_{\rm rep}|
+|\partial^\alpha(V-W_1-W_2)|
\le C_{I,K,\alpha,N}\min\{y_1,y_2\}^N \quad\mbox{ on }I\times(K\cap D).
\end{equation}
If both pairs $(f_i,W_i)$ are $T$-periodic in
time, then so are $p,V,\Phi_{\rm rep}$.
\end{proposition}

The three hypotheses have the following roles.
Condition~\ref{cond:merging-amplitude} ensures that neither wave
vanishes on the residual support of the other.
Condition~\ref{cond:merging-potential} states that each potential
decays faster than any power of the corresponding envelope, which
allows the repair coefficients to absorb the divisions.
Condition~\ref{cond:merging-chain} is the quantitative lower bound on
the finite chain near the zeros of the sum; it is the nondegeneracy
condition that makes the repair possible.

\begin{center}\small
\begin{tabular}{@{}p{.14\textwidth}p{.40\textwidth}p{.40\textwidth}@{}}
\toprule Condition&Role&Verified in\\\midrule
\ref{cond:merging-amplitude}&nonvanishing on the residual support&\eqref{ce:incoming-amplitude-bounds}\\
\ref{cond:merging-potential}&arbitrary decay relative to the envelope&\eqref{ce:incoming-paid}\\
\ref{cond:merging-chain}&quantitative chain bound near zeros&\eqref{ce:incoming-unweighted}\\\bottomrule
\end{tabular}
\end{center}
\begin{proof}[Proof of Proposition~\ref{ce:merging-criterion}]
The proof has four steps: control the divisions on each source, derive
the flatness of the correction coefficients, extend the repair across
both interfaces, and form the merged pair. The first two steps use the
decay and the nondegeneracy of the chain. The last two use the repair
lemma and the positivity of the multiplier.

\smallskip
\noindent\textbf{Step 1. Control the divisions on each source.}
Fix compact sets $I\subset\R$, $K\subset U$ and an index
$i\in\{1,2\}$, and put $j=3-i$. All constants in this proof may
depend on these compact sets. The functions $q,\mathcal C_l,E_i$ and
all their derivatives are bounded there. On
$\operatorname{supp}(\chi_*W_i)$,
\[
|q|\le b\delta\le bE_i,\qquad
|f_j|\ge(c-b)E_i,\qquad
E_j\ge\frac{c-b}{C}E_i.
\]
Thus
\[
\delta\ge\frac{c-b}{C+c-b}E_i\quad\hbox{on }\operatorname{supp}(\chi_*W_i).
\]
For any smooth nonvanishing function $h$ and $|\alpha|\ge1$,
\begin{equation}\label{ce:merging-reciprocal}
\partial^\alpha(h^{-1})
=\sum_{k=1}^{|\alpha|}(-1)^kh^{-k-1}
\sum_{\substack{\nu_1+\cdots+\nu_k=\alpha\\|\nu_s|\ge1}}
\frac{\alpha!}{\nu_1!\cdots\nu_k!}
\prod_{s=1}^k\partial^{\nu_s}h.
\end{equation}
Since $E_1+E_2\ge E_i$, this gives
\[
\begin{aligned}
\partial^\alpha\delta
&=\sum_{\alpha_1+\alpha_2+\alpha_3=\alpha}
\frac{\alpha!}{\alpha_1!\alpha_2!\alpha_3!}
(\partial^{\alpha_1}E_1)(\partial^{\alpha_2}E_2)
\partial^{\alpha_3}(E_1+E_2)^{-1},
\qquad
|\partial^\alpha\delta| \le C_\alpha E_i^{-M_\alpha}.
\end{aligned}
\]
On $\operatorname{supp}W_i\cap\{|q|\le b\delta\}$, the same comparison
$E_j\ge(c-b)E_i/C$ implies
\[
|\partial^\alpha\delta^{-1}|
+\left|\partial^\alpha(q/\delta)\right|
\le C_\alpha E_i^{-M_\alpha}.
\]
For $|\alpha|\ge1$, the chain rule reads
\[
\partial^\alpha\chi_*
=\sum_{k=1}^{|\alpha|}\frac{\vartheta^{(k)}(q/\delta)}{k!}
\sum_{\substack{\nu_1+\cdots+\nu_k=\alpha\\|\nu_s|\ge1}}
\frac{\alpha!}{\nu_1!\cdots\nu_k!}
\prod_{s=1}^k\partial^{\nu_s}(q/\delta).
\]
Since $\vartheta^{(k)}=0$ on $\{|s|>b\}$ for $k\ge1$,
\begin{equation}\label{ce:merging-cutoff-bound}
|\partial^\alpha\chi_*|\le C_\alpha E_i^{-M_\alpha}
\quad\hbox{on }\operatorname{supp}W_i.
\end{equation}
The order-zero bound is $|\chi_*|\le1$.
By \eqref{ce:criterion-nondegeneracy} and
\eqref{ce:merging-reciprocal},
\begin{equation}\label{ce:merging-chain-reciprocal}
|\partial^\alpha F_0^{-1}|
\le C_\alpha\sum_{k=0}^{|\alpha|}F_0^{-k-1}
\le C_\alpha E_i^{-\kappa(|\alpha|+1)}
\quad\hbox{on }\operatorname{supp}(\chi_*W_i).
\end{equation}

On $\operatorname{supp}((1-\chi_*)W_i)$, there is also a constant
$c_2>0$ such that $|q|\ge c_2E_i$. Indeed, if
$E_j<cE_i/(2C)$, then
\[
|q|\ge |f_i|-|f_j|\ge cE_i/2.
\]
If $E_j\ge cE_i/(2C)$, then
\[
|q|\ge a\delta\ge\frac{ac}{2C+c}E_i
\quad\hbox{on }\operatorname{supp}((1-\chi_*)W_i).
\]
Thus, on this source,
\begin{equation}\label{ce:merging-q-reciprocal}
|q|\ge\min\left\{\frac c2,\frac{ac}{2C+c}\right\}E_i,
\qquad
|\partial^\alpha q^{-1}|\le C_\alpha E_i^{-(|\alpha|+1)}.
\end{equation}

\medskip
\noindent\textbf{Step 2. Derive flatness of the correction
coefficients.}
Define
\begin{equation}\label{ce:merging-source-coefficients}
B_i\triangleq\frac{(1-\chi_*)W_if_j}{q},\qquad
\beta_{l,i}\triangleq-\frac{\chi_*W_if_j\mathcal C_l}{2F_0},\qquad
B\triangleq B_1+B_2,\quad \beta_l\triangleq \beta_{l,1}+\beta_{l,2},
\end{equation}
with the first quotient prescribed as zero at $q=0$ and the second at
$F_0=0$. These are smooth on $\R\times D$: the first numerator
vanishes on a neighborhood of every zero of $q$, and every zero of
$F_0$ has a neighborhood disjoint from
$\operatorname{supp}(\chi_*W_i)$ by
\eqref{ce:criterion-nondegeneracy}.
Where the respective denominators are nonzero,
\begin{equation}\label{ce:merging-source-derivatives}
\begin{aligned}
\partial^\alpha B_i
&=\sum_{\gamma+\eta+\nu=\alpha}
\frac{\alpha!}{\gamma!\eta!\nu!}
(\partial^\gamma W_i)(\partial^\eta f_j)
\partial^\nu\!\left(\frac{1-\chi_*}{q}\right),\\
\partial^\alpha \beta_{l,i}
&=-\frac12\sum_{\gamma+\eta+\nu=\alpha}
\frac{\alpha!}{\gamma!\eta!\nu!}
(\partial^\gamma W_i)(\partial^\eta f_j)
\partial^\nu\!\left(\frac{\chi_*\mathcal C_l}{F_0}\right).
\end{aligned}
\end{equation}
Combining \eqref{ce:merging-source-derivatives} with
\eqref{ce:merging-cutoff-bound}, \eqref{ce:merging-chain-reciprocal},
and \eqref{ce:merging-q-reciprocal} yields
\begin{equation}\label{ce:merging-source-intermediate-bound}
|\partial^\alpha B_i|+|\partial^\alpha \beta_{l,i}|
\le C_\alpha E_i^{-M_\alpha}
\Bigl(\max_{|\gamma|\le|\alpha|}|\partial^\gamma W_i|\Bigr)
\sum_{|\eta|\le|\alpha|}|\partial^\eta f_j|.
\end{equation}
Bound \eqref{ce:merging-source-intermediate-bound} is used on the
relevant source support; off these supports the derivatives of the
source vanish.
Choose an integer $P\ge M_\alpha+1$ in
\eqref{ce:criterion-potential-decay}. Then
\begin{equation}\label{ce:merging-source-decay}
\begin{aligned}
|\partial^\alpha B_i|+|\partial^\alpha \beta_{l,i}|
&\le C_{\alpha,P}E_i^{P-M_\alpha}
\sum_{|\eta|\le|\alpha|}|\partial^\eta f_j|
\le C_{I,K,\alpha}E_i
\sum_{|\eta|\le|\alpha|}|\partial^\eta f_j|.
\end{aligned}
\end{equation}
Both $E_i$ and $f_j$ extend smoothly by zero across their respective
interfaces. Taylor's formula in the normal variables therefore gives,
for every $N\in\mathbb N_0$,
\begin{equation}\label{ce:merging-interface-flatness}
E_i\le C_Ny_i^N,\qquad
\sum_{|\eta|\le|\alpha|}|\partial^\eta f_j|\le C_{\alpha,N}y_j^N.
\end{equation}
The constants are taken on compact neighborhoods of $I\times K$.
Using \eqref{ce:merging-source-coefficients},
\eqref{ce:merging-source-decay}, and
\eqref{ce:merging-interface-flatness},
\begin{equation}\label{ce:criterion-flatness}
\begin{aligned}
|\partial^\alpha B|+\sum_{l=0}^r|\partial^\alpha \beta_l|
&\le\sum_{i=1}^2\left(|\partial^\alpha B_i|
+\sum_{l=0}^r|\partial^\alpha \beta_{l,i}|\right)
\le C_\alpha\sum_{i=1}^2E_i
\sum_{|\eta|\le|\alpha|}|\partial^\eta f_{3-i}|\le C_{\alpha,N}y_1^Ny_2^N\\
&
=C_{\alpha,N}\min\{y_1,y_2\}^N\max\{y_1,y_2\}^N
\le C_{I,K,\alpha,N}\min\{y_1,y_2\}^N
\end{aligned}
\end{equation}
on $I\times(K\cap D)$.

\medskip
\noindent\textbf{Step 3. Extend the repair across both interfaces.}
For $\psi\in\{B,\beta_0,\ldots,\beta_r\}$, define
\[
\widetilde\psi(t,x)=
\begin{cases}
\psi(t,x),&x\in D,\\
0,&x\in U\setminus D.
\end{cases}
\]
Let $z_0=(t_0,x_0)\in\R\times(\partial D\cap U)$ and
$h=(h_t,h_x)$. If $x_0+h_x\in D$, then
\[
0<\min_i y_i(x_0+h_x)\le|h_x|\le|h|.
\]
By \eqref{ce:criterion-flatness} with $N=2$,
\[
\frac{\left|\widetilde{\partial^\alpha\psi}(z_0+h)
-\widetilde{\partial^\alpha\psi}(z_0)\right|}{|h|}
\le C\frac{|h_x|^2}{|h|}\le C|h|\longrightarrow0.
\]
For $x_0+h_x\notin D$, the numerator is zero. Hence, inductively,
\[
\partial^\alpha\widetilde\psi=\widetilde{\partial^\alpha\psi},\qquad
\widetilde\psi\in C^\infty(\R\times U),\qquad
\partial^\alpha\widetilde\psi=0\quad\hbox{on }\R\times(U\setminus D).
\]
Continue to write $B,\beta_l$ for these extensions.

On $\R\times D$, set
\[
N\triangleq\chi_*R,\qquad G\triangleq-W_1-W_2+B.
\]
By \eqref{ce:merging-data},
\[
\Box q=R-(W_1+W_2)q
=\chi_*R+\bigl(B-W_1-W_2\bigr)q=N+Gq.
\]
Moreover,
\[
\operatorname{supp}N\subset
\operatorname{supp}(\chi_*W_1)\cup\operatorname{supp}(\chi_*W_2),
\qquad \beta_l=-\frac{N\mathcal C_l}{2F_0}\quad(F_0>0).
\]
Apply the explicit construction of Lemma~\ref{ce:repair}. Use the
zero-extended coefficients $\beta_l$ in \eqref{ce:repair-exponent} to
define $\Phi_{\rm rep},{H_{\rm rep}}$ on all of $\R\times U$.
They are smooth and flat outside $D$, and
\begin{equation}\label{ce:merging-repair-identity}
2D_q\Phi_{\rm rep}=-\chi_*R+q{H_{\rm rep}}\quad\hbox{on }\R\times D.
\end{equation}

\medskip
\noindent\textbf{Step 4. Form the merged pair.}
Set
\[
p\triangleq e^{\Phi_{\rm rep}} q,\qquad
V\triangleq W_1+W_2-B-{H_{\rm rep}}-\Box\Phi_{\rm rep}-\mathcal B(\mathrm{d}\Phi_{\rm rep},\mathrm{d}\Phi_{\rm rep}).
\]
On $\R\times D$, \eqref{ce:merging-repair-identity} gives
\[
\begin{aligned}
\Box p+Vp
&=e^{\Phi_{\rm rep}}\{R-Bq+2D_q\Phi_{\rm rep}-{H_{\rm rep}}q\}\\
&=e^{\Phi_{\rm rep}}\{R-(1-\chi_*)R-\chi_*R+q{H_{\rm rep}}-{H_{\rm rep}}q\}=0.
\end{aligned}
\]
On $\R\times(U\setminus D)$,
\[
(B,\Phi_{\rm rep},{H_{\rm rep}})=(0,0,0),\qquad
(p,V)=(q,W_1+W_2),\qquad
\Box p+Vp=R=0,
\]
and
\[
(p,V)=
\begin{cases}
(f_i,W_i),&y_{3-i}\le0<y_i,\\
(0,0),&y_1\le0,\ y_2\le0.
\end{cases}
\]
On each compact neighborhood, \eqref{ce:repair-exponent} gives
\begin{equation}\label{ce:merging-exponent-flatness}
\begin{aligned}
|\partial^\alpha\Phi_{\rm rep}|+|\partial^\alpha {H_{\rm rep}}|+|\partial^\alpha\Box\Phi_{\rm rep}|
&\le C_\alpha\sum_{l=0}^r
\sum_{|\gamma|\le|\alpha|+r+2}|\partial^\gamma \beta_l|
\le C_{I,K,\alpha,N}\min\{y_1,y_2\}^N.
\end{aligned}
\end{equation}
Also,
\begin{equation}\label{ce:merging-potential-flatness}
\begin{aligned}
|\partial^\alpha\mathcal B(\mathrm d\Phi_{\rm rep},\mathrm d\Phi_{\rm rep})|
&\le\sum_{\eta\le\alpha}\binom\alpha\eta
|\mathrm d\partial^\eta\Phi_{\rm rep}|\,
|\mathrm d\partial^{\alpha-\eta}\Phi_{\rm rep}|
\le C_{I,K,\alpha,N}\min\{y_1,y_2\}^N,\\
|\partial^\alpha(V-W_1-W_2)|
&\le|\partial^\alpha B|+|\partial^\alpha {H_{\rm rep}}|
+|\partial^\alpha\Box\Phi_{\rm rep}|
+|\partial^\alpha\mathcal B(\mathrm d\Phi_{\rm rep},\mathrm d\Phi_{\rm rep})|\\
&\le C_{I,K,\alpha,N}\min\{y_1,y_2\}^N.
\end{aligned}
\end{equation}
Estimates \eqref{ce:merging-exponent-flatness} and
\eqref{ce:merging-potential-flatness} prove
\eqref{ce:merging-correction-flatness}.
For $T$-periodic input pairs, $\delta(t+T,x)=\delta(t,x)$ and
\[
\chi_*(t+T,x)
=\vartheta\!\left(\frac{q(t+T,x)}{\delta(t+T,x)}\right)
=\chi_*(t,x).
\]
The defining formulas then give
\[
(B,\beta_l,\Phi_{\rm rep},{H_{\rm rep}},p,V)(t+T,x)
=(B,\beta_l,\Phi_{\rm rep},{H_{\rm rep}},p,V)(t,x).
\]
\end{proof}

\subsection{Completion to a positive static solution}\label{ce:activation-section}

This subsection joins a given solution $f$ to a positive static
function $A$. The strategy is to interpolate between the two, make a
finite-dimensional perturbation in the transition layer so that the
finite-chain denominator is positive at the zeros of the perturbed
function, and then apply Lemma~\ref{ce:repair} to correct the wave
equation. The main result is
Proposition~\ref{ce:activation-criterion}, which states sufficient
conditions under which the local pair can be completed to the constant
solution.
\begin{proposition}[Completion to a positive static solution]\label{ce:activation-criterion}
Let $n\in\mathbb N$, $T>0$, let $U\subset\R^n$ be open, and let
$X\in C^\infty(U;\R)$ be a spatial level function. Assume that, for
some $\sigma>0$,
\begin{equation}\label{ce:completion-compact-band}
\{x\in U:-\sigma\le X(x)\le1+\sigma\}\Subset U.
\end{equation}
Let $f,W\in C^\infty(\R\times U;\R)$ be $T$-periodic in time and
satisfy $\Box f+Wf=0$, and let $A\in C^\infty(U;(0,\infty))$.
Suppose in addition that
\begin{equation}\label{ce:completion-entrance-hypothesis}
\mathcal B(\mathrm{d}f,\mathrm{d}f)\ne0
\quad\text{on }\{(t,x)\in\R\times U:X(x)=0,\ f(t,x)=0\}.
\end{equation}
Then there exist $\varepsilon\in(0,1/2)$ and real smooth $T$-periodic
functions $p,V$ on $\R\times U$ such that
\[
\begin{cases}
\ds\Box p+Vp=0\qquad \text{in }\R\times U,\\
\ns\ds (p,V)=(f,W)\quad\text{where }X<\varepsilon,\\[2pt]
\ns\ds (p,V)=(A,\Delta A/A)\quad\text{where }X>1-\varepsilon.
\end{cases}
\]
\end{proposition}
\begin{proof}[Proof of Proposition~\ref{ce:activation-criterion}]
We identify $T$-periodic functions of $t$ with functions on
$\R/T\mathbb Z$. Condition \eqref{ce:completion-compact-band} makes
\[
K\triangleq(\R/T\mathbb Z)\times\{x\in U:0\le X(x)\le1\}
\]
compact. If $K=\varnothing$, the disjoint open sets $\{X<0\}$ and
$\{X>1\}$ cover $U$. The direct definition
\[
(p,V)=\begin{cases}
(f,W),&X<0,\\
(A,\Delta A/A),&X>1
\end{cases}
\]
satisfies the conclusion with $\varepsilon=1/4$. Hence assume
$K\ne\varnothing$. Until Step~\ref{step:completion-3}, the
construction takes place on $(\R/T\mathbb Z)\times U$, and supports
are taken relative to this space.

\medskip
\noindent\textbf{Step 1. Interpolate between $f$ and $A$.}
\proofsteplabel{step:completion-1}{1}
At an entrance zero, \eqref{ce:completion-entrance-hypothesis} says
$f_t^2-|\nabla_x f|^2\ne0$. We first extend this inequality to a thin
layer. The set
\[
\{(t,x)\in(\R/T\mathbb Z)\times U:X(x)=0,\ f(t,x)=0\}
\]
is compact. Condition \eqref{ce:completion-entrance-hypothesis} on
this set gives a sufficiently small $\delta>0$ such that
\begin{equation}\label{ce:completion-entrance-persistence}
\mathcal B(\mathrm{d}f,\mathrm{d}f)\ne0
\quad\hbox{on }\{f=0,\ |X|\le4\delta\}.
\end{equation}
Indeed, otherwise there would be points $z_k$ in the compact set
$(\R/T\mathbb Z)\times\{-\sigma\le X\le1+\sigma\}$ with
\[
f(z_k)=0,\qquad X(z_k)\longrightarrow0,\qquad
\mathcal B(\mathrm{d}f,\mathrm{d}f)(z_k)=0.
\]
Passing to a subsequence, $z_k\to z_*$, with
\[
X(z_*)=0,\qquad f(z_*)=0,\qquad
\mathcal B(\mathrm df,\mathrm df)(z_*)=0,
\]
contrary to \eqref{ce:completion-entrance-hypothesis}.
Take $\delta<\min\{\sigma/8,1/16\}$.
Write $\chi_{\rm cut}$ for the smooth step $b$ in
\eqref{ce:cutoff-functions}; thus $\chi_{\rm cut}(s)=0$ for $s\le0$
and $\chi_{\rm cut}(s)=1$ for $s\ge1$. Define
\begin{equation}\label{ce:completion-interpolation}
\begin{cases}
\ds \theta_{\rm int}\triangleq\chi_{\rm cut}\bigl((X-4\delta)/(1-8\delta)\bigr),\\
\ns\ds q_0\triangleq(1-\theta_{\rm int})f+\theta_{\rm int} A,\\
\ns\ds W_0\triangleq(1-\theta_{\rm int})W+\theta_{\rm int}\,\Delta A/A.
\end{cases}
\end{equation}
Here $1-8\delta>0$ because $\delta<1/16$.
Thus $(q_0,W_0)=(f,W)$ on $\{X\le4\delta\}$ and
$(q_0,W_0)=(A,\Delta A/A)$ on $\{X\ge1-4\delta\}$.
The cutoff depends only on $X(x)$, so these functions are $T$-periodic.
Interpolation preserves the endpoint data. The wave equation may still
have a residual. Step~\ref{step:completion-3} will compute and correct
it after the perturbation has been chosen.

Compactness and \eqref{ce:completion-entrance-persistence} give
$c_0,\eta_0>0$ such that
\begin{equation}\label{ce:completion-endpoint-bounds}
\begin{aligned}
&|\mathcal B(\mathrm{d}q_0,\mathrm{d}q_0)|\ge2c_0
\qquad(0\le X\le2\delta,\ |q_0|\le2\eta_0),\\
&A\ge2c_0\qquad(1-2\delta\le X\le1).
\end{aligned}
\end{equation}
For the first assertion, failure would give points $z_k$ in the
compact set $K\cap\{0\le X\le2\delta\}$ such that
\[
|q_0(z_k)|\longrightarrow0,\qquad
|\mathcal B(\mathrm{d}q_0,\mathrm{d}q_0)(z_k)|\longrightarrow0.
\]
At a limit point, $q_0=f=0$ and
$\mathcal B(\mathrm df,\mathrm df)=0$ with $0\le X\le2\delta$,
contradicting \eqref{ce:completion-entrance-persistence}.
The second assertion follows from $A>0$ on the compact set
$K\cap\{1-2\delta\le X\le1\}$.
The first bound in \eqref{ce:completion-endpoint-bounds} will keep the
quadratic form nonzero at perturbed entrance zeros. The second will
keep the exit strictly positive, so there will be no zeros to check
there. Only the middle layer remains for Step~\ref{step:completion-2}.

\medskip
\noindent\textbf{Step 2. Make the finite chain nonzero at every zero
of the interpolation.}
\proofsteplabel{step:completion-2}{2}
To apply Lemma~\ref{ce:repair}, we perturb the middle layer so that
the denominator $F$ is positive near the zeros. Set $r=n$.

\smallskip
\noindent\textit{Auxiliary lemma (finite-dimensional perturbation).}
For every $\epsilon_1>0$, there is a real smooth function $q$ on
$(\R/T\mathbb Z)\times U$ such that
\begin{equation}\label{ce:completion-perturbation-claim}
\begin{cases}
\ds\operatorname{supp}(q-q_0)\Subset
(\R/T\mathbb Z)\times\{\delta<X<1-\delta\},\\
\ns\ds \|q-q_0\|_{C^1(K)}<\epsilon_1,\\
\ns\ds q=f\quad(X\le\delta),\qquad q=A\quad(X\ge1-\delta),\\
\ns\ds F\triangleq\sum_{\ell=0}^r
\bigl(D_q^\ell\mathcal B(\mathrm dq,\mathrm dq)\bigr)^2>0
\quad\hbox{on }K\cap\{q=0\}.
\end{cases}
\end{equation}
Because the summands are real squares, the positivity condition in
\eqref{ce:completion-perturbation-claim} is equivalent to
\[
z\in K,\quad q(z)=0\quad\Longrightarrow\quad
D_q^\ell\mathcal B(\mathrm dq,\mathrm dq)(z)\ne0
\quad\hbox{for at least one }\ell\in\{0,\ldots,r\}.
\]
The support condition preserves $f$ and $A$ at the two ends;
positivity of $F$ makes the repair in Step~\ref{step:completion-3}
possible.
Figure~\ref{fig:completion-regions} shows the roles of the intervals.

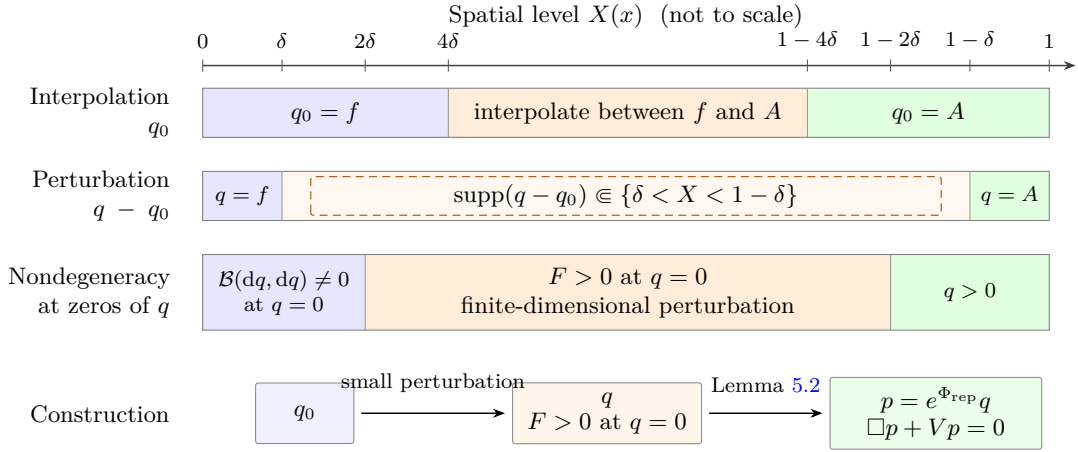
\begin{figure}[!htbp]
\centering
\begin{tikzpicture}[x=1cm,y=1cm,font=\footnotesize,
>={Stealth[length=1.6mm]}]
\node at (8.65,4.95) {Spatial level $X(x)$ \ (not to scale)};
\draw[->,black!75,line width=.55pt] (3.05,4.3)--(14.6,4.3);
\foreach \pos/\lab in {3.05/0,4.1/\delta,5.2/2\delta,6.3/4\delta,
11.05/1-4\delta,12.15/1-2\delta,13.2/1-\delta,14.25/1} {
\draw[black!65] (\pos,4.22)--(\pos,4.38);
\node[font=\scriptsize,anchor=south] at (\pos,4.4) {$\lab$};
}

\node[anchor=east,align=right,text width=2.6cm] at (2.75,3.675)
{Interpolation\\$q_0$};
\fill[blue!10] (3.05,3.35) rectangle (6.3,4.0);
\fill[orange!15] (6.3,3.35) rectangle (11.05,4.0);
\fill[green!12] (11.05,3.35) rectangle (14.25,4.0);
\draw[black!50] (3.05,3.35) rectangle (14.25,4.0);
\draw[black!40] (6.3,3.35)--(6.3,4.0) (11.05,3.35)--(11.05,4.0);
\node at (4.675,3.675) {$q_0=f$};
\node at (8.675,3.675) {interpolate between $f$ and $A$};
\node at (12.65,3.675) {$q_0=A$};

\node[anchor=east,align=right,text width=2.6cm] at (2.75,2.575)
{Perturbation\\$q-q_0$};
\fill[blue!10] (3.05,2.25) rectangle (4.1,2.9);
\fill[orange!6] (4.1,2.25) rectangle (13.2,2.9);
\fill[green!12] (13.2,2.25) rectangle (14.25,2.9);
\draw[black!50] (3.05,2.25) rectangle (14.25,2.9);
\draw[black!40] (4.1,2.25)--(4.1,2.9) (13.2,2.25)--(13.2,2.9);
\draw[densely dashed,orange!60!black,rounded corners=1pt]
(4.48,2.32) rectangle (12.82,2.83);
\node[font=\scriptsize] at (3.575,2.575) {$q=f$};
\node at (8.65,2.575)
{$\operatorname{supp}(q-q_0)\Subset\{\delta<X<1-\delta\}$};
\node[font=\scriptsize] at (13.725,2.575) {$q=A$};

\node[anchor=east,align=right,text width=2.6cm] at (2.75,1.3)
{Nondegeneracy\\at zeros of $q$};
\fill[blue!10] (3.05,.8) rectangle (5.2,1.8);
\fill[orange!15] (5.2,.8) rectangle (12.15,1.8);
\fill[green!12] (12.15,.8) rectangle (14.25,1.8);
\draw[black!50] (3.05,.8) rectangle (14.25,1.8);
\draw[black!40] (5.2,.8)--(5.2,1.8) (12.15,.8)--(12.15,1.8);
\node[font=\scriptsize,align=center] at (4.125,1.3)
{$\mathcal B(\mathrm dq,\mathrm dq)\ne0$\\at $q=0$};
\node[align=center] at (8.675,1.3)
{$F>0$ at $q=0$\\finite-dimensional perturbation};
\node[font=\scriptsize,align=center] at (13.2,1.3)
{$q>0$};

\node[anchor=east] at (2.75,-.3) {Construction};
\node[draw=black!45,fill=blue!6,rounded corners=1pt,
minimum width=1.3cm,minimum height=.8cm] at (4.4,-.3) {$q_0$};
\node[draw=black!45,fill=orange!8,rounded corners=1pt,
minimum width=2.5cm,minimum height=.8cm,align=center] at (8.4,-.3)
{$q$\\[-1pt]$F>0$ at $q=0$};
\node[draw=black!45,fill=green!8,rounded corners=1pt,
minimum width=2.8cm,minimum height=.8cm,align=center] at (12.75,-.3)
{$p=e^{\Phi_{\rm rep}} q$\\[-1pt]$\Box p+Vp=0$};
\draw[->,line width=.55pt] (5.15,-.3)--(7.05,-.3)
node[midway,above=4pt,font=\scriptsize] {small perturbation};
\draw[->,line width=.55pt] (9.75,-.3)--(11.25,-.3)
node[midway,above=4pt,font=\scriptsize] {Lemma~\ref{ce:repair}};
\end{tikzpicture}
\caption{Spatial organization of the completion construction.
The horizontal axis shows values of $X(x)$ on $0\le X\le1$,
schematically and not to scale. The middle row indicates where
perturbations are allowed, with support strictly inside the marked
interval. In the third row, the two end conditions follow from
Step~\ref{step:completion-1} and smallness of the perturbation, while the middle
condition is obtained by the finite-dimensional perturbation in Step~\ref{step:completion-2}.}
\label{fig:completion-regions}
\end{figure}

Appendix~\ref{app:perturbation-details} proves this lemma. The bad
finite-order derivative data have two types: $q=0$, $\mathrm dq=0$; or
$q=0$, $\mathrm dq\ne0$ with every chain term zero. With $r=n$, both
have codimension $n+2$, exceeding the space-time dimension $n+1$.
Finite-dimensional perturbations therefore avoid a measure-zero set of
bad parameters (Sard's theorem), while preserving the endpoints and
the $C^1$ smallness condition.

\medskip
\noindent\textbf{Step 3. Remove the residual without dividing by $q$
at its zeros.}
\proofsteplabel{step:completion-3}{3}
The function $q$ now has the required finite chain, but we have not
yet made it solve the equation with potential $W_0$. Set
$R\triangleq\Box q+W_0q$. We will split $R$ into a part near $q=0$,
to be removed by Lemma~\ref{ce:repair}, and a part away from $q=0$,
which can be absorbed directly into the potential by division by $q$.

First compute the residual. Since $\theta_{\rm int}$ and $A$ are
independent of time, $\Box\theta_{\rm int}=-\Delta\theta_{\rm int}$,
$\Box A=-\Delta A$, and $\Box f=-Wf$. The product rule applied to
\eqref{ce:completion-interpolation} gives
\[
\begin{cases}
\ds	\Box q_0
=-(1-\theta_{\rm int})Wf-\theta_{\rm int}\Delta A
-2\nabla\theta_{\rm int}\cdot\nabla(A-f)-(\Delta\theta_{\rm int})(A-f),\\
\ns\ds	W_0q_0
=(1-\theta_{\rm int})^2Wf
+\theta_{\rm int}(1-\theta_{\rm int})\bigl(WA+(\Delta A/A)f\bigr)+\theta_{\rm int}^2\Delta A.
\end{cases}
\]
Adding these expressions and using $q=q_0+\sum_{j=1}^J a_j\varphi_j$
yields
\begin{equation}\label{ce:completion-residual}
\begin{aligned}
R&=-2\nabla\theta_{\rm int}\cdot\nabla(A-f)-(\Delta\theta_{\rm int})(A-f)\\
&\quad+\theta_{\rm int}(1-\theta_{\rm int})(W-\Delta A/A)(A-f)
+\sum_{j=1}^J a_j(\Box\varphi_j+W_0\varphi_j).
\end{aligned}
\end{equation}
The factors $\nabla\theta_{\rm int}$, $\Delta\theta_{\rm int}$, and
$\theta_{\rm int}(1-\theta_{\rm int})$ vanish outside
$\{4\delta\le X\le1-4\delta\}$, and each $\varphi_j$ has compact
support in $\{\delta<X<1-\delta\}$. Thus
\[
\operatorname{supp}R\Subset
(\R/T\mathbb Z)\times\{\delta<X<1-\delta\}\subset K.
\]

Next choose the neighborhood of the zeros where the repair will act.
Compactness of $K$ and \eqref{ce:completion-chain-positive} give
$\eta,c_F>0$ such that
\[
F\ge c_F\quad\hbox{on }K\cap\{|q|\le2\eta\}.
\]
Indeed, failure would give a sequence $z_k\in K$ with
$|q(z_k)|\to0$ and $F(z_k)\to0$. A convergent subsequence would give
$q(z_*)=F(z_*)=0$, contrary to
\eqref{ce:completion-chain-positive}. Using the cutoff from
Step~\ref{step:completion-1}, define
\begin{equation}\label{ce:completion-residual-split}
\begin{aligned}
\chi&\triangleq\chi_{\rm cut}(2+q/\eta)\chi_{\rm cut}(2-q/\eta),\\
N&\triangleq\chi R,\qquad
G\triangleq-W_0+\frac{(1-\chi)R}{q}.
\end{aligned}
\end{equation}
Here $\chi=1$ when $|q|\le\eta$ and $\chi=0$ when $|q|\ge2\eta$.
Define the quotient in \eqref{ce:completion-residual-split} to be zero
at $q=0$. It is then identically zero on the open set
$\{|q|<\eta\}$, and is a quotient with nonzero denominator wherever
$q\ne0$. These definitions agree on their overlap, so $G$ is smooth.
They also give, including at $q=0$,
\begin{equation}\label{ce:completion-source-support}
\begin{cases}
\ds\left(\frac{(1-\chi)R}{q}\right)q =(1-\chi)R,\\
\ns\ds	N+Gq =\chi R-W_0q+(1-\chi)R=R-W_0q=\Box q,\\
\ns\ds \operatorname{supp}N \subset\operatorname{supp}R\cap\{|q|\le2\eta\}
\subset K\cap\{F\ge c_F\}.
\end{cases}
\end{equation}
We now regard $q,N,G$ as $T$-periodic functions on $\R\times U$.
Inclusion \eqref{ce:completion-source-support} ensures $F>0$ on
$\operatorname{supp}_{\R\times U}N$, so Lemma~\ref{ce:repair} applies
with $r=n$. Use its coefficients $\beta_\ell$ from
\eqref{ce:repair-b} and its functions $\Phi_{\rm rep},{H_{\rm rep}}$
from \eqref{ce:repair-exponent}.
In particular, ${H_{\rm rep}}$ is determined by the same coefficients:
\[
{H_{\rm rep}}=2\sum_{\ell=0}^r(-1)^\ell D_q^{\ell+1}\beta_\ell,
\qquad 2D_q\Phi_{\rm rep}=-N+q{H_{\rm rep}},
\]
where the second identity is \eqref{ce:repair-multiplier-identity}.
All these functions are smooth by
Corollary~\ref{ce:repair-properties}\ref{cond:repair-smooth}.
The output formula \eqref{ce:repair-output} now gives
\[
\begin{aligned}
&p=e^{\Phi_{\rm rep}} q,\qquad
V =-G-{H_{\rm rep}}-\Box\Phi_{\rm rep}-\mathcal B(\mathrm d\Phi_{\rm rep},\mathrm d\Phi_{\rm rep}),\\
&e^{-\Phi_{\rm rep}}(\Box p+Vp)=N+2D_q\Phi_{\rm rep}-q{H_{\rm rep}}=0.
\end{aligned}
\]
The last line follows from \eqref{ce:repair-conjugation-check}.

\medskip
\noindent\textbf{Step 4. Recover the prescribed functions and the
time period.}
\proofsteplabel{step:completion-4}{4}
We check that the repair leaves the two ends unchanged. On the open
set $E\triangleq\{X<\delta\}\cup\{X>1-\delta\}$,
\eqref{ce:completion-residual} gives $R=0$, and hence
\eqref{ce:completion-residual-split} gives $N=0$ and $G=-W_0$.
The definition \eqref{ce:repair-b} then gives
\[
\beta_\ell|_E=0,\qquad (D_q^a \beta_\ell)|_E=0
\quad(0\le\ell\le r,\ a\in\mathbb N_0).
\]
Thus \eqref{ce:repair-exponent} implies
$\Phi_{\rm rep}={H_{\rm rep}}=0$ on $E$; all coordinate derivatives
of $\Phi_{\rm rep}$ vanish there as well.
It follows that $(p,V)=(q,W_0)$ on $E$. By
\eqref{ce:completion-interpolation} and
\eqref{ce:completion-perturbation-claim},
\[
(p,V)=(f,W)\quad(X<\delta),\qquad
(p,V)=(A,\Delta A/A)\quad(X>1-\delta).
\]
Take $\varepsilon=\delta\in(0,1/2)$.

Finally, $q,R,\chi,N,G$ are all $T$-periodic. For
$\tau_T(t,x)\triangleq(t+T,x)$, the commutation rule
\eqref{ce:repair-translation-rule} gives, by iteration,
\[
D_q^a(h\circ\tau_T)=(D_q^a h)\circ\tau_T
\quad(a\in\mathbb N_0,\ h\in C^\infty(\R\times U)).
\]
Applying this rule in \eqref{ce:repair-chain}--\eqref{ce:repair-output}
gives
\[
(\mathcal C_\ell,F,\beta_\ell,\Phi_{\rm rep},{H_{\rm rep}},p,V)\circ\tau_T
=(\mathcal C_\ell,F,\beta_\ell,\Phi_{\rm rep},{H_{\rm rep}},p,V).
\]
Thus the constructed pair has the prescribed period.
\end{proof}

\medskip\noindent\textbf{Periodicity and derivative bounds.}
For $T$-periodic input pairs, Proposition~\ref{ce:merging-criterion}
preserves the period $T$; Proposition~\ref{ce:activation-criterion}
preserves the period in its hypotheses as just proved. In either case,
smoothness and periodicity give, for every compact $K\subset U$ and
$m\in\mathbb N_0$,
\begin{equation}\label{ce:criterion-periodic-bounds}
\begin{aligned}
\sup_{\R\times K}\sum_{|\alpha|\le m}
\bigl(|\partial^\alpha p|+|\partial^\alpha V|\bigr)
&=\max_{[0,T]\times K}\sum_{|\alpha|\le m}
\bigl(|\partial^\alpha p|+|\partial^\alpha V|\bigr)<\infty.
\end{aligned}
\end{equation} 

\section{The planar counterexample to observability under GCC}\label{ce:completion-section}

This section completes the proof of Theorem~\ref{thm:counter}. We have
already constructed the local oscillatory pairs
$(f_\gamma,V_\gamma)$ in Section~\ref{ce:front-section}, the GCC
observation geometry in Subsection~\ref{ce:geometry-section}, the
finite-chain repair in Lemma~\ref{ce:repair}, and the two general
gluing criteria in Propositions~\ref{ce:merging-criterion}
and~\ref{ce:activation-criterion}. It remains to verify these criteria
for the specific construction, to assemble the local pairs into a
global pair, and to impose the Dirichlet boundary condition.

The section is organized as follows. Subsection~\ref{ce:vertex-verification}
verifies the merging criterion at each vertex. This is the most
technical part: it requires a quantitative lower bound for the
finite-chain denominator near the cancellation set of the two edge
fronts. Subsection~\ref{ce:activation-verification} constructs the
spatial transition, verifies the completion criterion, and produces a
common static exit. Subsection~\ref{ce:assembly} assembles the local
pairs, glues them to an exterior Dirichlet pair, and proves the failure
of observability.

\subsection{Merging at the vertices}\label{ce:vertex-verification}

In this subsection we verify Proposition~\ref{ce:merging-criterion} at
a vertex. The two edge fronts are the real smooth pairs
$(f_0,V_0)$ and $(f_\gamma,V_\gamma)$ constructed in
Section~\ref{ce:front-section} on $|z|<1$, $|y|<Y_*$, where
$\gamma\triangleq24\log2$ and $\beta\triangleq2^{1/16}$.
Both pairs solve the corresponding wave equation, vanish for $y\le0$,
and are $2\pi$-periodic in time and flat at $y=0$. We place one pair
along each of the two edges and merge them where both normal
coordinates are positive.

Let $v\in\R^2$ be the vertex. Choose unit tangents $e_1,e_2$ and let
$n_i$ be the rotation of $e_i$ through $\pi/2$ counterclockwise,
$i\in\{1,2\}$. Suppose
$e_2=\cos\alpha_v\,e_1+\sin\alpha_v\,n_1$, where
$0<|\alpha_v|=\alpha_0\le\pi/32$. The normals point out of the
interior wedge. Choose $R>0$ and fixed tangential offsets
$\zeta_1,\zeta_2$ such that $R<\min\{|\zeta_1|,|\zeta_2|\}$.
For each integer $M$, set
\begin{equation}\label{ce:vertex-coordinates}
\begin{aligned}
&y_*\triangleq\beta^{-M},\qquad
C_M\triangleq B(v,Ry_*),\qquad y_i\triangleq n_i\cdot(x-v),\qquad
z_i\triangleq\zeta_i y_*+e_i\cdot(x-v).
\end{aligned}
\end{equation}
Here $R,\zeta_i$ are fixed, while $y_*=\beta^{-M}$ is the normal
scale that shrinks as $M$ increases. The quantity $y_i$ is the signed
normal distance from edge $i$, and $z_i$ is its tangential coordinate.
The set $C_M$ is the open disk of radius $Ry_*$ centered at $v$. The
interior wedge is $\{y_1\le0,y_2\le0\}$ and has angle
$\pi-\alpha_0$. For $M$ large enough that the coordinates remain in
the domain of the local solutions, place the two solutions in these
coordinates:
\[
\begin{aligned}
&(f_1,W_1)(t,x)\triangleq(f_0,V_0)(t,z_1(x),y_1(x)),\qquad (f_2,W_2)(t,x)\triangleq(f_\gamma,V_\gamma)(t,z_2(x),y_2(x)).
\end{aligned}
\]
\begin{proposition}[Vertex merging]\label{ce:incoming}
There exists an integer $M_0=M_0(\alpha_0,R,\zeta_1,\zeta_2)$, which
may also depend on the chosen family of local solutions, with the
following property. For every integer $M\ge M_0$, there is a real
smooth $2\pi$-periodic pair $(p,V)$ on a neighborhood of
$\R\times\overline{C_M}$ satisfying
\[
\Box p+Vp=0.
\]
On $C_M$, this pair agrees with the original edge solutions in the
single-positive regions and vanishes in the double-negative region:
\[
(p,V)=
\begin{cases}
(f_i,W_i),&y_{3-i}<0<y_i,\quad i\in\{1,2\},\\
(0,0),&y_1\le0,\ y_2\le0.
\end{cases}
\]
In the double-positive region, it has the form
\[
p=e^{\Phi_0}(f_1+f_2)
\]
for a smooth real function $\Phi_0$. Both $p$ and $V$ are flat at
$x=v$ for every $t$, and their derivatives of every fixed order are
bounded on $\R\times\overline{C_M}$.
\end{proposition}
\begin{figure}[H]
\centering
\begin{tikzpicture}[x=1cm,y=1cm,font=\small,
>={Stealth[length=1.6mm]}]
\pgfmathsetmacro{\tipheight}{3.4*tan(180/64)}
\fill[black!7]
(-3.4,-\tipheight)--(0,0)--(3.4,-\tipheight)
--(3.4,-1.85)--(-3.4,-1.85)--cycle;
\fill[blue!9!white]
(-3.4,\tipheight)--(0,0)--(3.4,\tipheight)
--(3.4,1.75)--(-3.4,1.75)--cycle;
\fill[orange!22]
(-3.4,-\tipheight)--(0,0)--(-3.4,\tipheight)--cycle;
\fill[green!17!white]
(3.4,-\tipheight)--(0,0)--(3.4,\tipheight)--cycle;
\draw[black!30,line width=.4pt] (-3.4,-1.85) rectangle (3.4,1.75);
\draw[line width=1.1pt]
(-3.4,-\tipheight)--(0,0)--(3.4,-\tipheight);
\draw[densely dashed,line width=.65pt]
(-3.4,\tipheight)--(0,0)--(3.4,\tipheight);
\fill (0,0) circle[radius=1.6pt];
\node[anchor=south west,inner sep=3pt] at (0,0) {$v$};
\node[font=\footnotesize,rotate=2.8125] at (-2.0,-.32) {edge $1$};
\node[font=\footnotesize,rotate=-2.8125] at (2.0,-.32) {edge $2$};
\draw[black!65,line width=.45pt]
(182.8125:.68) arc[start angle=182.8125,end angle=357.1875,radius=.68];
\node[font=\footnotesize] at (0,-.47) {$\pi-\alpha_0$};

\node at (0,1.42) {$y_1>0,\quad y_2>0$};
\node at (0,.99) {Repaired overlap};
\node at (0,.54) {$p=e^{\Phi_0}(f_1+f_2)$};

\node at (0,-.98) {Interior wedge: $y_1\le0,\ y_2\le0$};
\node at (0,-1.48) {$(p,V)=(0,0)$};

\node[font=\footnotesize,align=right,anchor=east] at (-3.75,.85)
{$y_1>0,\ y_2<0$\\[3pt]Edge $1$ solution\\[3pt]$(p,V)=(f_1,W_1)$};
\draw[->,line width=.5pt]
(-3.90,.25)--(-3.55,.25)--(-2.70,0);
\node[font=\footnotesize,align=left,anchor=west] at (3.75,.85)
{$y_2>0,\ y_1<0$\\[3pt]Edge $2$ solution\\[3pt]$(p,V)=(f_2,W_2)$};
\draw[->,line width=.5pt]
(3.90,.25)--(3.55,.25)--(2.70,0);
\end{tikzpicture}
\caption{Vertex merging near $v$, drawn to scale for
$\alpha_0=\pi/32$ (interior angle $31\pi/32$). Solid rays are the
incident edges; dashed rays extend their supporting lines. The case
$\alpha_v=-\alpha_0$ is shown. The narrow side sectors retain the
original edge pairs, while the exterior overlap carries the repaired
sum. The pairs match smoothly across the separating lines.}
\label{fig:vertex-merging}
\end{figure}
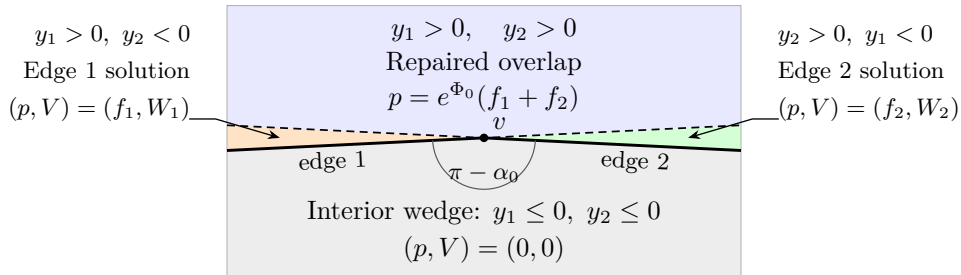

The sum $f_1+f_2$ already agrees with the required edge solution where
only one of the two solutions is nonzero. In the overlap
$\{y_1>0,y_2>0\}$, Proposition~\ref{ce:merging-criterion} supplies the
correction. To verify its three conditions, take the envelopes
\begin{equation}\label{ce:incoming-envelope-definition}
E_i(x)\triangleq
\begin{cases}
\displaystyle\sum_{j\ge4}
e^{-1/y_i-j^2-\gamma_i j}\rho(\beta^jy_i),&y_i>0,\\
0,&y_i\le0,
\end{cases}
\qquad \gamma_1=0,\quad\gamma_2=\gamma.
\end{equation}
We work on the slightly larger neighborhood
\[
R_+\triangleq\frac{R+\min_{i=1,2}|\zeta_i|}{2},\qquad
U\triangleq B(v,R_+\beta^{-M}),\qquad
D\triangleq U\cap\{y_1>0,y_2>0\}.
\]
Then $\overline{C_M}\subset U$ and $R_+<\min_i|\zeta_i|$.
All supports below are relative to $\R\times D$.
Constants may depend on the fixed local data and the two solution
families, but not on $M$.

The envelope $E_i$ measures the amplitude of $f_i$ without its
oscillatory factors. On $\{y_i>0\}$, it stays positive even at zeros
of $f_i$. For $y_i>0$, let $m_i$ be the integer satisfying
$\beta^{-(m_i+1)}<y_i\le\beta^{-m_i}$. For sufficiently large $M$,
both indices are at least $7$. Only the terms $j=m_i,m_i+1$ in
\eqref{ce:incoming-envelope-definition} can be nonzero, and
\[
\rho(\beta^{m_i}y_i)+\rho(\beta^{m_i+1}y_i)=1.
\]
Thus $E_i$ is a convex combination of the two positive amplitudes,
which gives
\begin{equation}\label{ce:envelope}
e^{-1/y_i-(m_i+1)^2-\gamma_i(m_i+1)}
\le E_i\le e^{-1/y_i-m_i^2-\gamma_i m_i}.
\end{equation}
The sum defining $E_i$ in \eqref{ce:incoming-envelope-definition} is
locally finite for $y_i>0$, and its derivatives satisfy
$|\partial^\alpha E_i|\le C_\alpha e^{-1/(2y_i)}$ near $y_i=0$.
Thus $E_i$ is smooth on $U$, positive where $y_i>0$, and flat where
$y_i=0$.

\medskip\noindent\textbf{The first two conditions of the merging
proposition.}
Equations~\eqref{ce:packet-envelope}, \eqref{ce:front-jets}, and
\eqref{ce:source-margin}, together with
$\operatorname{supp}W_i=\operatorname{supp}(\Box f_i)$, give
\begin{equation}\label{ce:incoming-amplitude-bounds}
|f_i|\le C_{\rm env}E_i\quad\hbox{on }\R\times U,\qquad
f_i\ge\tfrac58E_i\quad\hbox{on }\operatorname{supp}W_i.
\end{equation}
This verifies condition~\ref{cond:merging-amplitude} of
Proposition~\ref{ce:merging-criterion} with $c=5/8$.

On a neighborhood of $\operatorname{supp}W_i$ in $\R\times D$, we have
$W_i=-\Box f_i/f_i$. Estimate \eqref{ce:incoming-amplitude-bounds}
gives, for every fixed multi-index $\nu$,
\[
|\partial^\nu(f_i^{-1})|
\le C_\nu E_i^{-1}y_i^{-N_\nu}.
\]
Hence the product rule, \eqref{ce:reserve-estimate}, and
\eqref{ce:envelope} imply, for every fixed $\alpha\in\mathbb N_0^3$
and $P\in\mathbb N$,
\begin{equation}\label{ce:incoming-paid}
\begin{aligned}
E_i^{-P}|\partial^\alpha W_i|
&\le E_i^{-P}\sum_{\nu\le\alpha}\binom\alpha\nu
|\partial^{\alpha-\nu}(\Box f_i)|\,|\partial^\nu(f_i^{-1})|\\
&\le C_{\alpha,P}\exp\bigl\{-cm_i^3\beta^{m_i}
+C_{\alpha,P}m_i^2\beta^{m_i}\log(m_i+2)\bigr\}
\le C'_{\alpha,P}.
\end{aligned}
\end{equation}
Outside $\operatorname{supp}W_i$ all derivatives of $W_i$ vanish.
Thus condition~\ref{cond:merging-potential} of
Proposition~\ref{ce:merging-criterion} holds as well.

\medskip\noindent\textbf{The remaining condition: a lower bound for a
finite derivative chain.}
It remains to verify condition~\ref{cond:merging-chain} of
Proposition~\ref{ce:merging-criterion}. Choose
$\chi\in C^\infty(\R;[0,1])$ equal to one on $[-1/8,1/8]$
and zero outside $(-1/4,1/4)$, and set, on $\R\times D$,
\[
\begin{aligned}
&q\triangleq f_1+f_2,\qquad
\delta\triangleq\frac{E_1E_2}{E_1+E_2},\qquad
\chi_*\triangleq\chi(q/\delta),\\
&Q\triangleq\mathcal B(\mathrm dq,\mathrm dq),\qquad
\mathcal C_l\triangleq D_q^lQ\quad(0\le l\le3),\qquad
F_0\triangleq\sum_{l=0}^3\mathcal C_l^2.
\end{aligned}
\]
With $a=1/8$, $b=1/4$, and $\vartheta=\chi$, the remaining target is
\begin{equation}\label{ce:incoming-chain-target}
\begin{aligned}
F_0=\sum_{l=0}^3
\bigl(D_q^l\mathcal B(\mathrm dq,\mathrm dq)\bigr)^2
&\ge cE_i^{11}
\quad\hbox{on }\operatorname{supp}_{\R\times D}(\chi_*W_i),
\qquad i\in\{1,2\}.
\end{aligned}
\end{equation}
This is precisely condition~\ref{cond:merging-chain} of
Proposition~\ref{ce:merging-criterion} with $r=3$ and $\kappa=11$.
The calculations in Subsections~\ref{ce:small-owner}
and~\ref{ce:comparable} prove
\eqref{ce:incoming-chain-target} by separating the cases of a much
smaller envelope and comparable envelopes. Once
\eqref{ce:incoming-chain-target} is established,
Proposition~\ref{ce:merging-criterion} gives the pair shown in
Figure~\ref{fig:vertex-merging}. The exponent $11$ is a convenient
sufficient choice. To explain its source, set
$E_{\rm tot}\triangleq E_1+E_2$ and $m\triangleq\min(m_1,m_2)$. The
final normalization in Step~\ref{step:incoming-power} of the proof
below gives
\[
F_0\ge cE_{\rm tot}^{10}\beta^{-28m},
\qquad
E_{\rm tot}\le C\beta^{-28m}
\quad\Longrightarrow\quad
F_0\ge c'E_{\rm tot}^{11}.
\]

\medskip\noindent\textbf{Coordinates used in the calculation.}
For the calculation of $\mathcal C_0,\mathcal C_1,\mathcal C_2,
\mathcal C_3$, use the coordinates
\[
Z_i\triangleq z_i(v)=\zeta_i y_*,\qquad
u_i\triangleq t+z_i,\qquad v_i\triangleq t-z_i,
\qquad i\in\{1,2\}.
\]
Let $\mathsf J(x_1,x_2)\triangleq(-x_2,x_1)$ and define the bisector
frame by
\[
\begin{aligned}
&c_\alpha\triangleq\cos(\alpha_0/2),\quad
s_\alpha\triangleq\sin(\alpha_v/2),\quad
\mathfrak g\triangleq1-\cos\alpha_0=2s_\alpha^2,\quad e_*\triangleq(e_1+e_2)/(2c_\alpha),\quad
n_*\triangleq\mathsf J e_*.
\end{aligned}
\]
Projection onto the two edge frames gives
\begin{equation}\label{ce:vertex-affine}
\begin{aligned}
&z_1=Z_1+\frac{\cos\alpha_v\,y_1-y_2}{\sin\alpha_v},\qquad
z_2=Z_2+\frac{y_1-\cos\alpha_v\,y_2}{\sin\alpha_v},\\
&y_1-y_2=2s_\alpha\,e_*\cdot(x-v),\qquad
z_2-z_1=Z_2-Z_1+2s_\alpha\,n_*\cdot(x-v).
\end{aligned}
\end{equation}
Equivalently, the coordinate change and its inverse are
\[
\begin{aligned}
\binom{z_1-Z_1}{z_2-Z_2}
&=\frac1{\sin\alpha_v}
\begin{pmatrix}\cos\alpha_v&-1\\1&-\cos\alpha_v\end{pmatrix}
\binom{y_1}{y_2},\\
\binom{y_1}{y_2}
&=\frac1{\sin\alpha_v}
\begin{pmatrix}-\cos\alpha_v&1\\-1&\cos\alpha_v\end{pmatrix}
\binom{z_1-Z_1}{z_2-Z_2}.
\end{aligned}
\]
The sign of $\alpha_v$ allows either labelling of the two edge types.
Since $|e_i\cdot(x-v)|\le R_+\beta^{-M}$ and
$R_+<\min_i|\zeta_i|$, the offsets $Z_i$ give
\begin{equation}\label{ce:midpoint-size}
\begin{cases}
(\min\{|\zeta_1|, |\zeta_2|\}-R_+)\beta^{-M}\le |z_i|
\le(\min\{|\zeta_1|, |\zeta_2|\}+R_+)\beta^{-M},\\[3pt]
|y_i|\le R_+\beta^{-M}.
\end{cases}
\end{equation}

The rest of this subsection is devoted to proving
\eqref{ce:incoming-chain-target}. We first record the signs and
envelope bounds on $\operatorname{supp}(\chi_*W_i)$, then split into
the two cases of a much smaller envelope and comparable envelopes.
\begin{lemma}[Signs and envelope bounds]\label{ce:source-owner}
Set
\begin{equation}\label{ce:incoming-supports}
\mathcal U_i\triangleq
\operatorname{supp}_{\R\times D}(\chi_*W_i),\qquad i\in\{1,2\}.
\end{equation}
Then $\mathcal U_1\cap\mathcal U_2=\varnothing$. On $\mathcal U_i$,
\begin{equation}\label{ce:owner-signs}
f_i\ge5E_i/8,\quad f_{3-i}\le-3E_i/8,\quad
|f_{3-i}|\le(C_{\rm env}+1/4)E_i,\quad E_i\le K_0E_{3-i},
\end{equation}
where $K_0\triangleq8C_{\rm env}/3+1$ and $C_{\rm env}$ is the constant
in \eqref{ce:incoming-amplitude-bounds}.
\end{lemma}
\begin{proof}[Proof of Lemma~\ref{ce:source-owner}]
By the support of $\chi$ and \eqref{ce:incoming-amplitude-bounds},
\[
\mathcal U_i\subset
\operatorname{supp}_{\R\times D}W_i
\cap\{|q|\le\delta/4\},\qquad
f_i\ge\tfrac58E_i\quad\hbox{on }\mathcal U_i.
\]
If a point belonged to both sets, then at that point
\[
\tfrac58(E_1+E_2)\le q\le |q|
\le\tfrac14\delta\le\tfrac1{16}(E_1+E_2),
\]
which is impossible since $E_1+E_2>0$ on $D$.
Thus $\mathcal U_1\cap\mathcal U_2=\varnothing$.

Fix $i\in\{1,2\}$ and work on $\mathcal U_i$.
Since $\delta\le E_i$ and $|f_i|\le C_{\rm env}E_i$,
\[
\begin{cases}
\ds f_{3-i}=q-f_i
\le\tfrac14\delta-\tfrac58E_i\le-\tfrac38E_i,\\
\ns\ds |f_{3-i}| \le |q|+|f_i|\le(C_{\rm env}+\tfrac14)E_i.
\end{cases}
\]
Finally,
\[
\tfrac38E_i\le |f_{3-i}|\le C_{\rm env}E_{3-i},\qquad
E_i\le\tfrac{8C_{\rm env}}{3}E_{3-i}\le K_0E_{3-i}.\qedhere
\]
\end{proof}

By Lemma~\ref{ce:source-owner}, at each point where
\eqref{ce:incoming-chain-target} must be checked there is a unique
$o\in\{1,2\}$ such that $(t,x)\in\mathcal U_o$.
Put $c=3-o$, the index of the other solution; thus
$(o,c)\in\{(1,2),(2,1)\}$. On $\mathcal U_o$, the solution $f_o$
produces the residual and is positive, while $f_c$ is negative:
\[
f_o\ge\tfrac58E_o,\qquad f_c\le-\tfrac38E_o.
\]
Let $\varepsilon_2\in(0,1)$ be a small constant to be chosen below.
The bound $E_o\le K_0E_c$ leaves two cases:
\[
E_o\le\varepsilon_2E_c,
\qquad\hbox{or}\qquad
\varepsilon_2E_c\le E_o\le K_0E_c.
\]
We treat them in Subsections~\ref{ce:small-owner}
and~\ref{ce:comparable}, respectively, and combine their bounds in
Subsection~\ref{ce:incoming-repair} to obtain
\eqref{ce:incoming-chain-target}.

If $E_o$ is much smaller, cancellation forces the oscillatory part of
$f_c$ near one of its zeros. If the envelopes are comparable, their
band indices and derivative scales are comparable instead.
The geometric data and the two solution families are fixed first.
The constants in Lemma~\ref{ce:small-alternatives} are then chosen in
its proof, and Lemma~\ref{ce:comparable-windows} uses the resulting
threshold $\varepsilon_2$. The integer $M$ is chosen last, large
enough for both lemmas. Each later increase of $M$ keeps these
constants fixed.

\medskip\noindent\textbf{Derivative scales and contractions.}
To estimate $\mathcal C_1,\mathcal C_2,\mathcal C_3$, we use the
following derivative bounds for the complete local solutions $f_i$,
established in Section~\ref{ce:front-section}. For
$a,b,h\in\mathbb N_0$, there is a constant $C=C(a,b,h)>0$ such that
\begin{equation}\label{ce:front-interface}
\begin{aligned}
&|\partial_{u_i}^a\partial_{v_i}^b\partial_{y_i}^h f_i|
\le C(a,b,h)E_i\lambda_i^a\tau_i^b k_i^h,\qquad
\lambda_i\triangleq2^{m_i},\quad k_i\triangleq y_i^{-2},\quad
\tau_i\triangleq k_i^2/\lambda_i,\\
&f_i\ge5E_i/8,\qquad f_{i,u_iu_i}\le-cE_i\lambda_i^2,\qquad
f_{i,u_iu_iu_iu_i}\ge cE_i\lambda_i^4.
\end{aligned}
\end{equation}
The last three inequalities in \eqref{ce:front-interface} hold near
$\operatorname{supp}W_i$. Thus the three differentiation directions
$u_i,v_i,y_i$ carry the scales $\lambda_i,\tau_i,k_i$, respectively;
in particular,
\[
\lambda_i\asymp\beta^{16m_i},\qquad
k_i\asymp\beta^{2m_i},\qquad
\tau_i\asymp\beta^{-12m_i},\qquad
\frac{k_i}{\lambda_i}\asymp\beta^{-14m_i}.
\]

In the chain identities \eqref{ce:chain-expansion}, write
$w\triangleq(\mathrm{d}q)^\sharp$, $H\triangleq\partial^2q$, and
$\mathfrak a\triangleq(Hw)^\sharp$.
In the edge coordinates $(u_i,v_i,y_i)$, vector subscripts denote
coordinate components; for example, $\zeta_u=\mathrm{d}u_i(\zeta)$ and
$\mathfrak a_{u_i}=\mathrm{d}u_i(\mathfrak a)$.
We write $\zeta^n$ for $n$ repetitions of $\zeta$ as tensor arguments.
When a frame is specified, $O(A,B,C)$ denotes componentwise bounds in
the order $u,v,y$. Short derivative subscripts on $f_i$ refer to its
own coordinates $(u_i,v_i,y_i)$.
For spacetime vectors $\zeta,\xi$ expressed in the $i$th edge frame
and every integer $n\ge1$, all contractions below are expanded by
\begin{equation}\label{ce:tensor-component-sum}
\partial^nf_i(\zeta_1,\ldots,\zeta_n)
=\sum_{\eta_1,\ldots,\eta_n\in\{u_i,v_i,y_i\}}
(\partial_{\eta_1}\cdots\partial_{\eta_n}f_i)
\prod_{j=1}^n\mathrm d\eta_j(\zeta_j).
\end{equation}
In particular, for repeated arguments this becomes
\begin{equation}\label{ce:tensor-repeated-sum}
\partial^nf_i(\zeta^n)
=\sum_{a+b+h=n}\frac{n!}{a!b!h!}
(\partial_{u_i}^a\partial_{v_i}^b\partial_{y_i}^h f_i)
\zeta_{u_i}^a\zeta_{v_i}^b\zeta_{y_i}^h.
\end{equation}
Applying \eqref{ce:front-interface} to the finite sums in
\eqref{ce:tensor-component-sum} and \eqref{ce:tensor-repeated-sum}
gives
\begin{equation}\label{ce:tensor-contractions}
\begin{cases}
|\partial^nf_i(\zeta^n)|
\le C(n)E_i(\lambda_i|\zeta_u|+\tau_i|\zeta_v|+k_i|\zeta_y|)^n,\\[3pt]
|\partial^nf_i(\xi,\zeta^{n-1})|
\le C(n)E_i(\lambda_i|\xi_u|+\tau_i|\xi_v|+k_i|\xi_y|)
(\lambda_i|\zeta_u|+\tau_i|\zeta_v|+k_i|\zeta_y|)^{n-1},
\end{cases}\;\, i\in\{1,2\}.
\end{equation}
Write $U_i\triangleq f_{i,u_i}$, $H_i\triangleq f_{i,v_i}$,
$Y_i\triangleq f_{i,y_i}$, and
$J_*\triangleq\nabla q\cdot e_*-c_\alpha q_t$.
The identities are
\begin{equation}\label{ce:lorentz-identities}
\begin{aligned}
&J_*=s_\alpha(Y_1-Y_2)-2c_\alpha(H_1+H_2),\\
&(D_qu_1)(D_qu_2)=J_*^2+(\mathfrak g/2)Q,\qquad
Q=(D_qu_i)(D_qv_i)-(D_qy_i)^2.
\end{aligned}
\end{equation}
To check the signs in the first identity, the two edge frames are
\[
\begin{aligned}
e_1&=c_\alpha e_*-s_\alpha n_*,\quad
e_2 =c_\alpha e_*+s_\alpha n_*,\\
n_1 &=s_\alpha e_*+c_\alpha n_*,\quad
n_2 =-s_\alpha e_*+c_\alpha n_*.
\end{aligned}
\]
Since $f_{i,t}=U_i+H_i$ and $f_{i,z_i}=U_i-H_i$, their separate
contributions are
\[
\begin{aligned}
\nabla f_1\cdot e_*-c_\alpha f_{1,t}
&=c_\alpha(f_{1,z_1}-f_{1,t})+s_\alpha Y_1
=-2c_\alpha H_1+s_\alpha Y_1,\\
\nabla f_2\cdot e_*-c_\alpha f_{2,t}
&=c_\alpha(f_{2,z_2}-f_{2,t})-s_\alpha Y_2
=-2c_\alpha H_2-s_\alpha Y_2.
\end{aligned}
\]
Adding these equalities gives the first line of
\eqref{ce:lorentz-identities}. For the product identity, put
\[
T_*\triangleq q_t,\qquad X_*\triangleq\nabla q\cdot e_*,\qquad
N_*\triangleq\nabla q\cdot n_*.
\]
Then
\[
\begin{aligned}
&D_qu_1=T_*-c_\alpha X_*+s_\alpha N_*,\qquad
D_qu_2=T_*-c_\alpha X_*-s_\alpha N_*,\\
&(D_qu_1)(D_qu_2)
=(T_*-c_\alpha X_*)^2-s_\alpha^2N_*^2 =(X_*-c_\alpha T_*)^2+s_\alpha^2(T_*^2-X_*^2-N_*^2)
=J_*^2+(\mathfrak g/2)Q.
\end{aligned}
\]
Finally, in either orthonormal edge frame,
\[
\begin{aligned}
(D_qu_i)(D_qv_i)-(D_qy_i)^2
&=(q_t-\nabla q\cdot e_i)(q_t+\nabla q\cdot e_i)
-(\nabla q\cdot n_i)^2 =q_t^2-|\nabla q|^2=Q.
\end{aligned}
\]

The remainder of this subsection is organized as follows. In
Subsection~\ref{ce:small-owner} we prove
\eqref{ce:incoming-chain-target} when one envelope is much smaller
than the other. In Subsection~\ref{ce:comparable} we treat the case
of comparable envelopes. In Subsection~\ref{ce:incoming-repair} we
combine the two cases and complete the proof of
Proposition~\ref{ce:incoming}.

\subsubsection{A much smaller amplitude envelope}\label{ce:small-owner}
We prove \eqref{ce:incoming-chain-target} on $\mathcal U_o$ when
the solution producing the residual has a much smaller envelope.
Fix $o\in\{1,2\}$, put $c\triangleq3-o$, and assume
$
E_o\le\varepsilon_2E_c$ and $(o,c)\in\{(1,2),(2,1)\}$, 
where the threshold $\varepsilon_2\in(0,1)$ will be chosen in
Lemma~\ref{ce:small-alternatives}.
Here $E_c$ is the amplitude scale of $f_c$, but the smaller
solution $f_o$ can still control its higher derivatives through the
factor $\lambda_o$. The second-derivative scales give three ranges, using, respectively, the contribution of $f_c$ to $\mathcal C_1$, the fourth
derivative of $f_o$ in $\mathcal C_3$, and the negative $uu$ derivatives of
both solutions.

\par\noindent\textbf{Hierarchy lemma.} In the range $E_o\le\varepsilon_2 E_c$, sufficiently large $M$ satisfies \eqref{ce:small-hierarchy}.\par
The band definition and the size of $U$ give, on $D$,
\[
\beta^{-(m_i+1)}<y_i\le R_+\beta^{-M}
\quad\Longrightarrow\quad m_i\ge M-C_*,
\qquad i\in\{1,2\},
\]
where $C_*>0$ depends only on $R_+$. Thus both band indices tend
uniformly to infinity as $M\to\infty$.
The envelope bounds \eqref{ce:envelope} imply
\begin{equation}\label{ce:small-hierarchy}
\begin{cases}
m_o\ge m_c-1,\\[3pt]
\frac{E_o}{E_c}\left(
\frac{\lambda_o}{\lambda_c}+\frac{k_o}{k_c}
+\frac{\tau_o}{\tau_c}\right)\le C\varepsilon_2+o(1).
\end{cases}
\end{equation}
Indeed, $m_o\le m_c-2$ and \eqref{ce:envelope} would give
\[
\log(E_o/E_c)\ge(1-\beta^{-1})\beta^{m_c}-O(m_c+m_o)
\longrightarrow+\infty,
\]
contrary to $\log(E_o/E_c)\le\log\varepsilon_2<0$.
For $m_o\ge m_c+2$,
\begin{equation}\label{ce:high-owner-index}
E_o/E_c\le
e^{-\beta^{m_o}+\beta^{m_c+1}+(m_c+1)^2-m_o^2+\gamma_c(m_c+1)-\gamma_o m_o}.
\end{equation}
For every fixed $p,q\ge0$, the right-hand side of
\eqref{ce:high-owner-index} satisfies
\[
\sup_{m_o\ge m_c+2}
\frac{E_o}{E_c}\,2^{p(m_o-m_c)}\beta^{q(m_o-m_c)}
\longrightarrow0\qquad(m_c\longrightarrow\infty).
\]
For sufficiently large $m_c$, the upper bound is decreasing in
$m_o\ge m_c+2$, so its maximum occurs at $m_o=m_c+2$.
In fact, writing $j=m_o$, the change of its logarithm from $j$ to
$j+1$, after multiplication by the displayed frequency and normal
factors, is
\[
-(\beta-1)\beta^j-2j-1-\gamma_o+p\log2+q\log\beta<0.
\]

\hfill$\square$

To compare derivative scales, define
$
\epsilon\triangleq\frac{k_c}{\lambda_c}$ and
$\varrho\triangleq\frac{E_o\lambda_o^2}{E_c\lambda_c^2}.$
The ratio $\epsilon$ compares the normal derivative scale with the fast
frequency and tends to zero as $M\to\infty$.
The ratio $\varrho$ compares the second derivative scales of the two solutions.
The envelope comparison gives $\varrho\le1$. For every fixed
$c_*>0$, the condition $\varrho\ge c_*\epsilon^4$ also forces
$m_o\le m_c+1$ when $M$ is sufficiently large.
\begin{lemma}[Nondegeneracy for a smaller amplitude envelope]\label{ce:small-alternatives}
There exist constants
\[
\eta,\varepsilon_2\in(0,1),\qquad
c_*,q_*,\theta_Q,c(\eta),c>0,\qquad
K>1,
\]
depending only on the fixed local data and on the two solution
families, and independent of $M$, with the following property. For
every sufficiently large integer $M$, the alternatives below hold at
every point of the set
\[
\mathcal U_o\cap\{E_o\le\varepsilon_2E_c\}.
\]

\smallskip
\noindent\emph{Large time derivative.}
If $|q_t|\ge\eta E_c\lambda_c$, then
\begin{equation}\label{ce:small-large-gradient}
|Q|\ge c(\eta) E_c^2k_c^2
\quad\text{or}\quad
\mathcal C_1\le-c(\eta) E_oE_c^2\lambda_c^2\lambda_o^2.
\end{equation}

\smallskip
\noindent\emph{Small time derivative.}
If $|q_t|\le\eta E_c\lambda_c$, then there exists $\nu>0$ with
$\nu\asymp m_c^2\beta^{-m_c}$ such that
\[
|J_*|\ge cE_ck_c\nu.
\]
In this branch, the conclusion depends on the size of $\varrho$:
\begin{align}
\varrho\le c_*\epsilon^4 &:\quad
|Q|\ge q_*E_c^2k_c^2\nu^2
\quad\text{or}\quad
\mathcal C_1\ge cE_c^3k_c^4\nu,
\label{ce:small-effective-low}\\
c_*\epsilon^4\le \varrho\le K\epsilon^4 &:\quad
|Q|\ge q_*E_c^2k_c^2\nu^2
\quad\text{or}\quad
\mathcal C_3\ge cE_c^5k_c^4\lambda_c^4\nu^2,
\label{ce:small-effective-cross}\\
\varrho\ge K\epsilon^4 &:\quad
|Q|\ge\theta_QJ_*^2
\quad\text{or}\quad
\mathcal C_1\le-cE_ck_c^2J_*^2.
\label{ce:small-effective-high}
\end{align}
All alternatives hold for either labeling of the two fronts as
$(o,c)$ and for either sign of $\alpha_v$.
\end{lemma}
\begin{proof}[Proof of Lemma~\ref{ce:small-alternatives}]
The proof has five steps. Step~\ref{step:small-envelope-1} proves
\eqref{ce:small-large-gradient}. For the remaining steps we assume
$|q_t|\le\eta E_c\lambda_c$: Step~\ref{step:small-envelope-2} locates
the phase near $\pi$ and proves the lower bound for $J_*$.
Step~\ref{step:small-envelope-3} identifies the finer scale
$\xi^2\asymp\nu$ and proves \eqref{ce:small-effective-low}.
Steps~\ref{step:small-envelope-4} and~\ref{step:small-envelope-5}
prove \eqref{ce:small-effective-cross} and
\eqref{ce:small-effective-high}.

\medskip
\noindent\textbf{Step 1. The model derivatives and the large time
derivative.}\proofsteplabel{step:small-envelope-1}{1}

For \eqref{ce:small-large-gradient}, the quadratic form of $f_c$ is
negative. Small $Q$ requires compensation from the cross term with
$f_o$, forcing a large $D_qu_o$ so that the negative $uu$ derivative
of $f_o$ dominates $\mathcal C_1$. We first compute the leading
characteristic derivatives of $f_c$ and verify that its quadratic form
is negative whenever $|\sin\theta|$ is bounded below.

Set ${A_{c,0}}\triangleq A_{\gamma_c,m_c}(y_c)$ and
${A_{c,1}}\triangleq A_{\gamma_c,m_c+1}(y_c)$, so that
${A_{c,0}},{A_{c,1}}\ge0$ and ${A_{c,0}}+{A_{c,1}}=E_c$, and put
$\theta\triangleq\lambda_c u_c$. The leading characteristic
derivatives of $f_c$ are
\[
\begin{cases}
\ds	U^0 \triangleq-\lambda_c({A_{c,0}}\sin\theta+2{A_{c,1}}\sin2\theta),\\
\ns\ds	H^0 \triangleq\frac{k_c^2}{4\lambda_c}
\bigl({A_{c,0}}\sin\theta+({A_{c,1}}/2)\sin2\theta\bigr),\\
\ns\ds	Y^0 \triangleq k_c({A_{c,0}}\cos\theta+{A_{c,1}}\cos2\theta).
\end{cases}
\]
To verify the quadratic form, put $p=\cos\theta$. Then
\[
\begin{aligned}
\frac{4U^0H^0-(Y^0)^2}{k_c^2}
&=-(1-p^2)({A_{c,0}}+4{A_{c,1}}p)({A_{c,0}}+{A_{c,1}}p)
-[{A_{c,0}}p+{A_{c,1}}(2p^2-1)]^2\\
&=-{A_{c,0}}^2-{A_{c,1}}^2+{A_{c,0}}{A_{c,1}}(p^3-3p)\\
&=-({A_{c,0}}-{A_{c,1}})^2-{A_{c,0}}{A_{c,1}}(1+p)^2(2-p).
\end{aligned}
\]
Thus the coefficients of ${A_{c,0}}^2,{A_{c,1}}^2,{A_{c,0}}{A_{c,1}}$
are $-1,-1,p^3-3p$, and we obtain the exact identity
\begin{equation}\label{ce:negative-single-norm}
4U^0H^0-(Y^0)^2
=-k_c^2\bigl\{({A_{c,0}}-{A_{c,1}})^2
+{A_{c,0}}{A_{c,1}}(1+\cos\theta)^2(2-\cos\theta)\bigr\}.
\end{equation}
The errors in the complete solution are
$U_c-U^0=O(E_c\lambda_c\tau_c)$,
$H_c-H^0=O(E_c\tau_c(m_c^4\beta^{-m_c}+\tau_c))$, and
$Y_c-Y^0=O(E_ck_c(m_c^4\beta^{-m_c}+\tau_c))$.

We now use the hypothesis $|q_t|\ge\eta E_c\lambda_c$.
Equation~\eqref{ce:small-hierarchy} gives
$|U_c|\ge\eta E_c\lambda_c/2$, hence $|\sin\theta|\ge\eta/20$.
Equation~\eqref{ce:negative-single-norm} therefore gives
$\mathcal B(\mathrm{d}f_c,\mathrm{d}f_c)\le-c(\eta)E_c^2k_c^2$.
Write $\mathcal Q_i\triangleq\mathcal B(\mathrm{d}f_i,\mathrm{d}f_i)$
in this calculation. For some $a_0=a_0(\eta)>0$, decreasing
$\varepsilon_2$ and increasing $M$ give
\[
\mathcal Q_c\le-a_0E_c^2k_c^2,\qquad
|\mathcal Q_o|\le CE_o^2k_o^2\le\tfrac14a_0E_c^2k_c^2.
\]
If $|Q|<a_0E_c^2k_c^2/4$, then
\[
\tfrac12a_0E_c^2k_c^2
\le|Q-\mathcal Q_c-\mathcal Q_o|
=2|\mathcal B(\mathrm{d}f_o,\mathrm{d}f_c)|
\le CE_oE_c\lambda_o\lambda_c.
\]
Consequently,
\begin{equation}\label{ce:cross-required}
(E_o/E_c)\lambda_c\lambda_o/k_c^2\ge c(\eta).
\end{equation}
The corresponding vector component is
$D_qu_o=\mathfrak g U_c+(1+\cos\alpha_v)H_c\pm\sin\alpha_vY_c+2H_o$;
thus $|D_qu_o|\ge c(\eta)E_c\lambda_c$.
Its negative source Hessian contributes
$-cE_oE_c^2\lambda_c^2\lambda_o^2$. The mixed derivatives of $f_o$
have the additional factor $O(k_o/\lambda_o)$.
The errors from the Hessian of $f_c$, with zero, one, or two factors
from $\mathrm{d}f_o$, relative to this scale, are bounded by
$Ck_c^2/(\lambda_c\lambda_o)$, $Ck_c^2/(\lambda_c\lambda_o)$, and
$CE_o/E_c$, using \eqref{ce:cross-required}.
More explicitly, for $S=E_oE_c^2\lambda_c^2\lambda_o^2$,
\[
\begin{aligned}
\mathcal C_1&\le S\left[-c(\eta)+C\left(
\frac{k_o}{\lambda_o}+\frac{k_c^2}{\lambda_c\lambda_o}
+\frac{E_o}{E_c}\right)\right],\qquad\frac{k_c^2}{\lambda_c\lambda_o}
\le C(\eta)\frac{E_o}{E_c}\le C(\eta)\varepsilon_2.
\end{aligned}
\]
Choose $\varepsilon_2$ small relative to $c(\eta)$, and then $M$
large so that $k_o/\lambda_o$ is small. The bracket is at most
$-c(\eta)/2$, proving \eqref{ce:small-large-gradient}.

\medskip
\noindent\textbf{Step 2. The phase and normal scale for small time
derivative.}\proofsteplabel{step:small-envelope-2}{2}

From now on we assume $|q_t|\le\eta E_c\lambda_c$. We first show that
the two leading cosine terms of $f_c$ cancel, with
${A_{c,0}}\asymp{A_{c,1}}$ and $\theta$ near $\pi$. The difference
between the normal logarithmic rates of ${A_{c,0}}$ and ${A_{c,1}}$
then produces a negative normal derivative of size $E_ck_c\nu$, which
bounds $J_*$ away from zero.

The two small quantities are the leading value of $f_c$ and its fast
derivative:
\begin{equation}\label{ce:small-phase-identities}
\begin{cases}
\ds	{A_{c,0}}\cos\theta+{A_{c,1}}\cos2\theta
=q-f_o+O(E_c\tau_c),\\
\ns\ds	-\lambda_c\sin\theta({A_{c,0}}+4{A_{c,1}}\cos\theta)
=q_t-U_o-H_o-H_c+O(E_c\lambda_c\tau_c).
\end{cases}
\end{equation}
The nodal tube from Lemma~\ref{ce:source-owner} bounds $|q|$ by
$E_o/4$. Together with \eqref{ce:small-hierarchy}, equalities
\eqref{ce:small-phase-identities} give
\begin{equation}\label{ce:small-phase-errors}
\begin{cases}
\ds |{A_{c,0}}\cos\theta+{A_{c,1}}\cos2\theta| \le C(E_o+E_c\tau_c),\\
\ns\ds |\sin\theta({A_{c,0}}+4{A_{c,1}}\cos\theta)| \le C(\eta+\varepsilon_2+o(1))E_c.
\end{cases}
\end{equation}
Taking the errors in \eqref{ce:small-phase-errors} below $E_c/100$,
the identity
\[
{A_{c,0}}\cos\theta+{A_{c,1}}\cos2\theta
-\cos\theta({A_{c,0}}+4{A_{c,1}}\cos\theta)
=-{A_{c,1}}(1+2\cos^2\theta)
\]
implies
\begin{equation}\label{ce:balanced-phase}
|\sin\theta|<1/4,\qquad\cos\theta<-\sqrt{15}/4,\qquad
\cos2\theta>7/8,\qquad1/2<{A_{c,0}}/{A_{c,1}}<2.
\end{equation}
For an explicit elimination, write
$F={A_{c,0}}\cos\theta+{A_{c,1}}\cos2\theta$ and
$V={A_{c,0}}\sin\theta+2{A_{c,1}}\sin2\theta$.
If $|\sin\theta|\ge1/4$, then
\[
{A_{c,1}}\le|F|+4|V|,\qquad
{A_{c,0}}\le4|V|+4{A_{c,1}}\le4|F|+20|V|,
\]
so that
\[
E_c\le5|{A_{c,0}}\cos\theta+{A_{c,1}}\cos2\theta|
+24|{A_{c,0}}\sin\theta+2{A_{c,1}}\sin2\theta|<E_c,
\]
a contradiction. If instead $\cos\theta>0$, then
$F\ge(7/8)({A_{c,0}}+{A_{c,1}})$, again contrary to
\eqref{ce:small-phase-errors}. Finally, with $\delta_0=1/100$, the
bound for $F$ in \eqref{ce:small-phase-errors} gives
\[
\frac{7/8-\delta_0}{1+\delta_0}
\le\frac {A_{c,0}}{A_{c,1}}\le
\frac{1+\delta_0}{\sqrt{15}/4-\delta_0},
\]
which proves the last assertion of \eqref{ce:balanced-phase}.

By \eqref{ce:balanced-phase}, choose the representative $\xi$ of
$\theta-\pi$ modulo $2\pi$ with $|\xi|<\arcsin(1/4)$.
Denote the smooth step function $b(\cdot)$ in
\eqref{ce:cutoff-functions} by $\chi_{\rm cut}$, and put
\[
x\triangleq\frac{\beta^{m_c+1}y_c-1}{\beta-1},\qquad
\nu\triangleq\frac{\partial_{y_c}x}{k_c}
(1-\chi_{\rm cut}(x))(x^{-2}+(1-x)^{-2}).
\]
The two active cutoff factors are $\chi_{\rm cut}(x)$ and
$1-\chi_{\rm cut}(x)$, so
\begin{equation}\label{ce:small-cutoff-log-ratio}
\begin{aligned}
\log\frac{\chi_{\rm cut}(x)}{1-\chi_{\rm cut}(x)}
&=\log({A_{c,0}}/{A_{c,1}})-(2m_c+1+\gamma_c) =-x^{-1}+(1-x)^{-1},\\
x^{-1}-(1-x)^{-1}
&=2m_c+1+\gamma_c-\log({A_{c,0}}/{A_{c,1}})
=2m_c+O(1).
\end{aligned}
\end{equation}
The right side of \eqref{ce:small-cutoff-log-ratio} is positive, hence
$x<1/2$, hence $x^{-1}=2m_c+O(1)$ and $1/(3m_c)<x<1/m_c$ for large
$M$. Also $\chi_{\rm cut}(x)=O(e^{-2m_c})$ and
\[
\partial_{y_c}x=\frac{\beta^{m_c+1}}{\beta-1},\qquad
k_c\asymp\beta^{2m_c},\qquad
x^{-2}+(1-x)^{-2}\asymp m_c^2,
\]
which gives $\nu\asymp m_c^2\beta^{-m_c}$ directly from its
definition.
The logarithmic normal rates are exactly $k_c(1+\nu)$ and
$k_c[1-\nu\chi_{\rm cut}(x)/(1-\chi_{\rm cut}(x))]$. Their
derivatives satisfy
$|\partial_{y_c}^h(k_c\nu)|\le C(h)(k_c\nu)^{h+1}/m_c^h$ for
$1\le h\le3$, and $|k_c^{(h)}|/k_c^{h+1}\le C(h)y_c^h$.
By \eqref{ce:front-jets},
\[
\begin{aligned}
{A_{c,0}}'&=k_c(1+\nu){A_{c,0}},\qquad
{A_{c,1}}'=k_c\left(1-\nu
\frac{\chi_{\rm cut}(x)}{1-\chi_{\rm cut}(x)}\right){A_{c,1}},\\
Y_c&={A_{c,0}}'\cos\theta+{A_{c,1}}'\cos2\theta+O(E_ck_c\tau_c)\\
&=k_cf_c+k_c\nu\left[
{A_{c,0}}\cos\theta
-\frac{\chi_{\rm cut}(x)}{1-\chi_{\rm cut}(x)}
{A_{c,1}}\cos2\theta
\right]+O(E_ck_c\tau_c).
\end{aligned}
\]
On $\mathcal U_o$, the source margin and nodal tube imply
$f_o\ge5E_o/8$ and $f_c=q-f_o\le-3E_o/8$.
Also, \eqref{ce:balanced-phase} ensures that
\[
{A_{c,0}}\ge E_c/3,\qquad
{A_{c,0}}\cos\theta
-\frac{\chi_{\rm cut}(x)}{1-\chi_{\rm cut}(x)}
{A_{c,1}}\cos2\theta
\le-\frac{\sqrt{15}}{12}E_c.
\]
Consequently,
\[
\begin{aligned}
Y_c&\le-\frac38E_ok_c
-\frac{\sqrt{15}}{12}E_ck_c\nu+CE_ck_c\tau_c
\le-\frac{\sqrt{15}}{24}E_ck_c\nu,\\
\frac{\tau_c}{\nu}&\le Cm_c^{-2}\beta^{-11m_c}\longrightarrow0.
\end{aligned}
\]
Similarly, \eqref{ce:common-null} gives
\[
\begin{aligned}
Y_o&=k_of_o+
O\bigl(E_ok_o(m_o^4\beta^{-m_o}+\tau_o)\bigr) \ge\left[\frac58-C(m_o^4\beta^{-m_o}+\tau_o)\right]E_ok_o
\ge\frac12E_ok_o.
\end{aligned}
\]
By \eqref{ce:front-interface} and \eqref{ce:small-hierarchy},
\[
\frac{|H_1+H_2|}{E_ck_c\nu}
\le C\frac{E_c\tau_c+E_o\tau_o}{E_ck_c\nu}
\le C\frac{k_c}{\lambda_c\nu}
\le Cm_c^{-2}\beta^{-13m_c}\longrightarrow0.
\]
Thus \eqref{ce:lorentz-identities} implies
\begin{equation}\label{ce:small-Jstar-bound}
\begin{aligned}
|J_*|&\ge|s_\alpha|\,|Y_1-Y_2|-2c_\alpha|H_1+H_2| =|s_\alpha|(Y_o-Y_c)-2c_\alpha|H_1+H_2|\\
&\ge\left(\frac{\sqrt{15}}{24}|s_\alpha|-o(1)\right)E_ck_c\nu
\ge cE_ck_c\nu.
\end{aligned}
\end{equation}
Bound \eqref{ce:small-Jstar-bound} is independent of the sign of
$s_\alpha$.

\medskip
\noindent\textbf{Step 3. The small second-derivative ratio.}
\proofsteplabel{step:small-envelope-3}{3}

The positive $\mathcal C_1$ bound in \eqref{ce:small-effective-low},
under its small-$Q$ alternative, comes from the leading term
$(27/8){A_{c,0}}^3k_c^4\xi^2$. The quadratic form gives
$\xi^2\asymp\nu$, which controls the errors throughout
$\varrho\le K\epsilon^4$ and also applies in
Step~\ref{step:small-envelope-4}.

Assume $\varrho\le K\epsilon^4$ and
$|Q|\le q_*E_c^2k_c^2\nu^2$. For fixed $K$, set
$\overline {A_{c,1}}\triangleq-{A_{c,0}}\cos\theta/\cos2\theta$.
The source nodal tube and \eqref{ce:front-jets} give
\[
\frac{|{A_{c,1}}-\overline {A_{c,1}}|}{{A_{c,0}}}
\le C(E_o/E_c+\tau_c),
\]
since ${A_{c,0}}\asymp E_c$. We use $\overline {A_{c,1}}$ only in the
model, leaving the actual amplitudes unchanged. Freeze their
logarithmic normal rates to $k_c(1+\nu)$ and $k_c$, respectively.
Before substituting ${A_{c,1}}=\overline {A_{c,1}}$, the frozen
quadratic form, divided by $k_c^2$, is
\[
\begin{aligned}
&-({A_{c,0}}\sin\theta+2{A_{c,1}}\sin2\theta)
\bigl((1+\nu)^2{A_{c,0}}\sin\theta+({A_{c,1}}/2)\sin2\theta\bigr) -\bigl((1+\nu){A_{c,0}}\cos\theta+{A_{c,1}}\cos2\theta\bigr)^2\\
&=-(1+\nu)^2{A_{c,0}}^2-{A_{c,1}}^2
-2(1+\nu){A_{c,0}}{A_{c,1}}\cos\theta\cos2\theta -\bigl(\tfrac12+2(1+\nu)^2\bigr){A_{c,0}}{A_{c,1}}
\sin\theta\sin2\theta.
\end{aligned}
\]
Substituting ${A_{c,1}}=-{A_{c,0}}\cos\theta/\cos2\theta$ and dividing
by ${A_{c,0}}^2$ gives the normalized model
\begin{equation}\label{ce:small-norm-model}
\overline Q\triangleq-(1+\nu)^2-\frac{\cos^2\theta}{\cos^2(2\theta)}
+2(1+\nu)\cos^2\theta
+\left(\tfrac12+2(1+\nu)^2\right)
\frac{\cos\theta\sin\theta\sin2\theta}{\cos2\theta}.
\end{equation}
Equivalently,
\[
\begin{aligned}
\overline Q
&=-\frac{\sin^4\theta(2\cos^2\theta+1)}{\cos^2(2\theta)}
+2\nu\frac{\sin^2\theta(2\cos^2\theta+1)}{\cos(2\theta)} +\nu^2\left[-1
+\frac{4\cos^2\theta\sin^2\theta}{\cos(2\theta)}\right].
\end{aligned}
\]
To compare this model with the actual quadratic form, put
$h\triangleq k_c\nu$ and
$\sigma\triangleq\chi_{\rm cut}(x)/(1-\chi_{\rm cut}(x))$; primes
below denote $y_c$ derivatives. The exact amplitude rates satisfy
\begin{equation}\label{ce:small-amplitude-rates}
\begin{aligned}
&\frac{{A_{c,0}}'}{A_{c,0}}=k_c+h,\qquad
\frac{{A_{c,1}}'}{A_{c,1}}=k_c-h\sigma,\qquad\sigma=e^{-(2m_c+1+\gamma_c)}\frac {A_{c,0}}{A_{c,1}}
=O(e^{-2m_c}).
\end{aligned}
\end{equation}
For $L(x)\triangleq x^{-2}+(1-x)^{-2}$,
\[
\frac{h'}{h^2}
=\frac{\partial_xL}{(1-\chi_{\rm cut}(x))L^2}-\sigma=O(m_c^{-1}),
\qquad \sigma'=h\sigma(1+\sigma).
\]
Here we used $x\asymp m_c^{-1}$. Consequently,
\begin{equation}\label{ce:small-amplitude-second-rates}
\begin{aligned}
&\frac{{A_{c,0}}''}{A_{c,0}}=k_c^2(1+\nu)^2+k_c'+h',\qquad
\frac{{A_{c,1}}''}{A_{c,1}}=k_c^2+k_c'+O(k_c^2e^{-m_c}),\qquad
\frac{h'}{k_c^2}=O(\nu^2/m_c).
\end{aligned}
\end{equation}
The refined derivatives in \eqref{ce:front-extraction} and
\eqref{ce:normal-extraction}, without replacing amplitude derivatives
by their coarse bounds, yield
\[
\begin{aligned}
\mathcal B(\mathrm{d}f_c,\mathrm{d}f_c)
=&-({A_{c,0}}\sin\theta+2{A_{c,1}}\sin2\theta)
({A_{c,0}}''\sin\theta\\ &+\tfrac12{A_{c,1}}''\sin2\theta) -({A_{c,0}}'\cos\theta+{A_{c,1}}'\cos2\theta)^2
+O(E_c^2k_c^2\tau_c).
\end{aligned}
\]
The exact rates in
\eqref{ce:small-amplitude-rates}--\eqref{ce:small-amplitude-second-rates}
and $|{A_{c,1}}-\overline {A_{c,1}}|/{A_{c,0}}=O(E_o/E_c+\tau_c)$
leave, besides $\overline Q$, the common $k_c'$ term with factor
\[
\begin{aligned}
&({A_{c,0}}\sin\theta+2\overline {A_{c,1}}\sin2\theta)
({A_{c,0}}\sin\theta+(\overline {A_{c,1}}/2)\sin2\theta) =\frac{{A_{c,0}}^2\sin^4\theta(2\cos^2\theta+1)}
{\cos^2(2\theta)}
\end{aligned}
\]
by the definition of $\overline {A_{c,1}}$. Since $k_c'/k_c^2=-2y_c$,
this contributes $O(y_c\xi^4)$. The $h'$ term is
\[
-\frac{h'}{k_c^2}
\bigl(\sin\theta+2(\overline {A_{c,1}}/{A_{c,0}})\sin2\theta\bigr)
\sin\theta
=O\bigl((\nu^2/m_c)\xi^2\bigr).
\]
The errors from the second amplitude rate are $O(e^{-m_c})$.
Finally, $m_o\ge m_c-1$ and $\varrho\le K\epsilon^4$ imply
$E_o/E_c\le4K\epsilon^4$. The derivative bounds
\eqref{ce:front-interface} give
\[
\frac{|Q-\mathcal B(\mathrm{d}f_c,\mathrm{d}f_c)|}
{{A_{c,0}}^2k_c^2}
\le C\left[
\frac{E_o\lambda_o\lambda_c}{E_ck_c^2}
+\frac{E_o^2\lambda_o^2}{E_c^2k_c^2}\right]
\le C(K\epsilon^2+K^2\epsilon^6).
\]
The resulting comparison, with $C(K)>0$ independent of $M$, is
\begin{equation}\label{ce:small-norm-budget}
\begin{aligned}
\left|\frac{Q}{{A_{c,0}}^2k_c^2}-\overline Q\right|
\le C(K)\bigl[&y_c\xi^4+(\nu^2/m_c)\xi^2+\tau_c +e^{-m_c}+E_o/E_c+K\epsilon^2+K^2\epsilon^6\bigr].
\end{aligned}
\end{equation}
For fixed $K$, the last five terms are $o(\nu^2)$. In fact,
$\lambda_c=\beta^{16m_c}$, so
\[
\begin{aligned}
&\tau_c\asymp\beta^{-12m_c},\qquad
\epsilon\asymp\beta^{-14m_c},\qquad
\nu^2\asymp m_c^4\beta^{-2m_c},\\
&\frac{\tau_c+e^{-m_c}+E_o/E_c+K\epsilon^2+K^2\epsilon^6}{\nu^2}\\
&\le C(K)m_c^{-4}\bigl(
\beta^{-10m_c}+e^{-m_c}\beta^{2m_c}
+\beta^{-26m_c}+\beta^{-54m_c}+\beta^{-82m_c}\bigr)
\longrightarrow0.
\end{aligned}
\]
The first two terms still depend on the unknown phase $\xi$.
The exact model \eqref{ce:small-norm-model} and
\eqref{ce:small-norm-budget} imply, with constants independent of $M$,
\[
\begin{aligned}
\xi^4&\le C(K)(\nu\xi^2+\nu^2)+C(K)y_c\xi^4,\\
C(K)y_c\le\tfrac12
&\quad\Longrightarrow\quad
\xi^4\le2C(K)(\nu\xi^2+\nu^2),\\
\left(\frac{\xi^2}{\nu}\right)^2
&\le2C(K)\left(\frac{\xi^2}{\nu}+1\right)
\quad\Longrightarrow\quad \xi^2\le C(K)\nu.
\end{aligned}
\]
All errors in \eqref{ce:small-norm-budget} are then $o(\nu^2)$.
In this range the Taylor expansion is
\[
\overline Q=-\nu^2+6\nu\xi^2-3\xi^4
+O(\nu^2\xi^2+\nu\xi^4+\xi^6).
\]
Choose $q_*$ sufficiently small, using ${A_{c,0}}\asymp E_c$. Then
\[
\left|-1+6\frac{\xi^2}{\nu}
-3\left(\frac{\xi^2}{\nu}\right)^2\right|
\le Cq_*+o(1)<\tfrac12.
\]
Writing $T=\xi^2/\nu\ge0$, the inequality
$-1+6T-3T^2>-1/2$ gives the explicit bounds
$
1-\frac{\sqrt{30}}6<T<1+\frac{\sqrt{30}}6.
$
In particular $\xi^2/\nu$ is bounded away from zero, and
\begin{equation}\label{ce:small-phase-scale}
c\nu\le\xi^2\le C\nu.
\end{equation}

To compute the contribution of $f_c$ to $\mathcal C_1$ along its own
gradient, include the phase contribution from the coefficient with
$\ell=1$ and fix $k_c,\nu$ at the point under consideration:
\[
\Theta\triangleq\lambda_cu_c+
\frac{z_ck_c^2(1+\nu)^2}{2\lambda_c},\qquad
\Phi\triangleq2\lambda_cu_c+\frac{z_ck_c^2}{4\lambda_c}.
\]
In the following model calculation, regard the amplitudes and these
two phase arguments as independent variables. Define
\begin{equation}\label{ce:small-model-flow}
\begin{aligned}
&Q_0\triangleq-(1+\nu)^2{A_{c,0}}^2-{A_{c,1}}^2
-2(1+\nu){A_{c,0}}{A_{c,1}}\cos\Theta\cos\Phi\\
&\hspace{1.6cm}-\left[\tfrac12+2(1+\nu)^2\right]
{A_{c,0}}{A_{c,1}}\sin\Theta\sin\Phi,\\
&D{A_{c,0}}\triangleq-(1+\nu){A_{c,0}}
[(1+\nu){A_{c,0}}\cos\Theta+{A_{c,1}}\cos\Phi],\\
&D{A_{c,1}}\triangleq-{A_{c,1}}
[(1+\nu){A_{c,0}}\cos\Theta+{A_{c,1}}\cos\Phi],\\
&D\Theta\triangleq(1+\nu)^2{A_{c,0}}\sin\Theta
+\left(\tfrac14+(1+\nu)^2\right){A_{c,1}}\sin\Phi,\\
&D\Phi\triangleq\left(\tfrac14+(1+\nu)^2\right){A_{c,0}}\sin\Theta
+{A_{c,1}}\sin\Phi.
\end{aligned}
\end{equation}
The last four formulas in \eqref{ce:small-model-flow} give the
leading Lorentz gradient flow in
$({A_{c,0}},{A_{c,1}},\Theta,\Phi)$ after division by $k_c^2$, with
$\nu$ held fixed. They are the components of the first-order operator
\[
\begin{aligned}
&D=(D{A_{c,0}})\partial_{A_{c,0}}
+(D{A_{c,1}})\partial_{A_{c,1}}
+(D\Theta)\partial_\Theta+(D\Phi)\partial_\Phi,\\
&DQ_0=\sum_{z\in\{{A_{c,0}},{A_{c,1}},\Theta,\Phi\}}Q_{0,z}Dz.
\end{aligned}
\]
By the algebraic lemma in Appendix~\ref{app:intrinsic-details}, the
two phase contributions have quadratic coefficients $9/2$ and $-9/8$,
whose sum is $27/8$; at $\xi=0$ the amplitude terms start at
$O(\nu^3)$. Thus
\begin{equation}\label{ce:intrinsic-C1}
\begin{aligned}
&\frac1{{A_{c,0}}^3}\left.
\sum_{z\in\{{A_{c,0}},{A_{c,1}},\Theta,\Phi\}} Q_{0,z}Dz
\right|_{\substack{\Theta=\theta,\ \Phi=2\theta\\
{A_{c,1}}=-{A_{c,0}}\cos\theta/\cos2\theta}} =\frac{27}{8}\xi^2+O(\nu\xi^2+\xi^4+\nu^3).
\end{aligned}
\end{equation}
The connection with the complete Hessian is
\[
\begin{aligned}
&\Biggl|2\partial^2f_c
\bigl((\mathrm{d}f_c)^\sharp,(\mathrm{d}f_c)^\sharp\bigr)
-k_c^4\left.
\sum_{z\in\{{A_{c,0}},{A_{c,1}},\Theta,\Phi\}}Q_{0,z}Dz
\right|_{\substack{\Theta=\theta,\ \Phi=2\theta\\
{A_{c,1}}=-{A_{c,0}}\cos\theta/\cos2\theta}}\Biggr|\\
&\le CE_c^3k_c^4
\left(y_c+\nu^2/m_c+\tau_c+e^{-m_c}\right)
=o(E_c^3k_c^4\nu).
\end{aligned}
\]
Here the derivative terms omitted when freezing the rates contribute
$O(y_c+\nu^2/m_c)$ at the normalized scale. The phase shifts are
$O(\tau_c)$, and replacing ${A_{c,1}}$ by $\overline {A_{c,1}}$ costs
$O(E_o/E_c+\tau_c)=O(\tau_c)$, since $E_o/E_c\le4K\epsilon^4$.
The complete Hessian includes the $vv$ term
$-{A_{c,0}}^{(4)}\cos\theta/(16\lambda_c^2)$ from the coefficient with
$\ell=2$.
The vector $(\mathrm{d}f_c)^\sharp$ has components
$O(E_c\tau_c,E_c\lambda_c,E_ck_c)$ in $(u_c,v_c,y_c)$, so all six
contractions have scale $E_c^3k_c^4$.
For $\varrho\le c_*\epsilon^4$, the entire Hessian contribution of
$f_o$ is at most $Cc_*E_c^3k_c^4\nu$, and the errors from contracting
$\partial^2 f_c$ with $\mathrm{d}f_o$ are
$CE_c^3k_c^4(c_*\epsilon^2+c_*^2\epsilon^4)$.
By \eqref{ce:small-phase-scale}, the leading term is at least
$c_0E_c^3k_c^4\nu$ for some $c_0>0$. Dividing all errors by this
common scale gives
\[
\frac{\mathcal C_1}{E_c^3k_c^4\nu}
\ge c_0-Cc_*-o(1)
-C\left(c_*\frac{\epsilon^2}{\nu}
+c_*^2\frac{\epsilon^4}{\nu}\right).
\]
Choose $c_*$ so that $Cc_*\le c_0/4$. Since $\epsilon^2/\nu\to0$,
increasing $M$ makes the remaining errors at most $c_0/4$. This
proves \eqref{ce:small-effective-low}.

\medskip
\noindent\textbf{Step 4. The crossover range.}
\proofsteplabel{step:small-envelope-4}{4}

The $\mathcal C_3$ lower bound in \eqref{ce:small-effective-cross}
comes from the positive fourth derivative $f_{o,uuuu}$ on the source
support. Its contraction with four copies of $D_qu_o$ dominates the
other terms in \eqref{ce:chain-expansion}.

For $c_*\epsilon^4\le \varrho\le K\epsilon^4$, one has
$m_o=m_c+O(1)$, $|U_c|\asymp E_c\lambda_c\sqrt\nu$, and
$|D_qu_o|\ge cE_c\lambda_c\sqrt\nu$.
The positive fourth derivative of $f_o$ contributes at scale
\[
E_o\lambda_o^4E_c^4\lambda_c^4\nu^2
\asymp E_c^5k_c^4\lambda_c^4\nu^2.
\]
Indeed, \eqref{ce:small-phase-scale} gives
$|w|\le CE_c\lambda_c\sqrt\nu$ and
$|D_qu_o|\ge cE_c\lambda_c\sqrt\nu$; hence a $y_o$ or $v_o$
component in place of a $u_o$ component gives, respectively, a
factor $O(k_o/\lambda_o)$ or $O(\tau_o/\lambda_o)$.
In the coordinates of $f_c$,
\begin{equation}\label{ce:small-cross-components}
\begin{aligned}
&w=O(E_c\tau_c,E_c\lambda_c,E_ck_c),\qquad
\mathfrak a=O(E_c^2k_c^4/\lambda_c,E_c^2\lambda_c k_c^2,E_c^2k_c^3).
\end{aligned}
\end{equation}
For a vector $X$, write
$\mathcal L_i(X)=\lambda_i|X_{u_i}|+\tau_i|X_{v_i}|+k_i|X_{y_i}|$.
The component bounds \eqref{ce:small-cross-components} give
\[
\begin{aligned}
\mathcal L_c(w)&\le CE_ck_c^2,\qquad
\mathcal L_c(\mathfrak a)\le CE_c^2k_c^4,\\
\mathcal L_o(w)&\le CE_c\lambda_c\lambda_o\sqrt\nu,\qquad
\mathcal L_o(\mathfrak a)\le CE_c^2\lambda_c\lambda_ok_c^2.
\end{aligned}
\]
Using the finite sums in \eqref{ce:tensor-component-sum} gives
\[
\begin{aligned}
|\partial^4 f_c(w^4)|
&\le CE_c\mathcal L_c(w)^4\le CE_c^5k_c^8,\\
|\partial^3 f_c(\mathfrak a,w,w)|
&\le CE_c\mathcal L_c(\mathfrak a)\mathcal L_c(w)^2
\le CE_c^5k_c^8,\\
|\partial^2 f_c(\mathfrak a,\mathfrak a)|
&\le CE_c\mathcal L_c(\mathfrak a)^2\le CE_c^5k_c^8.
\end{aligned}
\]
Normalize by
$S_3=E_o\lambda_o^4E_c^4\lambda_c^4\nu^2
\asymp E_c^5k_c^4\lambda_c^4\nu^2$. The remaining estimates are
\[
\begin{aligned}
\frac{|\partial^3 f_o(\mathfrak a,w,w)|}{S_3}
&\le C\frac{k_c^2}{\lambda_c\lambda_o\nu}
\le C\frac{\epsilon^2}{\nu},\quad
\frac{|\partial^2 f_o(\mathfrak a,\mathfrak a)|}{S_3}
\le C\frac{k_c^4}{\lambda_c^2\lambda_o^2\nu^2}
\le C\frac{\epsilon^4}{\nu^2},\quad
\frac{E_c^5k_c^8}{S_3} \le C\frac{\epsilon^4}{\nu^2}.
\end{aligned}
\]
Thus $\mathcal C_3/S_3\ge c-C(k_o/\lambda_o+
\epsilon^2/\nu+\epsilon^4/\nu^2)>c/2$ for large $M$, which proves
\eqref{ce:small-effective-cross}.

\medskip
\noindent\textbf{Step 5. The large second-derivative ratio.}
\proofsteplabel{step:small-envelope-5}{5}

For \eqref{ce:small-effective-high}, the two negative pure $uu$ terms
have a sum containing $1+\varrho R_o^2$, which controls the mixed
terms uniformly in the component ratio $R_o$ defined below.

Finally let $\varrho\ge K\epsilon^4$. Both $f_o$ and $f_c$ have
negative $uu$ coefficients, the latter because
\[
f_{c,uu}=-\lambda_c^2
\bigl[{A_{c,0}}\cos\theta+4{A_{c,1}}\cos2\theta\bigr]
+O(E_c\lambda_c^2\tau_c).
\]
If $|Q|<\theta_QJ_*^2$, write $R_o\triangleq|D_qu_o|/|D_qu_c|$.
Then $|D_qu_o||D_qu_c|\asymp J_*^2$,
$|w|\le C|D_qu_c|(1+R_o)$, and
\[
|D_qy_c|^2\le C|D_qu_c|^2(1+R_o),\qquad
|D_qy_o|^2\le CR_o|D_qu_c|^2(1+R_o).
\]
To identify the five errors, the Hessian contraction in either edge
frame is
\[
\begin{aligned}
\partial^2f_i(w,w)
&=f_{i,uu}w_{u_i}^2+2f_{i,uv}w_{u_i}w_{v_i}
+2f_{i,uy}w_{u_i}w_{y_i} +f_{i,vv}w_{v_i}^2+2f_{i,vy}w_{v_i}w_{y_i}
+f_{i,yy}w_{y_i}^2.
\end{aligned}
\]
For $i=c$, \eqref{ce:front-interface} bounds the five derivative
coefficients, in the order $uv,uy,vv,vy,yy$, by
\[
CE_c\bigl(k_c^2,\lambda_ck_c,\tau_c^2,\tau_ck_c,k_c^2\bigr).
\]
Relative to the negative pure scale
$E_c\lambda_c^2|D_qu_c|^2(1+\varrho R_o^2)$, the five
mixed-derivative terms from $f_c$ are bounded by
\[
C\left\{\frac{\epsilon^2(1+R_o)}{1+\varrho R_o^2},
\frac{\epsilon\sqrt{1+R_o}}{1+\varrho R_o^2},
\frac{\epsilon^4(1+R_o)^2}{1+\varrho R_o^2},
\frac{\epsilon^3(1+R_o)^{3/2}}{1+\varrho R_o^2},
\frac{\epsilon^2(1+R_o)}{1+\varrho R_o^2}\right\}.
\]
For every $d\in[0,2]$, putting $z=\sqrt\varrho R_o$ gives
\[
\frac{(1+R_o)^d}{1+\varrho R_o^2}
\le\varrho^{-d/2}\frac{(1+z)^d}{1+z^2}
\le C(d)\varrho^{-d/2}.
\]
Their sum is therefore at most
$C(K^{-1/2}+K^{-1/4}+K^{-1}+K^{-3/4})$.
For the source terms the corresponding normalized factors are
\[
\begin{gathered}
\frac{\varrho}{1+\varrho R_o^2}
\left\{ 
\epsilon_o^2R_o(1+R_o),\quad
\epsilon_oR_o^{3/2}(1+R_o)^{1/2},\quad
\epsilon_o^4(1+R_o)^2,\quad
\epsilon_o^3R_o^{1/2}(1+R_o)^{3/2},\quad
\epsilon_o^2R_o(1+R_o)
\right\},\\
\epsilon_o=\frac{k_o}{\lambda_o}.
\end{gathered}
\]
Each rational factor is bounded uniformly for $0<\varrho\le1$ and
$R_o>0$, so their sum is $O(\epsilon_o)$.
Choose $K$ large and then $M$ large. The negative pure terms then
satisfy
\[
\begin{aligned}
\mathcal C_1
&\le-cE_c\lambda_c^2|D_qu_c|^2(1+\varrho R_o^2) \le-2cE_c\lambda_c^2\sqrt{\varrho}\,
|D_qu_oD_qu_c|
\le-cE_c\lambda_c^2\epsilon^2J_*^2
=-cE_ck_c^2J_*^2,
\end{aligned}
\]
where $\sqrt{\varrho}\ge\sqrt K\epsilon^2$ and
$|D_qu_oD_qu_c|\asymp J_*^2$. This is
\eqref{ce:small-effective-high}.

\smallskip
\noindent\textit{Simultaneous choice of the constants.}
The choices in Lemma~\ref{ce:small-alternatives} can be made
simultaneously: first take $\eta$ and $\varepsilon_2$ small enough
for Steps~\ref{step:small-envelope-1}
and~\ref{step:small-envelope-2}. Choose $\theta_Q$ small enough to
compare the product of the two $u$ components with $J_*^2$ in
Step~\ref{step:small-envelope-5}, and then choose $K$ to absorb its
mixed terms. With these fixed, choose $q_*$ for
\eqref{ce:small-phase-scale} and $c_*$ for the positive $\mathcal C_1$
bound in Step~\ref{step:small-envelope-3}. Finally choose one $M$
large enough for all the displayed errors, including those in
Step~\ref{step:small-envelope-4}.
\end{proof}

\subsubsection{Comparable amplitude envelopes and tangential phase shifts}\label{ce:comparable}
We now prove \eqref{ce:incoming-chain-target} on the remaining sets
$\mathcal U_o\cap\{\varepsilon_2E_c\le E_o\le K_0E_c\}$.
Set
\[
\begin{aligned}
&E_{\rm tot}\triangleq E_1+E_2,\quad m\triangleq\min(m_1,m_2),\quad
\lambda\triangleq2^m,\quad k\triangleq\max(y_1^{-2},y_2^{-2}),\quad\epsilon\triangleq k/\lambda,\quad
r\triangleq\frac{|D_qu_c|}{E_{\rm tot}k}.
\end{aligned}
\]
Here $E_{\rm tot}$ is the common amplitude scale, $\lambda$ the oscillation
frequency, and $k$ the scale of one normal derivative. Their ratio
$\epsilon=k/\lambda$ tends to zero as $M\to\infty$.
If $|Q|\ge\varepsilon_QE_{\rm tot}^2k^2$, then $\mathcal C_0=Q$ already gives a
nondegeneracy bound. We therefore consider small $Q$.
The parameter $r$ measures the $u_c$ component of $D_q$ relative to
$E_{\rm tot}k$. Small $Q$ will force the opposite component to have size $E_{\rm tot}k/r$;
thus $r$ records which of the two pure $u$ derivatives is amplified.
The following lemma states the four ranges used in this comparison.

\begin{lemma}[Nondegeneracy for comparable amplitudes]\label{ce:comparable-windows}
	There exist positive constants
	\[
	\varepsilon_Q,\ c,\ \varepsilon_1,\ \varepsilon_2',\
	c_{\rm low},\ c_0,\ K_1,\ C_{\max},
	\]
	depending only on the fixed local data, on the solution families, and
	on $\varepsilon_2$, and independent of $M$, with the following
	property. For every sufficiently large integer $M$, at every point of
	\[
	\mathcal U_o\cap\bigl\{\varepsilon_2E_c\le E_o\le K_0E_c,\
	|Q|<\varepsilon_QE_{\rm tot}^2k^2\bigr\}
	\]
	one has $0<r\le C_{\max}/\epsilon$, and the following alternatives
	hold.
	
	\begin{align}
		0<r\le1 &:\quad |\mathcal C_1|\ge cE_{\rm tot}^3\lambda^2k^2
		\quad\text{or}\quad \mathcal C_2\le-cE_{\rm tot}^4\lambda^4k^2;
		\label{ce:window-bounded}\\
		1\le r\le c_{\rm low}\epsilon^{-1/5} &:\quad
		|\mathcal C_1|\ge\varepsilon_1E_{\rm tot}^3\lambda^2k^2/r^2\quad\text{or}\quad
		|\mathcal C_2|\ge\varepsilon_2'E_{\rm tot}^4\lambda^4k^2/r^4\nonumber\\[-2pt]
		&\phantom{:}\quad\text{or}\quad \mathcal C_3\le-cE_{\rm tot}^5\lambda^6k^2/r^6;
		\label{ce:window-lower}\\
		c_0\epsilon^{-1/5}\le r\le K_1\epsilon^{-2/5} &:\quad
		|\mathcal C_1|>E_{\rm tot}^3\lambda^2k^2\epsilon^{2/5}\quad\text{or}\quad
		|\mathcal C_2|>E_{\rm tot}^4\lambda^4k^2\epsilon^{4/5};
		\label{ce:window-interface}\\
		\epsilon^{-2/5}\le r\le C_{\max}/\epsilon &:\quad
		|\mathcal C_1|>E_{\rm tot}^3\lambda^2k^2\epsilon^{2/5}\quad\text{or}\quad
		\mathcal C_3\ge cE_{\rm tot}^5\lambda^4k^4r^4.
		\label{ce:window-upper}
	\end{align}
	Here $c_0<c_{\rm low}$ and $K_1>1$, so the neighboring ranges of $r$
	overlap.
\end{lemma}
Before the proof, the mechanisms behind the four ranges can be read as follows:
\[
\begin{array}{c|l}
\text{range}&\text{conclusion when the preceding alternatives fail}\\\hline
\eqref{ce:window-bounded}&\text{a negative Hessian square gives }\mathcal C_2<0\\
\eqref{ce:window-lower}&\text{the }\mathcal C_1,\mathcal C_2\text{ balances give }\mathcal C_3<0\\
\eqref{ce:window-interface}&\text{the tangential phase shift excludes small }\mathcal C_1,\mathcal C_2\\
\eqref{ce:window-upper}&\text{the pure fourth derivative gives }\mathcal C_3>0.
\end{array}
\]
\begin{proof}[Proof of Lemma~\ref{ce:comparable-windows}]
The proof begins with \eqref{ce:comparable-scales} and \eqref{ce:comparable-flows}, which give common scales and the two characteristic components.
The bounds in \eqref{ce:envelope} imply
\begin{equation}\label{ce:comparable-envelope-log}
\begin{aligned}
m_2^2-(m_1+1)^2+\gamma m_2+\log(E_2/E_1)
&\le y_1^{-1}-y_2^{-1} \le(m_2+1)^2-m_1^2+\gamma(m_2+1)+\log(E_2/E_1).
\end{aligned}
\end{equation}
The bounded ratio $E_1/E_2$ makes $|\log(E_2/E_1)|\le C$.
If $m_1\le m_2-1$, the lower bound in \eqref{ce:comparable-envelope-log} and the definitions of
$m_i$ would give
\[
0>y_1^{-1}-y_2^{-1}\ge\gamma m_2-C>0.
\]
If $m_1\ge m_2+2$, the upper bound in \eqref{ce:comparable-envelope-log} would give
\[
(\beta-1)\beta^{m_2+1}
\le y_1^{-1}-y_2^{-1}
\le(\gamma-2)m_2+\gamma-3+C,
\]
which is impossible for large $M$. Hence $m_1=m_2$ or $m_2+1$.
Substituting $m_1=m_2$ and $m_1=m_2+1$ in \eqref{ce:comparable-envelope-log} yields
\begin{equation}\label{ce:comparable-scales}
m_1-m_2\in\{0,1\},\quad
(\gamma-4)m_2-C\le y_1^{-1}-y_2^{-1}\le(\gamma+2)m_2+C,
\quad y_1/y_2=1+O(m\beta^{-m}).
\end{equation}
The $y_1/y_2$ estimate in \eqref{ce:comparable-scales} follows from
\[
\left|\frac{y_1}{y_2}-1\right|
=\frac{|y_1^{-1}-y_2^{-1}|}{y_1^{-1}}
\le Cm\beta^{-m}.
\]
Consequently $\lambda_i\asymp\lambda$, $k_i\asymp k$,
$E_i\asymp E_{\rm tot}$, and $\tau_i\asymp k^2/\lambda$.
Equations~\eqref{ce:front-jets} and~\eqref{ce:common-null} give
\[
|Y_i-k_if_i|\le CE_{\rm tot}k(m^4\beta^{-m}+\tau_1+\tau_2)=o(E_{\rm tot}k),
\qquad |H_i|\le CE_{\rm tot}k\epsilon.
\]
On $\mathcal U_o$, Lemma~\ref{ce:source-owner} gives
$f_o\ge5E_o/8$ and $f_c\le-3E_o/8$; hence
\[
|Y_1-Y_2|=k_1|f_1|+k_2|f_2|+o(E_{\rm tot}k)\asymp E_{\rm tot}k,
\qquad
|J_*|=|s_\alpha(Y_1-Y_2)-2c_\alpha(H_1+H_2)|\asymp E_{\rm tot}k.
\]
Choose fixed constants $j_0,j_1>0$ such that
$j_0E_{\rm tot}k\le|J_*|\le j_1E_{\rm tot}k$. Recall that $\mathfrak g=1-\cos\alpha_0$ in the bisector frame.
Taking $\varepsilon_Q<j_0^2/\mathfrak g$ in \eqref{ce:lorentz-identities} gives
\[
\tfrac12j_0^2E_{\rm tot}^2k^2
\le J_*^2+\tfrac{\mathfrak g}{2}Q
=(D_qu_o)(D_qu_c)
\le(j_1^2+\tfrac{\mathfrak g}{2}\varepsilon_Q)E_{\rm tot}^2k^2.
\]
Dividing by $|D_qu_c|=E_{\rm tot}kr$ gives the second bound:
\begin{equation}\label{ce:comparable-flows}
\begin{aligned}
&(D_qu_o)(D_qu_c)\asymp E_{\rm tot}^2k^2>0,\qquad |D_qu_o|\asymp E_{\rm tot}k/r,\\
&|\mathrm{d}q|\le CE_{\rm tot}k(r+r^{-1}+1),\qquad0<r\le C_{\max}/\epsilon.
\end{aligned}
\end{equation}
Write $T=q_t$, $X=\nabla q\cdot e_*$,
$Y=\nabla q\cdot n_*$ and $V=(D_qu_1+D_qu_2)/2$.
Projection onto the bisector frame gives
\begin{equation}\label{ce:bisector-projection}
V=T-c_\alpha X,\quad J_*=X-c_\alpha T,\quad
T=\frac{V+c_\alpha J_*}{s_\alpha^2},\quad
X=\frac{J_*+c_\alpha V}{s_\alpha^2},\quad
Y=\frac{D_qu_1-D_qu_2}{2s_\alpha}.
\end{equation}
Projection formulas \eqref{ce:bisector-projection} prove the bound for $|\mathrm dq|$ in
\eqref{ce:comparable-flows}. The product is positive, so $r>0$.
The upper bound follows from
$|D_qu_c|\le C|\mathrm dq|\le CE_{\rm tot}\lambda$ and
$r=|D_qu_c|/(E_{\rm tot}k)\le C\lambda/k=C/\epsilon$.

\medskip\noindent\textbf{Step 1. The range $0<r\le1$.}
For \eqref{ce:window-bounded}, either the negative pure
$u_o$ Hessian term survives in $\mathcal C_1$, or its cancellation forces
a negative square in $\mathcal C_2$.
In the coordinates $(u_c,v_c,y_c)$, the components of $w$ are
$O(E_{\rm tot}kr,E_{\rm tot}k/r,E_{\rm tot}k)$; in the coordinates of $f_o$, each component is
$O(E_{\rm tot}k/r)$. Here $w_u=D_qu_i$, $w_v=D_qv_i$, and $w_y=D_qy_i$
in the indicated coordinates. The identity $w_y^2=w_uw_v-Q$
gives $|D_qy_c|\le CE_{\rm tot}k$.
Put $h_o\triangleq f_{o,uu}\le-cE_{\rm tot}\lambda^2$,
$h_c\triangleq f_{c,uu}$, and $S_1\triangleq E_{\rm tot}^3\lambda^2k^2/r^2$.
The negative pure contribution satisfies
$h_ow_{u_o}^2\le-cS_1$. The six contributions from $f_c$ obey
\[
\begin{aligned}
&\frac1{S_1}\bigl(|f_{c,uu}w_{u_c}^2|,
2|f_{c,uy}w_{u_c}w_{y_c}|,|f_{c,yy}w_{y_c}^2|, |f_{c,vv}w_{v_c}^2|,2|f_{c,vy}w_{v_c}w_{y_c}|,
2|f_{c,uv}w_{u_c}w_{v_c}|\bigr)\\ &\le C(r^4,\epsilon r^3,\epsilon^2r^2,\epsilon^4,
\epsilon^3r,\epsilon^2r^2).
\end{aligned}
\]
The mixed terms from $f_o$ are at most $C\epsilon S_1$.
Consequently
\[
\mathcal C_1\le-2cS_1+CS_1(r^4+\epsilon).
\]
Choose $r_0\in(0,1)$ with $Cr_0^4\le c/2$, then choose $M$
so that $C\epsilon\le c/2$. For $0<r\le r_0$, this gives
$\mathcal C_1\le-cS_1\le-cE_{\rm tot}^3\lambda^2k^2$.
For $r_0\le r\le1$, suppose the first alternative in
\eqref{ce:window-bounded} fails. After decreasing its threshold $c$ if
necessary, this forces $h_c\asymp E_{\rm tot}\lambda^2>0$ and
$h_c(D_qu_c)^2=-h_o(D_qu_o)^2[1+O(c+\epsilon)]$.
The pure covector $L\triangleq h_o(D_qu_o)\,\mathrm{d}u_o+h_c(D_qu_c)\,\mathrm{d}u_c$ has
\begin{equation}\label{ce:pure-negative-square}
\mathcal B(L,L)=2\mathfrak g h_oh_c(D_qu_o)(D_qu_c)\le-cE_{\rm tot}^4\lambda^4k^2.
\end{equation}
Put $R=Hw-L$. Since $r_0\le r\le1$, the mixed derivative bounds in \eqref{ce:front-interface}, together with \eqref{ce:comparable-flow-table}, give
\[
|R|\le CE_{\rm tot}^2\lambda k^2,\quad |L|\le CE_{\rm tot}^2\lambda^2k,\quad
|\partial^3q(w^3)|\le CE_{\rm tot}^4\lambda^3k^3.
\]
The $\mathcal C_2$ formula in \eqref{ce:chain-expansion} therefore yields
\[
\begin{aligned}
\mathcal C_2&=4\mathcal B(L,L)+2\partial^3q(w^3)
+8\mathcal B(L,R)+4\mathcal B(R,R)\\
&\le-cE_{\rm tot}^4\lambda^4k^2
+C(E_{\rm tot}^4\lambda^3k^3+E_{\rm tot}^4\lambda^2k^4)\\
&\le E_{\rm tot}^4\lambda^4k^2[-c+C(\epsilon+\epsilon^2)]
\le-c'E_{\rm tot}^4\lambda^4k^2
\end{aligned}
\]
for sufficiently large $M$. This proves \eqref{ce:window-bounded}.

\medskip\noindent\textbf{Step 2. The lower growth range.}
In the range of \eqref{ce:window-lower}, if $\mathcal C_1$ is below its threshold,
the two pure Hessian terms balance. If $\mathcal C_2$ is also below its
threshold, a third-derivative term is fixed by this balance and
contributes negatively to $\mathcal C_3$. All remaining contractions are bounded relative to the same scale.
For $r\ge1$, \eqref{ce:comparable-flows} and the definition of $r$ give
\[
\begin{aligned}
&|w_{u_o}|\le CE_{\rm tot}k/r,\qquad |w_{u_c}|=E_{\rm tot}kr,\\
&|w_{v_o}|+|w_{v_c}|+|w_{y_c}|
\le C|\mathrm dq|\le CE_{\rm tot}kr,\\
&|w_{y_o}|^2
=|w_{u_o}w_{v_o}-Q|
\le (CE_{\rm tot}k/r)(CE_{\rm tot}kr)+\varepsilon_QE_{\rm tot}^2k^2
\le CE_{\rm tot}^2k^2.
\end{aligned}
\]
Thus the components of $w$ in the two edge coordinates satisfy
\begin{equation}\label{ce:comparable-flow-table}
\begin{array}{c|ccc}
&w_u&w_v&w_y\\\hline
o&O(E_{\rm tot}k/r)&O(E_{\rm tot}kr)&O(E_{\rm tot}k)\\
c&O(E_{\rm tot}kr)&O(E_{\rm tot}kr)&O(E_{\rm tot}kr).
\end{array}
\end{equation}
To evaluate the contractions from these components, define
\begin{equation}\label{ce:weighted-component-norms}
\Lambda_i^{v,y}(\zeta)\triangleq
\epsilon^2|\zeta_{v_i}|+\epsilon|\zeta_{y_i}|,
\qquad
\Lambda_i(\zeta)\triangleq|\zeta_{u_i}|+\Lambda_i^{v,y}(\zeta).
\end{equation}
The component expansion \eqref{ce:tensor-component-sum},
with $\tau_i\asymp\epsilon^2\lambda$ and $k_i\asymp\epsilon\lambda$,
gives, for $1\le n\le4$,
\begin{equation}\label{ce:weighted-contractions}
\begin{aligned}
|\partial^nf_i(\zeta_1,\ldots,\zeta_n)|
&\le CE_{\rm tot}\lambda^n\prod_{j=1}^n\Lambda_i(\zeta_j),\\
\left|\partial^nf_i(\zeta_1,\ldots,\zeta_n)
-f_{i,u^n}\prod_{j=1}^n(\zeta_j)_{u_i}\right|
&\le CE_{\rm tot}\lambda^n\left\{
\prod_{j=1}^n\Lambda_i(\zeta_j)
-\prod_{j=1}^n|(\zeta_j)_{u_i}|\right\}.
\end{aligned}
\end{equation}
The difference on the right is a sum of nonnegative products, since
\[
\prod_{j=1}^n\Lambda_i(\zeta_j)-\prod_{j=1}^n|(\zeta_j)_{u_i}|
=\sum_{a=1}^n\Lambda_i^{v,y}(\zeta_a)
\prod_{j<a}|(\zeta_j)_{u_i}|\prod_{j>a}\Lambda_i(\zeta_j).
\]
For example, \eqref{ce:comparable-flow-table} yields
\[
\begin{aligned}
&\Lambda_o(w)\le CE_{\rm tot}k/r,\qquad
\Lambda_o^{v,y}(w)\le CE_{\rm tot}k\epsilon,\qquad\Lambda_c(w)\le CE_{\rm tot}kr,\qquad
\Lambda_c^{v,y}(w)\le CE_{\rm tot}k\epsilon r,
\end{aligned}
\]
where $\epsilon r\le C_{\max}$ was used in the $o$-frame.
In the range \eqref{ce:window-lower}, assume both the $\mathcal C_1$ and $\mathcal C_2$ alternatives
fail. Expanding $\mathcal C_1=2\partial^2q(w,w)$ and using
\eqref{ce:front-interface} gives
\begin{equation}\label{ce:lower-C1-errors}
\begin{aligned}
&|\mathcal C_1-2h_ow_{u_o}^2-2h_cw_{u_c}^2| \le CE_{\rm tot}\lambda^2\sum_{i=o,c}
\Lambda_i^{v,y}(w)\bigl(\Lambda_i(w)+|w_{u_i}|\bigr)
\le CE_{\rm tot}^3\lambda^2k^2\epsilon(r^{-1}+r^2),\\
&|h_cw_{u_c}^2+h_ow_{u_o}^2|
\le C(\varepsilon_1+\epsilon r^4)
\frac{E_{\rm tot}^3\lambda^2k^2}{r^2},\qquad
-h_ow_{u_o}^2\asymp\frac{E_{\rm tot}^3\lambda^2k^2}{r^2}.
\end{aligned}
\end{equation}
For small $\varepsilon_1$ and large $M$, \eqref{ce:lower-C1-errors} implies
\begin{equation}\label{ce:lower-balance}
h_c\asymp E_{\rm tot}\lambda^2/r^4,\qquad
h_c(D_qu_c)^2=-h_o(D_qu_o)^2[1+O(\varepsilon_1+\epsilon r^4)].
\end{equation}
We next compute the components of
$\mathfrak a_0\triangleq L^\sharp$ and
$\varrho_0\triangleq(Hw-L)^\sharp$.
For $\xi,\eta\in\{u,v,y\}$, put
$G_{ij}^{\xi\eta}\triangleq
\mathcal B(\mathrm d\xi_i,\mathrm d\eta_j)$ and
$\mathfrak s\triangleq e_c\cdot n_o$,
$\mathfrak j\triangleq1+\cos\alpha_0$.
Since $\mathrm du_i=(1,e_i)$, $\mathrm dv_i=(1,-e_i)$,
and $\mathrm dy_i=(0,n_i)$, the pairing matrices, with rows and
columns ordered as $u,v,y$, are
\begin{equation}\label{ce:frame-pairing-matrices}
G_{ii}=\begin{pmatrix}0&2&0\\2&0&0\\0&0&-1\end{pmatrix},
\qquad
G_{oc}=\begin{pmatrix}
\mathfrak g&\mathfrak j&\mathfrak s\\
\mathfrak j&\mathfrak g&-\mathfrak s\\
-\mathfrak s&\mathfrak s&-\cos\alpha_0
\end{pmatrix},\qquad G_{co}=G_{oc}^{\mathsf T}.
\end{equation}
For example, the first row is obtained from
\[
\begin{aligned}
G_{oc}^{uu}&=1-e_o\cdot e_c=\mathfrak g,\qquad
G_{oc}^{uv} =1+e_o\cdot e_c=\mathfrak j,\qquad
G_{oc}^{uy} =-e_o\cdot n_c=e_c\cdot n_o=\mathfrak s.
\end{aligned}
\]
The identity $e_o\cdot n_c=-e_c\cdot n_o$ follows from
$n_i=\mathsf J e_i$ and $\mathsf J^{\mathsf T}=-\mathsf J$.
Consequently, every component is given by the finite sums
\begin{equation}\label{ce:acceleration-component-sums}
\begin{aligned}
(\mathfrak a_0)_{\xi_i}
&=\sum_{j\in\{o,c\}}G_{ij}^{\xi u}h_jw_{u_j},\qquad
(\varrho_0)_{\xi_i}
=\sum_{j\in\{o,c\}}
\sum_{\substack{\eta,\zeta\in\{u,v,y\}\\(\eta,\zeta)\ne(u,u)}}
G_{ij}^{\xi\eta}f_{j,\eta\zeta}w_{\zeta_j}.
\end{aligned}
\end{equation}
For the pure part, the sums in \eqref{ce:acceleration-component-sums} read
\begin{equation}\label{ce:pure-acceleration-components}
\begin{aligned}
\bigl((\mathfrak a_0)_{u_o},(\mathfrak a_0)_{v_o},(\mathfrak a_0)_{y_o}\bigr)
&=\bigl(\mathfrak g h_cw_{u_c},\,
2h_ow_{u_o}+\mathfrak jh_cw_{u_c},\,-\mathfrak sh_cw_{u_c}\bigr),\\
\bigl((\mathfrak a_0)_{u_c},(\mathfrak a_0)_{v_c},(\mathfrak a_0)_{y_c}\bigr)
&=\bigl(\mathfrak g h_ow_{u_o},\,
\mathfrak jh_ow_{u_o}+2h_cw_{u_c},\,\mathfrak sh_ow_{u_o}\bigr),\\
|h_ow_{u_o}|&\le CE_{\rm tot}^2\lambda^2k/r,\qquad
|h_cw_{u_c}|\le CE_{\rm tot}^2\lambda^2k/r^3.
\end{aligned}
\end{equation}
For the remaining Hessian entries, \eqref{ce:front-interface} gives
the componentwise bound
\[
\left|\bigl(f_{i,\eta\zeta}
-h_i\mathbf1_{\{\eta=\zeta=u\}}\bigr)_{\eta,\zeta\in\{u,v,y\}}\right|
\le CE_{\rm tot}\lambda^2
\begin{pmatrix}
0&\epsilon^2&\epsilon\\
\epsilon^2&\epsilon^4&\epsilon^3\\
\epsilon&\epsilon^3&\epsilon^2
\end{pmatrix}.
\]
In the $o$-frame, multiplication by
$E_{\rm tot}k(r^{-1},r,1)^{\mathsf T}$ bounds the three covector coefficients by
\[
CE_{\rm tot}^2\lambda k^2
\begin{pmatrix}
1+\epsilon r\\
\epsilon/r+\epsilon^3r+\epsilon^2\\
r^{-1}+\epsilon^2r+\epsilon
\end{pmatrix}.
\]
In the $c$-frame the corresponding coefficients are bounded by
$CE_{\rm tot}^2\lambda k^2r(1,\epsilon,1)^{\mathsf T}$.
Put
$R_{j,\eta}=\sum_{\zeta}f_{j,\eta\zeta}w_{\zeta_j}
-h_jw_{u_j}\mathbf1_{\{\eta=u\}}$. Then
\begin{equation}\label{ce:lower-covector-row}
\begin{aligned}
(\varrho_0)_{u_c}
&=\mathfrak gR_{o,u}+\mathfrak jR_{o,v}
-\mathfrak sR_{o,y}+2R_{c,v},\\
R_{c,v}&=f_{c,vu}w_{u_c}+f_{c,vv}w_{v_c}+f_{c,vy}w_{y_c},\\
|R_{c,v}|&\le CE_{\rm tot}\lambda^2(\epsilon^2+\epsilon^4+\epsilon^3)E_{\rm tot}kr
\le CE_{\rm tot}^2\lambda k^2\epsilon r.
\end{aligned}
\end{equation}
The $u_c$ row of $G_{cc}$ in \eqref{ce:frame-pairing-matrices} is $(0,2,0)$, so the $c$-side covector in \eqref{ce:lower-covector-row} contributes only $CE_{\rm tot}^2\lambda k^2\epsilon r$ to
$(\varrho_0)_{u_c}$. The $o$-side covector in \eqref{ce:lower-covector-row} contributes at most
$CE_{\rm tot}^2\lambda k^2(1+\epsilon r)$ through $G_{co}$.
Thus, using $r\ge1$ and
$\epsilon r\le c_{\rm low}\epsilon^{4/5}\le C$,
\begin{equation}\label{ce:lower-remainder-components}
\begin{aligned}
| (\varrho_0)_{u_c}|&\le CE_{\rm tot}^2\lambda k^2(1+\epsilon r)
\le CE_{\rm tot}^2\lambda k^2,\\
| (\varrho_0)_{\xi_i}|&\le CE_{\rm tot}^2\lambda k^2r
\quad\bigl((i,\xi)\ne(c,u)\bigr).
\end{aligned}
\end{equation}
Combining \eqref{ce:pure-acceleration-components} and \eqref{ce:lower-remainder-components} gives
\begin{equation}\label{ce:lower-acceleration}
\begin{array}{c|ccc}
&u&v&y\\\hline
\mathfrak a_0\text{ in }o&O(E_{\rm tot}^2\lambda^2k/r^3)&O(E_{\rm tot}^2\lambda^2k/r)&O(E_{\rm tot}^2\lambda^2k/r^3)\\
\mathfrak a_0\text{ in }c&O(E_{\rm tot}^2\lambda^2k/r)&O(E_{\rm tot}^2\lambda^2k/r)&O(E_{\rm tot}^2\lambda^2k/r)\\
\varrho_0\text{ in }o&O(E_{\rm tot}^2\lambda k^2r)&O(E_{\rm tot}^2\lambda k^2r)&O(E_{\rm tot}^2\lambda k^2r)\\
\varrho_0\text{ in }c&O(E_{\rm tot}^2\lambda k^2)&O(E_{\rm tot}^2\lambda k^2r)&O(E_{\rm tot}^2\lambda k^2r).
\end{array}
\end{equation}
Since $\epsilon r^2\le c_{\rm low}^2\epsilon^{3/5}\le1$
for large $M$, the weighted component sums defined in \eqref{ce:weighted-component-norms} satisfy
\begin{equation}\label{ce:weighted-acceleration-bounds}
\begin{aligned}
\frac{\bigl(\Lambda_o(\mathfrak a_0),\Lambda_c(\mathfrak a_0),
\Lambda_o^{v,y}(\mathfrak a_0),\Lambda_c^{v,y}(\mathfrak a_0)\bigr)}
{E_{\rm tot}^2\lambda^2k}
&\le C(r^{-3},r^{-1},\epsilon r^{-3},\epsilon r^{-1}),\\
\frac{\bigl(\Lambda_o(\varrho_0),\Lambda_c(\varrho_0),
\Lambda_o^{v,y}(\varrho_0),\Lambda_c^{v,y}(\varrho_0)\bigr)}
{E_{\rm tot}^2\lambda^2k}
&\le C(\epsilon r,\epsilon,\epsilon^2r,\epsilon^2r).
\end{aligned}
\end{equation}
For example, the first row uses
\[
\Lambda_o^{v,y}(\mathfrak a_0)
\le CE_{\rm tot}^2\lambda^2k(\epsilon^2/r+\epsilon/r^3)
\le CE_{\rm tot}^2\lambda^2k\epsilon/r^3.
\]
Put $t_c\triangleq f_{c,uuu}$ and
$S_j\triangleq E_{\rm tot}^{j+2}\lambda^{2j}k^2/r^{2j}$ for $j\in\{2,3\}$.
These scales are the sizes of the pure Hessian products:
\[
|h_oh_cw_{u_o}w_{u_c}|\asymp
(E_{\rm tot}\lambda^2)(E_{\rm tot}\lambda^2/r^4)(E_{\rm tot}k/r)(E_{\rm tot}kr)=S_2,
\qquad
|h_o^2h_cw_{u_o}^2|\asymp S_3.
\]
The coefficients in the $\mathcal C_2$ identity follow from the product rule:
\[
D_qw=\mathfrak a,\quad
\mathcal C_2=D_q[2H(w,w)]
=2\partial^3q(w^3)+4H(\mathfrak a,w)
=2\partial^3q(w^3)+4\mathcal B(Hw,Hw).
\]
Since $Hw=L+(Hw-L)$, the last square is
\[
4\mathcal B(L,L)+8\mathcal B(L,Hw-L)
+4\mathcal B(Hw-L,Hw-L),
\quad \mathcal B(L,Hw-L)=L(\varrho_0).
\]
Using \eqref{ce:pure-negative-square} for the algebraic expression
of $\mathcal B(L,L)$ gives the complete identity
\begin{equation}\label{ce:lower-C2-identity}
\begin{aligned}
\mathcal C_2={}&2t_cw_{u_c}^3+8\mathfrak g h_oh_cw_{u_o}w_{u_c} +2\bigl\{\partial^3f_o(w^3)+\partial^3f_c(w^3)-t_cw_{u_c}^3\bigr\}
+8L(\varrho_0)+4\mathcal B(Hw-L,Hw-L).
\end{aligned}
\end{equation}
Removing the pure $c$-$uuu$ term uses the polynomial identity
\[
\Lambda_c(w)^3-|w_{u_c}|^3
=\Lambda_c^{v,y}(w)
\bigl(\Lambda_c(w)^2+\Lambda_c(w)|w_{u_c}|+|w_{u_c}|^2\bigr).
\]
Thus all the error terms in \eqref{ce:lower-C2-identity} satisfy
\begin{equation}\label{ce:lower-C2-errors}
\begin{aligned}
&\frac{|\partial^3f_o(w^3)+\partial^3f_c(w^3)-t_cw_{u_c}^3|}{S_2}
\le C\epsilon r^4(r^{-3}+\epsilon r^3),\\
&L(\varrho_0)
=h_ow_{u_o}(\varrho_0)_{u_o}+h_cw_{u_c}(\varrho_0)_{u_c},\qquad
\frac{|L(\varrho_0)|}{S_2}
\le Cr^4(\epsilon+\epsilon/r^3),\\
&\frac{|\mathcal B(Hw-L,Hw-L)|}{S_2}
=\frac{|(\varrho_0)_{u_c}(\varrho_0)_{v_c}-(\varrho_0)_{y_c}^2|}{S_2}
\le C\epsilon^2(r^5+r^6).
\end{aligned}
\end{equation}
Since $|\mathfrak g h_oh_cw_{u_o}w_{u_c}|\asymp S_2$ and
$|\mathcal C_2|<\varepsilon_2'S_2$, \eqref{ce:lower-C2-errors} gives
\begin{equation}\label{ce:lower-C2-balance}
\begin{aligned}
f_{c,uuu}(D_qu_c)^3&=-4\mathfrak g h_oh_c(D_qu_o)(D_qu_c) \bigl[1+O(\varepsilon_2'+\epsilon r^4+\epsilon^2r^6
+\epsilon r+\epsilon^2r^7)\bigr].
\end{aligned}
\end{equation}
In particular,
\[
|t_c|\le C\frac{S_2}{(E_{\rm tot}kr)^3}
=CE_{\rm tot}\lambda^3\epsilon^{-1}r^{-7}.
\]
Consequently the leading third-order acceleration term in $\mathcal C_3$ is
\begin{equation}\label{ce:lower-C3-leading}
\begin{aligned}
14f_{c,uuu}(\mathfrak a_0)_{u_c}(D_qu_c)^2
&=-56\mathfrak g^2h_o^2h_c(D_qu_o)^2 \bigl[1+O(\varepsilon_2'+\epsilon r^4+\epsilon^2r^6
+\epsilon r+\epsilon^2r^7)\bigr]\\
&\le-cE_{\rm tot}^5\lambda^6k^2/r^6.
\end{aligned}
\end{equation}
To compute every remaining contraction relative to $S_3$, note first that
\[
\frac{E_{\rm tot}\lambda^4(E_{\rm tot}k)^4}{S_3}=\epsilon^2r^6,\qquad
\frac{E_{\rm tot}\lambda^3(E_{\rm tot}^2\lambda^2k)(E_{\rm tot}k)^2}{S_3}=\epsilon r^6,\qquad
\frac{E_{\rm tot}\lambda^2(E_{\rm tot}^2\lambda^2k)^2}{S_3}=r^6.
\]
The fourth-order terms and the pure $c$-$uuu$ remainder therefore obey
\[
\begin{aligned}
&\frac{|\partial^4f_o(w^4)|+|\partial^4f_c(w^4)|}{S_3}
\le C\epsilon^2r^6(r^{-4}+r^4)
=C(\epsilon^2r^2+\epsilon^2r^{10}),\\
&\frac{|t_c(\varrho_0)_{u_c}w_{u_c}^2|}{S_3}
\le\frac{(CE_{\rm tot}\lambda^3\epsilon^{-1}r^{-7})
(CE_{\rm tot}^2\lambda^2k\epsilon)(E_{\rm tot}kr)^2}{S_3}
\le C\epsilon r.
\end{aligned}
\]
For all nonpure third-order terms, use
\[
\begin{aligned}
&\Lambda_i(\zeta)\Lambda_i(w)^2-|\zeta_{u_i}||w_{u_i}|^2 =\Lambda_i^{v,y}(\zeta)\Lambda_i(w)^2
+|\zeta_{u_i}|\Lambda_i^{v,y}(w)
\bigl(\Lambda_i(w)+|w_{u_i}|\bigr).
\end{aligned}
\]
Substitution of \eqref{ce:weighted-acceleration-bounds} into \eqref{ce:weighted-contractions} gives
\begin{equation}\label{ce:lower-third-errors}
\begin{aligned}
&\frac{|\partial^3f_o(\mathfrak a,w^2)|
+|\partial^3f_c(\mathfrak a,w^2)-t_c\mathfrak a_{u_c}w_{u_c}^2|}{S_3} \le C\epsilon r^6
\bigl(r^{-5}\!+\!\epsilon/r\!+\!\epsilon r\!+\!\epsilon^2r^3\bigr)
\le C(\epsilon r\!+\!\epsilon^2r^7\!+\!\epsilon^3r^9).
\end{aligned}
\end{equation}
The four terms in parentheses in \eqref{ce:lower-third-errors} arise from $\mathfrak a_0,\varrho_0$
in the $o$-frame and then the nonpure $\mathfrak a_0,\varrho_0$
terms in the $c$-frame, respectively.
For the Hessian terms, retain the exact cancellation
\[
\begin{aligned}
&h_o((\mathfrak a_0)_{u_o})^2+h_c((\mathfrak a_0)_{u_c})^2 =\mathfrak g^2h_oh_c\bigl[h_c(D_qu_c)^2+h_o(D_qu_o)^2\bigr],\\
&\left|h_o((\mathfrak a_0)_{u_o})^2+h_c((\mathfrak a_0)_{u_c})^2\right| \le C\frac{E_{\rm tot}^5\lambda^6k^2}{r^6}
(\varepsilon_1+\epsilon r^4),
\end{aligned}
\]
where the inequality follows from \eqref{ce:lower-balance}.
The component formula \eqref{ce:weighted-contractions}, with its pure $uu$ product removed, gives
\begin{equation}\label{ce:lower-Hessian-errors}
\begin{aligned}
&\frac{\left|\sum_{i=o,c}\bigl\{\partial^2f_i(\mathfrak a_0,\mathfrak a_0)
-h_i((\mathfrak a_0)_{u_i})^2\bigr\}\right|}{S_3}
\le Cr^6(\epsilon/r^6+\epsilon/r^2),\\
&\frac{|H(\mathfrak a_0,\varrho_0)|}{S_3}
\le Cr^6(\epsilon/r^2+\epsilon/r^5+\epsilon^2),\quad
\frac{|H(\varrho_0,\varrho_0)|}{S_3}
\le Cr^6(\epsilon^2r^2+\epsilon^2/r^4+\epsilon^3r).
\end{aligned}
\end{equation}
In the last two lines of \eqref{ce:lower-Hessian-errors}, the pure $c$-$uu$ products use
$|h_c|\le CE_{\rm tot}\lambda^2/r^4$ and contribute, respectively,
$C\epsilon/r^5$ and $C\epsilon^2/r^4$ inside parentheses.
All other products use the weighted component bounds in \eqref{ce:weighted-acceleration-bounds}.
Since
\[
H(\mathfrak a,\mathfrak a)
=H(\mathfrak a_0,\mathfrak a_0)
+2H(\mathfrak a_0,\varrho_0)+H(\varrho_0,\varrho_0),
\]
\eqref{ce:lower-third-errors} and \eqref{ce:lower-Hessian-errors} are summarized, relative to $S_3$, by
\begin{equation}\label{ce:lower-budget}
\begin{array}{c|c}
\partial^4f_c(w^4),\ \partial^4f_o(w^4)&C\cdot(\epsilon^2r^{10}+\epsilon^2r^2)\\
f_{c,uuu}(\varrho_0)_{u_c}(D_qu_c)^2&C\epsilon r\\
\text{other third-order terms}&C\cdot(\epsilon r+\epsilon^2r^7+\epsilon^3r^9)\\
H(\mathfrak a,\mathfrak a)&C\cdot(\varepsilon_1+\epsilon r^4+\epsilon^2r^8).
\end{array}
\end{equation}
The fixed coefficients $2,14,8$ in \eqref{ce:chain-expansion}
are included in $C$.
The upper endpoint $r\le c_{\rm low}\epsilon^{-1/5}$ was chosen
so that the largest fourth-order error is controlled by
\[
\epsilon^2r^{10}\le c_{\rm low}^{10},\quad
\epsilon r^4\le c_{\rm low}^4\epsilon^{1/5},\quad
\epsilon^2r^8\le c_{\rm low}^8\epsilon^{2/5},\quad
\epsilon^3r^9\le c_{\rm low}^9\epsilon^{6/5}.
\]
The remaining monomials in \eqref{ce:lower-C2-balance} and
\eqref{ce:lower-budget} have no larger upper bound. Thus their sum is at most
$C(\varepsilon_1+\varepsilon_2'+c_{\rm low}^{10}+\epsilon^{1/5})$.
First choose $\varepsilon_1,\varepsilon_2'$ and $c_{\rm low}$ small enough,
and then $M$ large enough, to make this sum less than half the coefficient
of the negative leading term in \eqref{ce:lower-C3-leading}. This proves~\eqref{ce:window-lower}.

\medskip\noindent\textbf{Step 3. The intermediate range and the tangential phase shift.}
For \eqref{ce:window-interface}, argue by contradiction. Small $\mathcal C_1$
and $\mathcal C_2$ would force $f_c$ to have a negative value but small first,
second, and third pure $u_c$ derivatives. The nonzero relative phase shift between the two active frequencies excludes these four conditions; the lower
bound for $|z_c|$ in \eqref{ce:midpoint-size} makes this shift larger
than every error retained below. In the coordinates of $f_c$ put
$\epsilon_c\triangleq k_c/\lambda_c$ and
$D_s\triangleq z_ck_c^2/\lambda_c$.
Equation~\eqref{ce:midpoint-size} gives,
uniformly for $m_c\ge M-C_*$,
\begin{equation}\label{ce:dispersion-ratios}
|D_s|/\epsilon_c\longrightarrow\infty,\quad
\epsilon_c^{6/5}/|D_s|\longrightarrow0,\quad
|D_s|/\epsilon_c^{4/5}\longrightarrow0,\quad
\epsilon_c^2/|D_s|\longrightarrow0
\qquad(M\longrightarrow\infty).
\end{equation}
These four ratios are respectively comparable to
$\beta^{-M+2m_c}$, $\beta^{M-24m_c/5}$, $\beta^{-M-4m_c/5}$,
and $\beta^{M-16m_c}$, with comparison constants depending only on the fixed local data.
Write ${A_{c,0}}\triangleq A_{\gamma_c,m_c}(y_c)$ and
${A_{c,1}}\triangleq A_{\gamma_c,m_c+1}(y_c)$ for its two amplitude coefficients;
thus ${A_{c,0}}+{A_{c,1}}=E_c$. The complete local solution has, for $0\le h\le3$,
\begin{equation}\label{ce:shifted-jets}
\begin{aligned}
&\lambda_c^{-h}\partial_{u_c}^hf_c
={A_{c,0}}\cos^{(h)}\theta+2^h{A_{c,1}}\cos^{(h)}\phi
+O\bigl(E_{\rm tot}[m_c^{20}\beta^{-m_c}|D_s|+D_s^2+\epsilon_c^2]\bigr),\\
&\theta\triangleq\lambda_c u_c+D_s/2,\qquad
\phi\triangleq2\lambda_c u_c+D_s/4,\qquad
d_s\triangleq\phi-2\theta=-3D_s/4.
\end{aligned}
\end{equation}
The term with $\ell=1$ is
$-z_cA_{\gamma_c,j}''\sin(\lambda_ju_c)/(2\lambda_j)$, where primes
denote derivatives in $y_c$. The amplitude estimate in
Lemma~\ref{ce:packet-bounds}\ref{cond:packet-amplitude} implies
\[
\sum_{j=m_c}^{m_c+1}|A_{\gamma_c,j}''-k_c^2A_{\gamma_c,j}|
\le CE_{\rm tot}k_c^2m_c^4\beta^{-m_c}.
\]
Consequently,
\[
\begin{aligned}
&\sum_{j=m_c}^{m_c+1}
\left[A_{\gamma_c,j}\cos(\lambda_ju_c)
-\frac{z_cA_{\gamma_c,j}''}{2\lambda_j}\sin(\lambda_ju_c)\right]\\
& ={A_{c,0}}\cos(\lambda_cu_c+D_s/2)
+{A_{c,1}}\cos(2\lambda_cu_c+D_s/4) +O\bigl(E_{\rm tot}[m_c^4\beta^{-m_c}|D_s|+D_s^2]\bigr).
\end{aligned}
\]
The coordinate dependence also satisfies
\[
\partial_{u_c}z_c=\tfrac12,\qquad
\lambda_c^{-1}\left|\partial_{u_c}
\left(\frac{z_ck_c^2}{2\lambda_j}\right)\right|
=\frac{k_c^2}{4\lambda_c\lambda_j}\le C\epsilon_c^2,
\quad j\in\{m_c,m_c+1\}.
\]
For $2\le l\le K_j$, each monomial has the form
$c\lambda_j^{-l}z_c^dA_{\gamma_c,j}^{(l+d)}$, with integer indices
$1\le d\le l$ and $l-d$ even.
Here $j\ge4$ is an active packet index.
For each fixed $l$, its normalized magnitude is bounded by
$CE_{\rm tot}(D_s^2+\epsilon_c^2)$. To estimate the sum over large $l$,
\eqref{ce:midpoint-size} and $m_c\ge M-C_*$ give
$|D_s|^2\ge c\beta^{-26m_c}$. In Lemma~\ref{ce:packet-bounds}\ref{cond:packet-tail},
fix $B>26\log\beta$ and choose the same threshold $L$ for the finitely
many derivative orders $0\le h\le3$ in \eqref{ce:shifted-jets}. The terms with $l\ge L$ are
$O(E_{\rm tot}e^{-Bm_c})=o(E_{\rm tot}|D_s|^2)$, so \eqref{ce:shifted-jets} holds for the full sum. The correction within each cell satisfies the stronger residual bound of Proposition~\ref{ce:reserve}.
Thus the error in~\eqref{ce:shifted-jets} is $o(E_{\rm tot}|D_s|)$.

\smallskip\noindent\textbf{Bad-set lemma (incompatible derivative conditions).}
For every fixed $c_7,C_7>0$ and sufficiently large $M$ depending
on these constants, no point in the stated range satisfies
\begin{equation}\label{ce:dispersion-exclusion}
f_c\le-c_7E_{\rm tot},\quad |f_{c,u}|\le C_7E_{\rm tot}k_c,\quad
|f_{c,uu}|\le C_7E_{\rm tot}\lambda_c^2\epsilon_c^{4/5},\quad
|f_{c,uuu}|\le C_7E_{\rm tot}\lambda_c^3\epsilon_c^{2/5}.
\end{equation}
We eliminate the even and odd derivative data separately. The matrices are Vandermonde-type systems generated by the squared frequencies $1,4$, with determinants $-3$ and $-6$. Thus derivatives of orders zero through three determine the cosine and sine coefficients of the two active frequencies.

To prove the claim, set
\[
\widehat f_h\triangleq\lambda_c^{-h}\partial_{u_c}^hf_c,\qquad
\eta_M\triangleq m_c^{20}\beta^{-m_c}|D_s|+D_s^2+\epsilon_c^2.
\]
By \eqref{ce:dispersion-ratios}, $\eta_M=o(|D_s|)$ and
$|D_s|=o(\epsilon_c^{4/5})$. Write the four errors in
\eqref{ce:shifted-jets} as $\mathcal E_h$, with $|\mathcal E_h|\le CE_{\rm tot}\eta_M$.
Putting $X={A_{c,0}}\cos\theta$ and $Y={A_{c,1}}\cos\phi$ gives
\[
\begin{pmatrix}\widehat f_0-\mathcal E_0\\\widehat f_2-\mathcal E_2\end{pmatrix}
=\begin{pmatrix}1&1\\-1&-4\end{pmatrix}
\begin{pmatrix}X\\Y\end{pmatrix},\qquad
X=\frac{4\widehat f_0+\widehat f_2-4\mathcal E_0-\mathcal E_2}{3},\quad
Y=-\frac{\widehat f_0+\widehat f_2-\mathcal E_0-\mathcal E_2}{3}.
\]
The bounds in \eqref{ce:dispersion-exclusion} therefore imply
$X\le-c'E_{\rm tot}$ and $Y\ge c'E_{\rm tot}$, so ${A_{c,0}},{A_{c,1}}\ge c'E_{\rm tot}$.
Likewise, putting $P_1={A_{c,0}}\sin\theta$ and $P_2={A_{c,1}}\sin\phi$ gives
\[
\begin{aligned}
\begin{pmatrix}\widehat f_1-\mathcal E_1\\\widehat f_3-\mathcal E_3\end{pmatrix}
&=\begin{pmatrix}-1&-2\\1&8\end{pmatrix}
\begin{pmatrix}P_1\\P_2\end{pmatrix},\quad
P_1 =-\frac{4\widehat f_1+\widehat f_3-4\mathcal E_1-\mathcal E_3}{3},\quad
P_2=\frac{\widehat f_1+\widehat f_3-\mathcal E_1-\mathcal E_3}{6},\\
|P_1|+|P_2|&\le CE_{\rm tot}(\epsilon_c+\epsilon_c^{2/5}+\eta_M)
\le CE_{\rm tot}\epsilon_c^{2/5}.
\end{aligned}
\]
The signs of $X,Y$ now give $\theta=\pi+\xi$ and
$\phi=2\xi+d_s$ modulo $2\pi$, where
$|\xi|\le C\epsilon_c^{2/5}$.
The second derivative then gives
\[
\begin{aligned}
|{A_{c,0}}-4{A_{c,1}}|&\le |{A_{c,0}}\cos\xi-4{A_{c,1}}\cos(2\xi+d_s)|
+CE_{\rm tot}(\xi^2+d_s^2)\\
&\le CE_{\rm tot}(\epsilon_c^{4/5}+\eta_M+\xi^2+d_s^2)
\le CE_{\rm tot}\epsilon_c^{4/5}.
\end{aligned}
\]
Thus $\delta_b\triangleq4{A_{c,1}}/{A_{c,0}}-1=O(\epsilon_c^{4/5})$.
Expanding $\widehat f_1={A_{c,0}}\sin\xi-2{A_{c,1}}\sin(2\xi+d_s)+\mathcal E_1$ gives
\[
\begin{aligned}
f_{c,u}/\lambda_c
&=-{A_{c,0}}d_s/2-{A_{c,0}}(4{A_{c,1}}/{A_{c,0}}-1)(\xi+d_s/2)\\
&\quad+O\Bigl(E_{\rm tot}\bigl[|\xi|^3+|d_s|\xi^2+d_s^2|\xi|+|d_s|^3 +m_c^{20}\beta^{-m_c}|D_s|+D_s^2+\epsilon_c^2\bigr]\Bigr).
\end{aligned}
\]
The absorption is explicit:
\[
\frac{|\widehat f_1+{A_{c,0}}d_s/2|}{E_{\rm tot}|D_s|}
\le C\left(\frac{\epsilon_c^{6/5}}{|D_s|}
+\epsilon_c^{4/5}+|D_s|+\frac{\eta_M}{|D_s|}\right)\longrightarrow0.
\]
Since ${A_{c,0}}\ge c'E_{\rm tot}$ and $|d_s|=3|D_s|/4$, this implies
$|\widehat f_1|\ge cE_{\rm tot}|D_s|$. But \eqref{ce:dispersion-exclusion} gives
$|\widehat f_1|\le C_7E_{\rm tot}\epsilon_c$, whereas
$|D_s|/\epsilon_c\to\infty$. This proves the claim, uniformly in time.

In the intermediate range, suppose that both alternatives in
\eqref{ce:window-interface} fail. For sufficiently large $M$, the
lower endpoint of this range ensures $r\ge1$. The first-derivative
bound follows from the opposite characteristic component. Indeed,
\eqref{ce:comparable-flows} gives $|D_qu_o|\le CE_{\rm tot}k/r$, and direct
projection yields
\[
D_qu_o=2H_o+\mathfrak g U_c
+(1+\cos\alpha_v)H_c\pm\sin\alpha_vY_c.
\]
By \eqref{ce:front-interface} and \eqref{ce:comparable-scales},
$|H_o|+|H_c|\le CE_{\rm tot}k\epsilon$ and $|Y_c|\le CE_{\rm tot}k$.
Since $\mathfrak g=1-\cos\alpha_0>0$ is fixed,
\begin{equation}\label{ce:interface-first-derivative-bound}
|U_c|\le C_{\alpha_0}E_{\rm tot}k(r^{-1}+\epsilon+1)
\le C_{\alpha_0}E_{\rm tot}k.
\end{equation}
The constant is independent of $M$, and \eqref{ce:interface-first-derivative-bound} does not require
smallness of $\mathcal C_1$ or $\mathcal C_2$. As $k\asymp k_c$, it also gives
$|f_{c,u_c}|\le CE_{\rm tot}k_c$, as required in
\eqref{ce:dispersion-exclusion}.

For the second- and third-derivative conditions in the claim, the small-$\mathcal C_1$ assumption and the complete Hessian
expansion \eqref{ce:lower-C1-errors} give
\begin{equation}\label{ce:interface-second-derivative-bound}
\begin{aligned}
|h_c|(E_{\rm tot}kr)^2&\le CE_{\rm tot}^3\lambda^2k^2
\{r^{-2}+\epsilon(r^{-1}+r^2)+\epsilon^{2/5}\},\\
\frac{|h_c|}{E_{\rm tot}\lambda^2}
&\le C\{r^{-4}+\epsilon r^{-3}+\epsilon+\epsilon^{2/5}r^{-2}\}
\le C\epsilon^{4/5},
\end{aligned}
\end{equation}
where $r\ge c_0\epsilon^{-1/5}$ was used in the last inequality of \eqref{ce:interface-second-derivative-bound}.
Put $R_H\triangleq Hw-L$. The mixed Hessian bounds in \eqref{ce:front-interface} and \eqref{ce:comparable-flow-table} give
$|R_H|\le CE_{\rm tot}^2\lambda k^2r$, whereas
$|L|\le CE_{\rm tot}^2\lambda^2k(r^{-1}+\epsilon^{4/5}r)$.
Thus
\[
\begin{aligned}
&\mathcal B(L,L)=2\mathfrak g h_oh_c(D_qu_o)(D_qu_c),\\
&\frac{|\mathcal B(Hw,Hw)|}{E_{\rm tot}^4\lambda^4k^2} \le C\bigl[\epsilon^{4/5}
+\epsilon(1+\epsilon^{4/5}r^2)+\epsilon^2r^2\bigr]
\le C\epsilon^{4/5},
\end{aligned}
\]
using $r\le K_1\epsilon^{-2/5}$ and
$\mathcal B(Hw,Hw)=\mathcal B(L,L)+2\mathcal B(L,R_H)+\mathcal B(R_H,R_H)$.
The nonpure third-order terms satisfy
\[
\frac{|\partial^3f_o(w^3)|+
|\partial^3f_c(w^3)-f_{c,uuu}w_{u_c}^3|}
{E_{\rm tot}^4\lambda^4k^2}
\le C\epsilon(r^{-3}+\epsilon r^3)
\le C\epsilon^{4/5}.
\]
Consequently, the identity
$\mathcal C_2=2\partial^3q(w^3)+4\mathcal B(Hw,Hw)$ gives
\[
2|f_{c,uuu}|(E_{\rm tot}kr)^3
\le |\mathcal C_2|+CE_{\rm tot}^4\lambda^4k^2\epsilon^{4/5}
\le CE_{\rm tot}^4\lambda^4k^2\epsilon^{4/5},\qquad
\frac{|f_{c,uuu}|}{E_{\rm tot}\lambda^3}
\le C\epsilon^{-1/5}r^{-3}\le C\epsilon^{2/5}.
\]
These are exactly the forbidden conditions~\eqref{ce:dispersion-exclusion},
proving~\eqref{ce:window-interface}.

\medskip\noindent\textbf{Step 4. The upper growth range.}
In the range of \eqref{ce:window-upper}, if $\mathcal C_1$ is small, the pure
second derivative of $f_c$ is small while $f_c$ stays negative.
The two active frequencies then force a positive fourth derivative.
Its contraction dominates the other terms in $\mathcal C_3$.
In this range, failure of the $\mathcal C_1$ alternative and the full
Hessian sum give
\[
\frac{|h_c|}{E_{\rm tot}\lambda^2}
\le C\bigl(r^{-4}+\epsilon+\epsilon^{2/5}r^{-2}\bigr)
\le C\epsilon,
\qquad r\ge\epsilon^{-2/5}.
\]
To obtain the fourth derivative of the complete solution, put
$P(z)\triangleq(z^2+\lambda_c^2)(z^2+4\lambda_c^2)
=z^4+5\lambda_c^2z^2+4\lambda_c^4$.
For an active packet $e^{\mathrm i\lambda_ju_c}b_j(z_c,y_c)$,
$\lambda_j\in\{\lambda_c,2\lambda_c\}$ and
$\partial_{u_c}z_c=1/2$ give the exact product formula
\[
P(\partial_{u_c})(e^{\mathrm i\lambda_ju_c}b_j)
=e^{\mathrm i\lambda_ju_c}
\sum_{p=1}^4\frac{P^{(p)}(\mathrm i\lambda_j)}{p!2^p}
\partial_{z_c}^pb_j.
\]
The four coefficients can also be read directly, with
$\mu=\lambda_j$, as
\[
\begin{aligned}
P(\partial_{u_c})(e^{\mathrm i\mu u_c}b_j)
=e^{\mathrm i\mu u_c}\bigl[
&\mathrm i\mu(5\lambda_c^2-2\mu^2)\partial_{z_c}b_j
+\tfrac14(5\lambda_c^2-6\mu^2)\partial_{z_c}^2b_j +\tfrac12\mathrm i\mu\partial_{z_c}^3b_j
+\tfrac1{16}\partial_{z_c}^4b_j\bigr].
\end{aligned}
\]
The term without derivatives of $b_j$ vanishes because
$P(\mathrm i\mu)=0$.
For a monomial $\lambda_j^{-\ell}z_c^d
A_{\gamma_c,j}^{(\ell+d)}$ in \eqref{ce:coefficient-formula},
each nonzero term, divided by $E_{\rm tot}\lambda^4$, is at most
\begin{equation}\label{ce:shifted-fixed-monomial}
C(\ell)\epsilon_c^{\ell-d+2p}|D_s|^{d-p}
\le C(\ell)\epsilon_c^2,
\qquad 1\le p\le d\le\ell,\quad |D_s|\le1.
\end{equation}
Apply \eqref{ce:shifted-fixed-monomial} to $\ell<L$, and use
Lemma~\ref{ce:packet-bounds}\ref{cond:packet-tail} through derivative order four with
$B>28\log\beta$ for $\ell\ge L$.
Estimate \eqref{ce:reserve-estimate} handles the periodic correction. Therefore
\[
|P(\partial_{u_c})f_c|
\le CE_{\rm tot}\lambda^4\epsilon_c^2
\le CE_{\rm tot}\lambda^4\tau_c.
\]
Expanding $P$ and using $|f_{c,uu}|\le CE_{\rm tot}\lambda_c^2\epsilon$
and $f_c\le-cE_{\rm tot}$ gives
\[
f_{c,uuuu}\ge 4cE_{\rm tot}\lambda_c^4
-CE_{\rm tot}\lambda_c^4(\epsilon+\tau_c)\ge c'E_{\rm tot}\lambda_c^4
\]
for large $M$. In particular,
\begin{equation}\label{ce:upper-fourth}
f_{c,uuuu}=-5\lambda_c^2f_{c,uu}-4\lambda_c^4f_c
+O(E_{\rm tot}\lambda^4\tau_c)\ge cE_{\rm tot}\lambda^4.
\end{equation}
The positive pure contraction satisfies
\[
f_{c,uuuu}w_{u_c}^4\ge cE_{\rm tot}\lambda^4(E_{\rm tot}kr)^4
=cE_{\rm tot}^5\lambda^4k^4r^4.
\]
The mixed Hessian sums in \eqref{ce:acceleration-component-sums}
do not use the lower-range bound on $h_c$. They give here
\[
|(\varrho_0)_{u_c}|\le CE_{\rm tot}^2\lambda k^2(1+\epsilon r),\qquad
|(\varrho_0)_{\xi_i}|\le CE_{\rm tot}^2\lambda k^2r.
\]
Using the pure components in \eqref{ce:pure-acceleration-components} and $\epsilon r\le C_{\max}$,
we obtain the stronger $u_c$ estimate and the other component bounds:
\begin{equation}\label{ce:upper-acceleration-components}
\begin{aligned}
|\mathfrak a_{u_c}|
&=|\mathfrak g h_ow_{u_o}+(\varrho_0)_{u_c}| \le CE_{\rm tot}^2\lambda^2k(r^{-1}+\epsilon+\epsilon^2r)
\le CE_{\rm tot}^2\lambda^2k/r,\\
|\mathfrak a_{\xi_i}|
&\le CE_{\rm tot}^2\lambda^2k(r^{-1}+\epsilon r)
\qquad\bigl((i,\xi)\ne(c,u)\bigr).
\end{aligned}
\end{equation}
In the last line of \eqref{ce:upper-acceleration-components}, the $h_cw_{u_c}$ terms are bounded by
$(CE_{\rm tot}\lambda^2\epsilon)(E_{\rm tot}kr)=CE_{\rm tot}^2\lambda^2k\epsilon r$.
In particular,
\begin{equation}\label{ce:upper-weighted-acceleration}
\begin{aligned}
\Lambda_c(\mathfrak a)&\le CE_{\rm tot}^2\lambda^2k
\left\{\frac{1+\epsilon+\epsilon^2}{r}
+(\epsilon^2+\epsilon^3)r\right\},\\
\Lambda_i(\mathfrak a)&\le CE_{\rm tot}^2\lambda^2k(r^{-1}+\epsilon r),
\qquad i\in\{o,c\}.
\end{aligned}
\end{equation}
Put $S_U\triangleq E_{\rm tot}^5\lambda^4k^4r^4$.
For $\alpha=(a,b,h)$, write
$f_{i,\alpha}\triangleq\partial_{u_i}^a\partial_{v_i}^b
\partial_{y_i}^hf_i$ and
$w_i^\alpha\triangleq w_{u_i}^aw_{v_i}^bw_{y_i}^h$.
The two fourth-order sums are then bounded by explicit polynomials:
\begin{equation}\label{ce:upper-fourth-sums}
\begin{aligned}
\frac{\displaystyle\sum_{\substack{|\alpha|=4\\\alpha\ne(4,0,0)}}
\frac{4!}{\alpha!}|f_{c,\alpha}w_c^\alpha|}{S_U}
&\le C\bigl\{(1+\epsilon+\epsilon^2)^4-1\bigr\}
\le C\epsilon,\\
\frac{\displaystyle\sum_{|\alpha|=4}\frac{4!}{\alpha!}
|f_{o,\alpha}w_o^\alpha|}{S_U}
&\le Cr^{-8}\bigl\{1+\epsilon r+(\epsilon r)^2\bigr\}^4
\le Cr^{-8}.
\end{aligned}
\end{equation}
For the third-order and Hessian sums, \eqref{ce:weighted-contractions} gives
\begin{equation}\label{ce:upper-third-hessian-sums}
\begin{aligned}
\frac{|\partial^3f_c(\mathfrak a,w^2)|}{S_U}
&\le\frac{CE_{\rm tot}\lambda^3\Lambda_c(\mathfrak a)\Lambda_c(w)^2}{S_U}
\le C\left\{\frac1{\epsilon r^3}+\frac{1+\epsilon}{r^3}
+\frac{\epsilon+\epsilon^2}{r}\right\},\\
\frac{|\partial^3f_o(\mathfrak a,w^2)|}{S_U}
&\le\frac{CE_{\rm tot}\lambda^3\Lambda_o(\mathfrak a)\Lambda_o(w)^2}{S_U}
\le C\left\{\frac1{\epsilon r^7}+\frac1{r^5}\right\},\\
\frac{|H(\mathfrak a,\mathfrak a)|}{S_U}
&\le\frac{C\sum_{i=o,c}E_{\rm tot}\lambda^2\Lambda_i(\mathfrak a)^2}{S_U}
\le C\frac{(r^{-1}+\epsilon r)^2}{\epsilon^2r^4} =C\left\{\frac1{\epsilon^2r^6}+\frac2{\epsilon r^4}+\frac1{r^2}\right\}.
\end{aligned}
\end{equation}
Including the fixed coefficients $2,14,8$ from
\eqref{ce:chain-expansion}, \eqref{ce:upper-fourth-sums} and \eqref{ce:upper-third-hessian-sums} give the following summary relative to $S_U$:
\begin{equation}\label{ce:upper-budget}
\begin{array}{c|c}
\text{nonpure fourth derivatives of }f_c&C\epsilon\\
\text{fourth derivatives of }f_o&Cr^{-8}\\
\text{third derivatives of }f_c&C\cdot((\epsilon r^3)^{-1}+r^{-3}+\epsilon/r)\\
\text{third derivatives of }f_o&C\cdot((\epsilon r^7)^{-1}+r^{-5})\\
H(\mathfrak a,\mathfrak a)&C\cdot((\epsilon^2r^6)^{-1}+(\epsilon r^4)^{-1}+r^{-2}).
\end{array}
\end{equation}
For $r\ge\epsilon^{-2/5}$, the potentially largest ratios are
\[
\frac1{\epsilon r^3}\le\epsilon^{1/5},\quad
\frac1{\epsilon^2r^6}\le\epsilon^{2/5},\quad
\frac1{\epsilon r^4}\le\epsilon^{3/5},\quad
r^{-2}\le\epsilon^{4/5}.
\]
Every term in \eqref{ce:upper-budget} is therefore bounded by
$C\epsilon^{1/5}$. It is therefore less than half the coefficient of
the positive fourth-order term for large $M$. This proves
\eqref{ce:window-upper} and completes the classification.
\end{proof}

The ranges in Lemmas~\ref{ce:small-alternatives}
and~\ref{ce:comparable-windows} include their endpoints and overlap
where they meet. The bounds of Lemmas~\ref{ce:small-alternatives} and~\ref{ce:comparable-windows} together give \eqref{ce:incoming-chain-target}.

\subsubsection{Combining the estimates}\label{ce:incoming-repair}
\begin{proof}[Proof of Proposition~\ref{ce:incoming}]
Conditions~\ref{cond:merging-amplitude}--\ref{cond:merging-potential} of
Proposition~\ref{ce:merging-criterion} were verified in
\eqref{ce:incoming-amplitude-bounds}--\eqref{ce:incoming-paid}.
It remains to combine Lemmas~\ref{ce:small-alternatives} and~\ref{ce:comparable-windows} to prove
\eqref{ce:incoming-chain-target}. Put
$\mathcal U\triangleq\mathcal U_1\cup\mathcal U_2$.

\medskip\noindent\textbf{Step 1. Put all derivative orders on the same amplitude scale.}\proofsteplabel{step:incoming-weights}{1}
At a point of $\mathcal U$, let $o$ be the unique index with
$(t,x)\in\mathcal U_o$, and put $c\triangleq3-o$.
Lemma~\ref{ce:source-owner} gives $E_o\le K_0E_c$, so
$E_o\le\varepsilon_2E_c$ and
$\varepsilon_2E_c\le E_o\le K_0E_c$ cover $\mathcal U$.
To use one denominator estimate in both cases, set
\[
E_{\rm tot}\triangleq E_1+E_2,\qquad m\triangleq\min(m_1,m_2).
\]
Under a constant rescaling $q\mapsto a q$, one has
$D_{a q}=aD_q$ and $Q(aq)=a^2Q(q)$, hence
\begin{equation}\label{ce:incoming-homogeneity}
\mathcal C_l(aq)=D_{a q}^{\,l}Q(aq)=a^{l+2}\mathcal C_l(q).
\end{equation}
Thus $\mathcal C_0,\mathcal C_1,\mathcal C_2,\mathcal C_3$ have different amplitude degrees.
Dividing $\mathcal C_l$ by a weight of amplitude degree $l$ puts all four
terms at degree two, without discarding their frequency factors. Set
\begin{equation}\label{ce:incoming-weights}
\mathcal W\triangleq E_1y_1^{-20}+E_2y_2^{-20},\qquad
F_{\mathrm w}\triangleq\sum_{l=0}^3\frac{\mathcal C_l^2}{\mathcal W^{2l}}.
\end{equation}
If $|m_1-m_2|\le1$, the band definitions give
$y_i^{-20}\asymp\beta^{20m}$ and hence
$\mathcal W\asymp E_{\rm tot}\beta^{20m}$.
For a larger index gap, label the larger band index by $h$ and the smaller
by $\ell$, and write $m_h=m+d$, $m_\ell=m$, with $d\ge2$.
By \eqref{ce:envelope} and $\gamma_i\in\{0,\gamma\}$,
\begin{equation}\label{ce:incoming-separated-envelope}
\frac{E_h}{E_\ell}\beta^{20d}
\le \exp\!\left[-(1-\beta^{-1})\beta^{m+d}
+\gamma(m+1)+20d\log\beta\right]\le1
\end{equation}
with \eqref{ce:incoming-separated-envelope} valid uniformly in $d\ge2$ for sufficiently large $m$.
Since $\beta^{20m_i}\le y_i^{-20}\le\beta^{20(m_i+1)}$,
\[
E_{\rm tot}\beta^{20m}\le\mathcal W
\le\beta^{20}\beta^{20m}(E_\ell+E_h\beta^{20d})
\le2\beta^{20}E_{\rm tot}\beta^{20m}.
\]
\medskip\noindent\textbf{Step 2. Insert the bounds from the two cases.}\proofsteplabel{step:incoming-cases}{2}
We divide each alternative by the weight for its derivative order. In the small-amplitude case, $E_{\rm tot}\asymp E_c$ and $m_c-1\le m\le m_c$.
When the $\mathcal C_1$ alternatives in \eqref{ce:small-effective-low} and
\eqref{ce:small-effective-high} are used, respectively, their normalized
bounds give
\[
\begin{aligned}
\frac{|\mathcal C_1|}{\mathcal W}
&\ge cE_{\rm tot}^2k_c^4\nu\beta^{-20m}
\asymp E_{\rm tot}^2\beta^{8m}(m^2\beta^{-m})\beta^{-20m}
=E_{\rm tot}^2m^2\beta^{-13m},\\
\frac{|\mathcal C_1|}{\mathcal W}
&\ge cE_{\rm tot}^2k_c^4\nu^2\beta^{-20m}
\asymp E_{\rm tot}^2\beta^{8m}(m^4\beta^{-2m})\beta^{-20m}
=E_{\rm tot}^2m^4\beta^{-14m}.
\end{aligned}
\]
Here $k_c\asymp\beta^{2m}$ and $\nu\asymp m^2\beta^{-m}$;
the second bound also uses $|J_*|\ge cE_{\rm tot}k_c\nu$.
For the remaining small-amplitude alternatives, the respective lower bounds are
\begin{equation}\label{ce:small-normalized-alternatives}
\begin{cases}
|Q| \ge cE_{\rm tot}^2\beta^{4m}\nu^2\ge cE_{\rm tot}^2,\\[2pt]
\frac{|\mathcal C_3|}{\mathcal W^3}
\ge cE_{\rm tot}^2\beta^{12m}\nu^2\ge cE_{\rm tot}^2,\\[2pt]
\frac{|\mathcal C_1|}{\mathcal W}
\ge cE_{\rm tot}^2 k_c^2\lambda_o\lambda_c\beta^{-20m}
\ge cE_{\rm tot}^2\beta^{16m}.
\end{cases}
\end{equation}
The last line of \eqref{ce:small-normalized-alternatives} concerns the large-time-derivative branch and uses
\eqref{ce:cross-required} and $\lambda_o\ge\lambda_c/2$.
In that branch the $Q$ alternative also gives $|Q|\ge cE_{\rm tot}^2\beta^{4m}$,
without using $\nu$.
For comparable amplitudes, every frequency factor is converted by
the same rule:
\begin{equation}\label{ce:incoming-frequency-conversion}
\begin{aligned}
|\mathcal C_l|\ge cE_{\rm tot}^{l+2}\lambda^p k^q\epsilon^\sigma r^d
\quad\Longrightarrow\quad
\frac{|\mathcal C_l|}{E_{\rm tot}^2\mathcal W^l}
&\ge c\beta^{(16p+2q-20l)m}\epsilon^\sigma r^d,
\end{aligned}
\end{equation}
because $\lambda\asymp\beta^{16m}$, $k\asymp\beta^{2m}$,
and $\mathcal W\asymp E_{\rm tot}\beta^{20m}$. In particular,
\begin{equation}\label{ce:incoming-frequency-powers}
\begin{aligned}
\lambda^2k^2\beta^{-20m} \asymp\beta^{16m},\quad
\lambda^4k^2\beta^{-40m} \asymp\beta^{28m},\quad
\lambda^6k^2\beta^{-60m} \asymp\beta^{40m}, \quad
\lambda^4k^4\beta^{-60m} \asymp\beta^{12m}.
\end{aligned}
\end{equation}
The conversions \eqref{ce:incoming-frequency-conversion} and \eqref{ce:incoming-frequency-powers} applied to \eqref{ce:window-bounded}--\eqref{ce:window-upper} give the following lower factors for
$\max_{0\le l\le3}|\mathcal C_l|/(E_{\rm tot}^2\mathcal W^l)$:
\begin{equation}\label{ce:incoming-normalized-factors}
\begin{array}{c|c}
\text{range}&\text{lower factor, up to a fixed positive constant}\\\hline
0<r\le1&\min\{\beta^{16m},\beta^{28m}\}\\
1\le r\le c_{\rm low}\epsilon^{-1/5}
&\min\{\beta^{16m}r^{-2},\beta^{28m}r^{-4},\beta^{40m}r^{-6}\}\\
c_0\epsilon^{-1/5}\le r\le K_1\epsilon^{-2/5}
&\min\{\beta^{16m}\epsilon^{2/5},\beta^{28m}\epsilon^{4/5}\}\\
\epsilon^{-2/5}\le r\le C_{\max}/\epsilon
&\min\{\beta^{16m}\epsilon^{2/5},\beta^{12m}r^4\}.
\end{array}
\end{equation}
For the second row of \eqref{ce:incoming-normalized-factors}, $r\le c_{\rm low}\epsilon^{-1/5}$ and
$\epsilon\asymp\beta^{-14m}$ give the three endpoint bounds
\[
\begin{aligned}
\beta^{16m}r^{-2}&\ge c\beta^{16m}\epsilon^{2/5}
\asymp\beta^{52m/5},\\
\beta^{28m}r^{-4}&\ge c\beta^{28m}\epsilon^{4/5}
\asymp\beta^{84m/5},\\
\beta^{40m}r^{-6}&\ge c\beta^{40m}\epsilon^{6/5}
\asymp\beta^{116m/5}.
\end{aligned}
\]
The two factors in the third row of \eqref{ce:incoming-normalized-factors} are likewise comparable to
$\beta^{52m/5}$ and $\beta^{84m/5}$. For the fourth row of \eqref{ce:incoming-normalized-factors},
\[
\beta^{12m}r^4\ge\beta^{12m}\epsilon^{-8/5}
\asymp\beta^{172m/5},\qquad
\beta^{16m}\epsilon^{2/5}\asymp\beta^{52m/5}.
\]
The first row of \eqref{ce:incoming-normalized-factors} has minimum $\beta^{16m}$, so all four row minima are
bounded below, respectively, by
\[
c\beta^{16m},\qquad c\beta^{52m/5},\qquad
c\beta^{52m/5},\qquad c\beta^{52m/5}.
\]
Outside the small-$Q$ range of Lemma~\ref{ce:comparable-windows},
$|Q|\ge\varepsilon_QE_{\rm tot}^2k^2\ge cE_{\rm tot}^2\beta^{4m}$.
Thus, on $\mathcal U$,
\[
\max_{0\le l\le3}\frac{|\mathcal C_l|}{\mathcal W^l}
\ge cE_{\rm tot}^2\beta^{-14m}.
\]
The common denominator estimate is therefore
\begin{equation}\label{ce:incoming-denominator}
\begin{aligned}
F_{\mathrm w}=\sum_{l=0}^3\frac{\mathcal C_l^2}{\mathcal W^{2l}}
&\ge\left(\max_{0\le l\le3}\frac{|\mathcal C_l|}{\mathcal W^l}\right)^2 \ge cE_{\rm tot}^4\beta^{-28m}\ge cE_{\rm tot}^4e^{-Cm}>0.
\end{aligned}
\end{equation}

\begin{center}\small
\begin{tabular}{@{}p{.25\textwidth}p{.68\textwidth}@{}}
\toprule Object&Amplitude degree and retained scale\\\midrule
$\mathcal C_l$&homogeneous of degree $l+2$, by the scaling identity \eqref{ce:incoming-homogeneity}\\
$\mathcal C_l/\mathcal W^l$&degree two; the weakest lower factor is $\beta^{-14m}$\\
$F_{\rm w}$&squaring gives $E_{\rm tot}^4\beta^{-28m}$ in \eqref{ce:incoming-denominator}\\
$F_0$&\eqref{ce:incoming-unweighting} costs $E_{\rm tot}^6$; \eqref{ce:incoming-amplitude-power} gives $E_{\rm tot}^{11}$\\\bottomrule
\end{tabular}
\end{center}
\medskip\noindent\textbf{Step 3. Remove the weights and obtain the power $11$.}\proofsteplabel{step:incoming-power}{3}
Shrink the vertex chart by increasing $M$ so that $E_{\rm tot}=E_1+E_2\le1$
and $y_i\le1$. Then
\[
\mathcal W\ge E_{\rm tot},\qquad
\mathcal W^{2l}\ge E_{\rm tot}^{2l}\ge E_{\rm tot}^6\quad(0\le l\le3),
\]
and \eqref{ce:incoming-denominator} gives
\begin{equation}\label{ce:incoming-unweighting}
F_0=\sum_{l=0}^3\mathcal C_l^2
=\sum_{l=0}^3\mathcal W^{2l}\frac{\mathcal C_l^2}{\mathcal W^{2l}}
\ge E_{\rm tot}^6F_{\mathrm w}\ge cE_{\rm tot}^{10}\beta^{-28m}.
\end{equation}
By \eqref{ce:envelope},
$E_{\rm tot}\le2e^{-\beta^m}\le C\beta^{-28m}$. Consequently,
\begin{equation}\label{ce:incoming-amplitude-power}
E_{\rm tot}^{10}\beta^{-28m}\ge C^{-1}E_{\rm tot}^{11}.
\end{equation}
Combining \eqref{ce:incoming-unweighting} and \eqref{ce:incoming-amplitude-power} proves
\begin{equation}\label{ce:incoming-unweighted}
F_0\ge cE_{\rm tot}^{11}\ge cE_i^{11}\quad\hbox{on }\mathcal U_i.
\end{equation}

\medskip\noindent\textbf{Step 4. Apply the merging proposition.}\proofsteplabel{step:incoming-apply}{4}
By the definition of $\mathcal U_i$ in \eqref{ce:incoming-supports},
\eqref{ce:incoming-unweighted} proves
\eqref{ce:incoming-chain-target}. Thus all the hypotheses of
Proposition~\ref{ce:merging-criterion} hold with
$a=1/8$, $b=1/4$, $c=5/8$, $\vartheta=\chi$, $r=3$, and $\kappa=11$.
The proposition gives a real smooth $2\pi$-periodic pair
$(p,V)$ on $\R\times U$, with $p=e^{\Phi_0}(f_1+f_2)$, satisfying
the wave equation and the required matching identities.

At $x=v$, the original pairs are flat. Estimate
\eqref{ce:merging-correction-flatness} and $y_i\le|x-v|$ on $D$
show that the corrections are flat there as well. Hence $p,V$ are
flat at $v$. Finally, smoothness, periodicity, and
$\overline{C_M}\subset U$ give the asserted bounds for every fixed
order of derivatives on $\R\times\overline{C_M}$.
\end{proof}

\subsection{Completion to the constant solution}\label{ce:activation-verification}
The pairs from Subsection~\ref{ce:vertex-verification} agree with the original edge solutions near the vertices. We now choose a smooth entrance curve outside each observation polygon and continue the pairs to the constant solution $(p,V)=(1,0)$.

The local solutions will be placed along the polygon edges after a spatial
rescaling. Choose $0<\delta<\delta_*$ as in
Proposition~\ref{ce:geometry}\ref{cond:geometry-escape}. By
Proposition~\ref{ce:geometry}\ref{cond:geometry-polygons}, the regular $64$-gons $P_j^0$
have a common positive edge length, denoted by $L_0$.
Choose $d_0>0$ smaller than $L_0\sin\alpha_0$,
the distance of each polygon from $\partial\Omega_0$, and the distances
between distinct polygons.

Recall $\beta=2^{1/16}$ from Section~\ref{ce:front-section}.
For $S>0$ and an integer $M\ge4$, define the reference normal scale
and the spatial scale by
\begin{equation}\label{ce:final-scales}
\begin{aligned}
&y_*\triangleq\beta^{-M},\qquad {r_{\rm geo}}\triangleq Sy_*,\qquad
\Omega\triangleq {r_{\rm geo}}\Omega_0,\qquad P_j\triangleq {r_{\rm geo}}P_j^0,
\qquad\omega\triangleq {r_{\rm geo}}\omega_0.
\end{aligned}
\end{equation}
We use the family of local solutions on $|z|<1$, with normal range
$|y|<Y_*$. The following conditions leave geometric room relative
to the reference normal scale while keeping the whole domain within that range.

\begin{assumption}[Scaling]\label{ce:geometric-scales}
The scale parameters satisfy
\begin{equation}\label{ce:geometric-scale-conditions}
\frac{32}{S}<\frac{d_0\sin\alpha_0}{10^4},\qquad
{r_{\rm geo}}\operatorname{diam}\Omega_0<\min\{1/100,Y_*\}.
\end{equation}
\end{assumption}

Choose $S$ to satisfy the first condition in \eqref{ce:geometric-scale-conditions} and keep it fixed throughout
this section. The second condition in \eqref{ce:geometric-scale-conditions} holds for every sufficiently large
$M$; the final common choice is made in Lemma~\ref{ce:global-local},
with the explicit scale restrictions collected in \eqref{ce:static-margins}.
These choices are compatible because increasing $M$
decreases ${r_{\rm geo}}=S\beta^{-M}$ without changing $S$.
The scaled pair $(\Omega,\omega)$ satisfies GCC, and
\eqref{ce:free-scaled} applies to it.

\medskip\noindent\textbf{The geometric lengths and the three kinds of charts.}
Fix reference widths
\begin{equation}\label{ce:reference-collar-scales}
\rho_g\triangleq d_0/100,\qquad
\frac{32}{S}<\delta_g<\frac{\rho_g\sin\alpha_0}{100},
\qquad \rho_p\triangleq {r_{\rm geo}}\rho_g,\quad\delta_p\triangleq {r_{\rm geo}}\delta_g.
\end{equation}
The interval for $\delta_g$ is nonempty because its upper endpoint is
$d_0\sin\alpha_0/10^4$. All reference data, including $S,\rho_g,\delta_g$,
are fixed before increasing $M$. The lengths used below are
\begin{center}
\begin{tabular}{c@{\qquad}l@{\qquad}l}
\text{symbol}&\text{size}&\text{role}\\\hline
$y_*$&$\beta^{-M}$&height of the principal packet\\
${r_{\rm geo}}$&$Sy_*$&spatial scale of the polygons\\
$\rho_g$&$d_0/100$&fixed reference geometric width\\
$\rho_p$&$S\rho_g y_*$&half the radius of a vertex chart\\
$\delta_p$&$S\delta_g y_*$&width of the exterior collar\\
$\ell$&$y_*^2$&height scale of the entrance level function
\end{tabular}
\end{center}
In particular,
\[
32y_*<\delta_p<\frac{\rho_p\sin\alpha_0}{100},\qquad
4\rho_p<\frac{{r_{\rm geo}}L_0}{10},\qquad \delta_p<\frac{{r_{\rm geo}}L_0}{10},
\]
and
\[
\frac{\ell}{y_*}=y_*\longrightarrow0,\qquad
\frac{M\ell}{\rho_p}=\frac{My_*}{S\rho_g}\longrightarrow0.
\]
Thus the entrance height fits well inside the exterior collar, while
the bends of the vertex curve fit inside a vertex chart.
We use the zero pair in the polygon interiors, the original local pair
near an edge, and the merged pair near a vertex. Step~\ref{step:activation-global-data}
will check that these three definitions agree on their overlaps.

Use the scaled polygons $P_j$ in \eqref{ce:final-scales}.
Orient each polygon clockwise.
On each polygon, label its edges by $i\in\{1,\ldots,64\}$.
On every whole edge choose its midpoint $c_i$ and unit tangent $e_i$
in that orientation, and set $n_i=\mathsf Je_i$, where
$\mathsf J(x_1,x_2)=(-x_2,x_1)$. Thus $n_i$ is the outward normal. Assign types $0$
and $\gamma=24\log2$ alternately around each polygon. The even
number of edges makes this assignment consistent. Use the same
midpoint, frame, and solution type at both endpoints of an edge.
At a vertex, relabel its two incident edges by $i\in\{1,2\}$ so that $\gamma_1=0$ and
$\gamma_2=\gamma$. Near that vertex define
\begin{equation}\label{ce:activation-edge-coordinates}
\begin{aligned}
&y_i\triangleq n_i\cdot(x-v),\quad z_i\triangleq e_i\cdot(x-c_i),\quad
Z_i\triangleq e_i\cdot(v-c_i),\\
&u_i\triangleq t+z_i,\quad v_i\triangleq t-z_i,\quad
(f_i,W_i)(t,x)\triangleq(f_{\gamma_i},V_{\gamma_i})(t,z_i,y_i).
\end{aligned}
\end{equation}
For $x,c_i,v\in\Omega$, \eqref{ce:geometric-scale-conditions} gives
\[
|e_i\cdot(x-c_i)|\le {r_{\rm geo}}\operatorname{diam}\Omega_0<1/100,\qquad
|n_i\cdot(x-v)|\le {r_{\rm geo}}\operatorname{diam}\Omega_0<Y_*.
\]
Thus the rescaled geometry stays inside the domains of the local solutions.
The edge coordinates agree with those in Subsection~\ref{ce:vertex-verification}:
\[
z_i=Z_i+e_i\cdot(x-v),\qquad
Z_i=\pm\frac{{r_{\rm geo}}L_0}{2}=\zeta_i y_*,\qquad
\zeta_i=\pm\frac{SL_0}{2}.
\]
Choose the vertex neighborhoods by $C_v\triangleq B(v,2\rho_p)$.
With $R\triangleq2S\rho_g$, they have the form $B(v,Ry_*)$, and
\[
R=2S\rho_g<SL_0/2=|\zeta_i|.
\]
These are exactly the coordinate and radius conditions required by
Proposition~\ref{ce:incoming}.
The larger closed balls of radius $4\rho_p$ are mutually disjoint
and contained in $\Omega$, and $4\rho_p<{r_{\rm geo}}L_0/10$.
Assumption~\ref{ce:geometric-scales} makes
$\rho_p/y_*=S\rho_g$ sufficiently large.
On each vertex chart, denote the pair $(p,V)$ supplied by
Proposition~\ref{ce:incoming} by $(p_{\rm in},V_{\rm in})$;
on an edge chart this notation denotes the original solution-potential pair.
Here and below, $C=C(S)>0$ may depend on $S$ and the fixed
reference geometry, but not on $M$.

\begin{proposition}[Completion to a constant solution]\label{ce:activation}
For every sufficiently large $M$, there are an open neighborhood $U$ of
the observation polygons, a function $X\in C^\infty(U;\R)$, and a real
smooth $2\pi$-periodic pair $(p,V)$ on $\R\times U$ such that
\[
\Box p+Vp=0,\qquad
(p,V)=(p_{\rm in},V_{\rm in})\quad(X\le0),\qquad
(p,V)=(1,0)\quad(X\ge1).
\]
Here $(p_{\rm in},V_{\rm in})$ is the smooth pair obtained by joining the
original edge pairs, the vertex pairs of Proposition~\ref{ce:incoming},
and the zero pair inside the polygons.
Moreover, $X<0$ on the polygons and
$\{-1/2\le X\le3/2\}\Subset U$.
The functions can be defined on slightly larger neighborhoods of
$\overline U$, and their derivatives of every fixed order are bounded
there uniformly in time.
\end{proposition}

\begin{proof}[Proof of Proposition~\ref{ce:activation}]
We apply Proposition~\ref{ce:activation-criterion} with $A\equiv1$.
We must construct a smooth spatial function $X$ such that
\begin{equation}\label{ce:activation-entrance-requirements}
\{-1/2\le X\le3/2\}\Subset U,\qquad
\mathcal B(\mathrm dp_{\rm in},\mathrm dp_{\rm in})<0
\quad\hbox{where }X=0,\ p_{\rm in}=0.
\end{equation}
We first prove the negativity in \eqref{ce:activation-entrance-requirements} along edge and vertex curves, then join
the local pairs and their level functions into global data.

\medskip\noindent\textbf{Step 1. A negative quadratic form for one local solution.}\proofsteplabel{step:activation-gradient}{1}
At an edge we will place the entrance curve close to the height
$y_*=\beta^{-M}$. 
Set
\begin{equation}\label{ce:activation-scales}
\ell\triangleq y_*^2,\quad
k\triangleq\ell^{-1},\quad\lambda\triangleq2^M,\quad
\epsilon\triangleq k/\lambda=\beta^{-14M},\quad
\tau\triangleq k^2/\lambda.
\end{equation}
Indeed, $\lambda=\beta^{16M}$, so $k=\beta^{2M}$ and
$\tau=\beta^{-12M}$. In the region
$y_i=y_*+O(M\ell)$, $|z_i|\le C(S)y_*$, write
\[
B_i\triangleq\begin{cases}
e^{-1/y_i-M^2-\gamma_iM}\rho(\beta^My_i),&y_i>0,\\
0,&y_i\le0,
\end{cases}\quad i\in\{1,2\},\qquad B_{\rm tot}\triangleq B_1+B_2.
\]
On an edge chart set $B_{\rm tot}\triangleq B_i$. Here
$B_i=A_{\gamma_i,M}(y_i)$ is the amplitude of the $M$th packet and
$k_i\triangleq y_i^{-2}$ satisfies $k_i/k=1+O(My_*)$.

\smallskip\noindent\textit{Claim.}
Uniformly on the region $y_i=y_*+O(M\ell)$, $|z_i|\le C(S)y_*$, each complete local solution satisfies
\[
\mathcal B(\mathrm df_i,\mathrm df_i)
=-B_i^2k_i^2+o(B_i^2k^2)
\le-\tfrac12 B_i^2k_i^2<0
\]
for sufficiently large $M$. The quadratic-form estimate in \eqref{ce:carrier-estimates} holds at every time and every oscillatory phase. 

\smallskip\noindent\textit{Amplitude and phase of the principal packet.}
For $j=M$, the coefficient of $e^{\mathrm i\lambda u_i}$ in
\eqref{ce:approximate-solution} is
\begin{equation}\label{ce:principal-coefficients}
\mathcal C_i\triangleq
\sum_{\nu=0}^{K_M}\lambda^{-\nu}H_{\gamma_i,M,\nu}
=B_i+\frac{\mathrm i z_iB_i''}{2\lambda}
-\frac{z_i^2B_i^{(4)}}{8\lambda^2}
+\sum_{\nu=3}^{K_M}\lambda^{-\nu}H_{\gamma_i,M,\nu}.
\end{equation}
Here $K_M$ is defined in \eqref{ce:front-parameters}. The first three
coefficients follow individually from the recursion in
\eqref{ce:approximate-solution}:
\[
\begin{aligned}
&H_{\gamma_i,M,0}=B_i,\qquad
H_{\gamma_i,M,1}
=\frac{\mathrm i}{2}\int_0^{z_i}B_i''\,\mathrm d\sigma
=\frac{\mathrm i z_i}{2}B_i'',\\
&H_{\gamma_i,M,2}
=\frac{\mathrm i}{2}\int_0^{z_i}
(\partial_\sigma^2+\partial_{y_i}^2)
\left(\frac{\mathrm i\sigma}{2}B_i''\right)\mathrm d\sigma
=-\frac{z_i^2}{8}B_i^{(4)}.
\end{aligned}
\]
Define
\begin{equation}\label{ce:principal-phase-definition}
A_i\triangleq|\mathcal C_i|,\qquad
\psi_i\triangleq\operatorname{Arg}\mathcal C_i\in(-\pi/2,\pi/2),\qquad
\theta_i\triangleq\lambda(t+z_i)+\psi_i,\qquad
r_i\triangleq f_i-A_i\cos\theta_i.
\end{equation}
Thus $B_i$ is the original packet amplitude, $\mathcal C_i$ includes its
finite correction series, $A_i,\psi_i$ are the resulting real amplitude
and phase correction, and $r_i$ contains the other packets and the
correction to the complete solution. We will show
$\operatorname{Re}(\mathcal C_i/B_i)>0$, which makes the branch of $\psi_i$ specified in \eqref{ce:principal-phase-definition} smooth.

With $U_i=f_{i,u_i}$, $H_i=f_{i,v_i}$, and $Y_i=f_{i,y_i}$,
these are
\[
U_i\approx-\lambda B_i\sin\theta_i,\qquad
H_i\approx\frac{B_i''}{4\lambda}\sin\theta_i,\qquad
Y_i\approx k_iB_i\cos\theta_i.
\]
Since $B_i''\approx k_i^2B_i$, the terms $4U_iH_i$ and $-Y_i^2$
give the negative $\sin^2\theta_i$ and $\cos^2\theta_i$ contributions,
respectively. To control the errors in the quadratic form, we first establish the following estimates:
\begin{equation}\label{ce:carrier-estimates}
\begin{aligned}
&A_i/B_i=1+o(1),\qquad
|\partial_{z_i}A_i|+|\partial_{y_i}A_i-k_iA_i|
\le C(S)B_ik(\epsilon^2/y_*+\epsilon^3),\\
&\partial_{z_i}\psi_i=\frac{B_i''}{2\lambda B_i}+o(\tau),
\qquad |\partial_{y_i}\psi_i|\le C\tau,\qquad
|r_i|+\lambda^{-1}|\mathrm{d}r_i|\le CB_i\epsilon^8,\\
&\mathcal B(\mathrm{d}f_i,\mathrm{d}f_i)=-B_i^2k_i^2+o(B_i^2k^2).
\end{aligned}
\end{equation}
The first line controls differentiation of the amplitude. In the
second line, $\partial_{z_i}\psi_i$ supplies the leading term in $H_i$,
while the other two estimates control the normal phase derivative and
the omitted terms. The last line of \eqref{ce:carrier-estimates} will follow from its first two lines. Every little-oh estimate is uniform in $t$ and on the
stated spatial sets.

\smallskip\noindent\textit{Derivatives of the coefficient series.}
The derivatives of $B_i$ come from differentiating its exponential and
cutoff. For example,
\[
(e^{-1/y_i})'=k_i e^{-1/y_i},\qquad
(e^{-1/y_i})''=(k_i^2-2y_i^{-3})e^{-1/y_i},\qquad
\partial_{y_i}\rho(\beta^My_i)=\beta^M\rho'(\beta^My_i).
\]
Here $\beta^My_i=1+O(My_*)$. Flatness of $\rho$ at its peak in
\eqref{ce:cutoff-functions}, followed by Leibniz's rule and
\eqref{ce:exponential-exact-derivative}, gives, for each fixed $r\ge0$,
\begin{equation}\label{ce:principal-amplitude-derivatives}
\frac{B_i^{(r)}}{B_i}=k_i^r[1+O(y_*)],
\qquad
\left|\partial_{y_i}\frac{B_i^{(r)}}{B_i}\right|
\le C(r)k_i^r/y_*.
\end{equation}
Indeed, every term containing a derivative of the peak cutoff, and
every neighboring packet with any fixed number of derivatives, has
relative size at most
$\exp[-c/(My_*)+C(r)M]$. This is smaller than any fixed power of
$\epsilon$. In particular,
\begin{equation}\label{ce:principal-envelope-comparison}
\frac{E_{\gamma_i}(y_i)}{B_i}
=1+O\!\left(\exp[-c/(My_*)+CM]\right)\le2.
\end{equation}

For integers $1\le p\le\nu$ with $\nu-p$ even, each monomial in
\eqref{ce:coefficient-formula} satisfies
\begin{equation}\label{ce:principal-monomial-bounds}
\begin{aligned}
&\frac{|\lambda^{-\nu}z_i^pB_i^{(\nu+p)}|}{B_i}
\le C(\nu)\lambda^{-\nu}|z_i|^pk^{\nu+p}
=C(\nu)\epsilon^{\nu-p}|z_i\tau|^p,\\
&\frac{|\partial_{z_i}(\lambda^{-\nu}z_i^pB_i^{(\nu+p)})|}{B_i\tau}
\le C(\nu)\epsilon^{\nu-p}|z_i\tau|^{p-1},\\
&\left|\partial_{y_i}\!\left(
\frac{\lambda^{-\nu}z_i^pB_i^{(\nu+p)}}{B_i}\right)\right|
\le\frac{C(\nu)}{y_*}\epsilon^{\nu-p}|z_i\tau|^p,
\qquad |z_i\tau|\le C(S)\epsilon/y_*\longrightarrow0.
\end{aligned}
\end{equation}
Estimate \eqref{ce:principal-monomial-bounds} treats any fixed finite set of indices $\nu$. For the
remaining indices choose $L\ge3$, independent of $M$, and write
\[
\mathcal C_{\rm tail}\triangleq
\sum_{\nu=L}^{K_M}\lambda^{-\nu}H_{\gamma_i,M,\nu}.
\]
For any prescribed $B_0>0$, \eqref{ce:coefficient-estimate} permits
$L$ large enough that this tail and its first derivatives, after
multiplication by $e^{\mathrm i\lambda u_i}$, are bounded by
$CB_i e^{-B_0M}$. To recover coefficient derivatives, use
\begin{equation}\label{ce:principal-tail-identities}
\begin{aligned}
&\partial_{z_i}\mathcal C_{\rm tail}
=e^{-\mathrm i\lambda u_i}(\partial_{z_i}-\partial_t)
(e^{\mathrm i\lambda u_i}\mathcal C_{\rm tail}),\\
&\partial_{y_i}\!\left(\frac{\mathcal C_{\rm tail}}{B_i}\right)
=\frac{\partial_{y_i}\mathcal C_{\rm tail}}{B_i}
-\frac{B_i'}{B_i}\frac{\mathcal C_{\rm tail}}{B_i}.
\end{aligned}
\end{equation}
Since $B_i'/B_i=O(k)$, \eqref{ce:principal-tail-identities} implies
\[
\frac{|\mathcal C_{\rm tail}|+|\partial_{z_i}\mathcal C_{\rm tail}|}{B_i}
+k^{-1}\left|\partial_{y_i}
\left(\frac{\mathcal C_{\rm tail}}{B_i}\right)\right|
\le C e^{-B_0M}.
\]
There is enough freedom in $B_0$ to absorb each fixed power loss:
\[
e^{-B_0M}k^p\tau^{-b}\epsilon^{-c}
=\exp\!\left[\bigl(-B_0+(2p+12b+14c)\log\beta\bigr)M\right]
\longrightarrow0
\]
for any fixed $p,b,c\ge0$ once
$B_0>(2p+12b+14c)\log\beta$.

Apply \eqref{ce:principal-monomial-bounds} to $\nu<L$ and include the tail in the errors.
The result is
\begin{equation}\label{ce:principal-normalized-expansions}
\begin{cases}
\mathcal C_i/B_i=1+\frac{\mathrm i z_iB_i''}{2\lambda B_i}
+O(\epsilon^2/y_*^2+\epsilon^3/y_*),\\[2pt]
\partial_{z_i}(\mathcal C_i/B_i)
=\frac{\mathrm i B_i''}{2\lambda B_i}
+O\bigl(\tau(\epsilon/y_*+\epsilon^2)\bigr),\\[2pt]
\partial_{y_i}(\mathcal C_i/B_i)
=\frac{\mathrm i z_i}{2\lambda}
\partial_{y_i}(B_i''/B_i)
+O\bigl(\epsilon^2/y_*^3+\epsilon^3/y_*^2\bigr).
\end{cases}
\end{equation}
To obtain amplitude and phase bounds, put $q_i\triangleq z_iB_i''/(2\lambda B_i)\in\R$ and
$D_i\triangleq\mathcal C_i/B_i-1-\mathrm i q_i$. Then
\begin{equation}\label{ce:principal-polar-identities}
\begin{aligned}
A_i/B_i&=\sqrt{(1+\operatorname{Re}D_i)^2
+(q_i+\operatorname{Im}D_i)^2},\\
\partial_x(A_i/B_i)
&=\frac{(1+\operatorname{Re}D_i)\operatorname{Re}\partial_xD_i
+(q_i+\operatorname{Im}D_i)(\partial_xq_i+\operatorname{Im}\partial_xD_i)}
{A_i/B_i},\\
\partial_x\psi_i
&=\frac{(1+\operatorname{Re}D_i)(\partial_xq_i+\operatorname{Im}\partial_xD_i)
-(q_i+\operatorname{Im}D_i)\operatorname{Re}\partial_xD_i}
{(A_i/B_i)^2},\qquad x\in\{z_i,y_i\}.
\end{aligned}
\end{equation}
In particular $q_i=O(\epsilon/y_*)$,
$\partial_{z_i}q_i=B_i''/(2\lambda B_i)=O(\tau)$, and
$\partial_{y_i}q_i=O(\tau)$. Substitution of
\eqref{ce:principal-normalized-expansions} into \eqref{ce:principal-polar-identities} gives
\[
\begin{cases}
\ds A_i/B_i-1 =O(\epsilon^2/y_*^2+\epsilon^3/y_*)=o(1),\\
\ns\ds B_i^{-1}|\partial_{z_i}A_i|
+B_i^{-1}|\partial_{y_i}A_i-k_iA_i|
\le Ck(\epsilon^2/y_*+\epsilon^3),\\
\ns\ds \partial_{z_i}\psi_i =\frac{B_i''}{2\lambda B_i}+o(\tau),
\qquad |\partial_{y_i}\psi_i|\le C\tau.
\end{cases}
\]
This proves the amplitude and phase estimates in
\eqref{ce:carrier-estimates}. The neighboring-packet comparison leading to \eqref{ce:principal-envelope-comparison}, together with \eqref{ce:front-correction-bound}, gives the remaining estimate
$|r_i|+\lambda^{-1}|\mathrm dr_i|\le CB_i\epsilon^8$.

\smallskip\noindent\textit{The quadratic form and its negative main term.}
The null coordinates in \eqref{ce:activation-edge-coordinates} satisfy
\[
t=\tfrac12(u_i+v_i),\qquad z_i=\tfrac12(u_i-v_i),\qquad
\partial_{u_i}=\tfrac12(\partial_t+\partial_{z_i}),\qquad
\partial_{v_i}=\tfrac12(\partial_t-\partial_{z_i}).
\]
The functions $A_i$ and $\psi_i$ are independent of $t$, whereas
$\partial_{u_i}\theta_i=\lambda+\tfrac12\psi_{i,z_i}$ and
$\partial_{v_i}\theta_i=-\tfrac12\psi_{i,z_i}$. Thus the full chain rule is
\begin{equation}\label{ce:principal-chain-rule}
\begin{cases}\ds
U_i=\tfrac12A_{i,z_i}\cos\theta_i
-A_i(\lambda+\tfrac12\psi_{i,z_i})\sin\theta_i+(r_i)_{u_i},\\
\ns\ds H_i =-\tfrac12A_{i,z_i}\cos\theta_i
+\tfrac12A_i\psi_{i,z_i}\sin\theta_i+(r_i)_{v_i},\\
\ns\ds Y_i =A_{i,y_i}\cos\theta_i-A_i\psi_{i,y_i}\sin\theta_i+(r_i)_{y_i}.
\end{cases}
\end{equation}
It follows that
\begin{equation}\label{ce:carrier-component-bounds}
\begin{aligned}
&U_i=-\lambda A_i\sin\theta_i
+O\bigl(B_ik[\epsilon|\sin\theta_i|+\epsilon^2/y_*+\epsilon^3]\bigr),\\
&H_i,\ Y_i-k_iA_i\cos\theta_i
=O\bigl(B_ik[\epsilon|\sin\theta_i|+\epsilon^2/y_*+\epsilon^3]\bigr),\\
&4U_iH_i=-B_i^2k_i^2\sin^2\theta_i+o(B_i^2k^2).
\end{aligned}
\end{equation}
To verify the last line, set
$T_i\triangleq A_{i,z_i}\cos\theta_i-A_i\psi_{i,z_i}\sin\theta_i$.
The exact derivatives in \eqref{ce:principal-chain-rule} give
\[
\begin{aligned}
4U_iH_i
&=2\lambda A_iA_{i,z_i}\sin\theta_i\cos\theta_i
-2\lambda A_i^2\psi_{i,z_i}\sin^2\theta_i -T_i^2+O(B_i^2\lambda^2\epsilon^8).
\end{aligned}
\]
Each error is small on the required scale, since
\[
\begin{aligned}
\lambda A_i|A_{i,z_i}|
&\le CB_i^2k^2(\epsilon/y_*+\epsilon^2)=o(B_i^2k^2),\\
T_i^2&\le CB_i^2k^2(\epsilon+\epsilon^2/y_*+\epsilon^3)^2
=o(B_i^2k^2),\\
B_i^2\lambda^2\epsilon^8&=B_i^2k^2\epsilon^6=o(B_i^2k^2),\\
2\lambda A_i^2\psi_{i,z_i}
&=\frac{A_i^2B_i''}{B_i}+o(B_i^2k^2)
=B_i^2k_i^2+o(B_i^2k^2).
\end{aligned}
\]
Also $Y_i=k_iA_i\cos\theta_i+o(B_ik)$. Consequently,
\[
\begin{aligned}
\mathcal B(\mathrm{d}f_i,\mathrm{d}f_i)
&=4U_iH_i-Y_i^2 =-B_i^2k_i^2\sin^2\theta_i-k_i^2A_i^2\cos^2\theta_i
+o(B_i^2k^2) =-B_i^2k_i^2+o(B_i^2k^2).
\end{aligned}
\]
This proves the claim and all of \eqref{ce:carrier-estimates}.

\medskip\noindent\textbf{Step 2. Choose the entrance curve near a vertex.}\proofsteplabel{step:activation-curve}{2}
We now construct a spatial graph $\eta=F_v(d)$ inside $C_v$ such that
\[
\mathcal B(\mathrm d(f_1+f_2),\mathrm d(f_1+f_2))<0
\quad\hbox{at every time on that graph},
\]
and whose two ends agree with the edge entrances. The construction
has three parts: opposite phases prevent cancellation of the gradient
in the center; one solution is small on either side; affine ends
match the edge curves. 

\smallskip\noindent\textit{Coordinates and the central antiphase curve.}
Write $e_2=\cos\alpha_v\,e_1+\sin\alpha_v\,n_1$, with
$|\alpha_v|=\alpha_0=\pi/32$, and set
\begin{equation}\label{ce:activation-bisector}
\begin{aligned}
&c_\alpha\triangleq\cos(\alpha_0/2),\quad s_\alpha\triangleq\sin(\alpha_v/2),\quad
\mathfrak g\triangleq1-\cos\alpha_0=2s_\alpha^2,\quad e_*\triangleq(e_1+e_2)/(2c_\alpha),\quad n_*\triangleq\mathsf Je_*.
\end{aligned}
\end{equation}
Here $c_\alpha>0$ and $|s_\alpha|=\sin(\alpha_0/2)>0$;
the sign of $s_\alpha$ records the orientation of the vertex.
Set $\xi\triangleq e_*\cdot(x-v)$, $\eta\triangleq n_*\cdot(x-v)$,
and $d\triangleq y_1-y_2$. Rotating the two edge frames into the $(e_*,n_*)$ frame of \eqref{ce:activation-bisector} gives
\begin{equation}\label{ce:activation-frames}
y_1=c_\alpha\eta+s_\alpha\xi,\quad y_2=c_\alpha\eta-s_\alpha\xi,\quad d=2s_\alpha\xi,\qquad
z_1=Z_1+c_\alpha\xi-s_\alpha\eta,\quad z_2=Z_2+c_\alpha\xi+s_\alpha\eta.
\end{equation}
In particular, $\eta=(y_1+y_2)/(2c_\alpha)$ and
$\xi=d/(2s_\alpha)$. Since $\theta_i=\lambda(t+z_i)+\psi_i$,
the common time term cancels:
\[
\theta_2-\theta_1
=[\lambda t+\lambda z_2+\psi_2]-[\lambda t+\lambda z_1+\psi_1]
=\lambda(z_2-z_1)+\psi_2-\psi_1.
\]
Thus one spatial curve fixes the relative phase for all $t$.
The coordinate derivatives and \eqref{ce:carrier-estimates} give
\[
\begin{aligned}
\partial_\eta(\psi_2-\psi_1)
&=s_\alpha(\psi_{2,z_2}+\psi_{1,z_1})
+c_\alpha(\psi_{2,y_2}-\psi_{1,y_1})=O(\tau),\\
\partial_d(\psi_2-\psi_1)
&=\frac{c_\alpha}{2s_\alpha}(\psi_{2,z_2}-\psi_{1,z_1})
-\tfrac12(\psi_{2,y_2}+\psi_{1,y_1})=O(\tau).
\end{aligned}
\]
Hence
\begin{equation}\label{ce:relative-phase}
\begin{aligned}
&\theta_2-\theta_1=\lambda(Z_2-Z_1+2s_\alpha\eta)+\psi_2-\psi_1,\\
&\partial_\eta(\theta_2-\theta_1)=2\lambda s_\alpha+O(\tau),\qquad
\partial_d(\theta_2-\theta_1)=O(\tau).
\end{aligned}
\end{equation}
Set
\[
d_L\triangleq-4\gamma M\ell/3,\qquad
d_R\triangleq-2\gamma M\ell/3.
\]
At $d=-\gamma M\ell$ choose one antiphase root $\eta_a(v)$ nearest to $y_*/c_\alpha$.
Indeed, \eqref{ce:relative-phase} gives
$|\partial_\eta(\theta_2-\theta_1)|\ge\lambda|s_\alpha|$
for large $M$. Consecutive antiphase roots are therefore separated
by at most $2\pi/(\lambda|s_\alpha|)$, so
$\eta_a=y_*/c_\alpha+O(\lambda^{-1})$.
Continue the same root by requiring
\[
\left.(\theta_2-\theta_1)\right|_{\eta=F_{\rm ph}(d)}
=(2N_v+1)\pi,\qquad N_v\in\Z,
\]
where $N_v$ is fixed for that vertex. The implicit function theorem and
\eqref{ce:relative-phase} give, on
$[d_L-M\ell,d_R+M\ell]$,
\[
F_{\rm ph}'
=-\frac{\partial_d(\theta_2-\theta_1)}
{\partial_\eta(\theta_2-\theta_1)}
=O(\tau/\lambda)=O(\epsilon^2),\qquad
F_{\rm ph}-\eta_a=O(M\epsilon^2\ell)=o(\lambda^{-1}).
\]
Here $\tau/\lambda=\epsilon^2$ and
$\lambda M\epsilon^2\ell=M\epsilon\to0$.
Moreover, $|d|=O(M\ell)$ and
$y_i=c_\alpha F_{\rm ph}(d)\pm d/2=y_*+O(M\ell)$.
The branch stays in the region of Step~\ref{step:activation-gradient},
so it continues across the whole stated interval.

\smallskip\noindent\textit{A baseline with affine ends and the final graph.}
Denote the flat step $b$ in \eqref{ce:cutoff-functions} by
$\chi_{\rm cut}$, as in \S\ref{ce:vertex-verification}.
Thus $\chi_{\rm cut}=0$ on $(-\infty,0]$,
$\chi_{\rm cut}=1$ on $[1,\infty)$, and
$\int_0^1\chi_{\rm cut}(q)\,\mathrm dq=1/2$.
We first define $H_{\rm base}$ to be horizontal in the center,
with slopes $1/(2c_\alpha)$ and $-1/(2c_\alpha)$ at its left
and right ends. Set $H_{\rm base}(d_L)=\eta_a$ and
\begin{equation}\label{ce:activation-base-graph}
H_{\rm base}'(d)=
\begin{cases}
[1-\chi_{\rm cut}(3+3d/(2\gamma M\ell))]/(2c_\alpha),&d<d_L,\\
0,&d_L\le d\le d_R,\\
-\chi_{\rm cut}(1+3d/(2\gamma M\ell))/(2c_\alpha),&d>d_R.
\end{cases}
\end{equation}
Flatness of $\chi_{\rm cut}$ at $0$ and $1$ makes the definition in \eqref{ce:activation-base-graph} smooth at $d_L$ and $d_R$. The arguments of $\chi_{\rm cut}$ run from $0$ to $1$ on
$[-2\gamma M\ell,d_L]$ and $[d_R,0]$. In particular,
\begin{equation}\label{ce:activation-base-end-heights}
\begin{aligned}
H_{\rm base}(d_L)-H_{\rm base}(-2\gamma M\ell)
&=\frac{\gamma M\ell}{3c_\alpha}
\int_0^1(1-\chi_{\rm cut}(q))\,\mathrm dq
=\frac{\gamma M\ell}{6c_\alpha},\\
H_{\rm base}(d_R)-H_{\rm base}(0)
&=\frac{\gamma M\ell}{3c_\alpha}
\int_0^1\chi_{\rm cut}(q)\,\mathrm dq
=\frac{\gamma M\ell}{6c_\alpha}.
\end{aligned}
\end{equation}
Hence the outer ends of both bending intervals have height
$\eta_a-\gamma M\ell/(6c_\alpha)$.
To replace the central horizontal segment by the antiphase curve,
choose $\chi\in C_c^\infty((d_L-M\ell,d_R+M\ell);[0,1])$
equal to one on $[d_L,d_R]$, with
$|\chi'|\le C/(M\ell)$, and set
\begin{equation}\label{ce:activation-graph}
F_v(d)\triangleq H_{\rm base}(d)+\chi(d)(F_{\rm ph}(d)-\eta_a),
\qquad X_v\triangleq\frac{\eta-F_v(d)}{\ell}.
\end{equation}
The product is extended by zero outside the interval where $F_{\rm ph}$
is defined. On $[d_L,d_R]$, $H_{\rm base}=\eta_a$ and $\chi=1$,
so $F_v=F_{\rm ph}$. At the affine ends, $\chi=0$, so
$F_v=H_{\rm base}$. The two matching errors are explicitly
\[
F_v-H_{\rm base}=\chi(F_{\rm ph}-\eta_a),\qquad
F_v'-H_{\rm base}'=\chi'(F_{\rm ph}-\eta_a)+\chi F_{\rm ph}'.
\]
They satisfy
\begin{equation}\label{ce:activation-graph-deviation}
|F_v-H_{\rm base}|\le CM\epsilon^2\ell,\qquad
|F_v'-H_{\rm base}'|\le C\epsilon^2.
\end{equation}
The function $X_v$ is the height difference $\eta-F_v(d)$, divided by
$\ell$; it vanishes exactly on this graph.
Its pieces and their respective reasons for negativity are
\[
\begin{aligned}
&d\in[d_L,d_R]: &&\theta_2-\theta_1=\pi\pmod{2\pi},\\
&d\in(-2\gamma M\ell,d_L)\cup(d_R,0):
&&\text{one principal amplitude is small},\\
&d\in(-\infty,-2\gamma M\ell]\cup[0,\infty):
&&\text{one complete solution is small}.
\end{aligned}
\]
\begin{figure}[H]
\centering
\begin{tikzpicture}[x=.85cm,y=.58cm,font=\scriptsize,>={Stealth[length=1.3mm]}]
\fill[gray!10] (-1,0)--(-1,2.5)--(1,2.5)--(1,0)--cycle;
\draw[->,thin] (-4.5,0)--(4.6,0) node[right] {$d$};
\draw[->,thin] (-4.45,0)--(-4.45,3.2) node[above] {$\eta$};
\draw[gray,line width=.9pt] (-4.2,.3)--(-3,1.15)
.. controls (-2.3,1.65) and (-1.8,1.93) .. (-1.4,1.97)
.. controls (-1.25,2) and (-1.15,2) .. (-1,2)
--(1,2).. controls (1.15,2) and (1.25,2) .. (1.4,1.97)
.. controls (1.8,1.93) and (2.3,1.65) .. (3,1.15)--(4.2,.3);
\draw[densely dashed] (-1.18,1.835)--(1.18,2.165);
\draw[line width=.8pt] (-4.2,.3)--(-3,1.15)
.. controls (-2.3,1.65) and (-1.8,1.93) .. (-1.4,1.97)
.. controls (-1.25,2) and (-1.15,1.84) .. (-1,1.86)
--(1,2.14).. controls (1.15,2.16) and (1.25,2) .. (1.4,1.97)
.. controls (1.8,1.93) and (2.3,1.65) .. (3,1.15)--(4.2,.3);
\foreach \x/\lab in {-3/{-2\gamma M\ell},-1/{d_L},1/{d_R},3/{0}}{
\draw[gray,dotted] (\x,0)--(\x,2.45);
\node[below,inner sep=1pt] at (\x,0) {$\lab$};}
\draw[gray,line width=.9pt] (-3.9,2.85)--(-3.3,2.85);
\node[right] at (-3.3,2.85) {$H_{\rm base}$};
\draw[densely dashed] (-.9,2.85)--(-.3,2.85);
\node[right] at (-.3,2.85) {$F_{\rm ph}$};
\draw[line width=.8pt] (2.1,2.85)--(2.7,2.85);
\node[right] at (2.7,2.85) {$F_v$};
\node at (-3.75,-.72) {affine};\node at (-2,-.72) {bend};
\node at (0,-.72) {antiphase};\node at (2,-.72) {bend};\node at (3.75,-.72) {affine};
\end{tikzpicture}
\caption{The vertex entrance curve (schematic). The final graph $F_v$ agrees
with $F_{\rm ph}$ in the center and with $H_{\rm base}$ at the affine ends.
The phase correction is exaggerated; its actual size is bounded by
\eqref{ce:activation-graph-deviation}.}
\label{fig:activation-entrance-curves}
\end{figure}
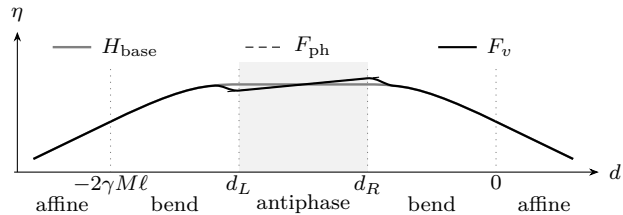

We prove negativity on these regions in turn; on the affine ends,
the estimate must continue until the smaller local solution disappears.

\smallskip\noindent\textit{Negativity in the center.}
Expanding the quadratic form in the original edge frames gives
\[
\begin{split}
\mathcal B(\mathrm{d}f_1,\mathrm{d}f_2)={}&
\mathfrak g(U_1U_2+H_1H_2)
+(1+\cos\alpha_v)(U_1H_2+H_1U_2)\\
&+\sin\alpha_v[(U_1-H_1)Y_2-Y_1(U_2-H_2)]-\cos\alpha_v Y_1Y_2.
\end{split}
\]
On the graph, $\theta=\theta_1$ and $\theta_2=\theta+\pi$ modulo $2\pi$.
By \eqref{ce:carrier-estimates},
\[
\begin{aligned}
&\mathcal B(\mathrm{d}(f_1+f_2),\mathrm{d}(f_1+f_2))\\
&= -(k_1^2A_1^2+k_2^2A_2^2)\sin^2\theta
-|k_1A_1n_1-k_2A_2n_2|^2\cos^2\theta -2\mathfrak g\lambda^2A_1A_2\sin^2\theta\\
&\quad +2\sin\alpha_v\lambda A_1A_2(k_2-k_1)\sin\theta\cos\theta +O(A_1A_2k^2\sin^2\theta)+o(B_{\rm tot}^2k^2).
\end{aligned}
\]
The elementary inequalities
\[
\begin{cases}
k_1^2A_1^2+k_2^2A_2^2\ge\tfrac12(k_1A_1+k_2A_2)^2,\\[3pt]
|k_1A_1n_1-k_2A_2n_2|^2
\ge\tfrac{1-\cos\alpha_0}{2}(k_1A_1+k_2A_2)^2
\end{cases}
\]
and the bounds $A_i\ge B_i/2$, $k_i\ge k/2$ for large $M$ imply
\begin{equation}\label{ce:activation-central-coercivity}
\begin{aligned}
&-(k_1^2A_1^2+k_2^2A_2^2)\sin^2\theta
-|k_1A_1n_1-k_2A_2n_2|^2\cos^2\theta \le-\frac{\mathfrak g}{2}(k_1A_1+k_2A_2)^2
\le-c_0B_{\rm tot}^2k^2, 
\end{aligned}
\end{equation}
where $c_0\triangleq\frac{\mathfrak g}{32}>0.$

For the mixed term, Young's inequality reads
\begin{equation}\label{ce:activation-mixed-young}
\begin{aligned}
&\left|2\sin\alpha_v\lambda A_1A_2(k_2-k_1)\sin\theta\cos\theta\right| \le\mathfrak g\lambda^2A_1A_2\sin^2\theta
+\frac{\sin^2\alpha_v}{\mathfrak g}
A_1A_2(k_2-k_1)^2\cos^2\theta.
\end{aligned}
\end{equation}
The remaining scale comparison is
\begin{equation}\label{ce:activation-normal-difference}
\frac{|k_2-k_1|}{k}
=\frac{|d|(y_1+y_2)y_*^2}{y_1^2y_2^2}
=O(My_*)\longrightarrow0.
\end{equation}
Since $k^2/\lambda^2=\epsilon^2\to0$, \eqref{ce:activation-central-coercivity}--\eqref{ce:activation-normal-difference} give
\begin{equation*}\label{ce:activation-central-negative}
\begin{aligned}
&\mathcal B(\mathrm{d}(f_1\!+f_2),\mathrm{d}(f_1\!+f_2)) \le-[c_0\!-C(My_*)^2\!-o(1)]B_{\rm tot}^2k^2
-(\mathfrak g\lambda^2\!-Ck^2)A_1A_2\sin^2\theta \le\!-\tfrac12c_0B_{\rm tot}^2k^2
\end{aligned}
\end{equation*}
for every sufficiently large $M$, with $c_0$ defined in \eqref{ce:activation-central-coercivity}.

\smallskip\noindent\textit{Negativity on the bends.}
Here we compare one solution with the other instead of adjusting phases.
On each side and its affine continuation, denote the larger and smaller
local solutions by $f_{\rm dom}$ and $f_{\rm minor}$, and their
principal amplitudes by $B_{\rm dom}$ and $B_{\rm minor}$.
Since $B_i=e^{-1/y_i-M^2-\gamma_iM}\rho(\beta^My_i)$ and
$\gamma_1=0$, $\gamma_2=\gamma$,
\begin{equation}\label{ce:activation-amplitude-log}
\log\frac{B_2}{B_1}
=-\gamma M+\frac1{y_1}-\frac1{y_2}
+\log\frac{\rho(\beta^My_2)}{\rho(\beta^My_1)}.
\end{equation}
The logarithm of the cutoff ratio in \eqref{ce:activation-amplitude-log} is $o(1)$ on the principal-packet region, and
$1/y_1-1/y_2=-d/(y_1y_2)$. Thus
\begin{equation}\label{ce:activation-bend-amplitude-ratio}
\log(B_2/B_1)=-\gamma M-d/(y_1y_2)+o(1).
\end{equation}
On both bends, $y_1y_2=y_*^2[1+O(My_*)]$.
The two ranges of $d$ give the separate estimates
\[
\begin{aligned}
d\le d_L=-4\gamma M\ell/3:
&\quad\log(B_1/B_2)\le-\gamma M/3+O(M^2y_*)+o(1),\\
d\ge d_R=-2\gamma M\ell/3:
&\quad\log(B_2/B_1)\le-\gamma M/3+O(M^2y_*)+o(1).
\end{aligned}
\]
The graph correction changes $1/y_i$ by
$O(kM\epsilon^2\ell)=O(M\epsilon^2)=o(1)$.
The smaller-to-larger amplitude ratio is consequently
$O(e^{-\gamma M/3+O(M^2y_*)})$. Our choice $\gamma=24\log2$ ensures
\[
\epsilon^8=2^{-7M},\qquad
\frac{e^{-\gamma M/3+O(M^2y_*)}}{\epsilon^8}
=2^{-M}e^{O(M^2y_*)}\longrightarrow0.
\]
Thus \eqref{ce:carrier-estimates} gives
$|\mathrm df_{\rm minor}|=o(B_{\rm dom}\lambda\epsilon^8)$
on the bends; the dominant solution retains the negative bound of
Step~\ref{step:activation-gradient}.

\smallskip\noindent\textit{Negativity on the affine ends.}
The smaller solution can now leave the principal-packet region, so we
estimate its full packet sum. By \eqref{ce:activation-base-end-heights}, on the right end, the minor normal height is at most
$c_\alpha\eta_a-\gamma M\ell/6
=y_*-\gamma M\ell/6+O(\lambda^{-1})$.
On the left end it is at most
$c_\alpha\eta_a-7\gamma M\ell/6
=y_*-7\gamma M\ell/6+O(\lambda^{-1})$.
Thus the minor height lies at least $cM\ell$ below $y_*$ in both cases.
Let $\gamma_{\rm minor}$ be the type of the smaller solution.
Write $y$ for the normal coordinate of the smaller solution. On the right
end the dominant solution has type $0$ and height
$y_{\rm dom}=c_\alpha\eta_a-\gamma M\ell/6$, with $y\le y_{\rm dom}$;
on the left it has type $\gamma$ and height
$y_{\rm dom}=c_\alpha\eta_a+5\gamma M\ell/6$, with
$y\le y_{\rm dom}-2\gamma M\ell$.
For $y>0$ all active packets of the smaller solution have indices $j\ge M$.
Since $\rho\le1$ and $\rho(\beta^My_{\rm dom})=1+o(1)$, their amplitudes obey
\[
\frac{A_{\gamma_{\rm minor},j}(y)}{B_{\rm dom}}
\le C\exp\{-\gamma M+CM^2y_*-(j^2-M^2)\}.
\]
For the left end this uses
$1/y-1/y_{\rm dom}\ge2\gamma M-CM^2y_*$; on the right,
$1/y-1/y_{\rm dom}\ge0$ and $\gamma_{\rm minor}j\ge\gamma M$.
The derivative bounds for $\rho$ and \eqref{ce:coefficient-estimate} imply,
for $r\in\{0,1\}$,
\[
\begin{aligned}
\frac{\max_{|\alpha|=r}|\partial^\alpha g_{\rm minor}|}
{B_{\rm dom}\lambda^r}
&\le C e^{-\gamma M+CM^2y_*}
\sum_{j\ge M}e^{-(j^2-M^2)}2^{r(j-M)}\beta^{2r(j-M)} \le C e^{-\gamma M+CM^2y_*},
\end{aligned}
\]
Here $g_{\rm minor}$ is the approximate solution \eqref{ce:approximate-sum}
of the same type in the same edge frame. The series satisfies
\[
\begin{aligned}
&\sum_{j\ge M}e^{-(j^2-M^2)}2^{r(j-M)}\beta^{2r(j-M)} =\sum_{h\ge0}
\exp\!\left[-h^2-h(2M-r\log2-2r\log\beta)\right]
\le\sum_{h\ge0}e^{-h}=\frac1{1-e^{-1}}
\end{aligned}
\]
for $2M\ge1+r\log2+2r\log\beta$. Likewise,
\[
\begin{aligned}
\frac{E_{\gamma_{\rm minor}}(y)}{B_{\rm dom}}
&=\sum_{j\ge M}\frac{A_{\gamma_{\rm minor},j}(y)}{B_{\rm dom}} \le Ce^{-\gamma M+CM^2y_*}\sum_{h\ge0}e^{-2Mh-h^2}
\le Ce^{-\gamma M+CM^2y_*}.
\end{aligned}
\]
For $\beta^{-(m+1)}<y\le\beta^{-m}$ we have $m\ge M$, and
\eqref{ce:reserve-estimate} gives
\[
\begin{aligned}
\frac{\max_{|\alpha|\le1}
|\partial^\alpha(f_{\rm minor}-g_{\rm minor})|}
{E_{\gamma_{\rm minor}}(y)}
&\le C e^{-cm^3\beta^m+Cm^2\beta^m\log(m+2)} \le C e^{-cM^3\beta^M/2}\longrightarrow0.
\end{aligned}
\]
Thus
\[
\frac{\max_{|\alpha|=r}|\partial^\alpha f_{\rm minor}|}
{B_{\rm dom}\lambda^r}
\le C e^{-\gamma M+CM^2y_*}=o(\epsilon^8).
\]
For $y\le0$ the smaller solution and all its derivatives vanish.
On the bends, the same derivative bound with $o(\epsilon^8)$ follows
from \eqref{ce:carrier-estimates} and the amplitude ratio \eqref{ce:activation-bend-amplitude-ratio}.
On all these pieces Step~\ref{step:activation-gradient} gives $|\mathrm{d}f_{\rm dom}|\le CB_{\rm dom}\lambda$
and $\mathcal B(\mathrm{d}f_{\rm dom},\mathrm{d}f_{\rm dom})
\le-cB_{\rm dom}^2k^2$. Consequently,
\[
\begin{aligned}
&\mathcal B(\mathrm{d}(f_{\rm dom}+f_{\rm minor}),
\mathrm{d}(f_{\rm dom}+f_{\rm minor}))\\
&=\mathcal B(\mathrm{d}f_{\rm dom},\mathrm{d}f_{\rm dom})
+2\mathcal B(\mathrm{d}f_{\rm dom},\mathrm{d}f_{\rm minor})
+\mathcal B(\mathrm{d}f_{\rm minor},\mathrm{d}f_{\rm minor})\\
&\le-cB_{\rm dom}^2k^2+o(B_{\rm dom}^2\lambda^2\epsilon^8)
\le-\tfrac c2B_{\rm dom}^2k^2,
\qquad \frac{\lambda^2\epsilon^8}{k^2}=\epsilon^6\longrightarrow0.
\end{aligned}
\]

We now pass from $f_1+f_2$ to the merged vertex solution.
In $\{y_1>0,y_2>0\}$, Proposition~\ref{ce:incoming} gives
$p_{\rm in}=e^{\Phi_0}(f_1+f_2)$. At a zero of $p_{\rm in}$,
\[
\begin{cases}
\ds\mathrm dp_{\rm in}
=e^{\Phi_0}[\mathrm d(f_1+f_2)+(f_1+f_2)\mathrm d\Phi_0]
=e^{\Phi_0}\mathrm d(f_1+f_2),\\
\ns\ds \mathcal B(\mathrm dp_{\rm in},\mathrm dp_{\rm in})
=e^{2\Phi_0}\mathcal B(\mathrm d(f_1+f_2),\mathrm d(f_1+f_2))<0.
\end{cases}
\]
If $y_i>0$ and $y_{3-i}\le0$, the same proposition gives
$p_{\rm in}=f_i$; the identity for first derivatives on
$y_{3-i}=0$ follows by smooth continuation. The positive height
of the dominant solution stays in $y_i=y_*+O(M\ell)$, so
Step~\ref{step:activation-gradient} gives
\[
\mathcal B(\mathrm dp_{\rm in},\mathrm dp_{\rm in})
=\mathcal B(\mathrm df_i,\mathrm df_i)<0.
\]
The graph never enters $\{y_1\le0,y_2\le0\}$: both heights
are $y_*+O(M\ell)>0$ in the center and bends, and the
dominant height has this form on the affine ends.
This proves the entrance condition on the entire vertex graph.

\medskip\noindent\textbf{Step 3. Match the local levels and construct the global data.}\proofsteplabel{step:activation-global-data}{3}
The vertex and edge definitions must agree on open overlaps. The heights of the affine ends of $F_v$ determine the entrance along each whole edge.

\smallskip\noindent\textit{Endpoint heights and interpolation along an edge.}
By \eqref{ce:activation-base-end-heights}, the two affine ends are
\[
F_v(d)=\begin{cases}
h_\gamma(v)/c_\alpha+d/(2c_\alpha),&d\le-2\gamma M\ell,\\
h_0(v)/c_\alpha-d/(2c_\alpha),&d\ge0,
\end{cases}
\]
where
\begin{equation}\label{ce:type-heights}
\begin{aligned}
&h_0(v)=c_\alpha\eta_a(v)-\gamma M\ell/6,\quad
h_\gamma(v)=c_\alpha\eta_a(v)+5\gamma M\ell/6,\quad
h_\gamma(v)-h_0(v)=\gamma M\ell.
\end{aligned}
\end{equation}
Since $\eta_a(v)=y_*/c_\alpha+O(\lambda^{-1})$ and $M\ell/y_*\to0$,
\[
h_0(v),h_\gamma(v)\in[y_*/2,3y_*/2]
\]
for large $M$. For an edge of fixed type $\gamma_i$, write
$
v_{i,\pm}\triangleq c_i\pm({r_{\rm geo}}L_0/2)e_i$ and
$h_{i,\pm}\triangleq h_{\gamma_i}(v_{i,\pm}).
$
The type-dependent terms in \eqref{ce:type-heights} cancel between its endpoints:
\[
h_{i,+}-h_{i,-}
=c_\alpha[\eta_a(v_{i,+})-\eta_a(v_{i,-})]=O(\lambda^{-1}).
\]
Define
\[
h_i(z)\triangleq h_{i,-}
+(h_{i,+}-h_{i,-})\chi_{\rm cut}\!\left(
\frac{z+{r_{\rm geo}}L_0/2-3\rho_p}{{r_{\rm geo}}L_0-6\rho_p}\right).
\]
The denominator is positive because $4\rho_p<{r_{\rm geo}}L_0/10$. The cutoff gives
\[
h_i(z)=\begin{cases}
h_{i,-},&z\le-{r_{\rm geo}}L_0/2+3\rho_p,\\
h_{i,+},&z\ge {r_{\rm geo}}L_0/2-3\rho_p.
\end{cases}
\]
Thus the height is constant near both endpoints, and always belongs to
$[y_*/2,3y_*/2]$. Set
\begin{equation}\label{ce:edge-level}
X_i\triangleq\frac{y_i-h_i(z_i)}{c_\alpha\ell}.
\end{equation}

\smallskip\noindent\textit{Exact matching on vertex--edge overlaps.}
Using the widths fixed before the proposition, define
\[
C_i\triangleq\left\{x\in\Omega:
|n_i\cdot(x-c_i)|<2\delta_p,\quad
|e_i\cdot(x-c_i)|<{r_{\rm geo}}L_0/2-\rho_p/2\right\}.
\]
The geometric gaps separate nonincident charts and distinct polygons;
adjacent edge charts do not intersect, and the vertex balls are disjoint.
On $C_v\cap C_i$, the distance of the edge projection from $v$ is at
least $\rho_p/2$. The other normal coordinate and their difference obey
\begin{equation}\label{ce:activation-overlap-radius}
y_{3-i}\le-\frac{\rho_p\sin\alpha_0}{2}+2\delta_p<0,\qquad
|d|\ge\frac{\rho_p\sin\alpha_0}{2}-4\delta_p
>\frac{\rho_p\sin\alpha_0}{4}>2\gamma M\ell.
\end{equation}
The last inequality in \eqref{ce:activation-overlap-radius} follows from $M\ell/\rho_p\to0$.
Also $|z_i-Z_i|<2\rho_p<3\rho_p$, so $h_i$ equals
$h_{\gamma_i}(v)$ throughout this overlap.
In the local labeling $\gamma_1=0$, $\gamma_2=\gamma$, the affine formulas give
\[
\begin{aligned}
d\ge0:\quad
X_v&=\frac{c_\alpha\eta+d/2-h_0(v)}{c_\alpha\ell}
=\frac{y_1-h_0(v)}{c_\alpha\ell}=X_1,\\
d\le-2\gamma M\ell:\quad
X_v&=\frac{c_\alpha\eta-d/2-h_\gamma(v)}{c_\alpha\ell}
=\frac{y_2-h_\gamma(v)}{c_\alpha\ell}=X_2.
\end{aligned}
\]
Because the other normal coordinate is negative, Proposition~\ref{ce:incoming}
also identifies the vertex pair with the original edge pair. We have proved
\begin{equation}\label{ce:common-chart-data}
X_v=X_i,\qquad(p_{\rm in},V_{\rm in})=(f_i,W_i)
\quad\hbox{on }C_v\cap C_i.
\end{equation}

\smallskip\noindent\textit{The pair and the level function on a common neighborhood.}
Let $P_1,\ldots,P_J$, $J\in\mathbb N$, denote the observation polygons and put
\[
U\triangleq\bigcup_{j=1}^J\{x:\operatorname{dist}(x,P_j)<\delta_p\},\qquad
U^+\triangleq\bigcup_{j=1}^J\{x:\operatorname{dist}(x,P_j)<3\delta_p/2\}.
\]
The geometric gaps give $\overline U\subset U^+\Subset\Omega$.
The polygon interiors, vertex balls and edge charts cover $\overline{U^+}$.
Indeed, for a point outside the polygons whose nearest point is a vertex
or is within $\rho_p$ of an edge endpoint,
\[
|x-v|\le\sqrt{\rho_p^2+(3\delta_p/2)^2}<2\rho_p.
\]
Otherwise its edge projection is more than $\rho_p$ from both endpoints,
and its normal distance is at most $3\delta_p/2<2\delta_p$, so it lies in $C_i$.
On portions inside a polygon, each edge pair and each vertex pair is
$(0,0)$, by the flat extension in Section~\ref{ce:front-section} and
Proposition~\ref{ce:incoming}, respectively.
The flat extensions in Section~\ref{ce:front-section} and Proposition~\ref{ce:incoming}, together with the matching identities \eqref{ce:common-chart-data}, give a smooth,
real, $2\pi$-periodic pair on $\R\times U^+$ satisfying
\[
\Box p_{\rm in}+V_{\rm in}p_{\rm in}=0.
\]

To extend the level over the polygon interiors, we first check that its
local values are less than $-2$ there. The derivative formula
\eqref{ce:activation-base-graph} gives
\[
\begin{aligned}
&\frac{\mathrm d}{\mathrm dd}\left(H_{\rm base}(d)
+\frac{|d|}{2c_\alpha}\right)\le0\quad(d<0),\\
&H_{\rm base}(d)+\frac{|d|}{2c_\alpha}
=\frac{h_0(v)}{c_\alpha}\quad(d\ge0).
\end{aligned}
\]
Hence, using \eqref{ce:activation-graph-deviation},
\[
H_{\rm base}(d)+\frac{|d|}{2c_\alpha}\ge\frac{h_0(v)}{c_\alpha},\qquad
F_v(d)+\frac{|d|}{2c_\alpha}\ge\frac{y_*}{2c_\alpha}.
\]
If $y_1,y_2\le0$, then $\eta\le-|d|/(2c_\alpha)$, and therefore
\[
X_v\le-\frac{y_*}{2c_\alpha\ell}<-2.
\]
Similarly, on a polygon portion of an edge chart,
\[
y_i\le0,\quad h_i\ge y_*/2
\quad\Longrightarrow\quad
X_i\le-\frac{y_*}{2c_\alpha\ell}<-2.
\]
Define a smooth scalar cutoff
\[
h_{\rm lev}(r)\triangleq-2+(r+2)\chi_{\rm cut}(r+2).
\]
It satisfies
\[
h_{\rm lev}(r)=r\ (r\ge-1),\quad
h_{\rm lev}(r)=-2\ (r\le-2),\quad
-2\le h_{\rm lev}(r)<-1\ (-2<r<-1).
\]
We can now define
Flatness of $\chi_{\rm cut}$ at $0,1$ matches every derivative of $h_{\rm lev}$ at $r=-2,-1$ with those of the constant $-2$ and the identity $r$, respectively. Elsewhere it is a product of smooth functions. Thus $h_{\rm lev}\in C^\infty(\R)$, and its compositions with the local level functions are smooth.
\begin{equation}\label{ce:activation-global-level}
X(x)\triangleq\begin{cases}
h_{\rm lev}(X_v(x)),&x\in C_v\cap U^+,\\
h_{\rm lev}(X_i(x)),&x\in C_i\cap U^+,\\
-2,&x\in\displaystyle\bigcup_{j=1}^J\operatorname{int}P_j.
\end{cases}
\end{equation}
The first two expressions agree by \eqref{ce:common-chart-data};
on polygon portions both equal $-2$. Thus $X\in C^\infty(U^+)$ and
$X=-2$ on every polygon. Moreover,
$
h_{\rm lev}^{-1}([-1/2,3/2])=[-1/2,3/2],
$
so the cutoff preserves the entrance and the entire transition band.

\smallskip\noindent\textit{Compactness of the transition band and the entrance condition.}
We prove that every point of the entrance graph is within $3y_*$ of
its polygon. On the central and bending pieces,
\[
|x-v|\le2y_*/c_\alpha+o(y_*)<3y_*.
\]
On an affine piece, let $\sigma\in[y_*/2,3y_*/2]$ be the dominant
normal height and set $D\triangleq|d|$. The point $x_0$ on that graph satisfies
\begin{equation}\label{ce:entrance-vertex-radius}
|x_0-v|^2=\frac{(2\sigma-D)^2}{4c_\alpha^2}
+\frac{D^2}{4s_\alpha^2}.
\end{equation}
The signed distance of its perpendicular projection along the incident
edge ray is
$
t_{\rm edge}=\frac{D-2s_\alpha^2\sigma}{2|s_\alpha|c_\alpha}.
$
If $t_{\rm edge}\ge0$, this projection lies on the edge and its normal
distance is $\sigma$. If $t_{\rm edge}\le0$, then
$0\le D\le2s_\alpha^2\sigma$; convexity in $D$ of \eqref{ce:entrance-vertex-radius} gives
\[
|x_0-v|^2\le
\max\left\{\frac{\sigma^2}{c_\alpha^2},\sigma^2\right\}
=\frac{\sigma^2}{c_\alpha^2}.
\]
Estimate \eqref{ce:entrance-polygon-distance} also applies to the graph points used in Step~\ref{step:local-constant-exit},
even if they lie slightly outside $C_v$: for large $M$,
\[
|x_0-v|\le2\rho_p+3\ell/2<{r_{\rm geo}}L_0/10.
\]
Thus the projection cannot pass the other endpoint. Consequently, on both affine pieces,
\begin{equation}\label{ce:entrance-polygon-distance}
\operatorname{dist}(x_0,P_j)\le\sigma/c_\alpha<3y_*.
\end{equation}
The edge graphs have normal height $h_i\le3y_*/2$, and their
projections lie on the corresponding edges, so \eqref{ce:entrance-polygon-distance} holds for the edge graphs as well.
All these graphs lie outside the polygons, because their dominant
normal height is positive.

For $-1/2\le X(x)\le3/2$, the cutoff is the identity on the local level.
The corresponding graph point is
\[
x_0=\begin{cases}
v+\xi e_*+F_v(2s_\alpha\xi)n_*,&\text{in a vertex chart},\\
x-(y_i-h_i(z_i))n_i,&\text{in an edge chart}.
\end{cases}
\]
Thus
\[
|x-x_0|\le\tfrac32\ell,\qquad
\operatorname{dist}\!\Big(x,\bigcup_jP_j\Big)
\le3y_*+\tfrac32\ell<4y_*<\delta_p/8.
\]
The band is relatively closed in $U^+$ and lies in the compact set
\[
\Big\{x:\operatorname{dist}\!\Big(x,\bigcup_jP_j\Big)
\le\delta_p/8\Big\}\Subset U.
\]
It is therefore compactly contained in $U$.
At entrance zeros in vertex charts, Step~\ref{step:activation-curve}
gives the required negative quadratic form. On an edge entrance,
\[
y_i=h_i(z_i)=y_*+O(M\ell),\qquad
|z_i|\le(SL_0/2)y_*.
\]
Step~\ref{step:activation-gradient} applies, so in every chart
\[
\mathcal B(\mathrm dp_{\rm in},\mathrm dp_{\rm in})<0
\quad\hbox{where }X=0,\ p_{\rm in}=0.
\]

\medskip\noindent\textbf{Step 4. Apply the completion proposition.}\proofsteplabel{step:activation-apply}{4}
We have now verified all hypotheses of Proposition~\ref{ce:activation-criterion}:
\[
\begin{cases}
\ds (p_{\rm in},V_{\rm in})\in C^\infty(\R\times U^+;\R^2),\qquad
\Box p_{\rm in}+V_{\rm in}p_{\rm in}=0,\\
\ns\ds (p_{\rm in},V_{\rm in})(t+2\pi,x)=(p_{\rm in},V_{\rm in})(t,x),\\
\ns\ds X\in C^\infty(U^+;\R),\qquad
\{-1/2\le X\le3/2\}\Subset U\Subset U^+,\\
\ns\ds \mathcal B(\mathrm dp_{\rm in},\mathrm dp_{\rm in})<0
\quad\hbox{on }\{X=0,\ p_{\rm in}=0\}.
\end{cases}
\]
Apply that proposition on $U^+$ with $(f,W)=(p_{\rm in},V_{\rm in})$,
$A\equiv1$ and $T=2\pi$; the final potential is $\Delta A/A=0$.
It gives a smooth periodic pair and $\varepsilon\in(0,1/2)$ such that
\[
\Box p+Vp=0,\qquad
(p,V)=(p_{\rm in},V_{\rm in})\ (X<\varepsilon),\qquad
(p,V)=(1,0)\ (X>1-\varepsilon).
\]
Since $\varepsilon>0$, these equalities imply the required matching on
$X\le0$ and $X\ge1$. On every polygon $X=-2$, so $(p,V)=(0,0)$ there.

Only finitely many vertices and estimates occur. Choose $M$ above all
thresholds in Steps~\ref{step:activation-gradient}--\ref{step:activation-global-data};
the final common choice is made in Lemma~\ref{ce:global-local}.
Restrict the constructed functions from $U^+$ to $U$. Since
\[
\overline{\mathcal N_{\delta_p/4}(\overline U)}\subset U^+,
\]
smoothness and $2\pi$-periodicity give, for every multi-index $\alpha$,
\[
\sup_{(t,x)\in\R\times\overline{\mathcal N_{\delta_p/4}(\overline U)}}
\left(|\partial^\alpha p|+|\partial^\alpha V|
+|\partial^\alpha X|\right)<\infty.
\]
This also gives the asserted bounds on $\R\times\overline U$.
\end{proof}

\subsection{Assembly and the Dirichlet boundary condition}
In the same neighborhood $U$ as in Proposition~\ref{ce:activation}, choose a compact core $K\Subset U$ outside which the completed local pair is $(1,0)$. Join it to a positive time-independent solution equal to one near $\overline U$ and zero at $\partial\Omega$.

\begin{lemma}[Assembly with a constant exterior state]\label{ce:global-local}
The local solution data constructed in Section~\ref{ce:front-section}
and the data in Propositions~\ref{ce:incoming} and~\ref{ce:activation}
can be chosen simultaneously. On the resulting scaled geometry there
are a compact set $K$ and an open set $U$, with
$K\Subset U\Subset\Omega$, and a smooth $2\pi$-periodic pair
$(p_{\rm loc},V_{\rm loc})$ on $\R\times U$ satisfying
\begin{equation}\label{ce:local-exit}
\Box p_{\rm loc}+V_{\rm loc}p_{\rm loc}=0,\qquad
(p_{\rm loc},V_{\rm loc})=\begin{cases}
(0,0),&x\in\omega,\\
(1,0),&x\in U\setminus K.
\end{cases}
\end{equation}
\end{lemma}

\begin{proof}[Proof of Lemma~\ref{ce:global-local}]
We choose one integer $M$ for all vertices and edges, then prove that
the pair from Proposition~\ref{ce:activation} is already constant on
$U\setminus K$.

\medskip\noindent\textbf{Step 1. Fix one choice of parameters for all local constructions.}\proofsteplabel{step:local-common-parameters}{1}
Fix the reference geometry, the local solution family with normal width
$Y_*$, and $S$ satisfying the first inequality in
\eqref{ce:geometric-scale-conditions}. The reference widths $\rho_g,\delta_g$
and the cutoff shapes in Proposition~\ref{ce:activation} are also fixed
before $M$ is chosen.
For a reference vertex $v^0$ and an incident edge midpoint $c_i^0$, set
\[
v={r_{\rm geo}}v^0,\qquad c_i={r_{\rm geo}}c_i^0,\qquad
R\triangleq2S\rho_g,\qquad
\zeta_i\triangleq S e_i\cdot(v^0-c_i^0),\quad i\in\{1,2\}.
\]
Since ${r_{\rm geo}}=Sy_*$ and the reference edge length is $L_0$,
\[
\begin{aligned}
Z_i&=e_i\cdot(v-c_i)={r_{\rm geo}} e_i\cdot(v^0-c_i^0)=\zeta_i y_*,\\
C_v&=B(v,2\rho_p)=B(v,Ry_*),\qquad
R=2S\rho_g<SL_0/2=|\zeta_i|.
\end{aligned}
\]
Thus the coordinate and radius hypotheses of Proposition~\ref{ce:incoming}
hold. Its signed angle satisfies $|\alpha_v|=\alpha_0=\pi/32$, and the incident types are $0$ and $\gamma$.
Choose
\[
M_{\rm merge}\triangleq\max_{v^0}M_0(\alpha_0,R,\zeta_1,\zeta_2),
\]
where the maximum is over the finitely many reference vertices, with
$\zeta_1,\zeta_2$ evaluated at each vertex. Then
$M\ge M_{\rm merge}$ permits all vertex constructions simultaneously.
Each whole edge retains its fixed midpoint, frame and type at both ends,
so these constructions preserve the same original edge pair.

Choose $C(S)$ to dominate the finitely many index losses $C_*$ in the
band comparison of Subsection~\ref{ce:small-owner}. At positive heights on the scaled
geometry this gives
\[
m_i\ge M-C_*\ge M-C(S),\qquad
M\ge C(S)+14\quad\Longrightarrow\quad m_i\ge14.
\]
Let $M_{\rm act}$ be a common threshold for the finitely many estimates in
Steps~\ref{step:activation-gradient}--\ref{step:activation-global-data} of
Proposition~\ref{ce:activation}. Increase this threshold, if needed, to
include the entrance-height bounds used in Step~\ref{step:local-constant-exit}
below. With the reference geometry and local solution family in Step~\ref{step:local-common-parameters} fixed, both
$M_{\rm merge}$ and $M_{\rm act}$ are independent of the final $M$.
\Needspace{10\baselineskip}
The threshold controls the following conditions.
\begin{center}\small
\begin{tabular}{@{}>{\raggedright\arraybackslash}p{.22\textwidth}>{\raggedright\arraybackslash}p{.71\textwidth}@{}}
\toprule Conditions controlled&Specific locations\\\midrule
first-derivative errors&\eqref{ce:principal-amplitude-derivatives}--\eqref{ce:principal-normalized-expansions} in Step~\ref{step:activation-gradient}\\
vertex curve and affine ends&\eqref{ce:relative-phase}--\eqref{ce:activation-graph-deviation} and \eqref{ce:activation-bend-amplitude-ratio}\\
level band and overlaps&compactness and entrance bounds after \eqref{ce:type-heights}--\eqref{ce:activation-global-level} in Step~\ref{step:activation-global-data}\\
constant-exit height margins&$F_v\le2y_*/c_\alpha$ and $y_*/2\le h_i\le3y_*/2$ in Step~\ref{step:local-constant-exit}\\\bottomrule
\end{tabular}
\end{center}
Impose also
\begin{equation}\label{ce:static-margins}
\begin{aligned}
&M\ge C(S)+14,\qquad
{r_{\rm geo}}\operatorname{diam}\Omega_0<\min\{1/100,Y_*\},\\
&2\ell<\min\{8c_\alpha-2/c_\alpha,\ 2/c_\alpha,\ 6/c_\alpha\}y_*,
\qquad 2\gamma M\ell<8y_*.
\end{aligned}
\end{equation}
Here $c_\alpha=\cos(\alpha_0/2)$ is defined in
\eqref{ce:activation-bisector}, and $\ell=y_*^2$ in
\eqref{ce:activation-scales}. The diameter bound keeps the geometry
inside the domains of the local solutions. The last line supplies the
three positive margins used in Step~\ref{step:local-constant-exit} and
puts $8y_*$ beyond the bending interval of the vertex graph.
All conditions can hold at once, because
\[
\begin{aligned}
{r_{\rm geo}}\operatorname{diam}\Omega_0 =S\operatorname{diam}\Omega_0\,\beta^{-M}\rightarrow0,\quad
\ell/y_* =\beta^{-M}\rightarrow0,\quad
M\ell/y_*=M\beta^{-M}\rightarrow0.
\end{aligned}
\]
Fix one integer $M\ge\max\{M_{\rm merge},M_{\rm act}\}$ satisfying
\eqref{ce:static-margins}, and apply Propositions~\ref{ce:incoming}
and~\ref{ce:activation} with this choice.

\medskip\noindent\textbf{Step 2. Define the compact core inside the same neighborhood $U$.}\proofsteplabel{step:local-regions}{2}
For the observation polygons $P_j$, $j\in\{1,\ldots,J\}$, set
\begin{equation}\label{ce:global-neighborhoods}
\begin{aligned}
&U_j\triangleq\{x:\operatorname{dist}(x,P_j)<\delta_p\},\qquad
K_j\triangleq\{x:\operatorname{dist}(x,P_j)\le8y_*\},\\
&U=\bigcup_{j=1}^J U_j,\qquad K\triangleq\bigcup_{j=1}^J K_j.
\end{aligned}
\end{equation}
This formula for $U$ is exactly the one in Step~\ref{step:activation-global-data}
of Proposition~\ref{ce:activation}. Recall also
\[
U^+=\bigcup_{j=1}^J\{x:\operatorname{dist}(x,P_j)<3\delta_p/2\},\qquad
\overline{\mathcal N_{\delta_p/4}(\overline U)}\subset U^+\Subset\Omega.
\]
The reference gaps exceed $d_0$, and
$32y_*<\delta_p<{r_{\rm geo}} d_0\sin\alpha_0/10^4$. Consequently,
\[
\begin{aligned}
\operatorname{dist}(K,\R^2\setminus U)&\ge\delta_p-8y_*>24y_*,\\
\operatorname{dist}(\overline U,\partial\Omega)&\ge {r_{\rm geo}}d_0-\delta_p>0,\\
\operatorname{dist}(\overline U_j,\overline U_k)&\ge {r_{\rm geo}}d_0-2\delta_p>0
\quad(j\ne k).
\end{aligned}
\]
Thus $K_j\Subset U_j$, $K\Subset U\Subset\Omega$, and the closures of
the $U_j$ are pairwise disjoint. The roles of these regions are:
\begin{center}
\small
\begin{tabular}{@{}p{.19\textwidth}p{.35\textwidth}p{.40\textwidth}@{}}
\toprule
Region&Location&Use in the construction\\
\midrule
$\omega$&polygon interiors&retain the zero pair\\
$K=\bigcup_jK_j$&closed $8y_*$ neighborhoods&contain the local corrections\\
$U=\bigcup_jU_j$&open $\delta_p$ neighborhoods&domain of the completed local pair\\
$U\setminus K$&between the two widths&prove the pair equals $(1,0)$\\
$\Omega\setminus K$&exterior of the compact core&construct the Dirichlet exterior pair\\
\bottomrule
\end{tabular}
\end{center}

\medskip\noindent\textbf{Step 3. Retain the zero pair on the observation set.}\proofsteplabel{step:local-zero}{3}
Denote the completed pair constructed on $U^+$ in
Proposition~\ref{ce:activation} by $(p,V)$, and define
\[
(p_{\rm loc},V_{\rm loc})\triangleq(p,V)\big|_{\R\times U}.
\]
It is the pair after completion; $(p_{\rm in},V_{\rm in})$ denotes the
pair supplied to that construction. Their relation is
\[
\Box p_{\rm loc}+V_{\rm loc}p_{\rm loc}=0,\qquad
(p_{\rm loc},V_{\rm loc})=(p_{\rm in},V_{\rm in})\quad(X\le0).
\]
Step~\ref{step:activation-global-data} of that proposition gives
\[
\omega=\bigcup_{j=1}^J\operatorname{int}P_j\subset\{X=-2\},\qquad
(p_{\rm in},V_{\rm in})\big|_{\R\times\omega}=(0,0).
\]
Therefore $(p_{\rm loc},V_{\rm loc})=(0,0)$ on $\R\times\omega$.

\medskip\noindent\textbf{Step 4. Prove that the pair is constant on $U\setminus K$.}\proofsteplabel{step:local-constant-exit}{4}
The matching conclusion of Proposition~\ref{ce:activation} reduces this
task to the spatial inclusion
$
U\setminus K\subset\{X>2\}.
$
Recall the local levels from \eqref{ce:activation-graph} and \eqref{ce:edge-level}:
\[
X_v=\frac{\eta-F_v(d)}{\ell},\qquad
X_i=\frac{y_i-h_i(z_i)}{c_\alpha\ell}.
\]
We first bound the entrance heights. By \eqref{ce:activation-base-graph},
\[
H_{\rm base}'\ge0\ (d<d_L),\qquad
H_{\rm base}=\eta_a\ (d_L\le d\le d_R),\qquad
H_{\rm base}'\le0\ (d>d_R).
\]
Hence $H_{\rm base}(d)\le\eta_a$ for every $d$, and
\eqref{ce:activation-graph-deviation} gives
\[
\begin{aligned}
F_v(d)&\le H_{\rm base}(d)+CM\epsilon^2\ell
\le\eta_a+CM\epsilon^2\ell =y_*/c_\alpha+O(\lambda^{-1})+CM\epsilon^2\ell
\le2y_*/c_\alpha.
\end{aligned}
\]
Also $y_*/2\le h_i\le3y_*/2$ by \eqref{ce:type-heights} and the
convex interpolation defining $h_i$.

Fix $x\in U_j\setminus K_j$. Since $P_j$ is convex, it has a unique
nearest point $x_P\in P_j$, and
\[
r\triangleq|x-x_P|=\operatorname{dist}(x,P_j)\in(8y_*,\delta_p).
\]
\Needspace{8\baselineskip}
The following three cases are disjoint and exhaustive:
\begin{center}
\small
\begin{tabular}{@{}cp{.49\textwidth}p{.37\textwidth}@{}}
\toprule
Case&Nearest point and chart&Quantity to estimate\\
\midrule
1&$x_P=v$ is a vertex&$X_v$, using $\eta\ge c_\alpha r$\\
2&$x_P$ is inside edge $i$, $x\in C_i$&$X_i$, using $y_i=r$\\
3&$x_P$ is inside edge $i$, $x\notin C_i$&$X_v$ at an endpoint, using one $y_i=r$\\
\bottomrule
\end{tabular}
\end{center}
Here ``inside edge'' means in its relative interior.

\smallskip\noindent\textit{Case 1: nearest point at a vertex.}
Since $r<\delta_p<2\rho_p$, we have $x\in C_v$. The normal cone at the
vertex is generated by the two outward normals, so
\[
\begin{aligned}
x-v&=a_1n_1+a_2n_2,\qquad a_1,a_2\ge0,\\
r&=|a_1n_1+a_2n_2|\le a_1+a_2,\\
\eta&=n_*\cdot(x-v)=c_\alpha(a_1+a_2)\ge c_\alpha r.
\end{aligned}
\]
The height bound and \eqref{ce:static-margins} imply
\[
X_v(x)>\frac{(8c_\alpha-2/c_\alpha)y_*}{\ell}>2.
\]

\smallskip\noindent\textit{Case 2: nearest point inside an edge, with $x\in C_i$.}
Orthogonal projection onto the edge gives
\[
x-x_P=rn_i,\qquad y_i=r,\qquad
X_i(x)=\frac{r-h_i(z_i)}{c_\alpha\ell}
>\frac{6y_*}{c_\alpha\ell}>2.
\]

\smallskip\noindent\textit{Case 3: nearest point inside an edge, with $x\notin C_i$.}
The normal bound $r<\delta_p<2\delta_p$ already satisfies the first
condition defining $C_i$. Thus its tangential condition fails: the
edge projection $x_P$ is within $\rho_p/2$ of an endpoint $v$. In particular,
\[
|x-v|\le\sqrt{(\rho_p/2)^2+\delta_p^2}<2\rho_p,
\]
so $x\in C_v$. Relabel its two incident edges by types $0,\gamma$,
as in \eqref{ce:activation-edge-coordinates}. One of $y_1,y_2$ equals $r$.
If both are nonnegative, then
\[
\eta=\frac{y_1+y_2}{2c_\alpha}\ge\frac r{2c_\alpha}
>\frac{4y_*}{c_\alpha},\qquad
X_v(x)>\frac{2y_*}{c_\alpha\ell}>2.
\]
If the other coordinate is negative, then
\[
|d|=|y_1-y_2|>r>8y_*>2\gamma M\ell.
\]
Thus the point lies beyond the bending interval. The affine formulas
computed in Step~\ref{step:activation-global-data} of
Proposition~\ref{ce:activation} give
\[
X_v(x)=\begin{cases}
(y_1-h_0(v))/(c_\alpha\ell),&d>0,\quad y_1=r,\\
(y_2-h_\gamma(v))/(c_\alpha\ell),&d<-2\gamma M\ell,\quad y_2=r.
\end{cases}
\]
In either case,
\[
X_v(x)\ge\frac{r-3y_*/2}{c_\alpha\ell}
>\frac{6y_*}{c_\alpha\ell}>2.
\]

The three exit cases in Step~\ref{step:local-constant-exit} give a local level greater than $2$ at every point of
$U\setminus K$. The scalar function $h_{\rm lev}$ in
\eqref{ce:activation-global-level} is the identity on $[-1,\infty)$,
so the global level also satisfies $X>2$ there. Therefore
\[
(p_{\rm loc},V_{\rm loc})=(1,0)\quad\hbox{on }\R\times(U\setminus K).
\]
Finally, Proposition~\ref{ce:activation} supplies smoothness on a
neighborhood of $\overline U$ and $2\pi$-periodicity. Thus, for every
fixed multi-index $\alpha$,
\[
\sup_{\R\times\overline U}
(|\partial^\alpha p_{\rm loc}|+|\partial^\alpha V_{\rm loc}|)
=\max_{[0,2\pi]\times\overline U}
(|\partial^\alpha p_{\rm loc}|+|\partial^\alpha V_{\rm loc}|)<\infty.
\]
\end{proof}

\paragraph{Parameter choices.}
Table~\ref{tab:counterexample-ledger} records the order of the choices.
The constants in Lemmas~\ref{ce:small-alternatives} and~\ref{ce:comparable-windows}
are included in the fixed local thresholds; the final integer $M$ is
chosen after all of them.
\begin{table}[H]
\centering
\caption{Parameter selection order for the counterexample in Theorem~\ref{thm:counter}.}
\label{tab:counterexample-ledger}
\small
\begin{tabular}{@{}p{.27\textwidth}p{.29\textwidth}p{.38\textwidth}@{}}
\toprule
Parameters&Fixed after&Purpose\\
\midrule
Local solution family; reference geometry
&Sections~\ref{ce:front-section} and~\ref{ce:geometry-section}
&Fix $Y_*,\beta,\gamma$ and $\Omega_0,P_j^0,d_0,\alpha_0,L_0$.\\
$S,\rho_g,\delta_g$
&Reference geometry
&Require $32/S<\rho_g\sin\alpha_0/100$ and choose $\delta_g$ as in \eqref{ce:reference-collar-scales}.\\
$R=2S\rho_g$, $\zeta_i=S e_i\cdot(v^0-c_i^0)$
&$S,\rho_g$ and the reference edge midpoints
&Fix the vertex data in Proposition~\ref{ce:incoming}.\\
$M_{\rm merge},M_{\rm act}$
&Fixed data in Step~\ref{step:local-common-parameters}
&Take the maximum vertex threshold and the completion threshold from Step~\ref{step:local-common-parameters}; both are independent of the final $M$.\\
$M$; then $y_*,{r_{\rm geo}},\ell,\rho_p,\delta_p$
&$M_{\rm merge},M_{\rm act},C(S)$
&Choose $M$ satisfying \eqref{ce:static-margins} and both thresholds; set the spatial lengths by \eqref{ce:final-scales} and \eqref{ce:activation-scales}.\\
\bottomrule
\end{tabular}
\end{table}

\begin{lemma}[Positive static completion at the Dirichlet boundary]\label{ce:boundary}
Let $n\ge1$, let $\Omega\subset\R^n$ be a bounded connected domain
with smooth boundary, and let $U\Subset\Omega$ be open. There exist
$0<\varepsilon<\operatorname{dist}(\overline U,\partial\Omega)/4$ and
$a_{\rm out}\in C^\infty(\overline\Omega)$, positive in $\Omega$ and zero
on $\partial\Omega$, such that $a_{\rm out}=1$ on
$\mathcal N_\varepsilon(\overline U)$ and
$V_{\rm out}\triangleq\Delta a_{\rm out}/a_{\rm out}$ extends smoothly
to $\overline\Omega$. Moreover,
$(a_{\rm out},V_{\rm out})=(\varphi_1,-\lambda_1)$ in a
neighborhood of $\partial\Omega$ inside $\Omega$, where
$\varphi_1>0$ is a first Dirichlet eigenfunction and $\lambda_1$
its eigenvalue.
\end{lemma}

\begin{proof}[Proof of Lemma~\ref{ce:boundary}]
We construct a positive function which is one near $\overline U$ and
is a Dirichlet eigenfunction near $\partial\Omega$. The eigenfunction
equation will make the resulting potential constant near the boundary.
Let $A_\Omega=-\Delta_D$, and choose a positive first eigenfunction
$\varphi_1\in C^\infty(\overline\Omega)$ with
\[
-\Delta\varphi_1=\lambda_1\varphi_1\quad\hbox{in }\Omega,
\qquad \varphi_1=0\quad\hbox{on }\partial\Omega,
\qquad \lambda_1=\min\sigma(A_\Omega)>0.
\]
Fix $0<\varepsilon<\operatorname{dist}(\overline U,\partial\Omega)/4$
and choose $\chi\in C_c^\infty(\Omega;[0,1])$ such that
\[
\chi=1\quad\hbox{on }\mathcal N_\varepsilon(\overline U),\qquad
\chi=0\quad\hbox{on }\Omega\setminus\mathcal N_{2\varepsilon}(\overline U).
\]
Define
\begin{equation}\label{ce:boundary-profile}
a_{\rm out}\triangleq\chi+(1-\chi)\varphi_1.
\end{equation}
At every $x\in\Omega$,
\[
a_{\rm out}(x)\ge\min\{1,\varphi_1(x)\}>0,
\qquad
a_{\rm out}=1\quad\hbox{on }\mathcal N_\varepsilon(\overline U).
\]
The product rule gives
\[
\begin{aligned}
\Delta a_{\rm out}
&=\Delta\chi-\varphi_1\Delta\chi
-2\nabla\chi\cdot\nabla\varphi_1+(1-\chi)\Delta\varphi_1 =(1-\varphi_1)\Delta\chi
-2\nabla\chi\cdot\nabla\varphi_1
-\lambda_1(1-\chi)\varphi_1.
\end{aligned}
\]
Consequently, the required potential in $\Omega$ is
\begin{equation}\label{ce:boundary-potential}
V_{\rm out}\triangleq
\frac{(1-\varphi_1)\Delta\chi
-2\nabla\chi\cdot\nabla\varphi_1
-\lambda_1(1-\chi)\varphi_1}
{\chi+(1-\chi)\varphi_1}.
\end{equation}
This quotient is smooth in the interior. To check it at the boundary,
put $C_\varepsilon\triangleq
\{x\in\Omega:\operatorname{dist}(x,\partial\Omega)<\varepsilon\}$.
For $x\in C_\varepsilon$,
\[
\operatorname{dist}(x,\overline U)
\ge\operatorname{dist}(\overline U,\partial\Omega)
-\operatorname{dist}(x,\partial\Omega)>3\varepsilon.
\]
Thus $\chi=0$ near $C_\varepsilon$, and
\[
(a_{\rm out},V_{\rm out})=(\varphi_1,-\lambda_1)
\quad\hbox{on }C_\varepsilon.
\]
The constant value $-\lambda_1$ extends $V_{\rm out}$ smoothly across
$\partial\Omega$; no division by the boundary value
$a_{\rm out}=0$ is needed. Finally, \eqref{ce:boundary-profile}--\eqref{ce:boundary-potential} give
\[
\begin{gathered}
a_{\rm out}|_{\partial\Omega}=0,\qquad
(a_{\rm out},V_{\rm out})=(1,0)
\quad\hbox{on }\mathcal N_\varepsilon(\overline U),\\
-\Delta a_{\rm out}+V_{\rm out}a_{\rm out}=0
\quad\hbox{in }\Omega.
\end{gathered}
\]
\end{proof}

\begin{proof}[Proof of Theorem~\ref{thm:counter}]
\medskip\noindent\textbf{Step 1. Glue on an open cover and verify the equation.}\proofsteplabel{step:counter-glue}{1}
Take $K\Subset U$ and $(p_{\rm loc},V_{\rm loc})$ from
Lemma~\ref{ce:global-local}, and apply Lemma~\ref{ce:boundary} to this
$U$. The sets and their matching data satisfy
\[
U\cup(\Omega\setminus K)=\Omega,\qquad
U\cap(\Omega\setminus K)=U\setminus K,
\]
\[
(p_{\rm loc},V_{\rm loc})=(1,0)
=(a_{\rm out},V_{\rm out})\quad\hbox{on }\R\times(U\setminus K).
\]
Both sets in the cover are open because $K$ is compact. We may therefore
define one smooth pair by
\begin{equation}\label{ce:global-pair}
(p,V)\triangleq
\begin{cases}
(p_{\rm loc},V_{\rm loc}),&x\in U,\\
(a_{\rm out},V_{\rm out}),&x\in\Omega\setminus K.
\end{cases}
\end{equation}
On these two open sets, respectively,
\[
\Box p+Vp=
\begin{cases}
\Box p_{\rm loc}+V_{\rm loc}p_{\rm loc}=0,&x\in U,\\
-\Delta a_{\rm out}+V_{\rm out}a_{\rm out}=0,&x\in\Omega\setminus K.
\end{cases}
\]
Since $\overline U\Subset\Omega$, the global pair equals the exterior
pair in a boundary collar. Lemma~\ref{ce:boundary} then yields
\[
p,V\in C^\infty(\R\times\overline\Omega),\qquad
p|_{\R\times\partial\Omega}=0,\qquad
p|_{\R\times\omega}=p_{\rm loc}|_{\R\times\omega}=0.
\]

\medskip\noindent\textbf{Step 2. Verify nontriviality, periodicity, and derivative bounds.}\proofsteplabel{step:counter-bounds}{2}
Choose an open ball $B$ whose closure lies in a nonempty boundary collar
where $p(t,x)=a_{\rm out}(x)=\varphi_1(x)>0$. For every $t\in\R$,
\[
\|p(t)\|_2^2\ge\int_Ba_{\rm out}(x)^2\,\mathrm dx
\ge |B|\min_{\overline B}a_{\rm out}^2>0.
\]
The local pair is $2\pi$-periodic and the exterior pair is independent
of time, so the glued pair is $2\pi$-periodic. Compactness now gives,
for each fixed spacetime multi-index $\alpha$,
\[
\sup_{\R\times\overline\Omega}
(|\partial^\alpha p|+|\partial^\alpha V|)
=\max_{[0,2\pi]\times\overline\Omega}
(|\partial^\alpha p|+|\partial^\alpha V|)<\infty.
\]
This proves \eqref{eq:counter-bounds}.

\medskip\noindent\textbf{Step 3. Identify the unique weak evolution.}\proofsteplabel{step:counter-weak}{3}
We next check the mild-solution assertion in the theorem. Recall
$A_\Omega=-\Delta_D$ and
$X_\Omega=L^2(\Omega)\times H^{-1}(\Omega)$ with the norm
\eqref{eq:state-norm}. Write $S(t)$ for the free wave group on this
space. In spectral form,
\[
S(t)\binom{a}{b}=
\begin{pmatrix}
\cos(t\sqrt{A_\Omega})a+A_\Omega^{-1/2}\sin(t\sqrt{A_\Omega})b\\
-A_\Omega^{1/2}\sin(t\sqrt{A_\Omega})a+\cos(t\sqrt{A_\Omega})b
\end{pmatrix}.
\]
\[
\|S(t)(a,b)\|_{X_\Omega}=\|(a,b)\|_{X_\Omega}.
\]
Define the bounded operator $K_V(t)(a,b)\triangleq(0,-V(t,\cdot)a)$.
Its norm satisfies
\[
\begin{aligned}
\|K_V(t)(a,b)\|_{X_\Omega}
&=\|A_\Omega^{-1/2}(V(t)a)\|_2 \le\lambda_1^{-1/2}\|V(t)\|_\infty\|a\|_2 \le\lambda_1^{-1/2}\|V(t)\|_\infty\|(a,b)\|_{X_\Omega}.
\end{aligned}
\]
For the smooth Dirichlet solution already constructed, put
$U_p(t)\triangleq(p(t),p_t(t))$. Its equation gives
\[
\begin{aligned}
&U_p'(t) =\mathcal A_\Omega U_p(t)+K_V(t)U_p(t),
\qquad\mathcal A_\Omega(a,b) \triangleq(b,-A_\Omega a),\\
&\frac{\mathrm d}{\mathrm dt}\bigl(S(-t)U_p(t)\bigr)
=S(-t)\bigl(U_p'(t)-\mathcal A_\Omega U_p(t)\bigr)
=S(-t)K_V(t)U_p(t).
\end{aligned}
\]
Integrating and multiplying by $S(t)$ yields the state Duhamel identity
and its first component:
\[
\begin{aligned}
&U_p(t) =S(t)U_p(0)+\int_0^tS(t-s)K_V(s)U_p(s)\,\mathrm ds,\\
&p(t) =\cos(t\sqrt{A_\Omega})p(0)
+A_\Omega^{-1/2}\sin(t\sqrt{A_\Omega})p_t(0) -\int_0^tA_\Omega^{-1/2}\sin((t-s)\sqrt{A_\Omega}) [V(s)p(s)]\,\mathrm ds.
\end{aligned}
\]
To verify uniqueness, let $U_1,U_2$ be two continuous $X_\Omega$-valued
solutions of the state identity on any interval $[t_0,T]$, with the
initial time changed to $t_0$. Set
$L_{t_0,T}\triangleq\lambda_1^{-1/2}
\max_{t_0\le s\le T}\|V(s)\|_\infty$. Unitarity and the bound on $K_V$
give
\begin{equation}\label{ce:forward-evolution-bound}
\begin{aligned}
D(t)&\triangleq\|U_1(t)-U_2(t)\|_{X_\Omega},\\
D(t)&\le D(t_0)+L_{t_0,T}\int_{t_0}^tD(s)\,\mathrm ds
\le e^{L_{t_0,T}(t-t_0)}D(t_0),\qquad t_0\le t\le T.
\end{aligned}
\end{equation}
The last inequality follows from Gronwall's inequality. With $t=-\tau$, the reversed free group
$S(-\tau)$ and perturbation $-K_V(-\tau)$ give the reverse-time version of \eqref{ce:forward-evolution-bound}:
\[
\|U_1(-\tau)-U_2(-\tau)\|_{X_\Omega}
\le e^{L_{-T,0}\tau}\|U_1(0)-U_2(0)\|_{X_\Omega},\qquad 0\le\tau\le T.
\]
Thus $U_p$ is the unique weak evolution for its initial state on every
finite interval.

\medskip\noindent\textbf{Step 4. Disprove observability for the fixed potential.}\proofsteplabel{step:counter-observation}{4}
Proposition~\ref{ce:geometry}\ref{cond:geometry-free} gives GCC and the
free observability estimate~\eqref{ce:free-scaled} for $(\Omega,\omega)$.
For our fixed $V$, Steps~\ref{step:counter-glue}--\ref{step:counter-bounds} give, for every $T>0$,
\[
\int_0^T\!\int_\omega|p(t,x)|^2\,\mathrm dx\,\mathrm dt=0,
\qquad
\|U_p(0)\|_{X_\Omega}^2\ge\|p(0)\|_2^2>0.
\]
An observability inequality would therefore imply the contradiction
\[
0<\|U_p(0)\|_{X_\Omega}^2
\le C_{T,V}\int_0^T\!\int_\omega|p(t,x)|^2\,\mathrm dx\,\mathrm dt=0.
\]
Hence observability fails for this potential at every finite time.
\end{proof}

\section{Counterexamples in higher dimensions}
\label{sec:general-dimensions}

This section proves Theorem~\ref{thm:counter-high-dim}. Fix $n\ge3$
throughout; all constants below may depend on this fixed dimension.

The proof has two independent geometric inputs and one analytic input.
On the analytic side, we lift the planar local solution of
Lemma~\ref{ce:global-local} to an $n$-dimensional pair which vanishes on
a ball and equals $(1,0)$ outside a compact set
(Proposition~\ref{mod:nd-module}). On the geometric side, we construct a
bounded star-shaped domain with real-analytic boundary such that every
generalized ray leaves a thin boundary collar within a uniform time
(Proposition~\ref{geo:escape}), and then place finitely many separated
interior balls which intercept every sufficiently long straight segment
(Lemma~\ref{asm:GCC}). The two sides are joined in
Subsection~\ref{sec:assembly}: the local pairs are glued to the exterior
Dirichlet solution supplied by Lemma~\ref{ce:boundary}, exactly as in the
planar case.

\begin{figure}[!htbp]
	\centering
	\begin{tikzpicture}[x=1cm,y=1cm,font=\small,
		box/.style={draw=black!55,rounded corners=2pt,align=center,
			text width=3.35cm,minimum height=1.05cm,inner sep=4pt},
		flow/.style={-{Stealth[length=2mm]},line width=.7pt}]
		\node[box,fill=orange!8] (dimples) at (0,0)
		{Boundary dimples\\Lemma~\ref{geo:concave-capture}};
		\node[box,fill=orange!8] (escape) at (4.4,0)
		{Rays reach the interior\\Proposition~\ref{geo:escape}};
		\node[box,fill=orange!8] (gcc) at (8.8,0)
		{Interior balls give GCC\\Lemma~\ref{asm:GCC}};
		\node[box,fill=blue!5] (local) at (0,-2)
		{Local zero ball\\Proposition~\ref{mod:nd-module}};
		\node[box,fill=green!7] (assembly) at (4.4,-2)
		{Dirichlet solution\\Theorem~\ref{thm:counter-high-dim}};
		\node[box,fill=blue!5] (boundary) at (8.8,-2)
		{Dirichlet completion\\Lemma~\ref{ce:boundary}};
		\draw[flow] (dimples)--(escape);
		\draw[flow] (escape)--(gcc);
		\draw[flow] (local)--(assembly);
		\draw[flow] (boundary)--(assembly);
		\draw[flow] (gcc.south)--(assembly.north east);
	\end{tikzpicture}
	\caption{Inputs to the higher-dimensional construction. The upper row
		chooses the observation geometry. The lower row supplies and joins the
		local solutions without changing their zero balls.}
	\label{fig:dimension-construction-roadmap}
\end{figure}
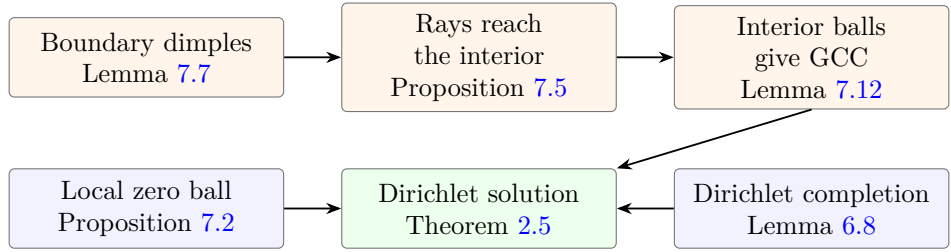

Figure~\ref{fig:dimension-construction-roadmap} records the dependencies
between these ingredients.

The section is organized as follows.
Subsection~\ref{sec:modules} constructs the higher-dimensional local
solutions, first lifting the planar pair and then removing the radial
first derivative by a positive multiplier.
Subsection~\ref{sec:geometry} builds the observation geometry: separated
balls intercepting straight segments, concave boundary regions reached by
every boundary geodesic, and the resulting escape from a thin collar.
Subsection~\ref{sec:assembly} glues the local and exterior pairs on an
open cover and proves Theorem~\ref{thm:counter-high-dim}.

\subsection{Construction of higher-dimensional local solutions}
\label{sec:modules}

The planar case is already available from Lemma~\ref{ce:global-local}.
We record it first, since the higher-dimensional construction starts from
it.
\begin{corollary}[A planar solution equal to one outside a compact set]
	\label{mod:planar}
	There exist a compact set $K_2\subset\R^2$, a nonempty open ball
	$O_2\subset\R^2$, and real-valued functions
	$u,W\in C^\infty(\R\times\R^2)$, both $2\pi$-periodic in time, with the
	following properties. First, the pair $(u,W)$ solves the planar wave
	equation on the whole spacetime：
	\begin{equation}\label{mod:planar-equation}
		\begin{cases}
			(\partial_t^2-\Delta_2)u+Wu=0 &\text{in }\R\times\R^2,\\ 
			\ns\ds 
			(u,W)=(1,0) &\text{on }\R\times(\R^2\setminus K_2),\\ 
			\ns\ds  u=0 &\text{on }\R\times O_2.
		\end{cases}
	\end{equation}
	Thus $u$ is the constant solution $1$ with zero potential outside the
	compact set $K_2$, while it vanishes identically on the observed
	cylinder over $O_2$. Second, every space-time derivative of $u$ and $W$
	is globally bounded on $\R\times\R^2$.
\end{corollary}

\begin{proof}[Proof of Corollary~\ref{mod:planar}]
	We build the global pair by extending the local pair of
	Lemma~\ref{ce:global-local} by the constant pair $(1,0)$ outside a
	compact set. Four points have to be checked: the two definitions agree
	on their overlap, so that the extension is smooth; the extension is
	$2\pi$-periodic; it solves the equation everywhere; and all its
	space-time derivatives are globally bounded.
	
	\medskip\noindent\textbf{Step 1. Define the extension.}
	Take $K\Subset U$ and $(u_{\rm loc},W_{\rm loc})$ from
	Lemma~\ref{ce:global-local}, and set $K_2=K$. Choose a nonempty ball
	$O_2$ contained in one of the observation polygons, so that
	$O_2\subset\omega$. Define
	\begin{equation}\label{mod:planar-extension}
		(u,W)(t,x)=
		\begin{cases}
			(u_{\rm loc},W_{\rm loc})(t,x),&x\in U,\\
			(1,0),&x\in\R^2\setminus K.
		\end{cases}
	\end{equation}
	Since $K\Subset U$, the sets $U$ and $\R^2\setminus K$ are open and
	cover $\R^2$. Their intersection is exactly $U\setminus K$.
	
	\medskip\noindent\textbf{Step 2. The definitions agree on the overlap.}
	By \eqref{ce:local-exit},
	\[
	(u_{\rm loc},W_{\rm loc})=(1,0)\qquad\text{on }\R\times(U\setminus K).
	\]
	This is precisely the value assigned by the second branch of
	\eqref{mod:planar-extension} on the overlap. Hence the two branches
	agree on $U\setminus K$, and the gluing lemma for smooth functions on
	an open cover shows that $(u,W)\in C^\infty(\R\times\R^2;\R^2)$.
	
	\medskip\noindent\textbf{Step 3. Periodicity and the equation.}
	The local pair $(u_{\rm loc},W_{\rm loc})$ is $2\pi$-periodic by
	Lemma~\ref{ce:global-local}, and the constant pair $(1,0)$ is trivially
	$2\pi$-periodic. Since the two branches agree on the overlap, the
	extension $(u,W)$ is $2\pi$-periodic on all of $\R\times\R^2$.
	
	For the equation, argue on the same two open sets. On $\R\times U$,
	\[
	(\partial_t^2-\Delta_2)u+Wu
	=(\partial_t^2-\Delta_2)u_{\rm loc}+W_{\rm loc}u_{\rm loc}=0
	\]
	by \eqref{ce:local-exit}. On $\R\times(\R^2\setminus K)$ we have
	$(u,W)=(1,0)$, hence
	\[
	(\partial_t^2-\Delta_2)u+Wu=0+0\cdot1=0.
	\]
	Since the two open sets cover $\R\times\R^2$, the equation holds
	everywhere.
	
	\medskip\noindent\textbf{Step 4. Global boundedness of the derivatives.}
	Fix a space-time multi-index $\alpha$. On $\R\times U$ the pair equals
	$(u_{\rm loc},W_{\rm loc})$, whose derivatives are bounded on
	$\R\times\overline U$ by Lemma~\ref{ce:global-local}; on
	$\R\times(\R^2\setminus K)$ the pair is the constant $(1,0)$, whose
	derivatives of order $|\alpha|\ge1$ vanish and whose value is bounded
	by $1$. Using that $U$ and $\R^2\setminus K$ cover $\R^2$, we obtain
	\[
	\sup_{\R\times\R^2}
	(|\partial^\alpha u|+|\partial^\alpha W|)
	\le 1+\max_{[0,2\pi]\times K_2}
	(|\partial^\alpha u|+|\partial^\alpha W|)<\infty.
	\]
	The maximum over $[0,2\pi]\times K_2$ is finite because
	$(u,W)\in C^\infty(\R\times\R^2;\R^2)$ and $K_2$ is compact, and it
	controls the supremum over all of $\R\times\R^2$ by
	$2\pi$-periodicity.
	
	\medskip
	Finally, on $\R\times O_2$ we have $O_2\subset\omega$, so the local
	pair vanishes there by \eqref{ce:local-exit}; this gives
	$u=0$ on $\R\times O_2$, as required in \eqref{mod:planar-equation}.
\end{proof}

We now lift this planar pair. Write the planar coordinates as $(\rho,z)$
and translate them so that $K_2\subset\{\rho>2\}$. For
$x=(y,z)\in\R^{n-1}\times\R$ put $\rho=|y|$ and set
\begin{equation}\label{mod:dimension-constants}
	m=n-2,\qquad \kappa_n=\frac{m(m-2)}4.
\end{equation}
The following proposition is the analytic input for the whole section.
\begin{proposition}[A local solution in dimension $n\ge3$]
	\label{mod:nd-module}
	Fix $n\ge3$. There exist radii $0<r<r_0<\infty$, a compact set
	$D\Subset B_{r_0}(0)$, and real-valued functions
	$P,Z\in C^\infty(\R\times\R^n)$, both $2\pi$-periodic in time, with the
	following properties. First, the pair $(P,Z)$ solves the
	$n$-dimensional wave equation on the whole spacetime,
	\begin{equation}\label{mod:module}
		\begin{cases}
			\ds(\partial_t^2-\Delta_n)P+ZP=0 &\text{in }\R\times\R^n,\\[2pt]
			\ns\ds P=0 &\text{on }\R\times B_r(0),\\[2pt]
			\ns\ds (P,Z)=(1,0) &\text{on }\R\times(\R^n\setminus D).
		\end{cases}
	\end{equation}
	Thus $P$ vanishes identically on the observed cylinder over the ball
	$B_r(0)$, while outside the larger compact set $D$ the pair is the
	constant solution $(1,0)$. Second, every space-time derivative of $P$
	and $Z$ is globally bounded on $\R\times\R^n$.
\end{proposition}
\begin{proof}[Proof of Proposition~\ref{mod:nd-module}]
	The construction has four steps, each addressing one obstruction.
	First, the lifted planar pair does not solve the $n$-dimensional
	equation, because the radial Laplacian produces a first-order term in
	$\rho$; we remove it by conjugating with a positive factor. Second, the
	conjugated pair equals $(f,-\kappa_n/\rho^2)$ outside the lifted
	compact set, not $(1,0)$; we change the factor smoothly to one where the
	planar solution equals one, so that the new potential can still be
	defined by division. Third, the construction is so far valid only away
	from the axis $\rho=0$; we extend it across the axis by the constant
	pair. Fourth, we locate a ball in the zero set and verify periodicity
	and global derivative bounds.
	
	\medskip\noindent\textbf{Step 1. Lift the planar equation away from the axis.}
	\proofsteplabel{step:module-lift}{1}
	The goal of this step is a pair $(P_0,Z_0)$, smooth for $\rho>0$, which
	solves the $n$-dimensional equation and reduces to the planar pair on
	each meridian. For $x=(y,z)\in\R^{n-1}\times\R$ with $\rho=|y|>0$, any
	smooth $H(\rho,z)$ satisfies
	\[
	\partial_{y_j}H=H_\rho\frac{y_j}{\rho},\qquad
	\sum_{j=1}^{n-1}\partial_{y_j}^2H
	=H_{\rho\rho}+\frac m\rho H_\rho .
	\]
	For $fu$ the coefficient of $u_\rho$ is $2f'+m\rho^{-1}f$. To eliminate
	it we need
	\[
	(\log f)'=-\frac{m}{2\rho},
	\qquad\text{so we take}\qquad f(\rho)=\rho^{-m/2}>0 .
	\]
	For this $f$,
	\[
	\frac{f'}f=-\frac{m}{2\rho},\qquad
	\frac{f''}f=\frac{m(m+2)}{4\rho^2},\qquad
	\frac{f''}f+\frac m\rho\frac{f'}f=-\frac{\kappa_n}{\rho^2},
	\]
	and therefore
	\begin{equation}\label{mod:conjugation}
		\Delta_n(fu)
		=f\left(u_{\rho\rho}+u_{zz}-\frac{\kappa_n}{\rho^2}u\right).
	\end{equation}
	Consequently the pair
	\begin{equation}\label{mod:raw-pair}
		P_0(t,y,z)=f(\rho)u(t,\rho,z),\qquad
		Z_0(t,y,z)=W(t,\rho,z)-\frac{\kappa_n}{\rho^2}
	\end{equation}
	satisfies $(\partial_t^2-\Delta_n)P_0+Z_0P_0=0$ for $\rho>0$, by
	\eqref{mod:planar-equation} and \eqref{mod:conjugation}. The coefficient
	$\kappa_n$ in \eqref{mod:conjugation} is exactly the value
	\eqref{mod:dimension-constants}; this is why the correction term
	$-\kappa_n/\rho^2$ appears in $Z_0$.
	
	\medskip\noindent\textbf{Step 2. Make the solution equal to one outside a compact set.}
	\proofsteplabel{step:module-exterior}{2}
	The pair $(P_0,Z_0)$ solves the equation, but outside the lifted $K_2$
	its value is $(f,-\kappa_n/\rho^2)$, not $(1,0)$. We keep $f$ near
	$K_2$ and change it smoothly to one outside a larger set. The change is
	made only where $u=1$, so the new solution stays positive there and its
	potential can still be defined by division.
	
	Choose bounded open sets and a cutoff with
	\[
	K_2\subset N\Subset U_2\Subset\{\rho>1\},\qquad
	\chi\in C_c^\infty(U_2;[0,1]),\qquad \chi=1\text{ on }N,
	\]
	and set, for $\rho>0$,
	\begin{equation}\label{mod:static-multiplier}
		b(\rho,z)=\exp\!\left(-\frac m2\chi(\rho,z)\log\rho\right)>0,
		\qquad P(t,y,z)=b(\rho,z)u(t,\rho,z).
	\end{equation}
	Then $b=f$ on $N$, while $b=1$ outside $\supp\chi$. Define
	\begin{equation}\label{mod:potential-cover}
		Z(t,y,z)=
		\begin{cases}
			W(t,\rho,z)-\dfrac{\kappa_n}{\rho^2},&(\rho,z)\in N,\\[6pt]
			\dfrac{b_{\rho\rho}+m\rho^{-1}b_\rho+b_{zz}}{b},
			&(\rho,z)\in\{\rho>0\}\setminus K_2.
		\end{cases}
	\end{equation}
	The two open sets cover the half-plane $\{\rho>0\}$. On their
	intersection $N\setminus K_2$ we have $u=1$, $W=0$, $b=f$, hence
	\[
	\frac{b_{\rho\rho}+m\rho^{-1}b_\rho+b_{zz}}b
	=\frac{f''+m\rho^{-1}f'}f
	=-\frac{\kappa_n}{\rho^2}
	=W-\frac{\kappa_n}{\rho^2},
	\]
	so the two smooth formulas agree there. On the lift of $N$ we have
	$(P,Z)=(P_0,Z_0)$, and the equation holds by Step~\ref{step:module-lift}.
	On the lift of $\{\rho>0\}\setminus K_2$ we have $P=b$, and the
	definition of $Z$ gives
	\[
	(\partial_t^2-\Delta_n)P+ZP
	=-\bigl(b_{\rho\rho}+m\rho^{-1}b_\rho+b_{zz}\bigr)
	+\frac{b_{\rho\rho}+m\rho^{-1}b_\rho+b_{zz}}b\,b=0.
	\]
	Thus the equation holds on both sets; by the smooth gluing on the open
	cover, it holds on $\{\rho>0\}$.
	
	\medskip\noindent\textbf{Step 3. Extend smoothly across the axis and locate the supports.}
	\proofsteplabel{step:module-axis}{3}
	The construction so far is valid only for $\rho>0$. We now check that it
	is already constant near the axis, so that it extends by the constant
	pair. If $0<|y|<1$, then $(\rho,z)\notin U_2\cup K_2$, so
	\[
	u=1,\quad W=0,\quad b=1,\quad
	b_\rho=b_{\rho\rho}=b_{zz}=0,\quad (P,Z)=(1,0).
	\]
	We may therefore define $(P,Z)=(1,0)$ on the whole tube $\{|y|<1\}$,
	including its axis. This agrees with \eqref{mod:static-multiplier} and
	\eqref{mod:potential-cover} on the open set $\{0<|y|<1\}$, and
	$(\partial_t^2-\Delta_n)1+0\cdot1=0$ there. Hence the extension is
	smooth and satisfies the equation across the axis.
	
	It remains to show that the pair is constant outside a compact set. Let
	$\pi(y,z)=(|y|,z)$ and $C_2=\supp\chi$. Since $C_2$ is compact and
	contained in $\{1<\rho\le L,\ |z|\le L\}$ for some finite $L$, the set
	\[
	D_0=\pi^{-1}(C_2)\subset\{|y|\le L,\ |z|\le L\}
	\]
	is closed and bounded, hence compact. Moreover,
	\[
	x\notin D_0\quad\Longrightarrow\quad
	(u,W)=(1,0),\quad b=1,\quad (P,Z)=(1,0).
	\]
	Thus the spatial supports of $P-1$ and $Z$, uniformly in time, are
	contained in the compact set $D_0$.
	
	\medskip\noindent\textbf{Step 4. Obtain a zero ball and verify periodicity.}
	\proofsteplabel{step:module-zero-ball}{4}
	The positive multiplier preserves the zero set, so it suffices to find a
	ball on which the planar solution $u$ vanishes. Choose
	$B_\eta((\rho_*,z_*))\Subset O_2$ and set
	$x_*=(\rho_*e_1,z_*)\in\R^{n-1}\times\R$. Since
	\[
	|\pi(y,z)-(\rho_*,z_*)|^2
	=(|y|-\rho_*)^2+(z-z_*)^2
	\le |y-\rho_*e_1|^2+(z-z_*)^2=|x-x_*|^2,
	\]
	we have $P=bu=0$ on $\R\times B_\eta(x_*)$. Take $r=\eta/2$ and choose
	$r_0>r$ with $D_0\Subset B_{r_0}(x_*)$; translating $x_*$ to the origin
	and putting $D=D_0-x_*$ gives \eqref{mod:module}.
	
	Since $b$ is independent of time and $u,W$ are $2\pi$-periodic, $P$ and
	$Z$ are $2\pi$-periodic as well. Finally, for every space-time
	multi-index $\alpha$, the pair equals $(1,0)$ outside $D$ and is smooth
	on $\R\times\R^n$, so
	\[
	\sup_{\R\times\R^n}
	(|\partial^\alpha P|+|\partial^\alpha Z|)
	\le 1+\max_{[0,2\pi]\times D}
	(|\partial^\alpha P|+|\partial^\alpha Z|)<\infty.
	\]
\end{proof}

The construction is invariant under simultaneous scaling of time and
space, which preserves the wave operator. For a center $c\in\R^n$ and a
scale $h>0$ define
\begin{equation}\label{mod:scaled}
	P_{h,c}(t,x)=P\!\left(\frac th,\frac{x-c}{h}\right),\qquad
	Z_{h,c}(t,x)=h^{-2}Z\!\left(\frac th,\frac{x-c}{h}\right).
\end{equation}
Writing $\tau=t/h$ and $\xi=(x-c)/h$, differentiation gives
\[
(\partial_t^2-\Delta_n)P_{h,c}+Z_{h,c}P_{h,c}
=h^{-2}(P_{\tau\tau}-\Delta_\xi P+ZP)(\tau,\xi)=0,
\]
and the spatial conditions become
\begin{equation}\label{mod:scaled-properties}
	P_{h,c}=0\text{ on }\R\times B_{hr}(c),\qquad
	\supp_x(P_{h,c}-1),\ \supp_xZ_{h,c}
	\subset c+hD\Subset B_{hr_0}(c).
\end{equation}
Thus the observation and support radii are $hr$ and $hr_0$, and the pair
equals $(1,0)$ near the support sphere, where it can be joined to the
exterior solution. For every $j\ge0$ and spatial multi-index $\alpha$,
\begin{equation}\label{mod:scaling-derivatives}
	\begin{cases}
		\ds\partial_t^j\partial_x^\alpha P_{h,c}(t,x)
		=h^{-j-|\alpha|}(\partial_t^j\partial_x^\alpha P)(\tau,\xi),\\
		\ns\ds\partial_t^j\partial_x^\alpha Z_{h,c}(t,x)
		=h^{-2-j-|\alpha|}(\partial_t^j\partial_x^\alpha Z)(\tau,\xi).
	\end{cases}
\end{equation}
Each fixed scaled solution has bounded derivatives and period
$2\pi h$. In the assembly we take $h=1/N_{\rm scale}$ with
$N_{\rm scale}\in\N$, for which
\[
P_{h,c}(t+2\pi,x)
=P(t/h+2\pi N_{\rm scale},(x-c)/h)=P_{h,c}(t,x),
\]
and similarly for $Z_{h,c}$, so the common period remains $2\pi$.

\subsection{Construction of observation geometry}
\label{sec:geometry}

Fix the radii $r,r_0$ from Proposition~\ref{mod:nd-module}. The
observation balls must satisfy two requirements: the larger balls of
radius $r_0$ must be separated, so that the local solutions can be
joined, and the smaller balls of radius $r$ must intercept all
sufficiently long straight segments. Both are supplied by the packing
theorem of Henk--Zong.

\subsubsection{Separated balls intercepting straight segments}

We use the following theorem twice below: once in dimension $n$, for the
observation balls, and once in dimension $n-1$, for the boundary
depressions. The two center sets are chosen independently.

\begin{theorem}[Henk--Zong~\cite{HZ}]
	\label{mod:HZ}
	For every integer $d\geq1$ and every $\varepsilon>0$, there are a
	full-rank lattice $\Lambda\subset\R^d$, a finite set $F\subset\R^d$,
	and a number $H_0<\infty$ such that the periodic center set
	$\mathcal C=F+\Lambda$ satisfies
	\begin{enumerate}
		\item $|c-c'|\geq2$ for all distinct $c,c'\in \mathcal C$;
		\item every Euclidean line segment of length $H_0$ meets an open
		ball $B_\varepsilon(c)$ with $c\in \mathcal C$.
	\end{enumerate}
\end{theorem}

For $d\geq2$ this is Theorem~1 of~\cite{HZ}; its proof proceeds by
dimension induction. The one-dimensional assertion also follows directly
by taking $\mathcal C=2\Z$ and $H_0=3$. In the source, the unit-ball
packing gives the separation, and the finite horizon bounds segments
avoiding the interiors of the small balls. Taking $H_0$ strictly larger
than that bound gives the open-ball assertion in Theorem~\ref{mod:HZ}.
Here $\varepsilon$ may be arbitrarily small.

\begin{corollary}[Separated centers for the observation balls]
	\label{mod:array}
	Fix the radii $r,r_0$ from Proposition~\ref{mod:nd-module}. Then there
	exist a periodic center set $\mathcal C_{r_0}\subset\R^n$ and a length
	$H<\infty$ with the following two properties. First, the centers are
	well separated:
	\begin{equation}\label{mod:separation}
		|c-c'|\geq4r_0\qquad(c,c'\in \mathcal C_{r_0},\ c\ne c').
	\end{equation}
	Second, every Euclidean line segment of length $H$ enters some open ball
	$B_r(c)$ with $c\in \mathcal C_{r_0}$.
	
	Both properties are stable under the scaling $c\mapsto hc$ with $h>0$.
	Consequently, for every $h>0$ the closed support balls
	$\overline{B_{hr_0}(c)}$, $c\in h\mathcal C_{r_0}$, are pairwise
	separated by at least $2hr_0$, and every segment of length $hH$ enters
	some $B_{hr}(c)$ with $c\in h\mathcal C_{r_0}$.
\end{corollary}

\begin{proof}[Proof of Corollary~\ref{mod:array}]
	We apply Theorem~\ref{mod:HZ} in dimension $d=n$ with
	$\varepsilon=r/(2r_0)$, and then rescale the resulting center set by
	$2r_0$. There are two things to verify: the unscaled separation and
	interception properties, and their behavior under the scaling
	$c\mapsto hc$.
	
	\medskip\noindent\textbf{Step 1. Unscaled separation and interception.}
	Let $\mathcal C$ and $H_0$ be the center set and length supplied by
	Theorem~\ref{mod:HZ}, and set
	\[
	\mathcal C_{r_0}=2r_0\mathcal C,\qquad H=2r_0H_0 .
	\]
	Since distinct centers of $\mathcal C$ are at distance at least $2$,
	distinct centers of $\mathcal C_{r_0}$ satisfy
	\[
	|2r_0c_0-2r_0c'_0|=2r_0|c_0-c'_0|\ge4r_0 .
	\]
	This is \eqref{mod:separation}.
	
	For the interception property, fix a unit vector $v$ and a segment
	$x+sv$, $0\le s\le H$. Dividing by $2r_0$ gives a segment of length
	$H_0$. By Theorem~\ref{mod:HZ} there are $s\in[0,H]$ and
	$c_0\in \mathcal C$ with
	\begin{equation}\label{mod:interception-unscaled}
		\left|\frac{x+sv}{2r_0}-c_0\right|<\varepsilon,
		\qquad |x+sv-2r_0c_0|<2r_0\varepsilon=r .
	\end{equation}
	Thus the segment enters the open ball $B_r(2r_0c_0)$ centered at a
	point of $\mathcal C_{r_0}$, as required.
	
	\medskip\noindent\textbf{Step 2. Scaling the centers by $h$.}
	Fix $h>0$. For distinct $c,c'\in \mathcal C_{r_0}$, the corresponding
	support balls satisfy
	\[
	\dist\bigl(\overline{B_{hr_0}(hc)},\overline{B_{hr_0}(hc')}\bigr)
	\ge h|c-c'|-2hr_0\ge2hr_0 ,
	\]
	where we used \eqref{mod:separation}. Hence the scaled support balls
	are pairwise separated.
	
	For the interception property, divide any segment of length $hH$ by
	$h$. This gives a segment of length $H$, so
	\eqref{mod:interception-unscaled} becomes
	\[
	|x+sv-hc|<hr
	\]
	for some $s\in[0,hH]$ and $c\in \mathcal C_{r_0}$. Thus every segment
	of length $hH$ enters some ball $B_{hr}(hc)$ with
	$c\in \mathcal C_{r_0}$.
	
	\medskip\noindent\textbf{Step 3. Relation between the two families of
		balls.}
	Since $r<r_0$, each observation ball $B_{hr}(hc)$ is contained in its
	support ball $B_{hr_0}(hc)$. By Step 1 the support balls at scale
	$h=1$ are disjoint, and by Step 2 the scaled support balls are
	separated by at least $2hr_0$. In particular, the observation balls
	are contained in pairwise disjoint support balls, which is exactly the
	configuration required for the local construction.
\end{proof}

\subsubsection{Making every boundary geodesic enter a concave region}

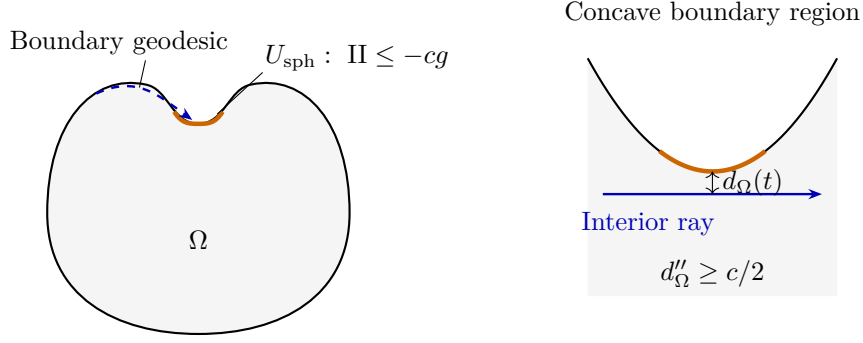
\begin{figure}[H]
	\centering
	\begin{tikzpicture}[x=1cm,y=1cm,font=\small,
		ray/.style={blue!70!black,line width=.9pt,-{Stealth[length=1.8mm]}},
		concave/.style={orange!80!black,line width=1.7pt}]
		\begin{scope}
			\path[fill=black!4,draw=black,line width=.8pt]
			(-2,0) .. controls (-2,1.2) and (-1.4,1.85) .. (-.7,1.7)
			.. controls (-.4,1.65) and (-.35,1.18) .. (0,1.18)
			.. controls (.35,1.18) and (.4,1.65) .. (.7,1.7)
			.. controls (1.4,1.85) and (2,1.2) .. (2,0)
			.. controls (2,-1.1) and (1.1,-1.6) .. (0,-1.6)
			.. controls (-1.1,-1.6) and (-2,-1.1) .. (-2,0);
			\draw[concave] (-.32,1.34)
			.. controls (-.22,1.18) and (-.12,1.18) .. (0,1.18)
			.. controls (.12,1.18) and (.22,1.18) .. (.32,1.34);
			\draw[ray,dashed] (-1.35,1.58)
			.. controls (-.82,1.83) and (-.53,1.56) .. (-.08,1.23);
			\node[align=center] at (-1.05,2.26) {Boundary geodesic};
			\draw[thin] (-.7,2.02)--(-.78,1.66);
			\node[anchor=west] at (.75,2.05) {$U_{\rm sph}:\ \mathrm{II}\le-cg$};
			\draw[thin] (.85,1.86)--(.25,1.31);
			\node at (0,-.35) {$\Omega$};
			\node at (0,-2.02) {(a) Every boundary geodesic reaches $U_{\rm sph}$.};
		\end{scope}
		\begin{scope}[shift={(6.8,0)}]
			\path[fill=black!4] (-1.65,2.05)
			plot[domain=-1.65:1.65,samples=50] (\x,{.55+.55*\x*\x})
			--(1.65,-1.1)--(-1.65,-1.1)--cycle;
			\draw[line width=.8pt]
			plot[domain=-1.65:1.65,samples=50] (\x,{.55+.55*\x*\x});
			\draw[concave]
			plot[domain=-.7:.7,samples=30] (\x,{.55+.55*\x*\x});
			\draw[ray] (-1.45,.25)--(1.45,.25);
			\draw[<->,thin] (0,.25)--node[right] {$d_\Omega(t)$}(0,.55);
			\node[anchor=south] at (0,2.35) {Concave boundary region};
			\node[blue!70!black,anchor=north] at (-.85,.15) {Interior ray};
			\node at (0,-.8) {$d_\Omega''\ge c/2$};
			\node at (0,-2.02) {(b) Strictly positive normal acceleration.};
		\end{scope}
	\end{tikzpicture}
	\caption{The mechanism for leaving a thin boundary collar.
		Boundary geodesics enter a concave region $U_{\rm sph}$ and remain there for a
		uniform positive time. A ray confined near the boundary has a nearby
		boundary geodesic; in $U_{\rm sph}$, the distance $d_\Omega(t)$ has strictly positive
		second derivative. Both panels are schematic; panel (b) shows a normal section.}
	\label{fig:high-dimensional-collar}
\end{figure}

Fix $n\ge3$ and put $d=n-1$. We construct an analytic boundary so that
every boundary geodesic spends a uniform positive time in a concave
region; the next subsubsection uses this property to prove that rays
escape a thin boundary collar. See Figure~\ref{fig:high-dimensional-collar}.

All ambient propagation uses the Euclidean metric of $\R^n$. A
generalized ray is the unit-speed spatial projection of a compressed
generalized bicharacteristic for the Dirichlet wave operator, with the
convention of \cite[Section~1C1, Definitions~1.3 and~1.4]{LRLT} and the
propagation framework of \cite{MS1,MS2}. The boundary compression
identifies incoming and outgoing normal covariables while retaining the
tangential covector.
\begin{proposition}[Uniform escape from a boundary collar]\label{geo:escape}
	Fix $n\ge3$. There exist a bounded star-shaped domain
	$\Omega\subset\R^n$ with real-analytic boundary and constants
	$\delta>0$ and $T_b<\infty$ with the following two properties. First,
	every generalized ray $x(t)$ in $\overline\Omega$ leaves a thin
	boundary collar within the uniform time $T_b$; more precisely,
	\begin{equation}\label{geo:escape-conclusion}
		\dist(x(t),\partial\Omega)>4\delta
		\qquad\text{for some }0<t<T_b.
	\end{equation}
	Second, the boundary has no infinite-order contact with a free
	Euclidean ray. The second property excludes infinite-order glancing,
	so that the free Dirichlet wave is observable on this geometry by the
	GCC theorem.
\end{proposition}
The convention for generalized rays covers transverse reflection,
glancing, gliding, and accumulation of reflections. The proof has three
geometric ingredients, which we develop in turn: concave boundary
regions reached by every boundary geodesic, uniform escape from a thin
collar, and the finite-order contact property of the analytic boundary.

\medskip\noindent\textbf{Radial graphs and their curvature.}
Let $\sigma$ be the round metric on $\Sphere^d$. We seek a positive
analytic radial function $R$ and an open region
$U_{\rm sph}\subset \Sigma_R(\Sphere^d)$, where
$\Sigma_R(\theta)=R(\theta)\theta$. For a radial surface write $g$ for
its induced metric and $n_{\rm in}$ for its inward unit normal. Define
the shape operator $\mathsf{W}$ and second fundamental form by
\[
D_vn_{\rm in}=-\mathsf{W}v,\qquad \mathrm{II}(v,w)=g(\mathsf{W}v,w).
\]
With this convention the unit sphere has $\mathrm{II}=\sigma$. Our
target is to find $c,T_g,\tau_0,m_0>0$ such that
\[
\mathrm{II}\le-cg\quad\hbox{on }U_{\rm sph},
\]
and every unit-speed boundary geodesic $q$ admits $t_0\in[0,T_g]$
such that
\[
\dist_g\bigl(q(t),\Sigma_R(\Sphere^d)\setminus U_{\rm sph}\bigr)\ge m_0
\quad(t_0\le t\le t_0+\tau_0).
\]
Here $T_g$ bounds the waiting time, $\tau_0$ is the length of the stay,
and $m_0$ is its distance margin from the edge of $U_{\rm sph}$.
Lemma~\ref{geo:concave-capture} establishes this property.

\medskip\noindent\textbf{Choosing chart regions with an interior margin.}
Choose disjoint chart regions which every geodesic reaches with enough
interior margin to extend the depressions smoothly by zero.

Set
\[
N_{\rm cap}=2(d+1),\qquad b=\frac1{2\sqrt{d+1}},
\]
and enumerate the vectors $\pm e_k$ of $\R^{d+1}$ as
$a_1,\ldots,a_{N_{\rm cap}}$. The open spherical caps
$
C_i=\{\theta\in\Sphere^d:a_i\cdot\theta>b\}
$
cover the sphere, since a unit vector has a coordinate of absolute value
at least $1/\sqrt{d+1}$. Each closed cap is contained in an open
hemisphere. In increasing index order, define pairwise disjoint open sets
\begin{equation}\label{geo:priority-partition}
	P_i=\{\theta:a_i\cdot\theta>b,\quad
	a_j\cdot\theta<b\text{ for }j<i\},
	\qquad
	\mathcal S=\bigcup_{i=1}^{N_{\rm cap}}\{\theta:a_i\cdot\theta=b\}.
\end{equation}
Empty sets $P_i$ are harmless. Every point outside $\mathcal S$ belongs
to exactly one $P_i$, namely that of the first cap which contains it.
Moreover $\partial P_i\subset \mathcal S$. Neither connectedness nor
smoothness of the boundary of $P_i$ is needed below.

\begin{lemma}[Geodesics away from the chart boundaries]\label{geo:uniform-depth}
	\begin{enumerate}
		\item[(1)] There is $\eta>0$ such that every round unit-speed geodesic
		on $\Sphere^d$ has, at some time in $[0,2\pi]$, distance at least
		$\eta$ from $\mathcal S$.
		\item[(2)] For every smooth metric $g$ sufficiently close to $\sigma$
		in $C^1$, every $g$ unit-speed geodesic has, in the same time
		interval, distance at least $\eta/2$ from $\mathcal S$, measured with
		$\sigma$.
	\end{enumerate}
\end{lemma}

\begin{proof}[Proof of Lemma~\ref{geo:uniform-depth}]
	The proof follows the two assertions of the lemma. Part (1) establishes
	the uniform margin for the round metric; Part (2) transfers it to
	$C^1$-small perturbations.
	
	\medskip\noindent\textbf{Part (1): A uniform margin for the round metric.}
	A round unit-speed geodesic is a great circle
	$\theta(t)=x\cos t+y\sin t$, with $x,y$ orthonormal. For a fixed $i$,
	its intersections with $\{a_i\cdot\theta=b\}$ solve
	\[
	(a_i\cdot x)\cos t+(a_i\cdot y)\sin t=b .
	\]
	If the amplitude $(a_i\cdot x)^2+(a_i\cdot y)^2$ is zero, there are no
	solutions, since $b>0$; otherwise this is a nonzero level of a single
	sinusoid, and there are at most two solutions per period. Hence no
	complete great circle lies in
	$\mathcal S=\bigcup_{i=1}^{N_{\rm cap}}\{a_i\cdot\theta=b\}$.
	
	Since the round unit tangent bundle is compact, we may define
	\[
	F(\zeta)=\max_{0\le t\le2\pi}\dist_\sigma(\theta_\zeta(t),\mathcal S).
	\]
	This function is continuous. It is strictly positive, because otherwise
	some geodesic would remain in $\mathcal S$ for all $t\in[0,2\pi]$, in
	particular its complete great circle would lie in $\mathcal S$,
	contrary to the previous paragraph. A continuous strictly positive
	function on a compact set has a positive minimum, so
	\[
	\eta\triangleq\min F>0 .
	\]
	This proves assertion (1).
	
	\medskip\noindent\textbf{Part (2): Stability under $C^1$-small perturbations.}
	Let $\varepsilon_g=\norm{g-\sigma}_{C^1}$. We compare a $g$-geodesic
	with a round geodesic over the fixed time interval $[0,2\pi]$. The
	comparison has three ingredients: the two metrics have close geodesic
	fields, the two trajectories start close, and the integral equations
	then keep them close for the whole interval.
	
	\smallskip\noindent\textit{Step 2.1: The geodesic fields are close.}
	In any of the fixed charts,
	\[
	\begin{cases}\ds
		g^{-1}-\sigma^{-1} =g^{-1}(\sigma-g)\sigma^{-1},\\
		\ns\ds\Gamma^k_{g,ab}\!-\Gamma^k_{\sigma,ab}\!
		=\!\frac12g^{kj}\bigl(\partial_a(g\!-\!\sigma)_{bj}
		+\partial_b(g\!-\!\sigma)_{aj}\!-\!\partial_j(g-\sigma)_{ab}\bigr)
		+\frac12(g^{kj}\!-\!\sigma^{kj})
		(\partial_a\sigma_{bj}\!+\partial_b\sigma_{aj}\!-\partial_j\sigma_{ab}),\\
		\ns\ds\norm{\Gamma_g-\Gamma_\sigma}_{C^0} \le C\varepsilon_g.
	\end{cases}
	\]
	The geodesic field is $X_g(\theta,v)=(v,-\Gamma_g[v,v])$ in
	coordinates. Since the Christoffel symbols differ by $O(\varepsilon_g)$
	and the velocities range over a fixed compact annulus, we obtain
	\[
	\norm{X_g-X_\sigma}_{C^0}\le C\varepsilon_g
	\]
	on that annulus. Thus the two geodesic fields differ by
	$O(\varepsilon_g)$ in the sup norm.
	
	\smallskip\noindent\textit{Step 2.2: The initial data are close.}
	Let $Z_g(t)=(\theta_g(t),v_g(t))$ be the $g$-geodesic trajectory with
	$|v_g(0)|_g=1$, and let $Z_\sigma$ start at the same point with initial
	velocity $v_g(0)/|v_g(0)|_\sigma$. Since
	$|v_g(0)|_g=1$ and $g$ is $C^1$-close to $\sigma$, the normalization
	factors differ by $O(\varepsilon_g)$, hence
	\[
	|Z_g(0)-Z_\sigma(0)|\le C\varepsilon_g .
	\]
	
	\smallskip\noindent\textit{Step 2.3: The trajectories stay close.}
	The integral equations, with $L_\sigma$ a Lipschitz constant for the
	fixed field $X_\sigma$, give
	\[
	|Z_g(t)-Z_\sigma(t)|
	\le C(1+t)e^{L_\sigma t}\varepsilon_g,
	\qquad 0\le t\le2\pi .
	\]
	Here $Z_g$ and $Z_\sigma$ are viewed in $\R^{d+1}\times\R^{d+1}$ with
	its Euclidean norm; on the fixed annulus this norm controls the
	spherical base-point distance. Thus, after decreasing
	$\varepsilon_g$ if necessary, the two base points are at spherical
	distance less than $\eta/2$ for all $t\in[0,2\pi]$.
	
	\smallskip\noindent\textit{Step 2.4: Transfer of the margin.}
	At a time $t_*\in[0,2\pi]$ where the round geodesic has distance at
	least $\eta$ from $\mathcal S$, the triangle inequality for
	$\dist_\sigma$ gives
	\[
	\dist_\sigma(\theta_g(t_*),\mathcal S)
	\ge\dist_\sigma(\theta_\sigma(t_*),\mathcal S)
	-\dist_\sigma(\theta_g(t_*),\theta_\sigma(t_*))
	\ge\eta-\eta/2=\eta/2 .
	\]
	Since $t_*$ lies in the same interval $[0,2\pi]$, the $g$-geodesic has
	distance at least $\eta/2$ from $\mathcal S$ at some time in
	$[0,2\pi]$. This proves assertion (2).
\end{proof}

The interception theorem concerns straight segments. To apply it on the
sphere we use coordinates that send round geodesics to straight lines:
for each cap take the central projection, also called the gnomonic chart,
\begin{equation}\label{geo:gnomonic}
	u=\phi_i(\theta)=\frac{\theta}{a_i\cdot\theta}-a_i
	\in a_i^\perp\simeq\R^d,
	\qquad D_i=\phi_i(P_i).
\end{equation}
The inverse and the line property follow directly from
\[
\phi_i^{-1}(u)=\frac{a_i+u}{\sqrt{1+|u|^2}},\qquad\phi_i\bigl(\Sphere^d\cap E\cap\{a_i\cdot\theta>0\}\bigr)
=\{u\in a_i^\perp:a_i+u\in E\},
\]
where $E$ is any two-dimensional linear plane meeting the hemisphere:
the set on the right is an affine line. Its parametrization need not have
constant speed; below we use Euclidean arclength along it. The chart is
smooth on $a_i\cdot\theta>b/2$. On the closed cap
$a_i\cdot\theta\ge b$ all coordinate comparison constants are uniform,
and $|u|\le\sqrt{1-b^2}/b$; in particular $D_i$ is bounded. The sets
\[
K_i=\{\theta\in\overline{P_i}:
\dist_\sigma(\theta,\mathcal S)\ge\eta/2\}
\]
are compact subsets of $P_i$. After omitting empty $K_i$, compactness
and finiteness give a single number $a>0$ such that
\begin{equation}\label{geo:chart-margin}
	\dist_{\R^d}\bigl(\phi_i(K_i),\R^d\setminus D_i\bigr)\ge a .
\end{equation}
The region $P_i$ has coordinate image $D_i$; $K_i$ is reached by
Lemma~\ref{geo:uniform-depth}, and the margin $a$ ensures that
intercepted balls are retained.

\medskip\noindent\textbf{Constructing the concave regions.}
We construct depressions whose bottoms have negative definite second
fundamental form and whose induced metric is a $C^1$-small perturbation
of the round metric.

For this second application of Theorem~\ref{mod:HZ}, take $d=n-1$ and
$\varepsilon=1/4$. It gives a periodic set $\mathcal C\subset\R^d$ and
$H_b<\infty$ such that distinct centers are at least $2$ apart and every
straight segment of length $H_b$ enters some open ball $B(c,1/4)$,
$c\in\mathcal C$. Here the balls locate boundary depressions; they are
unrelated to the observation balls of Corollary~\ref{mod:array}, where
$d=n$ and $\varepsilon=r/(2r_0)$. Fix $\mathcal C$ and $H_b$ before
choosing the width of the depressions.

Choose a smooth function $\chi:[0,\infty)\to[0,1]$ equal to one on
$[0,1/2]$ and to zero on $[3/4,\infty)$. Choose $A>0$ so large that in
all the fixed charts $2A I\ge4\sigma$, and then choose
$B_0>A(3/4)^2$. For a small number $h_b>0$, retain, in each $D_i$, the
centers $h_bc$, $c\in\mathcal C$, satisfying
\begin{equation}\label{geo:retained-dimples}
	\overline{B(h_bc,3h_b/4)}\subset D_i .
\end{equation}
There are finitely many, because $D_i$ is bounded and their mutual
distances are at least $2h_b$. In the $i$th chart define
\begin{equation}\label{geo:dimple-function}
	h_{i,c}(u)=h_b^2\chi\left(\frac{|u-h_bc|}{h_b}\right)
	\left(A\frac{|u-h_bc|^2}{h_b^2}-B_0\right),
	\qquad
	R_{h_b}=1+\sum_{i,c}h_{i,c}\circ\phi_i,
\end{equation}
with every summand extended by zero outside its chart support. The
scales are:
\begin{center}
	\small
	\begin{tabular}{@{}ll@{}}
		\hline
		$a$ & distance available before reaching a chart boundary,\\
		$h_b$ & width of a depression; its depth is of order $h_b^2$,\\
		$H_bh_b$ & straight-segment length needed to reach a ball of radius $h_b/4$,\\
		$h_b/2$ & radius of the bottom where negative curvature is required.\\
		\hline
	\end{tabular}
\end{center}
The depth is quadratic in the width: this makes the metric perturbation
small, while leaving the second derivatives large enough to change the
sign of the curvature. The cutoff is constant near zero, so the radial
expression is smooth at its center, and its zero extension is smooth at
the support boundary. For each fixed $h_b$ the supports are compactly
contained in their respective $P_i$. In one chart the support radii are
$3h_b/4$, whereas the centers are separated by $2h_b$; supports from
different charts belong to disjoint $P_i$. Thus at most one summand is
nonzero at any point. For a multi-index $\alpha$ in one chart,
\[
\partial_u^\alpha h_{i,c}(u)
=h_b^{2-|\alpha|}
\left.\partial_v^\alpha
\bigl[\chi(|v|)(A|v|^2-B_0)\bigr]\right|_{v=(u-h_bc)/h_b}.
\]
The coordinate changes belong to finitely many fixed charts, so the
increasing number of depressions does not enter the estimates
\begin{equation}\label{geo:radial-estimates}
	R_{h_b}=1+O(h_b^2),\qquad
	dR_{h_b}=O(h_b),\qquad
	\operatorname{Hess}_\sigma R_{h_b}=O(1).
\end{equation}
These estimates are global on the sphere, with constants independent of
sufficiently small $h_b$.

For a positive radial graph $\Sigma_R(\theta)=R(\theta)\theta$, its
induced metric and second fundamental form, pulled back to $\Sphere^d$,
are
\begin{align}
	g_R&=R^2\sigma+dR\otimes dR,
	\label{geo:radial-metric}\\
	\mathrm{II}_R
	&=\frac{R^2\sigma+2dR\otimes dR
		-R\operatorname{Hess}_\sigma R}
	{\sqrt{R^2+|\nabla_\sigma R|^2}}.
	\label{geo:radial-II}
\end{align}
The component calculation is given in Appendix~\ref{app:geo-radial-details}.

Write $g_{h_b}$ and $\mathrm{II}_{h_b}$ for the forms of $R_{h_b}$.
Formula \eqref{geo:radial-estimates} implies
\begin{equation}\label{geo:metric-close}
	\norm{g_{h_b}-\sigma}_{C^1}\le C h_b .
\end{equation}
Indeed, writing $R=R_{h_b}$, $R_a=\nabla_aR$ and
$R_{ab}=\nabla_a\nabla_bR$ for the covariant derivative of the round
metric,
\[
\begin{aligned}
	(g_{h_b}-\sigma)_{ab}
	&=(R^2-1)\sigma_{ab}+R_aR_b,\\
	\nabla_k(g_{h_b}-\sigma)_{ab}
	&=2RR_k\sigma_{ab}+R_{ka}R_b+R_aR_{kb},\\
	\norm{g_{h_b}-\sigma}_{C^0}
	&\le C\bigl(\norm{R-1}_{C^0}+\norm{dR}_{C^0}^2\bigr)
	\le Ch_b^2,\\
	\norm{\nabla(g_{h_b}-\sigma)}_{C^0}
	&\le C\bigl(\norm{R}_{C^0}\norm{dR}_{C^0}
	+\norm{\nabla^2R}_{C^0}\norm{dR}_{C^0}\bigr)
	\le Ch_b .
\end{aligned}
\]

On every closed bottom ball $|u-h_bc|\le h_b/2$ one has
\[
R_{h_b}=1+A|u-h_bc|^2-B_0h_b^2,
\qquad \partial_u^2R_{h_b}=2A I .
\]
Let $\Gamma_{\sigma,ab}^{k}$ be the Christoffel symbols of the round
metric in this chart. Its spherical Hessian satisfies, as an inequality
of quadratic forms,
\[
(\operatorname{Hess}_\sigma R_{h_b})_{ab}
=2A\delta_{ab}-2A\Gamma_{\sigma,ab}^{k}(u)(u-h_bc)_k
\ge4\sigma_{ab}-Ch_b\sigma_{ab}\ge3\sigma_{ab}
\]
for small $h_b$. Consequently
\[
R_{h_b}^2\sigma+2dR_{h_b}\otimes dR_{h_b}
-R_{h_b}\operatorname{Hess}_\sigma R_{h_b}
\le\bigl(-2+Ch_b^2\bigr)\sigma .
\]
Formula~\eqref{geo:radial-II}, with $g_{h_b}=\sigma+O(h_b^2)$, now gives
a constant $c_1>0$, independent of small $h_b$, with
\begin{equation}\label{geo:negative-bottom}
	\mathrm{II}_{h_b}\le-c_1g_{h_b}
	\qquad\text{on every }|u-h_bc|\le h_b/2 .
\end{equation}
Inequality \eqref{geo:negative-bottom} holds in all $d$ tangential
directions. Also $R_{h_b}>0$, so its radial graph is a smooth embedded
boundary of a bounded star-shaped domain.

\Needspace{6\baselineskip}
\medskip\noindent\textbf{Reaching a concave region and remaining inside it.}
The next lemma packages the preceding construction.

\begin{lemma}[Concave regions reached by every boundary geodesic]\label{geo:concave-capture}
	Fix $\kappa>0$. Then there is a positive real-analytic function
	$R:\Sphere^d\to(0,\infty)$ with the following two properties. First,
	the radial graph of $R$ has induced metric $g$ satisfying
	$\norm{g-\sigma}_{C^1}<\kappa$; that is, the graph is a
	$C^1$-small perturbation of the round sphere. Second, there are an open
	subset $U_{\rm sph}$ of the graph and constants $c,T_g,\tau_0,m_0>0$
	such that
	\begin{enumerate}
		\item[(1)] $\mathrm{II}\le-cg$ on $U_{\rm sph}$, so $U_{\rm sph}$
		is strictly concave;
		\item[(2)] every unit-speed boundary geodesic has an interval
		$[t_0,t_0+\tau_0]$, with $0\le t_0\le T_g$, entirely in
		$U_{\rm sph}$ and at distance at least $m_0$ from its complement.
	\end{enumerate}
	Thus $T_g$ bounds the waiting time before entering $U_{\rm sph}$,
	$\tau_0$ is the length of the stay, and $m_0$ is the distance margin
	from the edge of $U_{\rm sph}$.
\end{lemma}

\begin{proof}[Proof of Lemma~\ref{geo:concave-capture}]
	The construction proceeds in three stages. First, for the smooth radial
	graph $R_{h_b}$ we prove that every boundary geodesic hits a depression
	bottom within a uniform time. Second, we prove that it then stays in
	the bottom for a uniform positive time, with a uniform distance margin.
	Third, we replace $R_{h_b}$ by a real-analytic approximation which
	preserves both the strict concavity and these two time-margin
	properties. This yields the constants $c,T_g,\tau_0,m_0$ of the lemma.
	
	\medskip\noindent\textbf{Step 1. Every geodesic reaches a bottom ball.}
	\proofsteplabel{step:concave-hit}{1}
	\emph{Goal.} We show that every unit-speed geodesic of the smooth graph
	$R_{h_b}$ enters some chart ball $B(h_bc,3h_b/8)$ within a uniformly
	bounded time.
	
	Make $h_b$ small enough for Lemma~\ref{geo:uniform-depth} to apply to
	$g_{h_b}$. Every unit-speed $g_{h_b}$ geodesic then reaches a point
	$\theta_0\in K_i$ at a time $t_0\in[0,2\pi]$. Put
	$u_0=\phi_i(\theta_0)$. The coordinate speed of every such geodesic has
	uniform bounds
	\begin{equation}\label{geo:coordinate-speed}
		0<v_{\min}\le\left|\frac{du}{dt}\right|\le V<\infty
	\end{equation}
	as long as it remains in this chart; these constants hold in all the
	fixed cap charts and for all sufficiently small $h_b$.
	
	\smallskip\noindent\textit{Step 1.1: The geodesic stays close to a
		straight line.}
	Starting at $u_0$, parametrize the geodesic forward by Euclidean
	arclength $\ell$. If $H_b h_b<a/2$, the curve remains in $D_i$ for
	$0\le\ell\le H_b h_b$: before exiting it would need Euclidean
	displacement at least $a$ by \eqref{geo:chart-margin}, hence arclength
	at least $a$. Let $w=du/dt$ and $\tau=du/d\ell=w/|w|$. Then
	$\mathrm d\ell/\mathrm dt=|w|$ and
	\[
	\frac{\mathrm d\tau}{\mathrm d\ell}
	=\frac{(I-\tau\otimes\tau)w'}{|w|^2}
	=-\operatorname{proj}_{\tau^\perp}
	\Gamma_{g_{h_b}}(u)[\tau,\tau].
	\]
	The gnomonic property gives
	$\operatorname{proj}_{\tau^\perp}\Gamma_\sigma(u)[\tau,\tau]=0$.
	Subtracting this zero term and using \eqref{geo:metric-close} gives
	\begin{align}
		&\frac{d\tau}{d\ell}
		=-\operatorname{proj}_{\tau^\perp} (\Gamma_{g_{h_b}}-\Gamma_\sigma)[\tau,\tau],
		\qquad\left|\frac{d\tau}{d\ell}\right| \le C h_b, \label{geo:turning}\\
		&|u(\ell)-u_0-\ell\tau(0)|
		\le\frac C2h_b\ell^2,
		\qquad 0\le\ell\le H_b h_b.  \label{geo:line-error}
	\end{align}
	Estimate~\eqref{geo:line-error} follows from the explicit integral
	formula
	\[
	u(\ell)-u_0-\ell\tau(0)
	=\int_0^\ell\int_0^s\tau'(a)\,\mathrm da\,\mathrm ds,
	\qquad
	|u(\ell)-u_0-\ell\tau(0)|
	\le Ch_b\int_0^\ell\int_0^s1\,\mathrm da\,\mathrm ds
	=\frac C2h_b\ell^2 .
	\]
	Thus, on the relevant arclength scale, the geodesic is a small
	perturbation of the straight segment $u_0+\ell\tau(0)$.
	
	\smallskip\noindent\textit{Step 1.2: The straight segment is intercepted
		by a retained center.}
	The segment $u_0+\ell\tau(0)$, $0\le\ell\le H_bh_b$, has length
	$H_bh_b$. By the scaled interception property of $\mathcal C$, there
	are $c\in\mathcal C$ and $\ell_*\in[0,H_b h_b]$ with
	\begin{equation}\label{geo:virtual-hit}
		|u_0+\ell_*\tau(0)-h_bc|<h_b/4 .
	\end{equation}
	We must check that this center was retained in
	\eqref{geo:retained-dimples}; an arbitrary center of $\mathcal C$ need
	not have that property. For every $v\in\overline{B(h_bc,3h_b/4)}$,
	\[
	|v-u_0|\le|v-h_bc|+|h_bc-u_0|
	<\frac{3h_b}4+\left(H_b+\frac14\right)h_b=(H_b+1)h_b .
	\]
	By \eqref{geo:chart-margin}, this center therefore satisfies
	\eqref{geo:retained-dimples} if
	\begin{equation}\label{geo:small-scale-conditions}
		(H_b+1)h_b<a,\qquad \frac C2H_b^2h_b^3<h_b/8 .
	\end{equation}
	The first condition keeps the entire depression inside its chart; the
	second makes the geodesic-to-line error smaller than the spare radius
	$h_b/8$. Together with $H_bh_b<a/2$ from above, they can all be met by
	decreasing $h_b$, since $a,H_b,C$ are fixed.
	
	\smallskip\noindent\textit{Step 1.3: The actual geodesic hits the ball.}
	Combining \eqref{geo:line-error} at $\ell=\ell_*$ with
	\eqref{geo:virtual-hit} gives
	\[
	|u(\ell_*)-h_bc|
	\le|u(\ell_*)-u_0-\ell_*\tau(0)|
	+|u_0+\ell_*\tau(0)-h_bc|
	<\frac{h_b}8+\frac{h_b}4=\frac{3h_b}8 .
	\]
	By \eqref{geo:coordinate-speed}, the hitting time is at most
	\[
	T_g^{(0)}=2\pi+\frac{H_b h_b}{v_{\min}} .
	\]
	This proves that every geodesic reaches a bottom ball within the
	uniform time $T_g^{(0)}$.
	
	\medskip\noindent\textbf{Step 2. The geodesic remains in the bottom for a fixed time.}
	\proofsteplabel{step:concave-stay}{2}
	\emph{Goal.} We show that after the hit, the geodesic stays inside the
	open bottom ball for a uniform positive time, with a uniform spherical
	margin from the edge of the union of bottom balls.
	
	Let $t_*$ be the hitting time from Step~\ref{step:concave-hit} and set
	$\tau_0^{(0)}=h_b/(16V)$. For $t\in[t_*,t_*+\tau_0^{(0)}]$,
	\[
	|u(t)-h_bc|
	\le|u(t_*)-h_bc|+\int_{t_*}^t|u'(s)|\,\mathrm ds
	<\frac{3h_b}8+V\tau_0^{(0)}=\frac{7h_b}{16}<\frac{h_b}2 .
	\]
	A first-exit argument justifies using the speed estimate on this
	interval, since the closed bottom ball lies inside $D_i$. Thus the
	geodesic stays in the open bottom ball $B(h_bc,h_b/2)$ for the whole
	interval of length $\tau_0^{(0)}$.
	
	Let $U_{{\rm sph},0}\subset\Sphere^d$ be the union of the inverse
	images of all open bottom balls $B(h_bc,h_b/2)$. The fixed chart
	comparisons give $a_0>0$, independent of small $h_b$, such that
	\[
	\dist_\sigma(\theta(t),\Sphere^d\setminus U_{{\rm sph},0})
	\ge a_0\left(\frac{h_b}2-\frac{7h_b}{16}\right)
	=\frac{a_0h_b}{16}=:m_* .
	\]
	Thus the stay is not only long but also stays a uniform distance
	$m_*$ away from the edge of $U_{{\rm sph},0}$. This also covers a hit
	at an endpoint of the straight segment.
	
	\medskip\noindent\textbf{Step 3. Analytic approximation preserves these properties.}
	\emph{Goal.} We replace $R_{h_b}$ by a real-analytic $R$ close in
	$C^2$, so that the metric is $C^1$-close to $\sigma$, the bottom balls
	remain strictly concave, and the two time-margin properties of
	Steps~\ref{step:concave-hit}--\ref{step:concave-stay} persist.
	
	Fix $h_b$ satisfying the estimates of Steps~\ref{step:concave-hit}
	and~\ref{step:concave-stay}, and
	$\norm{g_{h_b}-\sigma}_{C^1}<\kappa/2$. For a positive radial function
	$\widetilde R$ near this fixed $R_{h_b}$, put
	$e_R=\norm{\widetilde R-R_{h_b}}_{C^2}$.
	
	\smallskip\noindent\textit{Step 3.1: Curvature and metric are stable.}
	The explicit formulas \eqref{geo:radial-metric}--\eqref{geo:radial-II}
	give
	\begin{equation}\label{geo:analytic-stability}
		\norm{g_{\widetilde R}-g_{h_b}}_{C^1}
		+\norm{\mathrm{II}_{\widetilde R}-\mathrm{II}_{h_b}}_{C^0}
		\le C_{h_b}e_R .
	\end{equation}
	In particular, for $e_R$ small,
	$\norm{g_{\widetilde R}-g_{h_b}}_{C^1}<\kappa/2$, hence
	$\norm{g_{\widetilde R}-\sigma}_{C^1}<\kappa$. Moreover, by
	\eqref{geo:negative-bottom} and the $C^0$ part of
	\eqref{geo:analytic-stability}, the bottom balls remain strictly
	concave:
	\[
	\mathrm{II}_{\widetilde R}\le-\frac{c_1}{2}g_{\widetilde R}
	\quad\text{on the closed bottom balls}.
	\]
	
	\smallskip\noindent\textit{Step 3.2: Geodesic trajectories are stable.}
	Let $Z_{\widetilde R}$ and $Z_{h_b}$ be the geodesic trajectories of
	the two metrics with the same initial point and normalized initial
	velocities. On $[0,T]$, where $T=T_g^{(0)}+\tau_0^{(0)}$, their
	integral equations give
	\begin{equation}\label{geo:analytic-flow-comparison}
		\begin{aligned}
			|Z_{\widetilde R}(t)-Z_{h_b}(t)|
			&\le C_{h_b}(1+t)e_R
			+L_{h_b}\int_0^t|Z_{\widetilde R}(s)-Z_{h_b}(s)|\,\mathrm ds\\
			&\le C_{h_b}(1+T)e^{L_{h_b}T}e_R=:C_{h_b,T}e_R .
		\end{aligned}
	\end{equation}
	Here $L_{h_b}$ belongs to the fixed smooth geodesic field of $g_{h_b}$.
	Estimate \eqref{geo:analytic-flow-comparison} includes the
	normalization of initial velocities. The approximation is chosen only
	after $h_b$ is fixed, so no bound uniform as $h_b\to0$ is needed for
	$L_{h_b}$ or $C_{h_b,T}$.
	
	\smallskip\noindent\textit{Step 3.3: The time-margin properties persist.}
	Choose $e_R$ so small that
	\begin{equation}\label{geo:analytic-neighborhood}
		\begin{gathered}
			\widetilde R>0,\qquad
			\norm{g_{\widetilde R}-g_{h_b}}_{C^1}<\kappa/2,
			\qquad C_{h_b,T}e_R<m_*/2 .
		\end{gathered}
	\end{equation}
	Then every perturbed trajectory remains inside $U_{{\rm sph},0}$ during
	the corresponding interval from Step~\ref{step:concave-stay}, with
	spherical margin at least $m_*/2$. Choose $c_{\rm met}>0$ so that
	$g_{\widetilde R}\ge c_{\rm met}^2\sigma$ for these approximations.
	Since $g_{\widetilde R}$ and $\sigma$ are uniformly equivalent on the
	compact sphere, the $g_{\widetilde R}$-distance to the complement of
	$U_{{\rm sph},0}$ satisfies, on that whole time interval,
	\[
	\dist_{g_{\widetilde R}}
	\bigl(\theta_{\widetilde R}(t),\Sphere^d\setminus U_{{\rm sph},0}\bigr)
	\ge c_{\rm met}\dist_\sigma
	\bigl(\theta_{\widetilde R}(t),\Sphere^d\setminus U_{{\rm sph},0}\bigr)
	\ge c_{\rm met}m_*/2 .
	\]
	
	\smallskip\noindent\textit{Step 3.4: Passage to a real-analytic surface.}
	Finite spherical harmonic sums approximate smooth functions in $C^2$.
	Choose such a sum $R$ within the $C^2$ neighborhood of $R_{h_b}$
	specified by \eqref{geo:analytic-neighborhood}. Then $R$ is real
	analytic and satisfies the required metric bound. With
	$U_{\rm sph}=\Sigma_R(U_{{\rm sph},0})$, one may take $c=c_1/2$,
	$T_g=T_g^{(0)}$, $\tau_0=\tau_0^{(0)}$, and $m_0=c_{\rm met}m_*/2$.
	Thus $T_g,\tau_0,m_0$ are independent of the chosen geodesic, once the
	surface and its scale have been fixed. Together with the strict
	concavity of Step 3.1, this gives all the assertions of the lemma.
\end{proof}

Fix the analytic radial function of Lemma~\ref{geo:concave-capture} and
set
\begin{equation}\label{geo:domain}
	\Omega=\{\rho\theta:\theta\in\Sphere^d,\ 0\le\rho<R(\theta)\}.
\end{equation}
Its boundary is the regular analytic radial graph
$\Sigma_R(\Sphere^d)$. Positivity and compactness give
\[
0<R_{\min}:=\min_{\Sphere^d}R\le R\le
R_{\max}:=\max_{\Sphere^d}R<\infty,\qquad
B_{R_{\min}}(0)\subset\Omega\subset B_{R_{\max}}(0).
\]
Every segment from $0$ to a point of $\Omega$ remains in $\Omega$, so
$\Omega$ is bounded, connected and star-shaped. The analytic
approximation in Lemma~\ref{geo:concave-capture} also gives the
following finite-contact property.

Let $G^\infty$ denote the glancing set corresponding to infinite-order
contact of free ambient rays with $\partial\Omega$.

\begin{lemma}[Finite-order boundary contact]\label{geo:finite-type}
	Every free Euclidean line tangent to $\partial\Omega$ has finite-order
	contact with it. In particular, the infinite-order glancing set
	$G^\infty$ is empty.
\end{lemma}

\begin{proof}[Proof of Lemma~\ref{geo:finite-type}]
	The function
	\begin{equation}\label{geo:analytic-defining-function}
		F(x)=|x|-R(x/|x|)
	\end{equation}
	is analytic on $\R^n\setminus\{0\}$, with
	$F(\rho\theta)=\rho-R(\theta)$ and $\partial_\rho F=1$. It is therefore
	a regular defining function for $\partial\Omega$. Suppose its
	restriction to a tangent line $x_0+tv$, $|v|=1$, vanishes to infinite
	order at $t=0$. Analyticity gives
	\[
	F(x_0+tv)=0,\qquad
	|x_0+tv|=R\left(\frac{x_0+tv}{|x_0+tv|}\right)\le R_{\max}
	\quad(t\in I),
	\]
	where $I$ is the component of $\{t:x_0+tv\ne0\}$ containing $0$. This
	component contains an unbounded sequence of times, whereas
	$|x_0+tv|\ge |t|-|x_0|\to\infty$ along that sequence. The contradiction
	proves finite-order contact.
	
	In the convention of \cite[Section~1C1]{LRLT}, infinite-order glancing
	means infinite-order contact of the free ambient ray with the boundary.
	For the Euclidean wave operator that ray is a straight line, so the
	analyticity argument for \eqref{geo:analytic-defining-function} gives
	$G^\infty=\varnothing$. This contact condition is distinct from the
	duration of a gliding interval along the boundary.
\end{proof}

\subsubsection{Escape from a thin boundary collar}

Fix the boundary and write $d_\Omega(t):=\dist(x(t),\partial\Omega)$;
$d'_{\Omega,-}$ and $d'_{\Omega,+}$ denote its left and right
derivatives. The contradiction has four steps:
\begin{enumerate}
	\item Assume $0\le d_\Omega\le\varepsilon$ on a fixed interval $[0,L]$.
	\item Bound the normal velocity by Lemma~\ref{geo:semiconvexity}, then
	compare with a unit-speed boundary geodesic $q_{\rm geo}$ by
	Lemma~\ref{geo:tangential-shadow}.
	\item Lemma~\ref{geo:concave-capture} puts the boundary projection in
	$U_{\rm sph}$, where Lemma~\ref{geo:concave-acceleration} gives the
	distributional second-derivative bound $d_\Omega''\ge c/2$.
	\item At the midpoint $m$ of a subinterval of length $\tau_0/2$,
	convexity gives $0\le d_\Omega(m)\le\varepsilon-c\tau_0^2/64<0$, having
	chosen $\varepsilon<c\tau_0^2/64$.
\end{enumerate}

\medskip\noindent\textbf{Normal and tangential equations.}
Write $s(x)=\dist(x,\partial\Omega)$ for the inner boundary distance.
In a sufficiently thin collar the nearest projection $q$ and $s$ are
smooth, and normal coordinates give
\begin{equation}\label{geo:normal-coordinates}
	x=q+s\,n_{\rm in}(q),\qquad
	ds_{\rm Eucl}^2=(ds)^2+g_s,\qquad
	g_s(v,w)=g((I-s\mathsf{W})v,(I-s\mathsf{W})w).
\end{equation}
Here $g$ is the metric of the fixed analytic boundary, and
$d=n-1$ continues to denote its dimension. Fix a collar width $\rho>0$
on which these coordinates are defined, and then fix $M\ge1$ with
\[
\sup_{0\le s(x)\le\rho}|\operatorname{Hess}s(x)|\le M .
\]
All later choices satisfy $0<\varepsilon<\rho$; in particular $M$ is
fixed before $\varepsilon$ and does not increase when the collar used by
a ray is narrowed.

In local boundary coordinates $y\in\R^d$, let $\eta$ be the tangential
covector and $\nu$ the normal covector. In physical time the unit-speed
spatial Hamiltonian is
\begin{equation}\label{geo:hamiltonian}
	\mathcal H(y,s,\eta,\nu)
	=\frac12\bigl(\nu^2+g_s^{ab}(y)\eta_a\eta_b\bigr)=\frac12 .
\end{equation}
On a free segment $d_\Omega'=s'=\nu$. At a transverse reflection,
\[
\nu_+=-\nu_->0,\qquad
(y_+,\eta_+)=(y_-,\eta_-),\qquad
d'_{\Omega,+}-d'_{\Omega,-}=2|\nu_-| .
\]
At every other time the defining differential equation of a generalized
bicharacteristic applies pointwise: the free Hamiltonian field is used
at diffractive points, and the boundary Hamiltonian field at the
remaining glancing points. These rules apply also at accumulation points
of reflections; see \cite[Definition~24.3.7 and the subsequent discussion
on p.~435]{HormanderIII}. Convexity of $s(x(t))+Mt^2/2$ controls both
one-sided normal velocities.

\begin{lemma}[Semiconvexity of the boundary distance]\label{geo:semiconvexity}
	If a generalized ray stays in the collar on an interval $I$, then
	\[
	t\longmapsto s(x(t))+\frac M2t^2
	\]
	is convex on $I$. In particular, if $0\le s(x(t))\le\varepsilon$ on
	$[0,L]$ and $a_\varepsilon=\sqrt{\varepsilon/M}<L/2$, then, with
	$d_\Omega(t)=s(x(t))$,
	\begin{equation}\label{geo:normal-speed}
		-\frac32\sqrt{M\varepsilon}
		\le d'_{\Omega,-}(t)\le d'_{\Omega,+}(t)
		\le\frac32\sqrt{M\varepsilon},
		\qquad a_\varepsilon\le t\le L-a_\varepsilon .
	\end{equation}
	The one-sided derivatives include both sides of a reflection.
\end{lemma}

\begin{proof}[Proof of Lemma~\ref{geo:semiconvexity}]
	On free segments $d_\Omega''=\operatorname{Hess}s(x',x')\ge-M$, and
	reflections give upward jumps of $d_\Omega'$. At glancing points
	$d_\Omega'=0$, so the slopes glue monotonically across the entire closed
	glancing set. Thus $d_\Omega+Mt^2/2$ is convex.
	Appendix~\ref{app:geo-semiconvex-details} verifies this also at
	accumulation points of reflections.
	
	For $a=a_\varepsilon$ and $a\le t\le L-a$, the forward and backward
	chord inequalities for $f=d_\Omega+Mt^2/2$ give
	\[
	\begin{aligned}
		d'_{\Omega,+}(t)&\le\frac{d_\Omega(t+a)-d_\Omega(t)}a+\frac{Ma}{2}
		\le\frac\varepsilon a+\frac{Ma}{2},\\
		d'_{\Omega,-}(t)&\ge\frac{d_\Omega(t)-d_\Omega(t-a)}a-\frac{Ma}{2}
		\ge-\frac\varepsilon a-\frac{Ma}{2}.
	\end{aligned}
	\]
	Together with the ordering of the one-sided derivatives, this yields
	\[
	-\frac{\varepsilon}{a_\varepsilon}-\frac{Ma_\varepsilon}{2}
	\le d'_{\Omega,-}(t)\le d'_{\Omega,+}(t)
	\le\frac{\varepsilon}{a_\varepsilon}+\frac{Ma_\varepsilon}{2}.
	\]
	The chosen value of $a_\varepsilon$ proves \eqref{geo:normal-speed}.
\end{proof}

The next lemma compares the tangential motion of a confined ray with a
unit-speed boundary geodesic, so that Lemma~\ref{geo:concave-capture}
applies.

\begin{lemma}[Comparison with boundary geodesics]\label{geo:tangential-shadow}
	For every generalized ray in the collar the tangential variables are
	$C^1$ and satisfy, in each boundary chart,
	\begin{equation}\label{geo:tangential-system}
		y'^a=g_s^{ab}(y)\eta_b,\qquad
		\eta'_a=-\frac12\partial_{y_a}g_s^{bc}(y)\eta_b\eta_c .
	\end{equation}
	If $0\le s(x(t))\le\varepsilon$ on $[0,L]$, their trajectory on any
	fixed time window in $[a_\varepsilon,L-a_\varepsilon]$ is
	$O(\varepsilon)$ close to the lift of a unit-speed boundary geodesic
	having the same initial point and normalized initial tangential covector.
	The constant depends only on the fixed boundary and the window length.
\end{lemma}

\begin{proof}[Proof of Lemma~\ref{geo:tangential-shadow}]
	In the interior \eqref{geo:tangential-system} is Hamilton's equation for
	\eqref{geo:hamiltonian}. The right-hand side does not contain $\nu$.
	At $s=0$ the free and boundary Hamiltonian fields have precisely the
	same tangential components. Hence the equations hold at each
	nonreflection time, including each glancing time, and $(y,\eta)$ is
	continuous across reflections.
	
	The continuous tangential equations on the two sides of reflections
	give $(y,\eta)\in C^1$, with no extra singular tangential motion on the
	closed glancing set. Appendix~\ref{app:geo-shadow-details} gives the
	difference-quotient argument at accumulation points.
	
	On the specified middle interval Lemma~\ref{geo:semiconvexity} and the
	energy identity give, on either side of a reflection when necessary,
	\[
	\nu^2\le\frac94M\varepsilon,\qquad
	1=\nu^2+g_s^{ab}\eta_a\eta_b
	=\nu^2+g_s(q',q').
	\]
	Since $0\le s\le\varepsilon$ and $g_s=g+O(\varepsilon)$,
	\[
	\begin{aligned}
		\bigl||\eta|_{g^{-1}}^2-1\bigr|
		&\le\nu^2+
		|(g^{-1}-g_s^{-1})^{ab}\eta_a\eta_b|
		\le C\varepsilon,\\
		\bigl||q'|_g^2-1\bigr|
		&\le\nu^2+|(g-g_s)(q',q')|
		\le C\varepsilon .
	\end{aligned}
	\]
	In particular,
	\begin{equation}\label{geo:tangential-energy}
		|\eta|_{g^{-1}}=1+O(\varepsilon).
	\end{equation}
	The covectors stay in a fixed compact annulus. Write $z=(y,\eta)$ and
	denote the right-hand side of \eqref{geo:tangential-system} by $V_s(z)$,
	with the convention $V_0:=V_s|_{s=0}$. For the fixed boundary,
	\begin{equation}\label{geo:tangential-vector-comparison}
		|V_s(z)-V_0(z)|\le C\varepsilon,
		\qquad |V_0(z)-V_0(\widetilde z)|\le C|z-\widetilde z|.
	\end{equation}
	At the initial time $t_*$ of the prescribed window let
	$\bar z(t_*)=(y(t_*),\eta(t_*)/|\eta(t_*)|_{g^{-1}})$, and let
	$\bar z'=V_0(\bar z)$. Denote by $\bar q$ its unit-speed boundary
	geodesic. By \eqref{geo:tangential-energy},
	$|z(t_*)-\bar z(t_*)|\le C\varepsilon$. The integral equations give,
	for $0\le t-t_*\le T$ within the middle interval,
	\begin{equation}\label{geo:shadow-gronwall}
		\begin{aligned}
			|z(t)-\bar z(t)|
			&\le C\varepsilon+C\varepsilon(t-t_*)
			+C\int_{t_*}^t|z(s)-\bar z(s)|\,\mathrm ds
			\le C(1+T)e^{CT}\varepsilon .
		\end{aligned}
	\end{equation}
	Coordinate changes have uniformly bounded derivatives on the compact
	overlaps of a finite atlas, so \eqref{geo:shadow-gronwall} extends over
	the fixed window (Appendix~\ref{app:geo-chart-details}). Hence
	$\dist_g(q(x(t)),\bar q(t))\le C_T\varepsilon$, where $C_T$ depends only
	on the fixed boundary and the window length.
\end{proof}

\medskip\noindent\textbf{Acceleration away from a concave boundary region.}

\begin{lemma}[Normal acceleration in concave regions]\label{geo:concave-acceleration}
	Suppose that on a time interval $J$ the projection $q(x(t))$ of a
	generalized ray lies in a region where $\mathrm{II}\le-cg$. If
	$0\le s(x(t))\le\varepsilon$ there and the squared tangential speed is
	$1+O(\varepsilon)$, then, for sufficiently small $\varepsilon$,
	\begin{equation}\label{geo:positive-acceleration}
		\frac{d^2}{dt^2}s(x(t))\ge\frac c2
		\qquad\text{on }J
	\end{equation}
	in the sense of distributions.
\end{lemma}

\begin{proof}[Proof of Lemma~\ref{geo:concave-acceleration}]
	On free segments,
	\begin{equation}\label{geo:free-normal-acceleration}
		d_\Omega''=-g\bigl(\mathsf{W}(I-d_\Omega\mathsf{W})q',q'\bigr),
		\qquad d_\Omega(t)=s(x(t)).
	\end{equation}
	This follows by differentiating $d_\Omega'=\langle n_{\rm in},x'\rangle$;
	see Appendix~\ref{app:geo-normal-details} for the calculation.
	
	Diagonalize the self-adjoint operator $\mathsf{W}$ in a $g$-orthonormal
	frame, with eigenvalues $\lambda_j\le-c$. Then
	\[
	d_\Omega''=\sum_{j=1}^d(-\lambda_j)(1-d_\Omega\lambda_j)(q'_j)^2
	\ge c|q'|_g^2\ge c/2,
	\]
	because $d_\Omega\ge0$ and the speed assumption, together with
	\eqref{geo:normal-coordinates}, gives $|q'|_g^2=1+O(\varepsilon)$. At a
	transverse reflection, $d'_{\Omega,+}-d'_{\Omega,-}=2|\nu_-|>0$. Such a
	reflection thus adds a positive point mass to $d_\Omega''$ and cannot
	spoil the lower bound.
	
	It remains to exclude an omitted glancing contribution.
	Formula~\eqref{geo:normal-coordinates} yields
	\[
	\partial_sg_s|_{s=0}=-2\mathrm{II},\qquad
	\partial_sg_s^{-1}|_{s=0}
	=2g^{-1}\mathrm{II}g^{-1}.
	\]
	At a glancing point $q'=g^{-1}\eta$ and $|q'|_g=1$. The normal Hamilton
	equation therefore gives
	\[
	\nu'=-\frac12\partial_sg_s^{ab}\eta_a\eta_b
	=-\mathrm{II}(q',q')\ge c.
	\]
	This is a strictly second-order diffractive point, where the defining
	rule is the free Hamiltonian field. If $t_0$ is such a contact time, the
	pointwise phase-space derivative gives, along nonreflection times,
	\begin{equation}\label{geo:diffractive-velocity}
		\nu(t)=A(t-t_0)+o(|t-t_0|),\qquad A>0.
	\end{equation}
	Thus $\nu$ is negative just before $t_0$ and positive just after it.
	There is no nearby transverse reflection on the left, since its outgoing
	free side has positive $\nu$; there is none on the right, since its
	incoming free side has negative $\nu$. Each side of a transverse
	reflection is a genuine free interval, so either reflection would
	contradict the signs in \eqref{geo:diffractive-velocity}. In this
	neighborhood all other contacts are also strictly diffractive. The
	trajectory therefore solves the smooth free Hamiltonian system
	throughout the neighborhood. Ordinary differential equation uniqueness
	shows that it is a free segment, with
	$d_\Omega(t_0)=d_\Omega'(t_0)=0$ and $d_\Omega''(t_0)>0$. In particular
	this contact is isolated.
	
	Every compact subinterval of $J$ now contains only finitely many
	contacts: an accumulation point would be glancing by continuity of the
	compressed trajectory and would contradict the local isolation implied
	by \eqref{geo:diffractive-velocity}. For any nonnegative
	$\varphi\in C_c^\infty(J)$ there are therefore only finitely many
	transverse reflection times $t_j$ in its support, and integration by
	parts on the free pieces gives
	\[
	\langle d_\Omega'',\varphi\rangle
	=\int_{J\setminus\{t_j\}}d_\Omega''(t)\varphi(t)\,\mathrm dt
	+\sum_j\bigl(d'_{\Omega,+}(t_j)-d'_{\Omega,-}(t_j)\bigr)\varphi(t_j)
	\ge\frac c2\int_J\varphi(t)\,\mathrm dt .
	\]
	Thus \eqref{geo:positive-acceleration} holds in distributions.
\end{proof}

We can now prove the escape proposition.

\begin{proof}[Proof of Proposition~\ref{geo:escape}]
	\noindent\textbf{Step 1. Fix the constants for the contradiction.}
	\proofsteplabel{step:escape-constants}{1}
	Take the analytic domain \eqref{geo:domain};
	Lemma~\ref{geo:finite-type} has established its finite-contact property.
	Fix $L>T_g+\tau_0+2$ and put $T=T_g+\tau_0$. We shall exclude
	$0\le d_\Omega\le\varepsilon$ on $[0,L]$ by finding a subinterval of
	length $\tau_0/2$ where $d_\Omega''\ge c/2$, and then applying its
	midpoint inequality.
	
	Let $C_T$ be the comparison constant in
	Lemma~\ref{geo:tangential-shadow} for a window of length $T$. Choose
	$0<\varepsilon<\rho$ below the threshold in
	Lemma~\ref{geo:concave-acceleration}, with
	\[
	a_\varepsilon<1,\qquad C_T\varepsilon<m_0/2,
	\qquad \varepsilon<c\tau_0^2/64 .
	\]
	The boundary and $L,T,M,C_T,c,\tau_0,m_0$ are fixed before this choice.
	Now suppose that a ray satisfies $0\le d_\Omega\le\varepsilon$ on
	$[0,L]$.
	
	\noindent\textbf{Step 2. Locate an interval inside the concave region.}
	\proofsteplabel{step:escape-capture}{2}
	Since $a_\varepsilon<1$ and $1+T<L-1<L-a_\varepsilon$, the window
	$[1,1+T]$ is contained in $[a_\varepsilon,L-a_\varepsilon]$.
	Lemma~\ref{geo:tangential-shadow} gives a unit-speed boundary geodesic
	$q_{\rm geo}$ with
	\[
	\dist_g(q(x(t)),q_{\rm geo}(t))\le C_T\varepsilon
	\qquad(1\le t\le1+T).
	\]
	Lemma~\ref{geo:concave-capture} supplies $t_0\in[0,T_g]$ such that for
	$1+t_0\le t\le1+t_0+\tau_0$,
	\[
	\begin{aligned}
		\dist_g(q(x(t)),\partial\Omega\setminus U_{\rm sph})
		&\ge\dist_g(q_{\rm geo}(t),\partial\Omega\setminus U_{\rm sph})
		-\dist_g(q(x(t)),q_{\rm geo}(t))
		\ge m_0-C_T\varepsilon>m_0/2 .
	\end{aligned}
	\]
	
	\noindent\textbf{Step 3. Apply the acceleration and midpoint estimates.}
	\proofsteplabel{step:escape-convexity}{3}
	Put
	\[
	\alpha=1+t_0+\tau_0/4,\qquad
	\beta=1+t_0+3\tau_0/4,\qquad
	J=[\alpha,\beta],\qquad \ell_0=\beta-\alpha=\tau_0/2 .
	\]
	On $J$ the projection lies in $U_{\rm sph}$. The normal-speed and energy
	estimates imply the speed hypothesis of
	Lemma~\ref{geo:concave-acceleration}, so $d_\Omega''\ge c/2$ there in
	distributions.
	
	The function $G(t)=d_\Omega(t)-ct^2/4$ is convex on $J$. At its midpoint
	$m=(\alpha+\beta)/2$,
	\[
	G(m)\le\frac{G(\alpha)+G(\beta)}2,\qquad
	\frac{\alpha^2+\beta^2}{2}-m^2=\frac{(\beta-\alpha)^2}{4}.
	\]
	Substitution gives
	\begin{equation}\label{geo:collar-convexity-contradiction}
		0\le d_\Omega(m)
		\le\frac{d_\Omega(\alpha)+d_\Omega(\beta)}2-\frac{c\ell_0^2}{16}
		\le\varepsilon-\frac{c\tau_0^2}{64}.
	\end{equation}
	Taking $\varepsilon<c\tau_0^2/64$ is a contradiction. All these choices
	concern the same fixed boundary and fixed time window, hence select a
	single positive $\varepsilon$ for every ray.
	
	\noindent\textbf{Step 4. Obtain the stated escape time.}
	\proofsteplabel{step:escape-time}{4}
	Set $\delta=\varepsilon/4$ and $T_b=L+1$. A ray starting with
	$s(x(0))>\varepsilon$ has the same property at a small positive time by
	continuity. Every other ray must attain $s>\varepsilon$ somewhere in
	$[0,L]$ by the contradiction in \eqref{geo:collar-convexity-contradiction};
	that time is positive and strictly smaller than $T_b$. This proves
	\eqref{geo:escape-conclusion}.
\end{proof}

\subsubsection{A finite family of interior balls satisfying GCC}
\label{sec:interior-observation-balls}

Every ray now reaches a point a fixed distance from the boundary. From
that point it travels along a straight segment long enough to meet one
of the scaled balls of Corollary~\ref{mod:array}. We make the scale
small enough that both this segment and the larger ball supporting the
corresponding local solution lie inside $\Omega$. These balls are the
observation balls; the boundary regions $U_{\rm sph}$ from
Lemma~\ref{geo:concave-capture} have forced escape from the collar
through Proposition~\ref{geo:escape}.

Fix the radii $r,r_0$ and the pair $(P,Z)$ from
Proposition~\ref{mod:nd-module}, and then fix $\mathcal C_{r_0},H$ from
Corollary~\ref{mod:array}. Independently, take the bounded connected
star-shaped domain $\Omega\subset\R^n$ with real-analytic boundary given
by Proposition~\ref{geo:escape}. Write
\[
d_\Omega(x)=\dist(x,\partial\Omega),\qquad x\in\overline\Omega .
\]
This proposition gives $\delta>0$ and $T_b<\infty$ such that every
compressed generalized unit-speed ray in $\overline\Omega$ has a time
$t_*\in(0,T_b)$ at which
\begin{equation}\label{asm:deep-point}
	d_\Omega(x(t_*))>4\delta .
\end{equation}

All of $r,r_0,H,\delta,T_b$ are now fixed. Choose an integer
\[
N_{\rm scale}>\frac{4\max\{r_0,H\}}\delta,\qquad N_{\rm scale}\ge1,
\]
so that, with $h=1/N_{\rm scale}$,
\begin{equation}\label{asm:small-scale}
	hr_0<\frac\delta4,\qquad hH<\frac\delta4 .
\end{equation}
The first inequality leaves room for the support ball; the second leaves
room for the straight segment of length $hH$. Both hold by increasing the
same integer $N_{\rm scale}$, without changing the domain or $\delta$.
The integer choice also preserves the common period $2\pi$, as proved
after \eqref{mod:scaling-derivatives}. The roles of the scales are
summarized below.
\begin{center}\small
	\begin{tabular}{@{}lll@{}}\toprule
		Object & Scale & Purpose\\\midrule
		$h=1/N_{\rm scale}$ & $h\max\{r_0,H\}<\delta/4$ & Controls both lengths\\
		$hH$ & $<\delta/4$ & Keeps the straight segment inside\\
		$hr_0$ & $<\delta/4$ & Keeps the support ball inside\\
		$\delta$ & Fixed before $h$ & Boundary-distance margin\\\bottomrule
\end{tabular}\end{center}

Retain the centers
\begin{equation}\label{asm:retained-centers}
	J=\{c\in h\mathcal C_{r_0}\cap\Omega:d_\Omega(c)\geq2\delta\},
	\qquad
	\omega=\bigcup_{c\in J}B_{hr}(c).
\end{equation}
For $c\in J$,
\[
B_{hr_0}(c)\subset\Omega,\qquad
|c-c'|\ge4hr_0\quad(c,c'\in J,\ c\ne c').
\]
The balls $B_{hr_0}(c)$ are disjoint. Comparing volumes for any finite
subset of $J$ therefore gives
\[
\#J\le\frac{|\Omega|}{|B_{hr_0}(0)|}<\infty .
\]
Moreover,
\[
\dist\bigl(\overline{B_{hr_0}(c)},\partial\Omega\bigr)
\geq2\delta-hr_0>\frac{7\delta}{4}\qquad(c\in J),
\]
and for distinct retained centers,
\[
\dist\bigl(\overline{B_{hr_0}(c)},
\overline{B_{hr_0}(c')}\bigr)
\ge |c-c'|-2hr_0\ge2hr_0 .
\]
Thus all support balls are compactly contained in $\Omega$ and
separated. The smaller observation balls are also separated and
$\omega\Subset\Omega$. Together with \eqref{mod:scaled-properties}, this
gives, for $c\ne c'$,
\begin{equation}\label{asm:disjoint-support-identities}
	(P_{h,c}-1)(P_{h,c'}-1)=0,
	\qquad Z_{h,c}(P_{h,c'}-1)=0 .
\end{equation}
These identities are used for the finite sums in Subsection~\ref{sec:assembly}.

\begin{lemma}[GCC for a finite family of interior balls]\label{asm:GCC}
	The set $J$ is nonempty, and $(\Omega,\omega)$ satisfies GCC for the
	standard generalized rays of the free Dirichlet wave equation. One may
	take the common observation time to be $T_{\mathrm{GCC}}=T_b+hH+1$.
\end{lemma}

\begin{proof}[Proof of Lemma~\ref{asm:GCC}]
	The proof has two tasks: keep the segment after $t_*$ inside $\Omega$,
	and show that the ball it meets has a retained center.
	Consider any generalized unit-speed ray, and choose $t_*$ as in
	\eqref{asm:deep-point}. For $0\le s\le hH$,
	\[
	\begin{aligned}
		|x(t_*+s)-x(t_*)|
		&\le\int_0^s|x'(t_*+\tau)|\,\mathrm d\tau=s\le hH,\\
		d_\Omega(x(t_*+s))&\ge d_\Omega(x(t_*))-|x(t_*+s)-x(t_*)|
		>4\delta-hH>0 .
	\end{aligned}
	\]
	The entire interval is in the interior, where the ray equation is
	$x''=0$. Hence
	\[
	x(t_*+s)=x(t_*)+sv,\qquad |v|=1,\qquad 0\le s\le hH .
	\]
	Corollary~\ref{mod:array} gives a time $s\in[0,hH]$ and a center
	$c\in h\mathcal C_{r_0}$ such that
	\[
	|x(t_*+s)-c|<hr .
	\]
	Consequently,
	\begin{equation}\label{asm:center-distance}
		|c-x(t_*)|<hH+hr<hH+hr_0<\frac\delta2 .
	\end{equation}
	The ball $B_{\delta/2}(x(t_*))$ lies inside $\Omega$ by
	\eqref{asm:deep-point}, so $c\in\Omega$. The distance function is
	one-Lipschitz, and hence
	\[
	d_\Omega(c)>4\delta-hH-hr>\frac{7\delta}{2}>2\delta .
	\]
	This center is retained in $J$, and the ray enters the corresponding
	open observation ball at time $t_{\rm hit}=t_*+s$, with
	\[
	0<t_{\rm hit}<T_b+hH<T_{\rm GCC}=T_b+hH+1 .
	\]
	The extra $1$ places the hitting time strictly below the common
	observation endpoint. Since the ball is open and the ray continuous, it
	is observed on a nonempty open time interval around the hitting time.
	Finally, choose any interior point and a unit initial direction to
	obtain a ray; \eqref{asm:center-distance} and $d_\Omega(c)>2\delta$ give
	a center in $J$, proving $J\ne\varnothing$.
	
	Proposition~\ref{geo:escape} includes every boundary contact; the
	interception argument is interior.
\end{proof}

\begin{remark}\label{geo:simplex-obstruction}
	A simplex partition does not have the same property in all dimensions.
	For example, in $E=\{x\in\R^6:\sum_i x_i=0\}\simeq\R^5$, the regular
	$5$-simplex with vertices $e_i-\frac16(1,\ldots,1)$ gives a radial
	partition of $\Sphere^4\subset E$. Its skeleton consists of points whose
	smallest coordinate is attained at least twice. The plane
	\[
	\{x\in E:x_1=x_2,\ x_3=x_4,\ x_5=x_6\}
	\]
	has dimension two, and its entire great circle lies in that skeleton. The
	nonzero cap levels in \eqref{geo:priority-partition} avoid this
	obstruction.
\end{remark}

\begin{remark}\label{geo:parameter-order}
	The choices have the order
	\[
	\begin{gathered}
		n;\quad (P_i,\eta,a);\quad (\mathcal C,H_b);\quad (\chi,A,B_0);
		\quad h_b;\quad R;\\
		(c,T_g,\tau_0,m_0,M);\quad L;\quad
		\varepsilon;\quad \delta=\varepsilon/4,\ T_b=L+1 .
	\end{gathered}
	\]
	In particular $H_b$ is fixed before $h_b$ is decreased, and the analytic
	surface $R$ is fixed before the collar width is chosen. No constants in
	the collar argument are required to be uniform as $h_b\to0$. The domain
	does not depend on the local waves or on the eventual observation time.
\end{remark}

\subsection{Assembly and the Dirichlet boundary condition}
\label{sec:assembly}

Use the domain $\Omega$, scale $h=1/N_{\rm scale}$, and center set $J$
from Subsubsection~\ref{sec:interior-observation-balls}. Join the local
and exterior Dirichlet solutions on an open cover of $\overline\Omega$,
then verify the equation, nontriviality, and failure of observability.
Write $\Delta=\Delta_n$ and $\Box=\partial_t^2-\Delta_n$.

\begin{proof}[Proof of Theorem~\ref{thm:counter-high-dim}]
	The proof has three steps. Step 1 glues the local solutions and the
	Dirichlet exterior solution into one smooth pair satisfying the wave
	equation on $\R\times\Omega$ and the two required zero conditions.
	Step 2 verifies that this pair is nontrivial, $2\pi$-periodic, and has
	bounded derivatives. Step 3 identifies it as the unique weak evolution
	of its initial state and uses its positivity in a boundary collar to
	disprove observability at every finite time.
	
	\medskip\noindent\textbf{Step 1. Join the local solutions to a Dirichlet exterior solution.}
	\proofsteplabel{step:highdim-glue}{1}
	Keep the domain in \eqref{geo:domain}, the scale in
	\eqref{asm:small-scale}, and the finite center set in
	\eqref{asm:retained-centers} fixed. By Lemma~\ref{asm:GCC}, the
	observation set $\omega=\bigcup_{c\in J}B_{hr}(c)$ satisfies GCC.
	Choose an open set $\mathcal U$ such that
	\[
	\bigcup_{c\in J}\overline{B_{hr_0}(c)}\Subset\mathcal U\Subset\Omega .
	\]
	Apply Lemma~\ref{ce:boundary} to $\mathcal U\Subset\Omega$. It gives a
	real smooth pair $(a_{\mathrm{out}},V_{\mathrm{out}})$ with
	\begin{equation}\label{asm:exterior-equation}
		\begin{gathered}
			-\Delta a_{\mathrm{out}}+V_{\mathrm{out}}a_{\mathrm{out}}=0,
			\qquad a_{\mathrm{out}}|_{\partial\Omega}=0,
			\qquad a_{\mathrm{out}}>0\quad\text{in }\Omega,\\
			(a_{\mathrm{out}},V_{\mathrm{out}})=(1,0)
			\quad\text{on a neighborhood of }\overline{\mathcal U}.
		\end{gathered}
	\end{equation}
	In a nonempty boundary collar this pair is
	$(\varphi_1,-\lambda_1)$, where $\varphi_1$ is a positive first
	Dirichlet eigenfunction and $\lambda_1>0$ its eigenvalue. Thus the
	exterior pair is a positive static solution which equals the constant
	$(1,0)$ on $\mathcal U$ and vanishes on $\partial\Omega$.
	
	For each $c\in J$, set $P_c=P_{h,c}$ and $Z_c=Z_{h,c}$ as in
	\eqref{mod:scaled}. By \eqref{mod:scaled-properties}, their equations
	and supports satisfy
	\begin{equation}\label{asm:local-equation-support}
		\Box(P_c-1)=-Z_cP_c,\qquad
		\supp_x(P_c-1)\cup\supp_xZ_c\Subset B_{hr_0}(c),
	\end{equation}
	the supports being taken over all times. To describe the regions on
	which we join the formulas, put
	\[
	D_c\triangleq c+hD\Subset B_{hr_0}(c),\qquad
	\mathcal U_c\triangleq B_{hr_0}(c),\qquad
	\mathcal U_0\triangleq\overline\Omega\setminus\bigcup_{c\in J}D_c .
	\]
	We now verify the four requirements for gluing: coverage, agreement on
	overlaps, smoothness, and the equation.
	
	\smallskip\noindent\textit{Coverage.}
	The sets $\mathcal U_0,\mathcal U_c$ are open relative to
	$\overline\Omega$. Since $D_c\subset\mathcal U_c$ and the larger balls
	are pairwise disjoint,
	\[
	\overline\Omega=\mathcal U_0\cup\bigcup_{c\in J}\mathcal U_c,
	\qquad\mathcal U_c\cap\mathcal U_{c'}=\varnothing\quad(c\ne c').
	\]
	Thus $\{\mathcal U_0\}\cup\{\mathcal U_c\}_{c\in J}$ is an open cover
	of $\overline\Omega$.
	
	\smallskip\noindent\textit{Agreement on overlaps.}
	The only possible overlaps are $\mathcal U_0\cap\mathcal U_c$. There,
	\eqref{mod:scaled-properties} and \eqref{asm:exterior-equation} give
	\[
	(P_c,Z_c)=(1,0)=(a_{\mathrm{out}},V_{\mathrm{out}})
	\quad\text{on }\R\times(\mathcal U_0\cap\mathcal U_c).
	\]
	
	\smallskip\noindent\textit{Smoothness.}
	Since the local pairs and the exterior pair agree on the overlaps, the
	gluing lemma for smooth functions on an open cover produces one smooth
	pair on $\R\times\overline\Omega$. Its finite-sum form is
	\begin{equation}\label{asm:global-formula}
		\begin{aligned}
			p(t,x)&=a_{\mathrm{out}}(x)+\sum_{c\in J}(P_c(t,x)-1),\\
			V(t,x)&=V_{\mathrm{out}}(x)+\sum_{c\in J}Z_c(t,x).
		\end{aligned}
	\end{equation}
	Both functions are real and smooth on $\R\times\overline\Omega$.
	
	\smallskip\noindent\textit{The equation.}
	The product $\Box p+Vp$ has three sources: the exterior pair, the
	diagonal local terms, and the mixed terms between different summands.
	The mixed terms vanish because, for every time,
	\[
	\begin{gathered}
		V_{\mathrm{out}}(P_c-1)=0,\qquad
		Z_c(a_{\mathrm{out}}-1)=0,\\
		Z_c(P_{c'}-1)=0\qquad(c,c'\in J,\ c\ne c').
	\end{gathered}
	\]
	Using \eqref{asm:exterior-equation} and
	\eqref{asm:local-equation-support}, expand \eqref{asm:global-formula}
	and split the double sum into diagonal and off-diagonal parts:
	\[
	\begin{aligned}
		\Box p+Vp
		&=-\Delta a_{\mathrm{out}}+V_{\mathrm{out}}a_{\mathrm{out}}
		-\sum_{c\in J}Z_cP_c
		+V_{\mathrm{out}}\sum_{c\in J}(P_c-1) +\sum_{c\in J}Z_ca_{\mathrm{out}}
		+\sum_{c,c'\in J}Z_c(P_{c'}-1)\\
		&=\sum_{c\in J}Z_c[-P_c+(P_c-1)+a_{\mathrm{out}}]
		+V_{\mathrm{out}}\sum_{c\in J}(P_c-1)
		+\sum_{\substack{c,c'\in J\\c\ne c'}}Z_c(P_{c'}-1)\\
		&=\sum_{c\in J}Z_c(a_{\mathrm{out}}-1)
		+V_{\mathrm{out}}\sum_{c\in J}(P_c-1)
		+\sum_{\substack{c,c'\in J\\c\ne c'}}Z_c(P_{c'}-1)=0 .
	\end{aligned}
	\]
	On $B_{hr_0}(c)$ the global pair is exactly $(P_c,Z_c)$; in
	particular, $p=0$ on $B_{hr}(c)$, hence on $\omega$. Near
	$\partial\Omega$ all local corrections vanish, so
	$(p,V)=(\varphi_1,-\lambda_1)$ and $p=0$ on $\partial\Omega$. We have
	therefore proved
	\begin{equation}\label{asm:global-equation}
		\begin{cases} 
			p_{tt}-\Delta p+Vp=0\quad\text{in }\R\times\Omega,\\
			p=0\quad\text{on }\R\times\partial\Omega,
			\qquad p=0\quad\text{on }\R\times\omega.
		\end{cases}
	\end{equation}
	
	\medskip\noindent\textbf{Step 2. Verify nontriviality, periodicity, and derivative bounds.}
	\proofsteplabel{step:highdim-bounds}{2}
	We check the three remaining properties of the pair $(p,V)$.
	
	\smallskip\noindent\textit{Nontriviality.}
	Choose a closed ball $\overline B\Subset\Omega$ contained in the
	boundary collar where
	\begin{equation}\label{asm:positive-collar}
		p(t,x)=\varphi_1(x)>0 .
	\end{equation}
	For every $t\in\R$,
	\[
	\|p(t)\|_{L^2(\Omega)}^2
	\geq\int_B\varphi_1(x)^2\,\mathrm dx
	\geq |B|\min_{\overline B}\varphi_1^2>0 .
	\]
	In particular, the initial state $(p(0),p_t(0))$ is nonzero.
	
	\smallskip\noindent\textit{Periodicity.}
	Each $P_c,Z_c$ has period $2\pi h=2\pi/N_{\rm scale}$, hence period
	$2\pi$; the exterior pair is independent of time. Therefore $p$ and
	$V$ are $2\pi$-periodic.
	
	\smallskip\noindent\textit{Derivative bounds.}
	Since $p,V$ are smooth and $2\pi$-periodic, every space-time derivative
	is bounded on $\R\times\overline\Omega$, and its supremum equals the
	maximum over one period. For every $k\geq0$ and spatial multi-index
	$\alpha$,
	\begin{equation}\label{asm:global-bounds}
		\begin{aligned}
			&\sup_{\R\times\overline\Omega}
			\left(|\partial_t^k\partial_x^\alpha p|
			+|\partial_t^k\partial_x^\alpha V|\right) =
			\max_{[0,2\pi]\times\overline\Omega}
			\left(|\partial_t^k\partial_x^\alpha p|
			+|\partial_t^k\partial_x^\alpha V|\right)<\infty .
		\end{aligned}
	\end{equation}
	
	All choices of $\Omega,\omega,p,V$ are now fixed. Let $T>0$ be
	arbitrary.
	
	\medskip\noindent\textbf{Step 3. Identify the weak evolution and disprove observability.}
	\proofsteplabel{step:highdim-observation}{3}
	We first show that the smooth solution is the unique weak evolution of
	its initial state in $X_\Omega$, and then use its positivity in the
	boundary collar to violate any observability inequality.
	
	\smallskip\noindent\textit{Step 3.1: The perturbation is bounded.}
	Let $A_\Omega=-\Delta$ with the Dirichlet condition, and let $S(t)$ be
	the free wave group on $X_\Omega=L^2(\Omega)\times H^{-1}(\Omega)$ with
	the spectral norm \eqref{eq:state-norm}. This group is unitary. For
	$K_V(t)(a,b)=(0,-V(t,\cdot)a)$,
	\begin{equation}\label{asm:bounded-perturbation}
		\begin{aligned}
			\|K_V(t)(a,b)\|_{X_\Omega}
			&=\|A_\Omega^{-1/2}(V(t)a)\|_{L^2(\Omega)}
			\leq\lambda_1^{-1/2}\|V(t)\|_{L^\infty(\Omega)}
			\|a\|_{L^2(\Omega)}\\
			&\leq\lambda_1^{-1/2}\|V(t)\|_{L^\infty(\Omega)}
			\|(a,b)\|_{X_\Omega}.
		\end{aligned}
	\end{equation}
	Put $M_T=\lambda_1^{-1/2}\max_{0\leq t\leq T}
	\|V(t)\|_{L^\infty(\Omega)}$.
	
	\smallskip\noindent\textit{Step 3.2: The Duhamel equation has a unique
		solution.}
	Given $U_0\in X_\Omega$, the iterates
	\begin{equation}\label{asm:picard-iteration}
		\begin{aligned}
			U^{(0)}(t)&=S(t)U_0,\\
			U^{(j+1)}(t)&=S(t)U_0+
			\int_0^t S(t-s)K_V(s)U^{(j)}(s)\,\mathrm ds
		\end{aligned}
	\end{equation}
	are uniformly summable by the successive-difference estimate in
	Appendix~\ref{app:geo-duhamel-details}. Its limit solves
	\begin{equation}\label{asm:Duhamel}
		U(t)=S(t)U_0+\int_0^t S(t-s)K_V(s)U(s)\,\mathrm ds .
	\end{equation}
	For two such solutions with initial states $U_0,\widetilde U_0$,
	\[
	\|U(t)-\widetilde U(t)\|_{X_\Omega}
	\leq\|U_0-\widetilde U_0\|_{X_\Omega}
	+M_T\int_0^t\|U(s)-\widetilde U(s)\|_{X_\Omega}\,\mathrm ds
	\leq e^{M_Tt}\|U_0-\widetilde U_0\|_{X_\Omega}.
	\]
	By Gronwall's inequality, the solution is unique for each initial
	state.
	
	\smallskip\noindent\textit{Step 3.3: The smooth solution is the weak
		evolution.}
	For the smooth Dirichlet solution in \eqref{asm:global-formula},
	$U(t)=(p(t),p_t(t))$ satisfies
	\[
	\frac{\mathrm d}{\mathrm dt}\bigl(S(-t)U(t)\bigr)
	=S(-t)K_V(t)U(t).
	\]
	Integration gives \eqref{asm:Duhamel} with $U_0=(p(0),p_t(0))$. By
	uniqueness, this smooth solution is the weak evolution of its initial
	state on $[0,T]$.
	
	\smallskip\noindent\textit{Step 3.4: Observability fails.}
	Equations \eqref{asm:global-equation} and \eqref{asm:positive-collar}
	give
	\begin{equation}\label{asm:observation-failure}
		\|(p(0),p_t(0))\|_{X_\Omega}^2>0,
		\qquad
		\int_0^T\!\int_\omega |p(t,x)|^2\,\mathrm dx\,\mathrm dt=0 .
	\end{equation}
	An observability inequality at time $T$ would bound the left side by a
	constant times the right side, which is impossible. Hence displacement
	observability fails for this fixed potential at every finite time.
\end{proof}

\begin{proof}[Proof of Corollary~\ref{cor:control}]
	The corollary has a positive part and a negative part. Step 1 derives
	exact controllability on the disk from the observability theorem by
	duality. Step 2 shows that, for the constructed potentials, exact
	controllability fails in every dimension $n\ge2$ and at every finite
	time.
	
	\medskip\noindent\textbf{Step 1. Exact controllability on the disk.}
	\proofsteplabel{step:control-disk}{1}
	Theorem~\ref{thm:disk}, applied to the adjoint equation, gives
	displacement observability in $L^2(\mathbb D)\times H^{-1}(\mathbb D)$
	for every $T>T_*$ and every bounded potential. The
	observability--controllability duality
	\cite{Lions,LionsBook,TucsnakWeiss} then gives exact controllability in
	$H_0^1(\mathbb D)\times L^2(\mathbb D)$ with
	$L^2((0,T)\times\omega)$ controls. This proves the positive assertion.
	
	\medskip\noindent\textbf{Step 2. Failure of exact controllability for the constructed potentials.}
	\proofsteplabel{step:control-obstruction}{2}
	The strategy is the following. The solution $p$ of the constructed
	invisible wave provides a nonzero continuous linear functional
	$\ell_T$ on the terminal energy space. The Green identity shows that
	$\ell_T$ annihilates every terminal state reachable from zero with an
	$L^2$ control. If some nonzero terminal state is not annihilated, it
	cannot be reached, and exact controllability fails. We now carry this
	out.
	
	\smallskip\noindent\textit{Step 2.1: The invisible wave is nonzero at
		every time.}
	For $n=2$, take $(\Omega,\omega,p,V)$ from Theorem~\ref{thm:counter};
	for $n\ge3$, take them from Theorem~\ref{thm:counter-high-dim},
	omitting the subscript $n$. In both constructions $p$ equals a positive
	time-independent Dirichlet eigenfunction in a nonempty boundary collar.
	Hence $\|p(t)\|_2>0$ for every $t\in\R$. Fix $T>0$; the domain,
	observation set, and potential are fixed before $T$ is chosen.
	
	\smallskip\noindent\textit{Step 2.2: The controlled solution is an
		energy solution.}
	For $f\in L^2((0,T)\times\omega)$, let $y$ solve
	\begin{equation}\label{asm:controlled}
		\begin{cases}
			y_{tt}-\Delta y+Vy=\mathbf1_\omega f
			&\text{in }(0,T)\times\Omega,\\
			y=0&\text{on }(0,T)\times\partial\Omega,\\
			y(0)=0,\quad y_t(0)=0&\text{in }\Omega.
		\end{cases}
	\end{equation}
	The energy-space perturbation is bounded because
	\begin{equation}\label{asm:energy-perturbation}
		\|V(t)a\|_{L^2(\Omega)}
		\leq\lambda_1^{-1/2}\|V(t)\|_{L^\infty(\Omega)}
		\|\nabla a\|_{L^2(\Omega)}
		\qquad(a\in H_0^1(\Omega)).
	\end{equation}
	Let $S_E(t)$ be the free wave group on
	$H_0^1(\Omega)\times L^2(\Omega)$. The zero-initial-state Duhamel
	equation is
	\begin{equation}\label{asm:energy-duhamel}
		Y(t)=\int_0^t S_E(t-s)
		\bigl(0,-V(s)y(s)+\mathbf1_\omega f(s)\bigr)\,\mathrm ds,
		\qquad Y(t)=(y(t),y_t(t)).
	\end{equation}
	Estimates \eqref{asm:energy-perturbation}--\eqref{asm:energy-duhamel}
	and the successive-difference bounds used for
	\eqref{asm:picard-iteration} give
	\[
	y\in C([0,T];H_0^1(\Omega))\cap C^1([0,T];L^2(\Omega)),
	\qquad y_{tt}\in L^2(0,T;H^{-1}(\Omega)).
	\]
	Thus both $y$ and the observer $p$ belong to the energy space, and the
	duality pairings below are legitimate.
	
	\smallskip\noindent\textit{Step 2.3: The Green identity kills the
		control pairing.}
	The weak product rule gives the following identity in distributions on
	$(0,T)$:
	\[
	\begin{aligned}
		\frac{\mathrm d}{\mathrm dt}\int_\Omega(y_tp-yp_t)\,\mathrm dx
		&=\langle y_{tt},p\rangle_{H^{-1},H_0^1}
		+\int_\Omega y_tp_t\,\mathrm dx
		-\int_\Omega y_tp_t\,\mathrm dx
		-\int_\Omega yp_{tt}\,\mathrm dx\\
		&=\langle\Delta y-Vy+\mathbf1_\omega f,p\rangle_{H^{-1},H_0^1}
		-\int_\Omega y(\Delta p-Vp)\,\mathrm dx\\
		&=-\int_\Omega\nabla y\cdot\nabla p\,\mathrm dx
		-\int_\Omega Vyp\,\mathrm dx+\int_\omega fp\,\mathrm dx
		+\int_\Omega\nabla y\cdot\nabla p\,\mathrm dx
		+\int_\Omega yVp\,\mathrm dx\\
		&=\int_\omega fp\,\mathrm dx .
	\end{aligned}
	\]
	The first line uses the duality pairing and $y,p\in H_0^1(\Omega)$,
	so both integrations by parts are justified for the energy solution
	and no boundary trace of $\nabla y$ is required. The last line uses
	that $f$ is supported in $(0,T)\times\omega$, hence the pairing with
	$p$ vanishes because $p=0$ on $\R\times\omega$. Thus
	\[
	\frac{\mathrm d}{\mathrm dt}\int_\Omega(y_tp-yp_t)\,\mathrm dx=0
	\quad\text{in }\mathcal D'(0,T).
	\]
	The integral on the left is continuous in time and has an integrable
	distributional derivative. Since $y(0)=y_t(0)=0$, we obtain
	\begin{equation}\label{asm:Green}
		\left[\int_\Omega(y_tp-yp_t)\,\mathrm dx\right]_0^T
		=\int_0^T\!\int_\omega fp\,\mathrm dx\,\mathrm dt=0 .
	\end{equation}
	
	\smallskip\noindent\textit{Step 2.4: A nonzero terminal functional
		vanishes on all reachable states.}
	Equation \eqref{asm:Green} shows that every terminal state reachable
	from zero belongs to the kernel of the continuous functional
	\begin{equation}\label{asm:terminal-functional}
		\ell_T(a,b)=\int_\Omega
		\bigl(b\,p(T)-a\,p_t(T)\bigr)\,\mathrm dx
		\quad\text{on }H_0^1(\Omega)\times L^2(\Omega).
	\end{equation}
	Its continuity follows from
	\[
	|\ell_T(a,b)|\leq
	\|b\|_2\|p(T)\|_2+
	\lambda_1^{-1/2}\|\nabla a\|_2\|p_t(T)\|_2 .
	\]
	The terminal state $(a,b)=(0,p(T))$ belongs to
	$H_0^1(\Omega)\times L^2(\Omega)$ and satisfies
	\[
	\ell_T(0,p(T))=\|p(T)\|_2^2>0
	\]
	by Step 2.1. Hence $(0,p(T))$ is not reachable from zero. This proves
	failure of exact internal controllability for every $T>0$ in every
	dimension $n\ge2$, and completes the proof of the negative assertions.
\end{proof}

\section*{Acknowledgements}

This work began more than ten years ago and made significant progress in
January 2024. An AI tool was used to assist with a substantial amount of
the technical computation in this paper: it was employed to carry out
large-scale algebraic and asymptotic calculations, to check the
consistency of the error budgets in the various regimes of the merging
argument, and to verify that the proposed lines of reasoning could be
completed as stated. The AI tool also helped polish the exposition of an
early version of the manuscript. All mathematical results, including the
counterexample, are due to the authors, who take full responsibility for
the content, for verifying every computation reported here, and for the
final exposition.

\section*{Appendix A. Calculations for the local packet construction}
\addcontentsline{toc}{section}{Appendix A. Calculations for the local packet construction}
\setcounter{equation}{0}
\renewcommand{\theequation}{A.\arabic{equation}}
\renewcommand{\theHequation}{appendix.A.\arabic{equation}}

This appendix gives three groups of calculations for the local packet construction: Gevrey bounds for recursive coefficients, angular means and the characteristic solver, and amplitude asymptotics with high-order tail bounds.

\subsection*{Packet coefficients and Gevrey estimates}
\phantomsection\proofsteplabel{app:Coefficients}{A}

\begingroup
\begin{proof}[Proof of Lemma~\ref{ce:coefficient-bounds}]

\medskip
\noindent\textbf{Step 1. Bound the cost of differentiating an amplitude.}\proofsteplabel{step:coefficient-amplitude}{1}
For $r\ge1$,
\begin{equation}\label{ce:exponential-exact-derivative}
\frac{\mathrm{d}^r}{\mathrm{d}y^r}e^{-1/y}=e^{-1/y}\sum_{h=1}^r(-1)^{r-h}
\frac{r!}{h!}\binom{r-1}{h-1}y^{-r-h}
\end{equation}
For $0<y<1$ and $1\le h\le r$,
\[
y^{-r-h}=y^{-2r}y^{r-h}\le y^{-2r},\qquad \frac1{h!}\le1.
\]
Thus
\begin{equation}\label{ce:exponential-derivative-bound}
\begin{aligned}
\left|\frac{\mathrm{d}^r}{\mathrm{d}y^r}e^{-1/y}\right|
&\le e^{-1/y}r!\sum_{h=1}^r
\frac1{h!}\binom{r-1}{h-1}y^{-r-h}\\
&\le e^{-1/y}r!y^{-2r}\sum_{h=1}^r\binom{r-1}{h-1}\\
&=e^{-1/y}2^{r-1}r!y^{-2r}
\le e^{-1/y}2^rr!y^{-2r}.
\end{aligned}
\end{equation}
Here $\sum_{h=1}^r\binom{r-1}{h-1}=2^{r-1}$ by the binomial theorem.
For $r=0$,
\[
\left|\frac{\mathrm{d}^0}{\mathrm{d}y^0}e^{-1/y}\right|
=e^{-1/y}=e^{-1/y}2^0\,0!\,y^0.
\]
Thus the final bound in \eqref{ce:exponential-derivative-bound} holds for every integer $r\ge0$.

Fix $m\ge 5$, $j\in\{m-1,m,m+1\}$, and $y\in J_m$.
Then $0<y<1$ and
\[
\beta^jy\le\beta^{m+1}\beta^{-(m-1)}=\beta^2.
\]
Let $C_0\ge1$ be the constant in the bound for $\rho$ in
Lemma~\ref{ce:cutoff-bounds}. For every integer $q\ge0$,
\[
\left|\frac{\mathrm{d}^q}{\mathrm{d}y^q}\rho(\beta^jy)\right|
=\beta^{jq}|\rho^{(q)}(\beta^jy)|
\le C_0^{q+1}(q!)^2\beta^{2q}y^{-q}.
\]
Also,
\[
A_{\gamma,j}^{(r)}(y)
=e^{-j^2-\gamma j}\sum_{q=0}^r\binom rq
\left(\frac{\mathrm{d}^{r-q}}{\mathrm{d}y^{r-q}}e^{-1/y}\right)
\beta^{jq}\rho^{(q)}(\beta^jy).
\]
Hence
\[
\begin{aligned}
|A_{\gamma,j}^{(r)}(y)|
&\le e^{-1/y-j^2-\gamma j}C_0y^{-2r}
\sum_{q=0}^r\binom rq 2^{r-q}(r-q)!(q!)^2
(C_0\beta^2y)^q\\
&\le e^{-1/y-j^2-\gamma j}C_0(r!)^2y^{-2r}
\sum_{q=0}^r2^{r-q}(C_0\beta^2)^q\\
&\le e^{-1/y-j^2-\gamma j}C^{r+1}(r!)^2y^{-2r}.
\end{aligned}
\]
Here we used $y^q\le1$,
\[
\binom rq(r-q)!(q!)^2=r!q!\le(r!)^2,
\qquad
\sum_{q=0}^r2^{r-q}(C_0\beta^2)^q\le(2+C_0\beta^2)^r,
\]
and may choose the single constant
\[
C\triangleq\max\left\{C_0,\,2+C_0\beta^2,\,
4+\frac{48\log2}{5}\right\}.
\]
Since $C_0$ is the fixed cutoff constant and $\beta=2^{1/16}$,
this $C$ is absolute and independent of $r$ and $Z$.
Finally, the lower bound \eqref{ce:band-envelope-lower} for $E_\gamma$ on $J_m$ gives
\[
\frac{e^{-1/y-j^2-\gamma j}}{E_\gamma(y)}
\le e^{(m+1)^2-j^2+\gamma(m+1-j)}
\le e^{4m+2\gamma}\le e^{Cm},
\]
because $j\ge m-1$, $m\ge5$, and $\gamma\in\{0,24\log2\}$.
Therefore,
\[
\frac{|A_{\gamma,j}^{(r)}(y)|}{E_\gamma(y)}
\le C^{r+1}(r!)^2y^{-2r}
\frac{e^{-1/y-j^2-\gamma j}}{E_\gamma(y)}
\le e^{Cm}C^{r+1}(r!)^2y^{-2r},
\]
which is \eqref{ce:growing-amplitude}.

\medskip
\noindent\textbf{Step 2. Solve the coefficient recurrence.}\proofsteplabel{step:coefficient-recurrence}{2}
Fix an integer $j\ge4$. We prove
\eqref{ce:coefficient-formula} by induction on $\ell$.
For $\ell=1$, the recursion gives
\[
H_{\gamma,j,1}(z,y)
=\frac{\mathrm i}{2}\int_0^z A_{\gamma,j}''(y)\,\mathrm{d}\sigma
=\frac{\mathrm i}{2}zA_{\gamma,j}''(y),
\]
which is \eqref{ce:coefficient-formula} for $\ell=1$.
Suppose that \eqref{ce:coefficient-formula} holds for some $\ell\ge1$.
For every integer $b\ge1$,
\[
\int_0^z\partial_\sigma^2\left(\frac{\sigma^b}{b!}\right)
\,\mathrm{d}\sigma
=\begin{cases}
0,&b=1,\\[2pt]
\displaystyle\frac{z^{b-1}}{(b-1)!},&b\ge2,
\end{cases}
\qquad
\int_0^z\frac{\sigma^b}{b!}\,\mathrm{d}\sigma
=\frac{z^{b+1}}{(b+1)!}.
\]
Thus \eqref{ce:approximate-solution} and the induction hypothesis give
\begin{equation}\label{ce:coefficient-induction-step}
\begin{aligned}
H_{\gamma,j,\ell+1}(z,y)
&=\left(\frac{\mathrm i}{2}\right)^{\ell+1}
\Biggl[
\sum_{\substack{2\le b\le\ell\\b\equiv\ell\pmod2}}
\frac b\ell\binom\ell{(\ell-b)/2}
\frac{z^{b-1}}{(b-1)!}A_{\gamma,j}^{(\ell+b)}(y)\\
&\qquad+
\sum_{\substack{1\le b\le\ell\\b\equiv\ell\pmod2}}
\frac b\ell\binom\ell{(\ell-b)/2}
\frac{z^{b+1}}{(b+1)!}A_{\gamma,j}^{(\ell+b+2)}(y)
\Biggr]\\
&=\left(\frac{\mathrm i}{2}\right)^{\ell+1}
\sum_{\substack{1\le a\le\ell+1\\a\equiv\ell+1\pmod2}}
\Biggl[
\frac{a-1}{\ell}\binom\ell{(\ell+1-a)/2}
+\frac{a+1}{\ell}\binom\ell{(\ell-1-a)/2}
\Biggr]\\
&\qquad\qquad\times
\frac{z^a}{a!}A_{\gamma,j}^{(\ell+1+a)}(y),
\end{aligned}
\end{equation}
where $\binom\ell q=0$ for integers $q<0$ or $q>\ell$.
For each index $a$ in the final sum of \eqref{ce:coefficient-induction-step}, set $q=(\ell+1-a)/2$.
Then $0\le q\le\ell$ and
\[
\begin{aligned}
\frac{a-1}{\ell}\binom\ell q
+\frac{a+1}{\ell}\binom\ell{q-1} &=
\frac{(a-1)(\ell+1-q)+(a+1)q}{\ell(\ell+1)}
\binom{\ell+1}q\\
& =
\frac{(a-1)(\ell+1)+2q}{\ell(\ell+1)}
\binom{\ell+1}q
=\frac a{\ell+1}\binom{\ell+1}q,
\end{aligned}
\]
since $2q=\ell+1-a$. Consequently,
\[
H_{\gamma,j,\ell+1}(z,y)
=\left(\frac{\mathrm i}{2}\right)^{\ell+1}
\sum_{\substack{1\le a\le\ell+1\\a\equiv\ell+1\pmod2}}
\frac a{\ell+1}\binom{\ell+1}{(\ell+1-a)/2}
\frac{z^a}{a!}A_{\gamma,j}^{(\ell+1+a)}(y).
\]
This proves \eqref{ce:coefficient-formula} for every $\ell\ge1$.
Every summand contains $z^a$ with $a\ge1$, so
$H_{\gamma,j,\ell}(0,y)=0$.

\medskip
\noindent\textbf{Step 3. Combine frequency decay with coefficient growth.}\proofsteplabel{step:coefficient-decay}{3}
Fix $r\ge0$, $m\ge 5$, $y\in J_m$,
$j\in\{m-1,m,m+1\}$, and $|z|\le Z+1$.
Let $1\le\ell\le K_j$ and
$\alpha=(\alpha_t,\alpha_z,\alpha_y)\in\mathbb N_0^3$ with
$|\alpha|\le r$. Formula \eqref{ce:coefficient-formula} gives the exact identity
\[
\begin{aligned}
&\partial^\alpha\bigl(e^{\mathrm i\lambda_ju}
\lambda_j^{-\ell}H_{\gamma,j,\ell}\bigr)\\
& =e^{\mathrm i\lambda_ju}\lambda_j^{-\ell}
\left(\frac{\mathrm i}{2}\right)^\ell
\sum_{\substack{1\le a\le\ell\\a\equiv\ell\pmod2}}
\frac a\ell\binom\ell{(\ell-a)/2} \Biggl[\sum_{p=0}^{\min\{a,\alpha_z\}}
\binom{\alpha_z}{p}(\mathrm i\lambda_j)^{\alpha_t+\alpha_z-p}
\frac{z^{a-p}}{(a-p)!}
A_{\gamma,j}^{(\ell+a+\alpha_y)}(y)\Biggr].
\end{aligned}
\]
For the indices in this sum,
\[
0\le\frac a\ell\binom\ell{(\ell-a)/2}\frac1{a!}
\le2^\ell\le(\ell+1)^\ell
\]
and
\[
\frac{a!}{(a-p)!}\le\ell^p\le\ell^r,\qquad
|z|^{a-p}\le(Z+2)^\ell,\qquad
\sum_{p=0}^{\min\{a,\alpha_z\}}\binom{\alpha_z}{p}\le2^r.
\]
Also, $\ell+a+\alpha_y\le2\ell+r$, so
\eqref{ce:growing-amplitude} and $0<y<1$ imply
\[
\begin{aligned}
\frac{|A_{\gamma,j}^{(\ell+a+\alpha_y)}(y)|}{E_\gamma(y)}
&\le e^{Cm}C^{\ell+a+\alpha_y+1}
[(\ell+a+\alpha_y)!]^2y^{-2(\ell+a+\alpha_y)}\\
&\le e^{Cm}C^{2\ell+r+1}[(2\ell+r)!]^2y^{-4\ell-2r}.
\end{aligned}
\]
There are at most $\ell$ choices of $a$, and
$\lambda_j^{\alpha_t+\alpha_z-p}\le\lambda_j^r$. Hence
\[
\begin{aligned}
\frac{|\partial^\alpha(e^{\mathrm i\lambda_ju}
\lambda_j^{-\ell}H_{\gamma,j,\ell})|}{E_\gamma}
&\le e^{Cm}C(r)^{\ell+1}\lambda_j^{r-\ell}
\ell^{r+1}(\ell+1)^\ell [(2\ell+r)!]^2y^{-4\ell-2r}.
\end{aligned}
\]
The factorial and polynomial factors satisfy
\begin{equation}\label{ce:factorial-polynomial-bound}
\begin{aligned}
[(2\ell+r)!]^2
&\le(2\ell+r)^{4\ell+2r}
\le2^{4\ell+2r}(\ell+r+2)^{4\ell+2r},\\
\ell^{r+1}(\ell+1)^\ell[(2\ell+r)!]^2
&\le2^{4\ell+2r}(\ell+r+2)^{5\ell+3r+1} \le e^{C(r)\ell}(\ell+r+2)^{5\ell}.
\end{aligned}
\end{equation}
The last inequality in \eqref{ce:factorial-polynomial-bound} uses
\[
\ell+r+2\le(r+3)2^{\ell-1},\qquad
(\ell+r+2)^{3r+1}
\le(r+3)^{3r+1}2^{(3r+1)(\ell-1)}
\le e^{C(r)\ell}.
\]
For the frequency and normal factors, $y^{-1}\le\beta^{m+1}$ and
$j\le m+1$ give
\[
\begin{aligned}
\log\bigl(\lambda_j^{r-\ell}y^{-4\ell-2r}\bigr)
&=(r-\ell)j\log2+(4\ell+2r)\log(1/y)\\
&\le-j\ell\log2+4m\ell\log\beta
+4\ell\log\beta +rj\log2+2r(m+1)\log\beta\\
&\le-j\ell\log2+4m\ell\log\beta
+C\ell+C(r)m.
\end{aligned}
\]
Moreover, $\ell\le K_j\le j^2\beta^j$ and $m\le j+1$ imply
\[
\ell+r+2\le j^2\beta^j+r+2\le(r+3)j^2\beta^j,
\]
and therefore
\[
\begin{aligned}
4m\ell\log\beta+5\ell\log(\ell+r+2)
&\le4(j+1)\ell\log\beta
+5\ell\bigl[j\log\beta+2\log j+\log(r+3)\bigr]\\
&\le9j\ell\log\beta+10\ell\log(j+2)+C(r)\ell,\\
-j\ell\log2+9j\ell\log\beta
&=-\frac7{16}j\ell\log2.
\end{aligned}
\]
Thus, for every fixed derivative order $r\ge0$,
\begin{align*}
\log\frac{|\partial^r(e^{\mathrm i\lambda_ju}
\lambda_j^{-\ell}H_{\gamma,j,\ell})|}{E_\gamma}
&\le-j\ell\log2+4m\ell\log\beta
+5\ell\log(\ell+r+2)+C(r)\ell+C(r)m\\
&\le-\tfrac7{16}j\ell\log2
+10\ell\log(j+2)+C(r)\ell+C(r)m.
\end{align*}
At zeros of the numerator, the corresponding exponential upper bound
holds trivially.
This proves \eqref{ce:coefficient-estimate}.
\end{proof}
\endgroup
\subsection*{Angular means and the Bessel characteristic solver}
\phantomsection\proofsteplabel{app:Characteristic}{A}

\begingroup
\begin{proof}[Proof of Lemma~\ref{ce:characteristic-solver}]
Regularize the angular mean at $q=0$, verify the equation by Fourier transform in $y$, then estimate the support and derivatives.

\noindent\textbf{Step 1: Angular derivative estimates.}
Differentiating $\sqrt q$ alone would introduce negative powers of
$q$. Integration by parts transfers these derivatives to $g$ and
gives a formula that remains regular at $q=0$.
For $g\in C^\infty(\R)$, $\nu,k\in\mathbb N_0$, and $q>0$,
\[
\frac{\mathrm d}{\mathrm d\theta}(1-\theta^2)^{\nu+1/2}
=-(2\nu+1)\theta(1-\theta^2)^{\nu-1/2}
\]
and
\[
\begin{aligned}
&\partial_q\int_{-1}^1(1-\theta^2)^{\nu-1/2}
g^{(2\nu+k)}(y+\sqrt q\,\theta)\,\mathrm d\theta\\
& =\frac1{2\sqrt q}\int_{-1}^1
\theta(1-\theta^2)^{\nu-1/2}
g^{(2\nu+k+1)}(y+\sqrt q\,\theta)\,\mathrm d\theta\\
& =\frac1{2(2\nu+1)}\int_{-1}^1
(1-\theta^2)^{\nu+1/2}
g^{(2\nu+k+2)}(y+\sqrt q\,\theta)\,\mathrm d\theta.
\end{aligned}
\]
The boundary terms vanish because $(1-\theta^2)^{\nu+1/2}=0$
at $\theta=\pm1$. Together with
\[
\frac{\nu!}{\pi(2\nu)!}\frac1{2(2\nu+1)}
=\frac{(\nu+1)!}{\pi(2\nu+2)!},
\]
this proves
\begin{equation}\label{ce:angular-derivative-formula}
\partial_q^\nu\partial_y^k\mathcal A_qg(y)
=\frac{\nu!}{\pi(2\nu)!}\int_{-1}^1
(1-\theta^2)^{\nu-1/2}
g^{(2\nu+k)}(y+\sqrt q\,\theta)\,\mathrm d\theta.
\end{equation}
Since
\[
\begin{aligned}
\int_{-1}^1(1-\theta^2)^{\nu-1/2}\,\mathrm d\theta
&=\pi\prod_{j=1}^{\nu}\frac{2j-1}{2j}
=\frac{\pi(2\nu)!}{4^\nu(\nu!)^2},\\
\frac{\nu!}{\pi(2\nu)!}\int_{-1}^1
(1-\theta^2)^{\nu-1/2}\,\mathrm d\theta
&=\frac{\nu!}{\pi(2\nu)!}
\frac{\pi(2\nu)!}{4^\nu(\nu!)^2}
=\frac1{4^\nu\nu!},
\end{aligned}
\]
\eqref{ce:angular-derivative-formula} extends continuously to $q=0$, with value
$g^{(2\nu+k)}(y)/(4^\nu\nu!)$. Taking absolute values proves
\eqref{ce:angular-bound}, including at $q=0$.

\medskip\noindent\textbf{Step 2: Verification of the solution formula.}
The change of variables $a=r\xi$, $b=s\eta$ in
\eqref{ce:causal-inverse} gives the fixed-square representation
\begin{equation}\label{ce:fixed-square-solver}
w(r,s,y)=\frac{rs}{4}\int_0^1\!\int_0^1
\mathcal A_{rs\xi\eta}
[h(r(1-\xi),s(1-\eta),\cdot)](y)\,\mathrm d\eta\,\mathrm d\xi.
\end{equation}
The derivative formula \eqref{ce:angular-derivative-formula} allows differentiation under this integral, including when
$r=0$ or $s=0$, and proves that $w$ is smooth up to the boundary.
It is real because $h$ and the averaging kernel are real.
The Bessel function gives the Fourier multiplier of the angular mean. We take the Fourier transform only in $y$, with convention
\[
\widehat g(\zeta)\triangleq\int_{\R}e^{-\mathrm i y\zeta}g(y)\,\mathrm dy,
\qquad
g(y)=\frac1{2\pi}\int_{\R}e^{\mathrm i y\zeta}\widehat g(\zeta)
\,\mathrm d\zeta,
\]
initially for $g\in C_c^\infty(\R)$; the real number $\zeta$ is the
frequency corresponding to $y$. In particular, integration by parts
gives $\widehat{-\partial_y^2g}(\zeta)=\zeta^2\widehat g(\zeta)$.
Define the order-zero Bessel function of the first kind by
\[
J_0(x)\triangleq\frac1\pi\int_{-1}^1
\frac{e^{\mathrm i x\theta}}{\sqrt{1-\theta^2}}\,\mathrm d\theta
=\frac1\pi\int_{-1}^1
\frac{\cos(x\theta)}{\sqrt{1-\theta^2}}\,\mathrm d\theta,
\qquad x\in\R.
\]
The imaginary part vanishes by oddness. To obtain its power series,
the substitution $\theta=\sin\varphi$ and integration by parts give,
for $k\ge0$,
\[
\frac1\pi\int_{-1}^1\frac{\theta^{2k}}{\sqrt{1-\theta^2}}
\,\mathrm d\theta
=\frac2\pi\int_0^{\pi/2}\sin^{2k}\varphi\,\mathrm d\varphi
=\prod_{j=1}^k\frac{2j-1}{2j}
=\frac{(2k)!}{4^k(k!)^2},
\]
where the empty product for $k=0$ is $1$.
The odd moments are zero, so termwise integration of the exponential
series yields
\begin{equation}\label{ce:bessel-series}
J_0(x)=\sum_{k=0}^{\infty}\frac{(-1)^k}{(k!)^2}
\left(\frac{x^2}{4}\right)^k.
\end{equation}
Next, Fubini's theorem and the change of variable
$x=y+\sqrt q\,\theta$ show directly that
\[
\begin{aligned}
\widehat{\mathcal A_qg}(\zeta)
&=\frac1\pi\int_{-1}^1\frac1{\sqrt{1-\theta^2}}
\left[\int_{\R}e^{-\mathrm i y\zeta}
g(y+\sqrt q\,\theta)\,\mathrm dy\right]\mathrm d\theta =\left[\frac1\pi\int_{-1}^1
\frac{e^{\mathrm i\zeta\sqrt q\,\theta}}{\sqrt{1-\theta^2}}
\,\mathrm d\theta\right]\widehat g(\zeta)
=J_0(\zeta\sqrt q)\widehat g(\zeta).
\end{aligned}
\]
Set $F_\zeta(q)\triangleq J_0(\zeta\sqrt q)$ for $q\ge0$.
The series \eqref{ce:bessel-series} gives
\[
F_\zeta(q)=\sum_{k=0}^{\infty}
\frac{(-\zeta^2/4)^kq^k}{(k!)^2},\qquad F_\zeta(0)=1.
\]
Differentiating this convergent power series term by term, we find
\begin{equation}\label{ce:bessel-ode}
\begin{aligned}
qF_\zeta''(q)+F_\zeta'(q)
&=\sum_{k=1}^{\infty}
\frac{k^2(-\zeta^2/4)^kq^{k-1}}{(k!)^2} =-\frac{\zeta^2}{4}\sum_{k=0}^{\infty}
\frac{(-\zeta^2/4)^kq^k}{(k!)^2}
=-\frac{\zeta^2}{4}F_\zeta(q).
\end{aligned}
\end{equation}
To use \eqref{ce:bessel-ode} in the two characteristic variables, put
\[
K_\zeta(a,b)\triangleq
\begin{cases}
\tfrac14 F_\zeta(ab),&a\ge0,\ b\ge0,\\
0,&\text{otherwise}.
\end{cases}
\]
Write $\delta_0$ for the Dirac mass at zero. In distributions,
$\partial_a\mathbf1_{a\ge0}=\delta_0(a)$ and
$a\delta_0(a)=0$, and similarly for $b$. Differentiating $K_\zeta$
therefore gives
\[
\begin{aligned}
(4\partial_a\partial_b+\zeta^2)K_\zeta
&=F_\zeta(0)\delta_0(a)\delta_0(b) +\mathbf1_{a\ge0}\mathbf1_{b\ge0}
\left[F_\zeta'(ab)+abF_\zeta''(ab)
+\frac{\zeta^2}{4}F_\zeta(ab)\right] =\delta_0(a)\delta_0(b).
\end{aligned}
\]
Taking the Fourier transform in $y$ in \eqref{ce:causal-inverse}
now gives
\[
\widehat w(r,s,\zeta)
=\frac14\int_0^r\int_0^s
F_\zeta(ab)\widehat h(r-a,s-b,\zeta)\,\mathrm db\,\mathrm da
=(K_\zeta*\widehat h)(r,s,\zeta),
\]
where $*$ denotes convolution in $(r,s)$ and $\widehat h$ is
extended by zero outside its rectangle in these variables. Hence,
in the interior of the rectangle,
\[
(4\partial_r\partial_s+\zeta^2)\widehat w
=\bigl[(4\partial_r\partial_s+\zeta^2)K_\zeta\bigr]
*\widehat h
=\widehat h.
\]
Fourier inversion proves
$(4\partial_r\partial_s-\partial_y^2)w=h$.
The characteristic conditions $w(0,s,y)=w(r,0,y)=0$ follow
directly from the zero integration ranges in \eqref{ce:causal-inverse}.

\medskip\noindent\textbf{Step 3: Dependence on the source and support.}
Writing out both integrals in \eqref{ce:causal-inverse}, we obtain
\[
w(r,s,y)=\frac1{4\pi}\int_0^r\int_0^s\int_{-1}^1
\frac{h(r-a,s-b,y+\sqrt{ab}\,\theta)}{\sqrt{1-\theta^2}}
\,\mathrm d\theta\,\mathrm db\,\mathrm da.
\]
Every source point in this integral has
\[
r'=r-a\in[0,r],\qquad s'=s-b\in[0,s],\qquad
|y'-y|=|\sqrt{ab}\,\theta|\le\sqrt{rs}\le\sqrt{dL}.
\]
If $y<y_- -\sqrt{dL}$ or $y>y_+ +\sqrt{dL}$, then every such
$y'$ lies outside $[y_-,y_+]$, so the integrand vanishes. This proves
\eqref{ce:causal-support}. All derivatives of $w$ vanish on the same
open set, so their supports satisfy the same inclusion.

\medskip\noindent\textbf{Step 4: Derivative estimates for the solution.}
The integration region in \eqref{ce:fixed-square-solver} is independent of $r,s$.
Its product and chain rules produce only polynomial coefficients in
$r,s,\xi,\eta$, so shrinking $d$ introduces no inverse power of $d$.
For every multiindex $\alpha$, differentiation gives a finite sum
of the form
\[
\begin{aligned}
\partial_{r,s,y}^\alpha w(r,s,y)
&=\sum_{\substack{0\le\ell\le|\alpha|\\
b,c,k\ge0,\ b+c+k\le|\alpha|}}
\int_0^1\!\int_0^1 P_{\alpha,\ell,b,c,k}(r,s,\xi,\eta)\\
&\quad\times\left.\partial_q^\ell\mathcal A_q
[(\partial_r^b\partial_s^c\partial_y^k h)
(r(1-\xi),s(1-\eta),\cdot)](y)\right|_{q=rs\xi\eta}
\,\mathrm d\eta\,\mathrm d\xi,
\end{aligned}
\]
where the coefficients $P_{\alpha,\ell,b,c,k}$ are polynomials.
They are bounded by constants depending only on $\alpha,D,L$ on
$0\le r\le d\le D$, $0\le s\le L$, and $0\le\xi,\eta\le1$.
By \eqref{ce:angular-bound}, each term requires at most
$b+c+k+2\ell\le3|\alpha|$ source derivatives, evaluated at points
with $0\le r'\le r$, $0\le s'\le s$, and
$|y'-y|\le\sqrt{rs\xi\eta}\le\sqrt{rs}$.
Integrating over the unit square and summing the finitely many
terms proves \eqref{ce:causal-derivative-bound}.
\end{proof}
\endgroup
\subsection*{Amplitude asymptotics and high-order tails}
\phantomsection\proofsteplabel{app:PacketTail}{A}

\begingroup
\begin{proof}[Proof of Lemma~\ref{ce:packet-bounds}]
We prove (i) and (ii) in Steps~\ref{step:packet-normal} and~\ref{step:packet-tail}, respectively.
In (i), the term $y^{-2r}A_{\gamma,j}$ is the leading contribution
from differentiating $e^{-1/y}$ $r$ times. In (ii), $L$ is fixed
after $r$ and $B$ are fixed; the actual truncation indices $K_j$
and the constructed solution remain unchanged.

\medskip
\noindent\textbf{Step 1: Separate the common exponential from the cutoffs.}\proofsteplabel{step:packet-normal}{1}
Fix an integer $r\ge1$. We first work on
$\beta^{-(m+1)}\le y\le\beta^{-m}$ for any integer $m\ge4$.
Only $A_{\gamma,m}$ and $A_{\gamma,m+1}$ and their derivatives can be
nonzero on this interval; see Figure~\ref{fig:packet-amplitudes}.
For $j=m,m+1$,
\[
\begin{aligned}
A_{\gamma,j}^{(r)}-y^{-2r}A_{\gamma,j}
&=e^{-j^2-\gamma j}\rho(\beta^jy)
\left(\frac{\mathrm{d}^r}{\mathrm{d}y^r}e^{-1/y}
-y^{-2r}e^{-1/y}\right)\\
&\quad+e^{-j^2-\gamma j}\sum_{q=1}^r\binom rq
\left(\frac{\mathrm{d}^{r-q}}{\mathrm{d}y^{r-q}}e^{-1/y}\right)
\frac{\mathrm{d}^q}{\mathrm{d}y^q}\rho(\beta^jy).
\end{aligned}
\]
We claim that
\begin{equation}\label{ce:exponential-derivative-remainder}
\left|\frac{\mathrm{d}^r}{\mathrm{d}y^r}e^{-1/y}
-y^{-2r}e^{-1/y}\right|
\le C(r)e^{-1/y}y^{-2r+1},
\end{equation}
and, for every integer $q\ge1$,
\begin{equation}\label{ce:weighted-cutoff-derivatives}
e^{-1/y}\sum_{j=m}^{m+1}e^{-j^2-\gamma j}
\left|\frac{\mathrm{d}^q}{\mathrm{d}y^q}\rho(\beta^jy)\right|
\le C(q)E_\gamma(y)m^{2q}\beta^{mq},\qquad q\ge1.
\end{equation}

Under \eqref{ce:exponential-derivative-remainder},
\eqref{ce:weighted-cutoff-derivatives}, and the derivative bound \eqref{ce:exponential-derivative-bound},
\begin{equation}\label{ce:normal-amplitude-half-band}
\begin{aligned}
\sum_{j\ge4}|A_{\gamma,j}^{(r)}-y^{-2r}A_{\gamma,j}|
&\le C(r)E_\gamma y^{-2r+1}
+C(r)E_\gamma\sum_{q=1}^ry^{-2(r-q)}m^{2q}\beta^{mq}\\
&=C(r)E_\gamma y^{-2r}
\left[y+\sum_{q=1}^rm^{2q}(\beta^my^2)^q\right].
\end{aligned}
\end{equation}
Since $y\le\beta^{-m}$, we have $\beta^my^2\le\beta^{-m}<1$.
Also $m^{2q}\le m^{2r}$ for $1\le q\le r$. Therefore,
\[
y+\sum_{q=1}^rm^{2q}(\beta^my^2)^q
\le\beta^{-m}+\sum_{q=1}^rm^{2r}\beta^{-mq}
\le(r+1)m^{2r}\beta^{-m}.
\]
Thus the two claims imply \eqref{ce:normal-amplitude} on
$[\beta^{-(m+1)},\beta^{-m}]$.

To prove \eqref{ce:exponential-derivative-remainder}, use the exact derivative formula
\eqref{ce:exponential-exact-derivative}: for $r\ge2$,
\[
\frac{\mathrm{d}^r}{\mathrm{d}y^r}e^{-1/y}-y^{-2r}e^{-1/y}
=e^{-1/y}\sum_{h=1}^{r-1}(-1)^{r-h}
\frac{r!}{h!}\binom{r-1}{h-1}y^{-r-h}.
\]
For $1\le h\le r-1$ and $0<y<1$, we have
\[
y^{-r-h}=y^{-2r+1}y^{r-h-1}\le y^{-2r+1},
\]
because $r-h-1\ge0$.
Since $1/h!\le1$,
\begin{equation}\label{ce:exponential-remainder-proof}
\begin{aligned}
\left|\frac{\mathrm{d}^r}{\mathrm{d}y^r}e^{-1/y}
-y^{-2r}e^{-1/y}\right|
&\le e^{-1/y}r!\sum_{h=1}^{r-1}
\frac1{h!}\binom{r-1}{h-1}y^{-r-h} \le e^{-1/y}r!y^{-2r+1}
\sum_{h=1}^{r-1}\binom{r-1}{h-1}\\
&=(2^{r-1}-1)r!\,e^{-1/y}y^{-2r+1}.
\end{aligned}
\end{equation}
The last equality in \eqref{ce:exponential-remainder-proof} follows from
\[
\sum_{h=1}^{r-1}\binom{r-1}{h-1}
=\sum_{p=0}^{r-2}\binom{r-1}{p}
=2^{r-1}-\binom{r-1}{r-1}=2^{r-1}-1.
\]
Thus \eqref{ce:exponential-derivative-remainder} holds for $r\ge2$;
in \eqref{ce:exponential-derivative-remainder} one may take $C(r)=2^{r-1}r!$, independently of $y$.
For $r=1$,
\[
\frac{\mathrm{d}}{\mathrm{d}y}e^{-1/y}-y^{-2}e^{-1/y}
=e^{-1/y}\frac1{y^2}-y^{-2}e^{-1/y}=0,
\]
so \eqref{ce:exponential-derivative-remainder} also holds in this case.

We next prove \eqref{ce:weighted-cutoff-derivatives}.
Set $\xi\triangleq(\beta^{m+1}y-1)/(\beta-1)$, so $0\le\xi\le1$.
The definitions give
\[
\rho(\beta^my)=b(\xi),\qquad
\rho(\beta^{m+1}y)=b(1-\xi)=1-b(\xi).
\]
Consequently,
\[
E_\gamma(y)=e^{-1/y-(m+1)^2-\gamma(m+1)}
\bigl[e^{2m+1+\gamma}b(\xi)+1-b(\xi)\bigr].
\]

We first estimate derivatives of $b$ relative to its size near the endpoints.
The denominator $a(\xi)+a(1-\xi)$ lies between $e^{-2}$ and $2$ on $[0,1]$,
so every fixed-order derivative of its reciprocal is bounded.
For $q\ge1$,
\[
b^{(q)}(\xi)=\sum_{p=0}^q\binom qp a^{(p)}(\xi)
\frac{\mathrm d^{q-p}}{\mathrm d\xi^{q-p}}
\bigl[a(\xi)+a(1-\xi)\bigr]^{-1},
\]
and the bound \eqref{ce:exponential-derivative-bound}, applied to $a(\xi)=e^{-1/\xi}$, implies
\[
\begin{aligned}
|b^{(q)}(\xi)|
&\le C(q)\sum_{p=0}^q|a^{(p)}(\xi)|
\le C(q)a(\xi)\xi^{-2q} \le C(q)b(\xi)\xi^{-2q},\qquad 0<\xi\le\tfrac12.
\end{aligned}
\]
Also, $b^{(q)}(\xi)=(-1)^{q+1}b^{(q)}(1-\xi)$, so
\[
|b^{(q)}(\xi)|\le C(q)(1-b(\xi))(1-\xi)^{-2q},
\qquad \tfrac12\le\xi<1.
\]
Moreover, directly from the definition of $b$,
\[
\frac{b(\xi)}{1-b(\xi)}
=\exp\left(-\frac1\xi+\frac1{1-\xi}\right).
\]

For $0<\xi\le1/2$, the denominator below is at least both
$e^{2m+1+\gamma}b(\xi)$ and $1-b(\xi)$, while
$e^{2m+1+\gamma}+1\le2e^{2m+1+\gamma}$. Hence
\begin{equation}\label{ce:cutoff-ratio-left}
\begin{aligned}
\frac{(e^{2m+1+\gamma}+1)|b^{(q)}(\xi)|}
{e^{2m+1+\gamma}b(\xi)+1-b(\xi)}
&\le C(q)\min\left\{1,
e^{2m+1+\gamma}\frac{b(\xi)}{1-b(\xi)}\right\}\xi^{-2q}\\
&\le C(q)\min\{1,e^{2m+3+\gamma-1/\xi}\}\xi^{-2q}.
\end{aligned}
\end{equation}
To bound the right side of \eqref{ce:cutoff-ratio-left}, put $t=\xi^{-1}$.
If $t\le4m+2\gamma+6$, then $t^{2q}\le C(q)m^{2q}$.
If $t>4m+2\gamma+6$, then
\[
e^{2m+3+\gamma-t}t^{2q}\le e^{-t/2}t^{2q}
\le\sup_{s>0}e^{-s/2}s^{2q}
=\left(\frac{4q}{e}\right)^{2q}.
\]
Thus the ratio is bounded by $C(q)m^{2q}$ on this half of $[0,1]$.
For $1/2\le\xi<1$,
\begin{equation}\label{ce:cutoff-ratio-right}
\begin{aligned}
\frac{(e^{2m+1+\gamma}+1)|b^{(q)}(\xi)|}
{e^{2m+1+\gamma}b(\xi)+1-b(\xi)}
&\le C(q)\frac{1-b(\xi)}{b(\xi)}(1-\xi)^{-2q} \le C(q)e^{2-1/(1-\xi)}(1-\xi)^{-2q}
\le C(q).
\end{aligned}
\end{equation}
Equations \eqref{ce:cutoff-ratio-left} and \eqref{ce:cutoff-ratio-right} extend continuously to $\xi=0,1$: the endpoint derivative
bounds tend to zero because $e^{-1/\xi}\xi^{-2q}\to0$ as
$\xi\downarrow0$, and similarly at $\xi=1$.
Thus, for $0\le\xi\le1$,
\begin{equation}\label{ce:weighted-b-ratio}
\frac{(e^{2m+1+\gamma}+1)|b^{(q)}(\xi)|}
{e^{2m+1+\gamma}b(\xi)+1-b(\xi)}
\le C(q)m^{2q},\qquad q\ge1.
\end{equation}

For $q\ge1$,
\[
\left|\frac{\mathrm{d}^q}{\mathrm{d}y^q}\rho(\beta^my)\right|
=\left|\frac{\mathrm{d}^q}{\mathrm{d}y^q}\rho(\beta^{m+1}y)\right|
=\left(\frac{\beta^{m+1}}{\beta-1}\right)^q|b^{(q)}(\xi)|.
\]
By \eqref{ce:weighted-b-ratio},
\[
\begin{aligned}
&e^{-1/y}\sum_{j=m}^{m+1}e^{-j^2-\gamma j}
\left|\frac{\mathrm{d}^q}{\mathrm{d}y^q}\rho(\beta^jy)\right|\\
& =e^{-1/y-(m+1)^2-\gamma(m+1)}
(e^{2m+1+\gamma}+1)
\left(\frac{\beta^{m+1}}{\beta-1}\right)^q|b^{(q)}(\xi)|\\
& =E_\gamma(y)\left(\frac{\beta^{m+1}}{\beta-1}\right)^q
\frac{(e^{2m+1+\gamma}+1)|b^{(q)}(\xi)|}
{e^{2m+1+\gamma}b(\xi)+1-b(\xi)}\\
& \le C(q)E_\gamma(y)
\left(\frac{\beta}{\beta-1}\right)^qm^{2q}\beta^{mq}
\le C(q)E_\gamma(y)m^{2q}\beta^{mq}.
\end{aligned}
\]
This is \eqref{ce:weighted-cutoff-derivatives}.

For the full band $J_m$ with $m\ge 5$, it remains to consider
$\beta^{-m}\le y\le\beta^{-(m-1)}$.
Here $m-1\ge4$, the only contributing functions are
$A_{\gamma,m-1}$ and $A_{\gamma,m}$, and \eqref{ce:normal-amplitude-half-band}, with $m$ replaced by $m-1$, reads
\[
\begin{aligned}
\sum_{j\ge4}|A_{\gamma,j}^{(r)}-y^{-2r}A_{\gamma,j}|
&\le C(r)E_\gamma y^{-2r}(m-1)^{2r}\beta^{-(m-1)} \le C(r)\beta E_\gamma y^{-2r}m^{2r}\beta^{-m}.
\end{aligned}
\]
Indeed,
\[
(m-1)^{2r}\beta^{-(m-1)}
\le\beta\,m^{2r}\beta^{-m}.
\]
Since $\beta$ is fixed, \eqref{ce:normal-amplitude} holds on all of $J_m$.

\medskip
\noindent\textbf{Step 2: Sum the terms after a fixed index.}\proofsteplabel{step:packet-tail}{2}
Fix an integer $r\ge0$ and $B>0$. Our goal is a bound
$E_\gamma e^{-Bm}$ for the sum over $\ell\ge L$, with $L$
independent of the band index $m$. The coefficient growth must be
included before summing; the factor $\lambda_j^{-\ell}$ alone
does not give the required estimate. Let $C_r\ge1$ be a fixed constant
such that \eqref{ce:coefficient-estimate} implies, for every
$|\alpha|=r$, $y\in J_m$, $j\in\{m-1,m,m+1\}$, and
$1\le\ell\le K_j$,
\[
\begin{aligned}
\frac{|\partial^\alpha(e^{\mathrm i\lambda_ju}
\lambda_j^{-\ell}H_{\gamma,j,\ell})|}{E_\gamma}
&\le\exp\left\{-\frac7{16}j\ell\log2
+10\ell\log(j+2)+C_r\ell+C_rm\right\}.
\end{aligned}
\]
The constant $C_r=C_r(Z)$ is independent of $j,m,\ell$, the spacetime
point, and the choice of $\gamma$.
Set $c_1\triangleq7\log2/64$. Since
\[
\lim_{m\to\infty}\frac{10\log(m+3)+C_r}{m}=0,
\]
there is an integer $m_r\ge5$ such that
\[
10\log(m+3)+C_r\le c_1m\qquad(m\ge m_r).
\]
For $m\ge m_r$ and $|j-m|\le1$, we have
$
j\ge m-1\ge m/2$ and $j+2\le m+3$, 
and hence
\[
\begin{aligned}
-\frac7{16}j\ell\log2+10\ell\log(j+2)+C_r\ell
&\le\left[-\frac7{32}m\log2
+10\log(m+3)+C_r\right]\ell\\
&\le(-2c_1m+c_1m)\ell=-c_1m\ell.
\end{aligned}
\]
Thus
\begin{equation}\label{ce:tail-coefficient-bound}
\frac{|\partial^\alpha(e^{\mathrm i\lambda_ju}
\lambda_j^{-\ell}H_{\gamma,j,\ell})|}{E_\gamma}
\le e^{C_rm-c_1m\ell}\qquad(m\ge m_r).
\end{equation}

Choose the integer
\[
L=L(r,B)\triangleq
\max\left\{K_{m_r}+1,
\left\lceil\frac{B+C_r}{c_1}\right\rceil\right\}.
\]
In particular,
\[
L>K_{m_r},\qquad c_1L-C_r\ge B.
\]
Both $m_r$ and $K_{m_r}$ depend only on $r$ and the fixed data, so
$L$ is independent of $m$ and the spacetime point.

On $J_m$, every derivative of $H_{\gamma,j,\ell}$ vanishes for
$j\notin\{m-1,m,m+1\}$. Therefore, for $|\alpha|=r$,
\begin{equation}\label{ce:tail-derivative-sum}
\begin{aligned}
&\partial^\alpha\operatorname{Re}
\sum_{j\ge4}e^{\mathrm i\lambda_ju}
\sum_{\ell=L}^{K_j}\lambda_j^{-\ell}H_{\gamma,j,\ell} =\operatorname{Re}\sum_{j=m-1}^{m+1}\sum_{\ell=L}^{K_j}
\partial^\alpha\bigl(e^{\mathrm i\lambda_ju}
\lambda_j^{-\ell}H_{\gamma,j,\ell}\bigr).
\end{aligned}
\end{equation}
For $m\ge m_r$, \eqref{ce:tail-derivative-sum} and \eqref{ce:tail-coefficient-bound} give
\begin{equation}\label{ce:tail-absolute-bound}
\begin{aligned}
&\frac1{E_\gamma}\left|\partial^\alpha\operatorname{Re}
\sum_{j\ge4}e^{\mathrm i\lambda_ju}
\sum_{\ell=L}^{K_j}\lambda_j^{-\ell}H_{\gamma,j,\ell}\right| \le\sum_{j=m-1}^{m+1}\sum_{\ell=L}^{K_j}
\frac{|\partial^\alpha(e^{\mathrm i\lambda_ju}
\lambda_j^{-\ell}H_{\gamma,j,\ell})|}{E_\gamma}\\
& \le3e^{C_rm}\sum_{\ell=L}^{\infty}e^{-c_1m\ell} =3e^{C_rm-c_1mL}
\sum_{q=0}^{\infty}(e^{-c_1m})^q
=\frac{3e^{-(c_1L-C_r)m}}{1-e^{-c_1m}} \le\frac{3}{1-e^{-c_1m_r}}e^{-Bm},
\end{aligned}
\end{equation}
where
\[
0<e^{-c_1m}\le e^{-c_1m_r}<1,\qquad
1-e^{-c_1m}\ge1-e^{-c_1m_r}>0.
\]

For the remaining integers $5\le m<m_r$, every contributing
index satisfies $4\le j\le m+1\le m_r$. Consequently,
\[
K_j=\lfloor j^2\beta^j\rfloor
\le\lfloor m_r^2\beta^{m_r}\rfloor=K_{m_r}<L,
\qquad
\sum_{\ell=L}^{K_j}\lambda_j^{-\ell}H_{\gamma,j,\ell}=0.
\]
Thus the tail and all its derivatives vanish on these bands.
The estimate in (ii) holds for every $m\ge 5$, with
$
C(r,B)=\frac{3}{1-e^{-c_1m_r}}.
$
\end{proof}
\endgroup

\section*{Appendix B: repair and finite-dimensional perturbation}
\addcontentsline{toc}{section}{Appendix B: repair and finite-dimensional perturbation}
\setcounter{equation}{0}
\renewcommand{\theequation}{B.\arabic{equation}}
\renewcommand{\theHequation}{appendix.B.\arabic{equation}}
\subsection*{B.1 Identities for the repair coefficients}\proofsteplabel{app:repair-details}{B.1}
\begingroup
\begin{proof}[Proof of Lemma~\ref{ce:repair}]
For smoothness at the zeros of $F$, take the open sets
\[
U_+=\{F>0\},\qquad U_0=O\setminus \operatorname{supp}_O N.
\]
The hypothesis $F>0$ on $\operatorname{supp}_O N$ implies $U_+\cup U_0=O$.
On $U_+$, the denominator $F$ is positive, so the quotient in \eqref{ce:repair-b} is real and smooth. On $U_0$,
the function $N$ vanishes identically. Consequently, the same quotient
is identically zero on $U_+\cap U_0$, and agrees there with the
zero function on $U_0$. These two smooth definitions therefore glue
to a real smooth function $\beta_l$ on $O$, with precisely the values
prescribed in \eqref{ce:repair-b}. In particular, each point of
$\{F=0\}$ has an open neighborhood on which every $\beta_l$ is zero;
no lower bound for $F$ away from such neighborhoods is needed.
Since $D_q$ is a smooth real vector field, all the functions defined
in \eqref{ce:repair-exponent}--\eqref{ce:repair-output}, as well as
$p=e^{\Phi_{\rm rep}} q$, are real and smooth on $O$.

On $U_+$, the definition of $F$ gives
\begin{equation}\label{ce:repair-coefficient-sum}
2\sum_{l=0}^r\mathcal C_l\beta_l
=-\frac{N}{F}\sum_{l=0}^r \mathcal C_l^2=-N.
\end{equation}
On $U_0$, both $N$ and every $\beta_l$ vanish, so \eqref{ce:repair-coefficient-sum} holds
throughout $O$. To differentiate the exponent, for each
$0\le l\le r$ write
\[
A_l\triangleq\sum_{a=0}^l(-1)^a\mathcal C_{l-1-a}D_q^a\beta_l,
\qquad \Phi_{\rm rep}=\sum_{l=0}^r A_l.
\]
Using $D_q\mathcal C_{l-1-a}=\mathcal C_{l-a}$ for $0\le a\le l$ and applying
the Leibniz rule term by term yields
\begin{equation}\label{ce:repair-telescoping}
\begin{aligned}
D_qA_l
&=\sum_{a=0}^l(-1)^a\mathcal C_{l-a}D_q^a\beta_l
+\sum_{a=0}^l(-1)^a\mathcal C_{l-1-a}D_q^{a+1}\beta_l\\
&=\sum_{a=0}^l(-1)^a\mathcal C_{l-a}D_q^a\beta_l
+\sum_{a=1}^{l+1}(-1)^{a-1}\mathcal C_{l-a}D_q^a\beta_l\\
&=\mathcal C_l\beta_l+(-1)^lqD_q^{l+1}\beta_l.
\end{aligned}
\end{equation}
In the second line of \eqref{ce:repair-telescoping}, the terms with $1\le a\le l$ cancel in
pairs; the surviving endpoints use $\mathcal C_{-1}=q$. When $l=0$, this
is simply $D_q(q\beta_0)=\mathcal C_0\beta_0+qD_q\beta_0$, so the argument also covers
$r=0$. Summing over $l$ and using the definition of ${H_{\rm rep}}$, we obtain
throughout $O$
\begin{equation}\label{ce:repair-multiplier-identity}
2D_q\Phi_{\rm rep}
=2\sum_{l=0}^r\mathcal C_l\beta_l
+2q\sum_{l=0}^r(-1)^lD_q^{l+1}\beta_l
=-N+q{H_{\rm rep}}.
\end{equation}

Differentiating the exponential gives
\[
\mathrm{d}(e^{\Phi_{\rm rep}})=e^{\Phi_{\rm rep}}\mathrm{d}\Phi_{\rm rep},\qquad
\Box(e^{\Phi_{\rm rep}})
=e^{\Phi_{\rm rep}}\bigl(\Box\Phi_{\rm rep}+\mathcal B(\mathrm{d}\Phi_{\rm rep},\mathrm{d}\Phi_{\rm rep})\bigr).
\]
The product rule for $\Box$ therefore implies
\begin{equation}\label{ce:repair-conjugation-check}
\begin{aligned}
\Box(e^{\Phi_{\rm rep}} q)
&=e^{\Phi_{\rm rep}}\Box q+2\mathcal B(\mathrm{d}(e^{\Phi_{\rm rep}}),\mathrm{d}q)
+q\Box(e^{\Phi_{\rm rep}})\\
&=e^{\Phi_{\rm rep}}\bigl[\Box q+2D_q\Phi_{\rm rep}
+q\{\Box\Phi_{\rm rep}+\mathcal B(\mathrm{d}\Phi_{\rm rep},\mathrm{d}\Phi_{\rm rep})\}\bigr]\\
&=e^{\Phi_{\rm rep}} q\bigl[G+{H_{\rm rep}}+\Box\Phi_{\rm rep}
+\mathcal B(\mathrm{d}\Phi_{\rm rep},\mathrm{d}\Phi_{\rm rep})\bigr].
\end{aligned}
\end{equation}
The last equality in \eqref{ce:repair-conjugation-check} uses $\Box q=N+Gq$ and \eqref{ce:repair-multiplier-identity}.
The choice of $V$ in \eqref{ce:repair-output} now gives
$\Box p+Vp=0$ on $O$. No step divides by $q$, so this verification
also holds at every zero of $q$.
\end{proof}
\endgroup
\subsection*{B.2 Proof of the finite-dimensional perturbation lemma}\proofsteplabel{app:perturbation-details}{B.2}
\begingroup
\emph{Proof of the claim \eqref{ce:completion-perturbation-claim}.}
We first treat the case where $\{2\delta\le X\le1-2\delta\}$ is
nonempty. 

\smallskip\noindent\textit{Excluded derivative data.}
At a point $z$, collect a real function's value and all its partial
derivatives through order $r+1$ (its $(r+1)$-jet, denoted by
$j_z^{r+1}h$). As $z$ varies over $(\R/T\mathbb Z)\times U$, these data
form the space $J^{r+1}((\R/T\mathbb Z)\times U;\R)$, called the
$(r+1)$-jet bundle. The task is to keep the derivative data of $q$
outside the following two sets.
For a scalar function $h$, write
$Q_h\triangleq\mathcal B(\mathrm dh,\mathrm dh)$ and use $D_h$ from
\eqref{ce:lorentz}. Define
\[
\begin{aligned}
\mathcal Z_0&\triangleq\{h=0,\ \mathrm dh=0\},\\
\mathcal Z_1&\triangleq\{h=0,\ \mathrm dh\ne0,\ D_h^\ell Q_h=0
\text{ for every }\ell\in\{0,\ldots,r\}\}.
\end{aligned}
\]
The first set consists of zeros where all first derivatives vanish;
then $Q_h=0$ and $D_h=0$ at that point, so the whole chain vanishes.
In the second set the first derivative is nonzero, but every chain
entry still vanishes. These alternatives exhaust failure of the
required condition:
\[
h(z)=0,\quad\sum_{\ell=0}^r(D_h^\ell Q_h(z))^2=0
\quad\Longleftrightarrow\quad
j_z^{r+1}h\in\mathcal Z_0\cup\mathcal Z_1.
\]

\smallskip\noindent\textit{Codimension of the excluded data.}
The equations $h=0$ and $\mathrm dh=0$ prescribe one value and $n+1$
first derivatives independently. Hence $\mathcal Z_0$ is a smooth
submanifold of codimension $n+2$.
For $\mathcal Z_1$, put $w\triangleq(\partial_t h,-\nabla_xh)$, so
$D_h=\sum_{\mu=0}^n w^\mu\partial_\mu$. The differential of $Q_h$
with respect to the first derivatives is
\begin{equation}\label{ce:completion-Q-first-derivative}
\mathrm d_{(\mathrm dh)}Q_h[\xi]=2\mathcal B(\mathrm dh,\xi),
\qquad
\frac{\partial Q_h}{\partial(\partial_\mu h)}=2w^\mu.
\end{equation}
It is nonzero when $\mathrm dh\ne0$. For the higher chain entries,
\[
D_h^\ell Q_h=2\partial^{\ell+1}h(w,\ldots,w)
+\text{a polynomial in }\partial^\alpha h,\ |\alpha|\le\ell,
\qquad 1\le\ell\le r.
\]
Indeed, $D_hQ_h=2\partial^2h(w,w)$, and differentiation of the leading
term gives
\begin{equation}\label{ce:completion-chain-induction}
\begin{cases}
\ds D_h[2\partial^{\ell+1}h(w,\ldots,w)]
=2\partial^{\ell+2}h(w,\ldots,w) +2(\ell+1)\partial^{\ell+1}h(D_hw,w,\ldots,w),\\
\ns\ds (D_hw)^\mu=\sum_{\nu=0}^n w^\nu\partial_\nu w^\mu.
\end{cases}
\end{equation}
The components of $D_hw$ contain at most second derivatives of $h$.
Thus the second term and the differentiated lower-order polynomial
involve derivatives of order at most $\ell+1$; only the first term
contains derivatives of order $\ell+2$. This proves the formula by
induction.

Choose $\mu$ such that $w^\mu\ne0$ and take the Jacobian minor
\[
\mathsf J_\mu\triangleq
\frac{\partial(h,Q_h,D_hQ_h,\ldots,D_h^rQ_h)}
{\partial(h,\partial_\mu h,\partial_\mu^2h,\ldots,\partial_\mu^{r+1}h)}.
\]
For $1\le\ell\le r$, the entry $D_h^\ell Q_h$ involves no derivative
above order $\ell+1$. Its coefficient of the pure derivative $\partial_\mu^{\ell+1}h$ is
$2(w^\mu)^{\ell+1}$. Together with \eqref{ce:completion-Q-first-derivative} and \eqref{ce:completion-chain-induction},
this gives the lower triangular matrix, shown schematically as
\[
\mathsf J_\mu=
\begin{pmatrix}
1&0&0&\cdots&0\\
0&2w^\mu&0&\cdots&0\\
0&*&2(w^\mu)^2&\cdots&0\\
\vdots&\vdots&\vdots&\ddots&\vdots\\
0&*&*&\cdots&2(w^\mu)^{r+1}
\end{pmatrix},
\]
where the stars are entries below the diagonal. Consequently,
\[
\det\mathsf J_\mu
=\prod_{j=1}^{r+1}2(w^\mu)^j
=2^{r+1}(w^\mu)^{(r+1)(r+2)/2}\ne0.
\]
The $r+2$ equations defining $\mathcal Z_1$ therefore have independent
differentials on $\{\mathrm dh\ne0\}$, so $\mathcal Z_1$ is a smooth
submanifold of codimension $r+2$. With $r=n$ we obtain
\[
\operatorname{codim}\mathcal Z_0
=\operatorname{codim}\mathcal Z_1=n+2>n+1.
\]

\smallskip\noindent\textit{Finite perturbations that specify all the derivative data.}
The middle layer is compact and lies strictly inside the permitted
perturbation region. Choose finitely many coordinate neighborhoods
$V_j,V'_j$, maps $\kappa_j:V'_j\to\R^{n+1}$, and cutoffs $\chi_j$ such
that
\begin{equation}\label{ce:completion-perturbation-cover}
\begin{aligned}
&V_j\Subset V'_j\Subset
(\R/T\mathbb Z)\times\{\delta<X<1-\delta\},
\qquad j\in\{1,\ldots,L\},\\
&\chi_j\in C_c^\infty(V'_j;[0,1]),\qquad
\chi_j=1\quad\hbox{on }V_j,\\
&(\R/T\mathbb Z)\times\{2\delta\le X\le1-2\delta\}
\subset\bigcup_{j=1}^L V_j.
\end{aligned}
\end{equation}
Define smooth functions by zero extension:
\begin{equation}\label{ce:completion-perturbation-monomials}
\varphi_{j,\alpha}(z)\triangleq
\begin{cases}
\chi_j(z)\kappa_j(z)^\alpha,&z\in V'_j,\\
0,&z\notin V'_j,
\end{cases}
\qquad\alpha\in\mathbb N_0^{n+1},\quad |\alpha|\le r+1.
\end{equation}
Reindex the finite pairs $(j,\alpha)$ as $1,\ldots,J$, writing the
corresponding functions as $\varphi_1,\ldots,\varphi_J$.
Fix $z_0\in V_j$ and arbitrary real numbers $\xi_\alpha$, $|\alpha|\le r+1$. Set
\[
P(v)\triangleq\sum_{|\alpha|\le r+1}
\frac{\xi_\alpha}{\alpha!}(v-\kappa_j(z_0))^\alpha,
\qquad v\in\R^{n+1}.
\]
Since $\chi_j=1$ near $z_0$,
\begin{equation}\label{ce:completion-jet-realization}
\partial_{\kappa_j}^{\alpha}[\chi_jP(\kappa_j)](z_0)=\xi_\alpha.
\end{equation}
The zero extension of this function is a linear combination of the
chosen perturbations, because
\[
\chi_jP(\kappa_j)
=\sum_{|\alpha|\le r+1}\frac{\xi_\alpha}{\alpha!}
\sum_{\beta\le\alpha}\binom{\alpha}{\beta}
[-\kappa_j(z_0)]^{\alpha-\beta}\varphi_{j,\beta}.
\]
Thus their linear combinations can specify the value and all
coordinate derivatives through order $r+1$ at $z_0$. A change from
$\kappa_j$ to the spacetime coordinates is invertible on these data,
so they also specify the derivatives used in the chain.

\smallskip\noindent\textit{Choose small parameters by Sard's theorem.}
For $a\in\R^J$, set
\[
q_a\triangleq q_0+\sum_{j=1}^J a_j\varphi_j,
\]
and consider
\[
\mathcal J:\R^J\times\bigcup_{j=1}^L V_j
\longrightarrow J^{r+1}((\R/T\mathbb Z)\times U;\R),
\qquad\mathcal J(a,z)\triangleq j_z^{r+1}q_a.
\]
At fixed $z$, parameter derivatives span all value and derivative
components by \eqref{ce:completion-perturbation-monomials} and \eqref{ce:completion-jet-realization}. The base component of
$\mathcal J(a,z)$ is $z$. Its full differential is therefore surjective.
For $i\in\{0,1\}$, the inverse image
$\mathcal M_i\triangleq\mathcal J^{-1}(\mathcal Z_i)$ is consequently
a smooth submanifold, if nonempty, and
\[
\dim\mathcal M_i=J+(n+1)-(n+2)=J-1.
\]
The projection $\pi_i:\mathcal M_i\to\R^J$, $\pi_i(a,z)=a$, has rank
at most $J-1$ everywhere. Sard's theorem gives
\[
\bigl|\pi_0(\mathcal M_0)\cup\pi_1(\mathcal M_1)\bigr|_{\R^J}=0.
\]
Thus almost every parameter avoids both excluded sets at every point
of $\bigcup_jV_j$. We can also make the parameter arbitrarily small.
Indeed, with $u_a\triangleq q_a-q_0$,
\[
\|u_a\|_{C^1(K)}
\le\sum_{j=1}^J|a_j|\,\|\varphi_j\|_{C^1(K)}\rightarrow0
\quad\hbox{as }a\to0,
\]
and
\[
\begin{aligned}
&\|\mathcal B(\mathrm dq_a,\mathrm dq_a)
-\mathcal B(\mathrm dq_0,\mathrm dq_0)\|_{C^0(K)} \le2\|\mathrm dq_0\|_{C^0(K)}\|\mathrm du_a\|_{C^0(K)}
+\|\mathrm du_a\|_{C^0(K)}^2\rightarrow0.
\end{aligned}
\]
For the prescribed $\epsilon_1>0$, the following conditions therefore
hold on an open neighborhood of $a=0$:
\begin{equation}\label{ce:completion-small-parameters}
\begin{aligned}
&\|q_a-q_0\|_{C^1(K)}<\epsilon_1,\qquad
\|q_a-q_0\|_{C^0(K)}<\min\{\eta_0,c_0\},\\
&\|\mathcal B(\mathrm dq_a,\mathrm dq_a)
-\mathcal B(\mathrm dq_0,\mathrm dq_0)\|_{C^0(K)}<c_0.
\end{aligned}
\end{equation}
An open neighborhood cannot be contained in the two measure-zero
images. Choose $a$ in this neighborhood outside those images and put
$q=q_a$.

\smallskip\noindent\textit{Verify the condition on the whole compact layer $K$.}
At the entrance, if $q=0$ and
$0\le X\le2\delta$, then \eqref{ce:completion-endpoint-bounds} and
\eqref{ce:completion-small-parameters} give
\begin{equation}\label{ce:completion-left-zero-check}
\begin{aligned}
|q_0|&=|q_0-q|<\eta_0,\\
|\mathcal B(\mathrm dq,\mathrm dq)|
&\ge|\mathcal B(\mathrm dq_0,\mathrm dq_0)|
-|\mathcal B(\mathrm dq,\mathrm dq)
-\mathcal B(\mathrm dq_0,\mathrm dq_0)|
>2c_0-c_0=c_0.
\end{aligned}
\end{equation}
Thus the zeroth chain entry is nonzero there. At the other end,
$q_0=A$ and
\begin{equation}\label{ce:completion-right-zero-check}
q\ge A-|q-q_0|>2c_0-c_0=c_0
\quad\hbox{on }\{1-2\delta\le X\le1\},
\end{equation}
so there are no zeros to check. In the middle region, which is covered
by $\bigcup_jV_j$, the chains from \eqref{ce:repair-chain} satisfy
\begin{equation}\label{ce:completion-middle-zero-check}
q(z)=0,\quad\sum_{\ell=0}^r\mathcal C_\ell(z)^2=0
\quad\Longrightarrow\quad j_z^{r+1}q\in\mathcal Z_0\cup\mathcal Z_1.
\end{equation}
This was excluded by the parameter choice. Since the three regions checked in \eqref{ce:completion-left-zero-check}--\eqref{ce:completion-middle-zero-check} cover $K$, we obtain
\begin{equation}\label{ce:completion-chain-positive}
F\triangleq\sum_{\ell=0}^r\mathcal C_\ell^2>0
\quad\hbox{on }K\cap\{q=0\}.
\end{equation}
The supports of the finitely many $\varphi_j$ also give
\[
\begin{aligned}
&\operatorname{supp}(q-q_0)\Subset
(\R/T\mathbb Z)\times\{\delta<X<1-\delta\},\\
&q=f\quad(X\le\delta),\qquad q=A\quad(X\ge1-\delta).
\end{aligned}
\]
Together with \eqref{ce:completion-small-parameters}, these facts
prove the claim in the nonempty case.

\smallskip\noindent\textit{An empty middle layer.}
If $\{2\delta\le X\le1-2\delta\}=\varnothing$, take $q=q_0$ and
$J=0$. Then
\[
\operatorname{supp}(q-q_0)=\varnothing,\qquad
\|q-q_0\|_{C^1(K)}=0<\epsilon_1.
\]
By \eqref{ce:completion-interpolation}, $q_0=f$ where $X<2\delta$
and $q_0=A>0$ where $X>1-2\delta$. Hence
\[
K\cap\{q_0=0\}\subset\{0\le X<2\delta,\ f=0\},\qquad
\mathcal B(\mathrm dq_0,\mathrm dq_0)\ne0
\quad\hbox{on }K\cap\{q_0=0\}.
\]
The claim follows in this case as well; all sums indexed by
$j=1,\ldots,J$ below are empty.

\endgroup
\section*{Appendix C: algebra at a vertex}
\addcontentsline{toc}{section}{Appendix C: algebra at a vertex}
\setcounter{equation}{0}
\renewcommand{\theequation}{C.\arabic{equation}}
\renewcommand{\theHequation}{appendix.C.\arabic{equation}}
\subsection*{Algebraic lemma: the coefficient $27/8$}\proofsteplabel{app:intrinsic-details}{C}
We use the variables and first-order operator $D$ from \eqref{ce:small-model-flow} to prove \eqref{ce:intrinsic-C1}.
\begingroup
The partial derivatives are
\[
\begin{aligned}
Q_{0,{A_{c,0}}}&=-2(1+\nu)^2{A_{c,0}}-2(1+\nu){A_{c,1}}\cos\Theta\cos\Phi
-\bigl(\tfrac12+2(1+\nu)^2\bigr){A_{c,1}}\sin\Theta\sin\Phi,\\
Q_{0,{A_{c,1}}}&=-2{A_{c,1}}-2(1+\nu){A_{c,0}}\cos\Theta\cos\Phi
-\bigl(\tfrac12+2(1+\nu)^2\bigr){A_{c,0}}\sin\Theta\sin\Phi,\\
Q_{0,\Theta}&=2(1+\nu){A_{c,0}}{A_{c,1}}\sin\Theta\cos\Phi
-\bigl(\tfrac12+2(1+\nu)^2\bigr){A_{c,0}}{A_{c,1}}\cos\Theta\sin\Phi,\\
Q_{0,\Phi}&=2(1+\nu){A_{c,0}}{A_{c,1}}\cos\Theta\sin\Phi
-\bigl(\tfrac12+2(1+\nu)^2\bigr){A_{c,0}}{A_{c,1}}\sin\Theta\cos\Phi.
\end{aligned}
\]
The phase terms satisfy
\[
\begin{aligned}
D\Theta
&=\tfrac12[(1+\nu)^2{A_{c,0}}\sin\Theta+({A_{c,1}}/2)\sin\Phi]
+\tfrac12(1+\nu)^2[{A_{c,0}}\sin\Theta+2{A_{c,1}}\sin\Phi],\\
D\Phi
&=[(1+\nu)^2{A_{c,0}}\sin\Theta+({A_{c,1}}/2)\sin\Phi]
+\tfrac14[{A_{c,0}}\sin\Theta+2{A_{c,1}}\sin\Phi].
\end{aligned}
\]
For the coefficient in $DQ_0$, set $\nu=0$, $\theta=\pi+\xi$ and
${A_{c,1}}=-{A_{c,0}}\cos\theta/\cos2\theta$. Then
\begin{equation}\label{ce:intrinsic-phase-expansion}
\begin{aligned}
& D{A_{c,0}}=D{A_{c,1}} =0,\qquad
{A_{c,1}}={A_{c,0}}\bigl(1+\tfrac32\xi^2+O(\xi^4)\bigr),\\
&Q_{0,\Theta} =3{A_{c,0}}^2\xi+O({A_{c,0}}^2\xi^3),\qquad
D\Theta =\tfrac32{A_{c,0}}\xi+O({A_{c,0}}\xi^3),\\
&Q_{0,\Phi} =-\tfrac32{A_{c,0}}^2\xi+O({A_{c,0}}^2\xi^3),\qquad
D\Phi =\tfrac34{A_{c,0}}\xi+O({A_{c,0}}\xi^3),\\
&{A_{c,0}}^{-3}Q_{0,\Theta}D\Theta =\tfrac92\xi^2+O(\xi^4),\qquad
{A_{c,0}}^{-3}Q_{0,\Phi}D\Phi =-\tfrac98\xi^2+O(\xi^4).
\end{aligned}
\end{equation}
At $\xi=0$ and arbitrary small $\nu$, the same substitution gives
\begin{equation}\label{ce:intrinsic-amplitude-expansion}
\begin{aligned}
& Q_{0,{A_{c,0}}}=-2(1+\nu)\nu {A_{c,0}},\qquad D{A_{c,0}} =(1+\nu)\nu {A_{c,0}}^2,\\
& Q_{0,{A_{c,1}}}=2\nu {A_{c,0}},\qquad D{A_{c,1}} =\nu {A_{c,0}}^2,\\
&{A_{c,0}}^{-3}(Q_{0,{A_{c,0}}}D{A_{c,0}}+Q_{0,{A_{c,1}}}D{A_{c,1}}) =-2\nu^2\bigl((1+\nu)^2-1\bigr)
=-4\nu^3-2\nu^4.
\end{aligned}
\end{equation}
The substituted expression is smooth in $\nu$ and even in $\xi$,
because $\cos2\theta>7/8$. Expansions \eqref{ce:intrinsic-phase-expansion} and \eqref{ce:intrinsic-amplitude-expansion} therefore give \eqref{ce:intrinsic-C1}.

\endgroup

\section*{Appendix D: boundary geometry and rays in a thin layer}
\addcontentsline{toc}{section}{Appendix D: boundary geometry and rays in a thin layer}
\setcounter{equation}{0}
\renewcommand{\theequation}{D.\arabic{equation}}
\renewcommand{\theHequation}{appendix.D.\arabic{equation}}
\subsection*{D.1 Metric and curvature of radial graphs}\phantomsection\proofsteplabel{app:geo-radial-details}{D.1}
Write $J=(R^2+|\nabla_\sigma R|^2)^{1/2}$.
For tangent vectors $v,w$ on the sphere,
\[
\begin{aligned}
d\Sigma_R(v)&=Rv+dR(v)\theta,\qquad
n_{\rm in}=\frac{-R\theta+\nabla_\sigma R}{J},\\
\langle d\Sigma_R(v),d\Sigma_R(w)\rangle
&=R^2\sigma(v,w)+dR(v)dR(w),\\
D_v(-R\theta+\nabla_\sigma R)
&=-2dR(v)\theta-Rv+\nabla_v^\sigma\nabla_\sigma R,\\
-\langle D_vn_{\rm in},d\Sigma_R(w)\rangle
&=\frac{R^2\sigma(v,w)+2dR(v)dR(w)
-R\operatorname{Hess}_\sigma R(v,w)}J.
\end{aligned}
\]
Here the derivative of $J^{-1}$ contributes zero because
$\langle -R\theta+\nabla_\sigma R,d\Sigma_R(w)\rangle=0$.
The derivation of \eqref{geo:radial-metric}--\eqref{geo:radial-II} holds in every dimension.

\subsection*{D.2 Semiconvexity across boundary contacts}\phantomsection\proofsteplabel{app:geo-semiconvex-details}{D.2}
We justify the estimate for $d'_{\Omega,\pm}$ across all boundary contacts. Let $\mathcal B$ be the
transverse reflection times and $K$ the glancing times.  The compressed variables $(y,d_\Omega,\eta)$ are continuous, and
\[
K=\{t:d_\Omega(t)=0,\ |\eta(t)|_{g^{-1}}=1\}
\]
is closed.  Every finite accumulation point of $\mathcal B$ belongs to
$K$.  Indeed it lies on the boundary by continuity, and a boundary
point which is not a transverse reflection is glancing.  A transverse
reflection is isolated and has free intervals on its two sides.
Consequently every compact subinterval of $I\setminus K$ contains
only finitely many reflections.

On a free segment,
$d_\Omega''=\operatorname{Hess}s(x',x')\ge-M$; a transverse reflection changes
$d_\Omega'$ from $-b$ to $b$ with $b>0$.  At $k\in K$ the defining pointwise
Hamiltonian equation gives $d_\Omega'(k)=0$.  This is a full time derivative
even if the uncompressed derivative away from reflections is initially
interpreted relative to their complement.  To see the extension, if
$b_j\in\mathcal B$ tends to $k$, choose by continuity
$z_j\notin\mathcal B$ so close to $b_j$ that
\[
|z_j-b_j|=o(|b_j-k|),\qquad
|d_\Omega(z_j)-d_\Omega(b_j)|=o(|b_j-k|).
\]
In fact,
\[
\frac{d_\Omega(b_j)-d_\Omega(k)}{b_j-k}
=\frac{d_\Omega(z_j)-d_\Omega(k)}{z_j-k}\frac{z_j-k}{b_j-k}
+\frac{d_\Omega(b_j)-d_\Omega(z_j)}{b_j-k}\longrightarrow0.
\]
Thus the derivative extends to every approaching sequence.

Put $f(t)=d_\Omega(t)+Mt^2/2$.  On each component of $I\setminus K$ it is
convex, by its free second derivative and upward reflection jumps.
At an endpoint $k\in K$, the full derivative is $f'(k)=Mk$, so the
endpoint slopes of the convex component tend to $Mk$.  These pieces
glue convexly across the entire closed set $K$.  More explicitly, if
$u<v$ lie in different complementary components and $b\le a$ are the
right endpoint of the first and the left endpoint of the second, then
\[
f'_+(u)\le Mb\le Ma\le f'_-(v).
\]
Points of $K$ themselves have slope $Mt$, and slopes within a single
component are ordered by convexity.  Hence slopes are nondecreasing
globally.  This argument places no countability or interval assumption
on $K$, and excludes a negative singular part of $d_\Omega''$ on that set.

\subsection*{D.3 Differentiability of tangential variables}\phantomsection\proofsteplabel{app:geo-shadow-details}{D.3}

Let $z$ be continuous and have relative derivative $V$ off the set
$\mathcal B$ of isolated transverse reflection times, where $V$ is
continuous. For $k\notin\mathcal B$ and $b_j\in\mathcal B$ tending
to $k$, choose $s_j\notin\mathcal B$ with
\[
|s_j-b_j|=o(|b_j-k|),\qquad
|z(s_j)-z(b_j)|=o(|b_j-k|).
\]
The relative derivative then gives
\[
\frac{z(b_j)-z(k)}{b_j-k}
=\frac{z(s_j)-z(k)}{s_j-k}\frac{s_j-k}{b_j-k}
+\frac{z(b_j)-z(s_j)}{b_j-k}\longrightarrow V(k).
\]
Thus the full derivative at $k$ is $V(k)$. At an isolated
$b\in\mathcal B$, integration on either adjacent interval gives
\[
\frac{z(b+h)-z(b)}h
=\frac1h\int_b^{b+h}V(s)\,\mathrm ds\longrightarrow V(b)
\qquad(h\longrightarrow0^\pm).
\]
Consequently $z\in C^1$ and $z'=V$ on the entire interval. Apply
this argument to $z=(y,\eta)$ and the continuous right-hand side of
\eqref{geo:tangential-system}. It also proves that there is no
additional singular tangential motion on a closed glancing set.

\subsection*{D.4 Continuing the comparison through a finite atlas}\phantomsection\proofsteplabel{app:geo-chart-details}{D.4}
The base curve of $\bar z$ is a unit-speed boundary geodesic.
To pass between charts, choose a finite atlas whose smaller chart
domains still cover the compact boundary. The larger domains have
a fixed positive margin around the smaller ones. Since the base
velocities are uniformly bounded, a fixed time step $t_{\rm ch}>0$
keeps both curves in one larger chart, once $\varepsilon$ is small.
A window of length $T$ therefore requires at most
\[
K_T=1+\left\lceil T/t_{\rm ch}\right\rceil
\]
chart intervals. Smooth coordinate changes on their compact overlaps
have bounded derivatives. If $e_j$ is the coordinate error at the
$j$th chart transition, the integral estimate and coordinate change
give, for a fixed $B_T\ge1$,
\[
e_{j+1}\le B_T(e_j+\varepsilon),\qquad e_0\le C\varepsilon,
\qquad
e_j\le C\varepsilon\sum_{\ell=0}^jB_T^\ell
\le C_T\varepsilon\quad(0\le j\le K_T).
\]
Applying \eqref{geo:shadow-gronwall} within each chart interval then gives
\[
\dist_g(q(x(t)),\bar q(t))\le C_T\varepsilon,
\qquad t_*\le t\le t_*+T,
\]
where $\bar q$ is the base curve of $\bar z$, and $C_T$ depends only
on the boundary and $T$. Neither continuity of $\nu$ nor a bound
on the number of reflections enters the integral comparison based on \eqref{geo:shadow-gronwall}.

\subsection*{D.5 Normal acceleration in boundary coordinates}\phantomsection\proofsteplabel{app:geo-normal-details}{D.5}
On a free segment $x''=0$, so
\[
\begin{aligned}
x'&=(I-d_\Omega\mathsf{W})q'+d_\Omega'n_{\rm in}(q),\quad
d_\Omega'=\langle n_{\rm in}(q),x'\rangle,\quad
d_\Omega'' =\langle D_{q'}n_{\rm in},x'\rangle
=-g\bigl(\mathsf{W}q',(I-d_\Omega\mathsf{W})q'\bigr).
\end{aligned}
\]

\subsection*{D.6 Successive differences for the bounded perturbation}\phantomsection\proofsteplabel{app:geo-duhamel-details}{D.6}
With the notation of \eqref{asm:bounded-perturbation}--\eqref{asm:picard-iteration}, the iterates satisfy, first for $j=0$ and then by induction,
\[
\begin{aligned}
\|U^{(1)}(t)-U^{(0)}(t)\|_{X_\Omega}
&\le M_Tt\|U_0\|_{X_\Omega},\\
\|U^{(j+1)}(t)-U^{(j)}(t)\|_{X_\Omega}
&\le M_T\int_0^t
\|U^{(j)}(s)-U^{(j-1)}(s)\|_{X_\Omega}\,\mathrm ds\\
&\le \|U_0\|_{X_\Omega}
\frac{M_T^{j+1}}{j!}\int_0^t s^j\,\mathrm ds
=\|U_0\|_{X_\Omega}\frac{(M_Tt)^{j+1}}{(j+1)!}
\quad(j\ge1).
\end{aligned}
\]
The series of successive differences is uniformly summable because
$\sum_{j\ge0}(M_TT)^{j+1}/(j+1)!=e^{M_TT}-1$.

\end{document}